%% file: A_Thesis.tex
\documentclass{dcthesis}

\committee[F. Jon Kull, Ph.D.]{}{}{}{}
\school{Dartmouth College}{Hanover, New Hampshire}
\degree{Doctor of Philosophy}
\field{Mathematics}

\usepackage{amsfonts, amssymb, amsthm, amsmath, braket, hyperref, mathtools, comment, mathrsfs, enumitem}
\usepackage[usenames,dvipsnames, table]{xcolor}
\usepackage[all]{xy}
\usepackage[T1]{fontenc}
\usepackage{tikz-cd}
\usepackage[nameinlink]{cleveref}

\newtheorem{theorem}[equation]{Theorem}

\newtheorem{corollary}[equation]{Corollary}
\newtheorem{cor}[equation]{Corollary}

\newtheorem{lemma}[equation]{Lemma}
\newtheorem{proposition}[equation]{Proposition}
\newtheorem{prop}[equation]{Proposition}
\newtheorem{example}[equation]{Example}

\theoremstyle{definition}
\newtheorem{definition}[equation]{Definition}

\theoremstyle{remark}
\newtheorem{remark}[equation]{Remark}

\numberwithin{equation}{section}

\newcommand{\defi}[1]{\textsf{#1}} %

\newenvironment{enumalg}
{\begin{enumerate}}
{\end{enumerate}}

\newcommand{\jvtable}[3]
{\begin{equation} \label{#1}\addtocounter{equation}{1} \notag
\begin{gathered}
#2 \\
\text{
\parbox[c]{4.5in}{\centering Table \ref*{#1}: #3}}
\end{gathered}
\end{equation}
}

\newcommand{\jvfigure}[3]
{\begin{equation} \label{#1}\addtocounter{equation}{1} \notag
\begin{gathered}
#2 \\
\text{
\parbox[c]{4.5in}{\centering Figure \ref*{#1}: #3}}
\end{gathered}
\end{equation}
}

\newcommand{\Z}{\mathbb{Z}}
\newcommand{\PP}{\mathbb{P}}

\newcommand{\bbC}{\mathbb{C}}
\newcommand{\bbP}{\mathbb{P}}
\newcommand{\bbQ}{\mathbb{Q}}
\newcommand{\bbR}{\mathbb{R}}
\newcommand{\bbF}{\mathbb{F}}
\newcommand{\bbZ}{\mathbb{Z}}

\newcommand{\repr}[1]{\gamma_{#1}}

\newcommand{\calE}{\mathcal{E}}

\newcommand{\calR}[1]{\mathcal{R}_{#1}}
\newcommand{\calT}[1]{\mathcal{T}_{#1}}

\newcommand{\R}[1]{R_{#1}}
\newcommand{\Q}[1]{Q_{#1}}
\newcommand{\Qly}[1]{Q_{#1,\leq Y}}

\newcommand{\scrE}{\mathscr{E}}
\newcommand{\scrEX}{\mathscr{E}_{\leq X}}
\newcommand{\scrO}{\mathscr{O}}

\newcommand{\Qalg}{\bbQ^{\textup{al}}}
\newcommand{\bbQalg}{\Qalg}

\newcommand{\prm}{^\prime}
\newcommand{\parent}[1]{\left(#1\right)}

\newcommand{\inv}{^{-1}}
\newcommand{\GL}{\mathrm{GL}}
\newcommand{\abs}[1]{\left|#1\right|}
\DeclareMathOperator{\Gal}{Gal}
\DeclareMathOperator{\ord}{ord}
\DeclareMathOperator{\Nm}{Nm}
\DeclareMathOperator{\repart}{Re}

\DeclareMathOperator{\Area}{area}
\DeclareMathOperator{\len}{len}

\DeclareMathOperator{\lcm}{lcm}
\DeclareMathOperator{\disc}{disc}

\newcommand{\ceil}[1]{\left\lceil #1 \right\rceil}
\newcommand{\floor}[1]{\left\lfloor #1 \right\rfloor}

\newcommand{\gcdP}[3]{P_{#1, #2}^{(#3)}}
\newcommand{\gcdQ}[3]{Q_{#1, #2}^{(#3)}}
\newcommand{\gcdm}[3]{m_{#1, #2}^{(#3)}}

\newenvironment{enumalph}
{\begin{enumerate}}
	{\end{enumerate}}

\newenvironment{enumroman}
{\begin{enumerate}}
	{\end{enumerate}}

\usepackage{colonequals}
\usepackage{hyperref}
\hypersetup{colorlinks=true,urlcolor=blue,citecolor=blue,linkcolor=blue}

\renewcommand{\set}[1]{\left\{#1\right\}}

\newcommand{\heaviside}{\theta}

\DeclareMathOperator{\height}{ht}
\DeclareMathOperator{\hht}{\height}
\DeclareMathOperator{\twistheight}{twht}
\DeclareMathOperator{\twht}{\twistheight}
\newcommand{\rawheight}{H}

\DeclareMathOperator{\Res}{Res}
\DeclareMathOperator{\res}{res}

\newcommand{\psmod}[1]{~(\textup{\text{mod}}~{#1})}

\newcommand{\A}[1]{A_{#1}}
\newcommand{\B}[1]{B_{#1}}
\newcommand{\tA}[1]{A\prm_{#1}}
\newcommand{\tB}[1]{B\prm_{#1}}
\newcommand{\C}[1]{C_{#1}}

\newcommand{\calcc}[1]{\mathcal{C}_{#1}}

\newcommand{\f}[1]{f_{#1}}
\newcommand{\g}[1]{g_{#1}}
\newcommand{\h}[1]{h_{#1}}
\newcommand{\hII}[1]{h_{#1, \textup{II}}}
\newcommand{\hIII}[1]{h_{#1, \textup{III}}}
\newcommand{\CII}[1]{C_{#1, \textup{II}}}
\newcommand{\CIII}[1]{C_{#1, \textup{III}}}

\newcommand{\tf}[1]{f\prm_{#1}}
\newcommand{\tg}[1]{g\prm_{#1}}

\newcommand{\degB}[1]{d\parent{#1}}
\newcommand{\twerror}[1]{e\parent{#1}}

\newcommand{\upperratio}[1]{\kappa_{#1}}

\newcommand{\cusps}[1]{\mathscr{C}_{#1}}

\DeclareMathOperator{\mindefect}{md}
\newcommand{\md}{\mindefect}
\DeclareMathOperator{\twistdefect}{tmd}
\newcommand{\tmd}{\twistdefect}

\newcommand{\twistEX}{\scrE^{{\rm tw}}_{\leq X}}
\newcommand{\twistE}{\scrE^{\textup{tw}}}

\newcommand{\twistNad}[1]{N_{#1}^{\textup{tw}}}
\newcommand{\twNad}[1]{\twistNad{#1}}

\newcommand{\twistNeq}[1]{\widetilde{N}_{#1}^{\textup{tw}}}
\newcommand{\twNeq}[1]{\twistNeq{#1}}

\newcommand{\twistNeqly}[1]{\widetilde{N}_{#1, \leq y}^{\rm tw}}
\newcommand{\twistNeqgy}[1]{\widetilde{N}_{#1, > y}^{\rm tw}}
\newcommand{\NQad}[1]{N_{#1}}
\newcommand{\NQeq}[1]{\widetilde{N}_{#1}}
\newcommand{\twistNadly}[1]{N_{#1, \leq y}^{\rm tw}}
\newcommand{\twistNadgy}[1]{N_{#1, > y}^{\rm tw}}
\newcommand{\twNadly}[1]{\twistNadly{#1}}
\newcommand{\twNadgy}[1]{\twistNadgy{#1}}

\newcommand{\cM}[1]{M_{#1}}

\newcommand{\cT}[1]{T_{#1}}

\newcommand{\tcT}[1]{\widetilde{T}_{#1}}

\newcommand{\tcalT}[1]{\widetilde{\mathcal{T}}_{#1}}

\newcommand{\LQeq}[1]{\widetilde{L}_{#1}}

\newcommand{\twistheq}[1]{\Delta \twistNeq {#1}}%
\newcommand{\twisthad}[1]{\Delta \twistNad {#1}}%
\newcommand{\hQad}[1]{\Delta \NQad {#1}}%
\newcommand{\hQeq}[1]{\Delta \NQeq {#1}}%

\newcommand{\twistLad}[1]{L^{\rm tw}_{#1}}

\newcommand{\twistLadR}[1]{L_{#1,\rm{rem}}^{\rm tw}}
\newcommand{\LQad}[1]{L_{#1}}

\newcommand{\twistLeq}[1]{\widetilde{L}^{\rm tw}_{#1}}

\newcommand{\constad}[1]{c_{#1}}
\newcommand{\constadprm}[1]{c\prm_{#1}}
\newcommand{\consteq}[1]{\widetilde{c}_{#1}}
\newcommand{\consteqprm}[1]{\widetilde{c}\prm_{#1}}
\newcommand{\twconst}[1]{c_{#1}^{\textup{tw}}}
\newcommand{\twconstad}[1]{c_{#1}^{\textup{tw}}}
\newcommand{\twconsteq}[1]{\widetilde{c}_{#1}^{\textup{tw}}}
\newcommand{\twconstprm}[1]{c^{\textup{tw}\prime}_{#1}}

\newcommand{\kfree}[1]{S_{#1}}

\title{Counting elliptic curves with a cyclic $m$-isogeny over $\bbQ$}
\author{Grant Molnar}
\date{April 24, 2023}
\field{Mathematics}
\degree{Doctor of Philosophy}
\committee{John Voight}{Asher Auel}{Robert Lemke Oliver}{Carl Pomerance}

\begin{document}

\frontmatter

\maketitle

\chapter*{Abstract}
\addcontentsline{toc}{section}{Abstract}
Using methods from analytic number theory, for $m > 5$ and for $m = 4$, we obtain asymptotics with power-saving error terms for counts of elliptic curves with a cyclic $m$-isogeny up to quadratic twist over the rational numbers. For $m > 5$, we then apply a Tauberian theorem to achieve asymptotics with power-saving error for counts of elliptic curves with a cyclic $m$-isogeny up to isomorphism over the rational numbers.

\chapter*{Preface}
\addcontentsline{toc}{section}{Preface}

Throughout my academic and personal life, I have been blessed with the support and guidance of many wonderful people. I would like to begin by honoring my mother and father, Wanda and Steven Molnar, for giving me a firm foundation to stand on and for fostering my love of exploration and enquiry. My parents have given me roots and wings throughout my life, and I am who I am because of them.

In the same breath, I recognize my wife, Brianna Molnar, for her heartfelt support, encouragement, and understanding. I like to think of myself as eloquent, but even perfectly placed words cannot capture or convey my affection and appreciation for her. Brianna is brave, and creative, and strong, and true. I love Brianna, and I have learned a great deal from her example.

My thanks go out to my thesis advisor, John Voight, whose expertise, enthusiasm, and example have been instrumental in teaching me to do mathematics gracefully. I came to Dartmouth College to work with John, and I cannot imagine a better advisor than he is. John's patient mentorship and constructive feedback have been essential in guiding my work: were it not for John, this thesis would not exist.

I likewise thank the members of my thesis committee, John Voight, Asher Auel, Robert Lemke Oliver, and Carl Pomerance, for taking the time to read and refine my work.

I am grateful to Eran Assaf, Jesse Elliott, Mits Kobayashi, David Lowry-Duda, Tristan Phillips, and Rakvi for their helpful comments and discussions, which have greatly improved the quality of my thesis.

I am deeply delighted to have made this academic voyage alongside my cohort: Lizzie Buchanan, Steve Fan, Richard Haburcak, and Alexander Wilson. Their support and encouragement have been invaluable, and I could not have asked for better colleagues. 

Nor are they the only fellow-travelers whom I now salute. Along with Steve Fan, I thank Jennifer Molnar and Taylor Petty for commiserating with me about the challenges of thesis work. Both their griping and goodwill kept me going. I also want to extol the 2021-2022 Dartmouth Graduate Student Council Executive Board and the Dartmouth Graduate Student Council Ad Hoc Healthcare Committee, for their work and for our friendship. Most especially, I thank Elizabeth Bien, Keighley Rockcliffe, and Kelly Cantwell. Their energy and dedication to improving the graduate student experience have been inspiring, and I am privileged to have worked alongside them.

I laud the tireless and ofttimes thankless work of the Baker-Berry Library staff, who processed and accommodated hundreds of my requests for books and papers during the half-decade I have been at Dartmouth. Finally, I would like to acknowledge the Gridley Fund for Graduate Mathematics and the Simons Collaboration (grant 550029, to John Voight) for their financial support, which enabled me to pursue my research.

\input{Corrections_Note.tex}

\tableofcontents

\mainmatter

\include{Ch1_Introduction}

\include{Ch2_Background}

\include{Ch3_Preliminaries}

\include{Ch4_m7}

\include{Ch5_m10_25}

\include{Ch6_m13}

\include{Ch7_m4_6_8_9_12_16_18}

\include{Ch8_m11_14_15_17_19_21_27_37_43_67_163}

\backmatter

\bibliographystyle{amsplain}
\addcontentsline{toc}{chapter}{References}
\bibliography{isogenySources}

\end{document}

%% file: Corrections_Note.tex
\chapter*{Note on this corrected edition}

It's been over three years since I defended my dissertation, and in the last few months I have finally circled back to doing the sort of follow-on research that every academic feels an urge to pick up almost as soon as his thesis is set down. As I tried to build on my isogeny research, I discovered typos, proof gaps, and the occasional error. I've made a dedicated effort to uproot every one of these little troublemakers, but if I have missed a couple then I hope you will be forgiving.

The two most consequential corrections concern the M\"obius-sieve argument for large twist-minimality defect and the $m=13$ argument in Chapter~6. Both required repairs to the error terms, though the main terms are unchanged. I have also corrected a small number of numerical constants, finite tables, and cross-references.

One deliberate addition goes beyond correction: I complete the $13$-adic local computation left unfinished in the original thesis, determining $Q_{13}(13)$ and permitting a numerical estimate of the twist-counting constant $c_{13}^{\mathrm{tw}}$.

A page-by-page record of the changes appears in the \href{https://grantmolnar.com/pdfs/Counting_elliptic_curves_with_a_cyclic_m_isogeny_over_Q_errata.pdf}{separate errata}.

Thanks for reading!

\begin{flushright}
Grant Molnar\\
September 24, 2026
\end{flushright}

%% file: Ch1_Introduction.tex
\chapter{Introduction}\label{Chapter: Introduction}

In \cref{Section: Motivation and setup}, we briefly motivate the study of the count of elliptic curves with a cyclic $m$-isogeny over $\bbQ$, and establish the notation necessary to state our main results. In \cref{Section: Results}, we report on our asymptotics for elliptic curves equipped with or admitting a cyclic $m$-isogeny over $\bbQ$, first up to $\bbQ$-isomorphism, and then up to twist equivalence. In \cref{Section: Our approach}, we sketch our approach for proving these asymptotics. Finally, in \cref{Section: Contents}, we outline the structure for the remainder of this thesis.

\section{Motivation and setup}\label{Section: Motivation and setup}

In this section, we outline the scope of our results and give impetus for studying the count of elliptic curves with a cyclic $m$-isogeny over $\bbQ$ whose na\"ive height is less than or equal to $X$ (readers unfamiliar with elliptic curves or cyclic isogenies are directed to \cref{Section: A primer on elliptic curves}). We then set up the notation necessary to state our results.

\subsection*{Motivation}

Elliptic curves have been an object of fascination for number theorists and geometers for over a century. Much effort has gone into developing tools to understand the behavior of particular elliptic curves. But in the last twenty years, there has been an explosion of interest %
in the statistical behavior of families of elliptic curves \cite{Bektemirov-Mazur-Stein, Bhargava-Shankar, Boggess-Sankar, Bruin-Najman, Cullinan-Kenney-Voight, Dabrowski-Pomykala, Duke, Fuchs-vonKanel-Wustholz, Grant, Harron-Snowden, Molnar-Voight, Phillips1, Pizzo-Pomerance-Voight, Pomerance-Schaefer, Poonen, Watkins, Young}. For more on the history of the arithmetic statistics of elliptic curves, see \cref{Section: History}.

In this thesis, we recount and strengthen arguments made in \cite{Molnar-Voight} to estimate the number of elliptic curves equipped with (or admitting) a cyclic $7$-isogeny over $\bbQ$. We then go further, and adapt the methods of \cite{Molnar-Voight} to estimate the number of elliptic curves equipped with (or admitting) a cyclic $m$-isogeny for 
\[
    m \in \set{6, 7, 8, 9, 10, 12, 13, 16, 18, 25}.
\]
Our work improves on prior results \cite{Molnar-Voight, Phillips1}, and gives entirely new asymptotics for $m \in \set{10, 13, 25}$. Finally, we give asymptotics for the number of elliptic curves equipped with (or admitting) a cyclic $m$-isogeny for 
\[
m \in \set{11, 14, 15, 17, 19, 21, 27, 37, 43, 67, 163},
\]
where the associated compactified modular curve $X_0(m)$ has nonzero genus. 

In conjunction with earlier work, for all $m \in \bbZ_{>0}$ except $m = 5$, we establish asymptotics for the number of elliptic curves equipped with or admitting a cyclic $m$-isogeny over $\bbQ$ (see \cref{Section: Results} and \cref{Section: Previous results}).

These asymptotics are a natural area of study for several reasons. Concretely, individual elliptic curves can be rather delicate objects to work with, and it is interesting to ask what behavior is typical of elliptic curves: for instance, \cref{Section: Results} and \cref{Section: Previous results} together show that almost all elliptic curves $E/\bbQ$ have no cyclic isogenies besides the trivial automorphisms $\pm 1: E \xrightarrow{\sim} E$. As a consequence, counting elliptic curves by height up to (cyclic) isogeny yields the same asymptotics as counting elliptic curves by height up to $\bbQ$-isomorphism.

In addition, the modular curve $Y_0(m) \subseteq X_0(m)$ is a moduli space for elliptic curves equipped with a cyclic $m$-isogeny, so counting elliptic curves equipped with a cyclic $m$-isogeny over $\bbQ$ may serve as an example and prototype for counting elliptic curves with other level structure (see \cite{Cullinan-Kenney-Voight, PhillipsThesis}). More abstractly, counting elliptic curves with a cyclic $m$-isogeny requires overtly or implicitly grappling with the ``stackiness'' of modular curves like $X_0(7)$ and therefore has implications for the Batyrev–Manin conjecture and even larger questions in arithmetic geometry \cite{Ellenberg-Satriano-Zureick-Brown}.

\subsection*{Setup}

In this subsection, we establish notation which will be necessary to state our main theorems. We then recall notation for asymptotics from analytic number theory.

We begin by setting up a fragment of the theory of elliptic curves. Every elliptic curve $E$ over $\bbQ$ has a unique minimal Weierstrass model of the form
\begin{equation}\label{Equation: Weierstrass equation}
E : y^2 = x^3 + A x + B,
\end{equation}
where $A, B \in \bbZ$, $4A^3 + 27B^2 \neq 0$, and for every prime $\ell$ we have $\ell^4 \nmid A$ or $\ell^6 \nmid B$. If the minimal model for $E$ is given by \eqref{Equation: Weierstrass equation}, we define the \defi{(na\"ive) height} of $E$ to be
\begin{equation}\label{Equation: Intro height}
\hht(E) \colonequals \max(4 \abs{A}^3, 27 \abs{B}^2).
\end{equation}
Let $\scrE$ be the set of elliptic curves over $\bbQ$ in their minimal model, and let 
\begin{equation}
\scrEX \colonequals \{E \in \scrE : \hht(E) \leq X\}.
\end{equation}
For $m \in \bbZ_{>0}$ and $E, E\prm \in \scrE$, a \defi{cyclic $m$-isogeny} $\phi : E \to E\prm$ is a morphism of elliptic curves such that $\ker \phi \subseteq E(\bbQalg)$ is a cyclic group of order $m$ (the unfamiliar reader may peruse \cref{Section: A primer on elliptic curves} for more information). Here, as usual, $\bbQalg$ denotes the algebraic closure of $\bbQ$. In this thesis, all isogenies are defined over $\bbQ$ unless otherwise indicated. An \defi{unsigned isogeny} is an isogeny up to postcomposition by $\pm 1$. 

We define
\begin{equation}\label{Equation: counts of m-isogenies}
\begin{aligned}
\NQeq m (X) &\colonequals \# \set{(E, \phi) : E \in \scrEX \ \text{and} \ \phi : E \to E\prm \ \textup{an unsigned cyclic $m$-isogeny}}, \\
\NQad m (X) &\colonequals \# \set{E \in \scrEX : \textup{$E$ admits a cyclic $m$-isogeny}},
\end{aligned}
\end{equation}
where as usual the cyclic $m$-isogeny $\phi$ is defined over $\bbQ$, and $E\prm \in \scrE$.

Let $E$ have a Weierstrass model
\begin{equation}
E : y^2 = x^3 + A x + B,
\end{equation}
which is not necessarily minimal, i.e., we might have $d^4 \mid A$ and $d^6 \mid B$ for some $d > 1$. For $c \in \bbQ^\times$, the \defi{quadratic twist} $E^{(c)}$ of $E$ via $c$ is defined by the Weierstrass equation
\begin{equation}\label{Equation: Defining E^(c)}
E^{(c)} : y^2 = x^3 + c^2 A x + c^3 B,
\end{equation}
and we have a $\bbQalg$-isomorphism $E \xrightarrow{\sim} E^{(c)}$ given by $(x, y) \mapsto (c x, c^{3/2} y)$. We say $E, E\prm \in \scrE$ are \defi{twist equivalent} if $E\prm = E^{(c)}$ for some $c \in \bbQ^\times$; if $j(E) \neq 0, 1728$, then $E, E\prm \in \scrE$ are twist equivalent if and only if they are $\bbQalg$-isomorphic (see \Cref{Corollary: Isomorphisms are quadratic twists for j not 0 and 1728}). 

We let $\twistE$ denote the set of elliptic curves over $\bbQ$ up to twist equivalence. We define the \defi{twist height} of $E$ to be
\begin{equation}\label{Equation: twist height}
\twistheight(E) \colonequals \min\set{\hht(E\prm) : E\prm \in \scrE \ \text{is twist equivalent to} \ E},
\end{equation}
and we let 
\begin{equation}
\twistEX \colonequals \{E \in \twistE : \twistheight(E) \leq X\}.
\end{equation}
Twist equivalence preserves (cyclic) isogenies (\Cref{Corollary: isogenies are stablized by quadratic twist}), so it is natural to define 
\begin{equation}\label{Equation: twist counts of m-isogenies}
\begin{aligned}
\twistNeq m (X) &\colonequals \# \set{(E, \phi) : E \in \twistEX \ \text{and} \ \phi : E \to E\prm \ \textup{an unsigned cyclic $m$-isogeny}}, \\
\twistNad m (X) &\colonequals \# \set{E \in \twistEX : \textup{$E$ admits a cyclic $m$-isogeny}}
\end{aligned}
\end{equation}
(as above, the $m$-isogeny $\phi$ is defined over $\bbQ$, and $E\prm \in \scrE$). The functions defined in \eqref{Equation: counts of m-isogenies} and \eqref{Equation: twist counts of m-isogenies} are the main objects of study in this thesis.

We adopt the following notations from analytic number theory. For eventually positive functions $f, g, h : \bbR_{>0} \to \bbR$, we write
\begin{equation}\label{Equation: big oh}
f(X) = g(X) + O(h(X))
\end{equation}
for $X \geq X_0$ if there is a constant $C$ such that for all $X \geq X_0$ we have
\begin{equation}\label{Equation: big-oh bound}
\abs{f(X) - g(X)} < C h(X).
\end{equation}
If we write \eqref{Equation: big oh} without specifying an $X_0$, then \eqref{Equation: big-oh bound} holds for all $X$ sufficiently large. If $f(X) = O(g(X))$, we also may write $f(X) \ll g(X)$.

Similarly, we write
\begin{equation}\label{Equation: little oh}
f(X) = g(X) + o(h(X))
\end{equation}
if
\begin{equation}
\lim_{X \to \infty} \abs{\frac{f(X) - g(X)}{h(X)}} = 0.
\end{equation}
If $g(X) = h(X)$ in \eqref{Equation: little oh}, we write $f(X) \sim g(X)$. As the notation suggests, this is an equivalence relation.

Finally, we write
\begin{equation}
f(X) \asymp g(X)
\end{equation}
if there are constants $C_1, C_2 \in \bbR_{>0}$ such that 
\begin{equation}
C_1 f(X) < g(X) < C_2 f(X)
\end{equation}
for all $X$ sufficiently large.

\section{Results}\label{Section: Results}

In this section, we give asymptotics for $\NQeq m (X)$ and $\NQad m (X)$ for $m > 5$. We then give asymptotics for $\twistNeq m (X)$ and $\twistNad m (X)$ for $m > 5$ and for $m = 4$. To our knowledge, for $m \in \set{7, 10, 13, 25}$, even the order of growth for $\NQad m (X)$ and $\NQeq m (X)$ was previously unknown (see \cite[Remark 4.2]{Boggess-Sankar}).

\subsection*{Main results}

In this subsection, we present asymptotics for $\NQeq m (X)$ and $\NQad m (X)$ for all $m$ such that the compactified modular curve $X_0(m)$ has genus $0$ and $m > 5$, i.e., for
\[
m \in \set{6, 7, 8, 9, 10, 12, 13, 16, 18, 25}.
\]
We then present results for all $m$ such that $X_0(m)$ has nonzero genus and the noncompactified modular curve $Y_0(m)$ has $Y_0(m)(\bbQ) \neq \emptyset$, i.e., for
\[
m \in \set{11, 14, 15, 17, 19, 21, 27, 37, 43, 67, 163}.
\]

\begin{theorem}\label{Intro Theorem: asymptotic for NQ(X) for 5 < m <= 9} 
    Let $m \in \set{6, 7, 8, 9}$. Then there are effectively computable constants $\consteq m$, $\consteqprm m$, $\constad m$, and $\constadprm m$ such that for any $\epsilon > 0$, we have
    \begin{equation}
    \NQeq m (X) = \consteq m X^{1/6} \log X + \consteqprm m X^{1/6} + O(X^{38/291 + \epsilon})
    \end{equation}
    and
    \begin{equation}
    \NQad m (X) = \constad m X^{1/6} \log X + \constadprm m X^{1/6} + O(X^{38/291 + \epsilon})
    \end{equation}
    for $X \geq 1$. The implicit constant depends on $m$ and $\epsilon$.
\end{theorem}

\begin{theorem}\label{Intro Theorem: asymptotic for NQ(X) for m > 9}
    Let $m \in \set{10, 12, 13, 16, 18, 25}$. Then there are effectively computable constants $\consteq m$ and $\constad m$ such that for any $\epsilon > 0$, we have
    \begin{equation}
    \NQeq m (X) = \consteq m X^{1/6} + O(X^{1/8 + \epsilon})
    \end{equation}
    and
    \begin{equation}
    \NQad m (X) = \constad m X^{1/6} + O(X^{1/8 + \epsilon})
    \end{equation}
    for $X \geq 1$. The implicit constant depends on $m$ and $\epsilon$.
\end{theorem}

In both \Cref{Intro Theorem: asymptotic for NQ(X) for 5 < m <= 9}  and \Cref{Intro Theorem: asymptotic for NQ(X) for m > 9}, the constants $\consteq m$ and $\constad m$ are positive, but the constants $\consteqprm m$ and $\constadprm m$ need not be: for instance, when $m = 7$, we have
\begin{equation}
\consteqprm 7 = \constadprm 7 \approx -0.16.
\end{equation}

\Cref{Intro Theorem: asymptotic for NQ(X) for 5 < m <= 9} summarizes results given in \Cref{Theorem: asymptotic for NQ(X) for m = 7}, \Cref{Theorem: asymptotic for twN(X) for m of genus $0$}, and \Cref{Corollary: asymptotic for twN(X) for m of genus $0$}; likewise, \Cref{Intro Theorem: asymptotic for NQ(X) for m > 9} summarizes results given in \Cref{Theorem: asymptotic for NQ(X) for m = 10 and 25}, \Cref{Theorem: asymptotic for NQ(X) for m = 13}, \Cref{Theorem: asymptotic for twN(X) for m of genus $0$}, and \Cref{Corollary: asymptotic for twN(X) for m of genus $0$}. \Cref{Intro Theorem: asymptotic for NQ(X) for 5 < m <= 9} and \Cref{Intro Theorem: asymptotic for NQ(X) for m > 9} extend and strengthen results in the literature \cite{Phillips1}; see \cref{Section: History} for more details.

The cases $m \in \set{7, 10, 13, 25}$ are of special interest, because their associated modular curves $X_0(m)$ have multiple elliptic points. Consequently, for these $m$, the elliptic surfaces that parameterize elliptic curves equipped with a cyclic $m$-isogeny have points of additive reduction.

We now present the asymptotics for $\NQeq m (X) = \NQad m (X)$ when $X_0(m)$ has nonzero genus and $Y_0(m) \subseteq X_0(m)$ has $Y_0(m)(\bbQ) \neq \emptyset$.

\begin{theorem}[\Cref{Theorem: asymptotic for NQ(X) for m of nonzero genus} and \Cref{Theorem: asymptotic for NQ(X) for m of nonzero genus given RH}]\label{Intro Theorem: asymptotic for NQ(X) for m of nonzero genus}
	Let 
    \[
    m \in \set{11, 14, 15, 17, 19, 21, 27, 37, 43, 67, 163}.
    \]
    Then there is a positive, effectively computable constant $\constad m$ such that
	\begin{equation}
	\NQeq m (X) = \NQad m (X) = \constad m X^{1/6} + o\parent{X^{1/12}}.
	\end{equation}
	If the Riemann hypothesis holds, then for any $\epsilon > 0$, we may replace the error term with $O\parent{X^{11/210 + \epsilon}}$ for $X \geq 1$, with the implicit constant now depending on $\epsilon$.
\end{theorem}
\Cref{Intro Theorem: asymptotic for NQ(X) for m of nonzero genus} is almost immediate from Walfisz's and Liu's estimates for counts of squarefree integers \cite{Liu, Walfisz}, both of which utilize zero-free regions for the Riemann zeta function. In fact, Walfisz's estimate gives a slight improvement on the $o(X^{1/12})$ error in \Cref{Intro Theorem: asymptotic for NQ(X) for m of nonzero genus}. The equality $\NQeq m (X) = \NQad m (X)$ in \Cref{Intro Theorem: asymptotic for NQ(X) for m of nonzero genus} is exact: no elliptic curve over $\bbQ$ admits more than one cyclic $m$-isogeny when $X_0(m)(\bbQ) < \infty$. The constants $\constad m$ are given explicitly in Table \ref{table:constantsfornonzerogenus}.%

The cases $m \in \set{1, 2, 3, 4, 5}$ were handled by previous authors, and we report them in \cref{Section: Previous results}; in the case $m = 5$, only the order of growth for $\NQeq 5 (X)$ and $\NQad 5 (X)$ is given (\Cref{Intro Theorem: asymptotics for N(X) = 5}). For all $m > 5$ not addressed by \Cref{Intro Theorem: asymptotic for NQ(X) for 5 < m <= 9}, \Cref{Intro Theorem: asymptotic for NQ(X) for m > 9}, and \Cref{Intro Theorem: asymptotic for NQ(X) for m of nonzero genus}, Mazur's theorem on isogenies (\Cref{Theorem: Kenku-Mazer Theorem}) implies $\NQeq m (X) = \NQad m (X) = 0$ identically. Thus for all $m \neq 5$, we have asymptotics with power-saving error for $\NQeq m (X)$ and $\NQad m (X)$. Of course, much work remains to be done in counting elliptic curves with a cyclic $m$-isogeny over global fields.

\subsection*{Twist results}

In this subsection, we present asymptotics for $\twistNeq m (X)$ and $\twistNad m (X)$ for all $m$ such that $X_0(m)$ has genus $0$ and $m \not\in \set{2, 3, 5}$. These asymptotics are stepping stones to the results in the previous subsection, but we view them as natural and interesting in their own right. 

\begin{theorem}\label{Intro Theorem: asymptotic for twN(X) for m of genus $0$}
	Let $m \in \set{1, 4, 6, 7, 8, 9, 10, 12, 13, 16, 18, 25}$. Then there are positive, effectively computable constants $\twconstad m$ and $\twconsteq m$ such that for all $\epsilon > 0$, we have
	\begin{equation}
	\twistNeq m (X) = \twconsteq m X^{1/\degB m} + O\parent{X^{1/\twerror m + \epsilon}}
	\end{equation}
        and
        \begin{equation}
	\twistNad m (X) = \twconstad m X^{1/\degB m} + O\parent{X^{1/\twerror m + \epsilon}}
	\end{equation}
	for $X \geq 1$. The exponents $\degB m$ and $\twerror m$ are given in Table \ref{table:degB and error} below. The implicit constants depend on $m$ and $\epsilon$.
\end{theorem}

It turns out that $\twconsteq m = \twconstad m$ when $m \in \set{5, 6, 7, 9, 13, 18, 25}$, and that $\twconsteq m = 2 \twconstad m$ when $m \in \set{4, 8, 12, 16}$ (see \Cref{Corollary: Elliptic curves for which twNeq(X) = twNad(X)} and \Cref{Lemma: Difference between twNeq and twNad}).

\jvtable{table:degB and error}{
\rowcolors{2}{white}{gray!10} \begin{tabular}{c | c c c c c c c c c c c c} 
$m$ & $1$ & $4$ & $6$ & $7$ & $8$ & $9$ & $10$ & $12$ & $13$ & $16$ & $18$ & $25$ \\
\hline\hline
$\degB m$ & $6/5$ & $3$ & $6$ & $6$ & $6$ & $6$ & $12$ & $12$ & $12$ & $12$ & $18$ & $18$ \\
$\twerror m$ & $2$ & $6$ & $12$ & $12$ & $12$ & $12$ & $18$ & $24$ & $51/4$ & $24$ & $36$ & $30$ \\
\end{tabular}
}{Exponents for asymptotics of $\twistNeq m (X)$ and $\twistNad m (X)$}

\Cref{Intro Theorem: asymptotic for twN(X) for m of genus $0$} summarizes the asymptotics given by \eqref{Equation: Count of elliptic curves up to twist}, \Cref{Theorem: asymptotic for twN(X) for m = 7}, \Cref{Theorem: asymptotic for twN(X) for m = 10 and 25}, \Cref{Theorem: asymptotic for twN(X) for m = 13}, and \Cref{Theorem: asymptotic for twN(X) for m of genus $0$}. We have coarsened the error terms of these theorems slightly for uniformity and clarity of exposition.

We do not have asymptotics for $\twistNeq 2 (X)$, $\twistNad 2 (X)$, $\twistNeq 3 (X)$, $\twistNad 3 (X)$, $\twistNeq 5 (X)$, or $\twistNad 5 (X)$ (see \Cref{Remark: Why not m = 5}). To our knowledge, no prior work has been done on counting elliptic curves with a cyclic $m$-isogeny over global fields up to quadratic twist.

\section{Our approach}\label{Section: Our approach}

In this section, we sketch our methodology for proving \Cref{Intro Theorem: asymptotic for NQ(X) for 5 < m <= 9}, \Cref{Intro Theorem: asymptotic for NQ(X) for m > 9}, \Cref{Intro Theorem: asymptotic for NQ(X) for m of nonzero genus}, and \Cref{Intro Theorem: asymptotic for twN(X) for m of genus $0$}. We revisit and expand on this sketch in \cref{Section: Our approach revisited} below.

Choose $m$ so that $X_0(m)$ is of genus $0$. Our approach to proving \Cref{Intro Theorem: asymptotic for NQ(X) for 5 < m <= 9}, \Cref{Intro Theorem: asymptotic for NQ(X) for m > 9}, and \Cref{Intro Theorem: asymptotic for twN(X) for m of genus $0$} proceeds through five main steps.
\begin{enumerate}[label=(\arabic*)]
    \item We employ the modular curve $X_0(m)$ to establish a parameterization for the family of elliptic curves equipped with a cyclic $m$-isogeny. We obtain two polynomials $\f m (t)$ and $\g m (t)$ such that up to twist equivalence, each elliptic curve with a cyclic $m$-isogeny may be written in the form
    \begin{equation}\label{Equation: parameterizing E with a cyclic m isogeny}
    E : y^2 = x^3 + \f m (t) x + \g m (t)
    \end{equation}
    for some $t \in \bbQ$ (\Cref{Lemma: parameterizing m-isogenies}).
    \item Writing $t = a/b$ and homogenizing \eqref{Equation: parameterizing E with a cyclic m isogeny}, the elliptic curves equipped with a cyclic $m$-isogeny up to twist equivalence are parameterized by coprime pairs $(a, b) \in \bbZ^2$. We apply the Principle of Lipschitz (\Cref{Theorem: Principle of Lipschitz}), which asserts that the number of lattice points within a region can be approximated by the area of that region, to derive a coarse estimate for the count of Weierstrass equations of the form $E : y^2 = x^3 + A x + B$ that arise from \eqref{Equation: parameterizing E with a cyclic m isogeny} in this fashion (\Cref{Corollary: Estimates for lattice counts}). 
    \item The discrepancy between the coefficients of the model \eqref{Equation: parameterizing E with a cyclic m isogeny} and the twist height of $E$ may be quite large. However, the set of pairs $(a, b)$ for which this ``twist minimality defect'' (see \eqref{Equation: twist defect powers}) is divisible by a given integer $e$ can be expressed as a finite union of sublattices of $\bbZ^2$, and therefore can be estimated using step 2. In addition, a single elliptic curve $E$ may occur more than once in the estimates given above, since $da / db = a/b$. Recall that the M\"obius sieve is a method that applies the inclusion-exclusion principle to the prime factorizations of integers. Using two M\"obius sieves (see \Cref{Lemma: fundamental sieve} and for example \Cref{Lemma: asymptotic for M(X; e) for m = 7}), one for each of these issues, we write $\twistNeq m (X)$ in terms of the estimates obtained in the previous step. This gives us \Cref{Intro Theorem: asymptotic for twN(X) for m of genus $0$} for $\twistNeq m (X)$.
    \item We bound the difference between $\twistNeq m (X)$ and $\twistNad m (X)$ or $\twistNeq m (X)$ and $2 \twistNad m (X)$ (\Cref{Corollary: Elliptic curves for which twNeq(X) = twNad(X)} and \Cref{Lemma: Difference between twNeq and twNad}). For each proper divisor $n$ of $m$, there is a modular curve parameterizing elliptic curves equipped with a pair of cyclic $m$-isogenies whose kernels have intersection of order $n$. If $m \in \set{4, 8, 12, 16}$ and $n = m/2$, this modular curve is $X_0(m)$ itself: this is why $\twconsteq m = 2 \twconstad m$ for these $m$. Otherwise, when the modular curves indexed by $m$ and $n$ are of genus $0$, we emulate our first step (see also \cite[Theorem 3.3.1]{Cullinan-Kenney-Voight}) to bound this contribution to $\twistNeq m (X) - \twistNad m (X)$ or to $\twistNeq m (X) - 2 \twistNad m (X)$, depending on whether $4 \mid m$. When the modular curves are of genus greater than 1, Faltings's theorem \cite{Faltings1, Faltings2} assures us that we get a contribution of at most $O(1)$ to $\twistNeq m (X) - \twistNad m (X)$ or $\twistNeq m (X) - 2 \twistNad m (X)$. When the modular curves are of genus 1, the contribution is still $O(1)$ by inspection. We obtain \Cref{Intro Theorem: asymptotic for twN(X) for m of genus $0$} in its entirety.
    \item Finally, we wish to estimate $\NQeq m (X)$ and $\NQad m (X)$ using our estimates for $\twistNeq m (X)$ and $\twistNad m (X)$. We first use \Cref{Intro Theorem: asymptotic for twN(X) for m of genus $0$} to establish a half-plane of convergence for the height zeta functions $\twistLeq m (s)$ and $\twistLad m (s)$ associated to $\twistNeq m (X)$ and $\twistNad m (X)$ (see for example \Cref{Corollary: twL(s) has a meromorphic continuation for m = 7}). The height zeta functions $\twistLeq m (s)$ and $\twistLad m (s)$ are closely related to the height zeta functions $\LQeq m (s)$ and $\LQad m (s)$ of $\NQeq m (X)$ and $\NQad m (X)$ (\Cref{Theorem: relationship between twistL(s) and L(s)}). 
    
    Recall that a Tauberian theorem establishes a connection between the asymptotics of a sequence and the analytic behavior of an associated function. We use Landau's Tauberian theorem (\Cref{Theorem: Landau's Tauberian theorem}) as presented by Roux \cite{Roux} to express the asymptotics of $\NQeq m (X)$ and $\NQad m (X)$ in terms of residues of $X^s \LQeq m (s)/s$ and $X^s \LQad m (s)/s$. This gives us \Cref{Intro Theorem: asymptotic for NQ(X) for 5 < m <= 9} and \Cref{Intro Theorem: asymptotic for NQ(X) for m > 9}.
\end{enumerate} 

Our key innovations occur in step 3, where we address the discrepancy between the size of a model and the twist of the associated elliptic curve, and in step 5, where we use our estimates for $\twNeq m (X)$ and $\twNad m (X)$ to obtain estimates for $\NQeq m (X)$ and $\NQad m (X)$.

\begin{remark}
    In conjunction with earlier work, we obtain asymptotics for $\NQeq m (X)$ and $\NQad m (X)$ for all $m \neq 5$. In the case $m = 5$, steps 1 and 2 of our approach go through without obstruction, and we are able to set up the sieves in step 3 as well. However, in this case our sieve tells us only that
    \begin{equation}
        \twNeq 5 (X), \twNad 5 (X) = O(X^{1/6} \log X),
    \end{equation}
    rather than giving us a power-saving asymptotic for $m = 5$. On a technical level, this occurs because the polynomials appearing in \eqref{Equation: parameterizing E with a cyclic m isogeny} are of too low degree. If we were able to obtain an asymptotic for $\twNeq 5 (X)$ with power-saving error, we could follow steps 4 and 5 of our approach without further obstruction. See \Cref{Remark: Why not m = 5} for more details.
\end{remark}

\begin{remark}
    The techniques of this thesis may be adapted to count the fibers of an elliptic surface
    \begin{equation}
    E(t) : y^2 = x^3 + f(t) x + g(t)
    \end{equation}
    over $\bbP^1$ according to their twist heights, provided that the geometry of this elliptic surface is sufficiently similar to the elliptic surfaces we obtain in \cref{Section: Parameterizing elliptic curves with a cyclic m-isogeny}.
\end{remark}

If the compactified modular curve $X_0(m)$ has nonzero genus, and the noncompactified modular curve $Y_0(m) \subseteq X_0(m)$ has $Y_0(m)(\bbQ) \neq \emptyset$, we can use the classical enumeration of the rational points on $X_0(m)$ (see Table \ref{table:nonzerogenus} or \cite{Mazur1}) in conjunction with Walfisz's count of squarefree integers (\Cref{Theorem: Estimating squarefree(X)}) and Liu's contingent refinement (\Cref{Theorem: Estimating squarefree(X) given RH}) to deduce \Cref{Intro Theorem: asymptotic for NQ(X) for m of nonzero genus}.

\section{Contents}\label{Section: Contents}

In this section, we summarize the remainder of the thesis.

In \cref{Chapter: Background}, we review the fundamentals of elliptic curves and relay the history and context for our problem. In \cref{Chapter: Preliminaries}, we gather several results from geometry and analysis for later use. We close the chapter by expanding on \cref{Section: Our approach} to give a more detailed outline of the argument we will use to prove our main theorems.
In \cref{Chapter: m = 7}, we apply the material from \cref{Chapter: Preliminaries} to prove \Cref{Intro Theorem: asymptotic for NQ(X) for 5 < m <= 9} and \Cref{Intro Theorem: asymptotic for twN(X) for m of genus $0$} when $m = 7$, establishing the asymptotics for $\twistNeq 7 (X) = \twistNad 7 (X)$ and $\NQeq 7 (X) = \NQad 7 (X)$. In \cref{Chapter: m = 10 and 25}, we adapt the methods of \cref{Chapter: m = 7} to prove \Cref{Intro Theorem: asymptotic for NQ(X) for m > 9} and \Cref{Intro Theorem: asymptotic for twN(X) for m of genus $0$} when $m = 10, 25$, establishing the asymptotics for these $\twistNeq m (X) = \twistNad m (X)$ and $\NQeq m (X) = \NQad m (X)$.
In \cref{Chapter: m = 13}, we elaborate on the methods of \cref{Chapter: m = 7} and \cref{Chapter: m = 10 and 25} to prove \Cref{Intro Theorem: asymptotic for NQ(X) for m > 9} and \Cref{Intro Theorem: asymptotic for twN(X) for m of genus $0$} when $m = 13$, establishing the asymptotics for $\twistNeq {13} (X) = \twistNad {13} (X)$ and $\NQeq {13} (X) = \NQad {13} (X)$.  
In \cref{Chapter: other m of genus $0$}, we prove \Cref{Intro Theorem: asymptotic for NQ(X) for 5 < m <= 9}, \Cref{Intro Theorem: asymptotic for NQ(X) for m > 9}, and \Cref{Intro Theorem: asymptotic for twN(X) for m of genus $0$} for all other $m > 5$ for which $X_0(m)$ has genus $0$, establishing the asymptotics for these $\twistNeq m (X)$, $\twistNad m (X)$, $\NQeq m (X)$, and $\NQad m (X)$, as well as for $\twistNeq 4(X)$ and $\twistNad 4 (X)$.
In \cref{Chapter: m of nonzero genus}, we prove \Cref{Intro Theorem: asymptotic for NQ(X) for m of nonzero genus}, establishing asymptotics for $\NQeq m (X) = \NQad m (X)$ when $X_0(m)$ is of nonzero genus.

%% file: Ch2_Background.tex
\chapter{Historical background}\label{Chapter: Background}

In \cref{Section: A primer on elliptic curves}, we review the basic theory of elliptic curves and their isogenies.
In \cref{Section: History}, we survey the history of the arithmetic statistics of elliptic curves, with a slant towards our main problem. %
In \cref{Section: Previous results}, we recall several results from the literature which, taken in conjunction with the results of \cref{Section: Results}, yield an almost complete picture of the asymptotics for counts of elliptic curves over $\bbQ$ with a cyclic $m$-isogeny.

This chapter may be skimmed or skipped in its entirety by readers familiar with both elliptic curves and the history of our problem.

\section{A primer on elliptic curves}\label{Section: A primer on elliptic curves}

In this section, we briefly review and motivate the study of elliptic curves, and highlight several important results from the discipline. We do not aim for generality in this section, and freely assume that our elliptic curves are defined over a field of characteristic $0$ or even over $\bbQ$ whenever convenient. The curious reader can learn more from any of a number of standard references \cite{Koblitz, Silverman2, Silverman, Silverman-Tate} (see also \cite[Chapter IV.4]{Hartshorne}).

An \defi{elliptic curve} is a pair $(E, O)$, where $E$ is a nonsingular projective curve of genus 1, and $O \in E$. We typically write $E$ for the elliptic curve $(E, O)$, and elide the description of this distinguished point. If $K$ is a field, the elliptic curve $E$ is \defi{defined over} $K$, written $E/K$, if $E$ is defined over $K$ as a curve and $O \in E(K)$. In this thesis, we shall be interested in elliptic curves over the rational numbers $\bbQ$, so we typically take $K = \bbQ$ to be the field of rational numbers.

Elliptic curves are natural objects of study for many reasons. We record several of them now. 

The genus $g \in \bbZ_{\geq 0}$ of a curve provides a coarse measure of its geometric and arithmetic complexity, and it is well-known that a curve has genus $0$ if and only if it is birationally equivalent to the projective line $\bbP^1$ \cite[Example IV.1.3.5]{Hartshorne}. We therefore have a thorough geometric understanding of genus $0$ curves, and genus 1 curves like elliptic curves are the natural next case to study.

Arithmetic geometry concerns itself in large part with determining the rational points of a variety $V/\bbQ$. If a curve of genus $g = 0$ has any rational points, then it has infinitely many rational points. On the other hand, Faltings proved the following remarkable theorem for curves of genus greater than 1.

\begin{theorem}[Faltings's Theorem]\label{Theorem: Faltings Theorem}
	Let $C/\bbQ$ be a nonsingular algebraic curve of genus $g$. If $g > 1$, then $C$ has at most finitely many rational points.
\end{theorem}

\begin{proof}
	Faltings \cite{Faltings1, Faltings2}.
\end{proof}

Algebraic curves of genus 1 lie on the boundary of these two cases: they can have no rational points, finitely many rational points, or infinitely many rational points.

There are ample other reasons to be fascinated with elliptic curves. For instance, the modularity theorem \cite[Preface]{Diamond-Shurman} asserts that every elliptic curve arises from a modular form. An elliptic curve $E$ also induces an adelic Galois representation via its Tate modules \cite[Section III.7]{Silverman}, and gives rise to an $L$-function via its $\bbF_p$-points \cite[Appendix C.16]{Silverman}. These four objects--elliptic curves, modular forms, Galois representations, and $L$-functions--are all interrelated via the Langlands program, which remains a lively area of research \cite{Caraiani-Emerton-Gee, Clozel, Gelbart, Glazunov, Langlands}. The Birch and Swinnerton-Dyer Conjecture, one of the seven famed millennium problems, relates the $L$-function of an elliptic curve to its set of rational points \cite{Claymath}.

Although elliptic curves are rather abstract objects as we have defined them, every elliptic curve may be concretely realized as the zero set of a Weierstrass equation, as the following theorem shows.

\begin{theorem}[{\cite[Proposition 3.1]{Silverman}}]\label{Theorem: Weierstrass equation for an elliptic curve}
	Let $E$ be an elliptic curve defined over a field $K$ of characteristic $0$.
	\begin{enumalph}
		\item There exist functions $x, y \in K(E)$ such that the map
		\begin{equation} 	\begin{aligned}
		\phi &: E \to \bbP^2 \\
		\phi &= (x : y : 1)
		\end{aligned}\end{equation} 
		gives a $K$-isomorphism of $E/K$ onto a curve given by a (simplified) Weierstrass equation
		\begin{equation}\label{Equation: Simplified Weierstrass equation}
		E : y^2 = x^3 + A x + B
		\end{equation}
		with coefficients $A, B \in K$, and such that $\phi(O) = (0 : 1 : 0) \in \bbP^2$.
		\item Any two Weierstrass equations as in (a) are related by a linear change of variables of the form
		\begin{equation}
		(x, y) \mapsto (u^2 x, u^3 y)
		\end{equation}
		with $u \in K^\times$. This change of variables yields the model
		\begin{equation}
		E^{(u^2)} : y^2 = x^3 + u^4 A x + u^6 B
		\end{equation}
		\item Conversely, every nonsingular cubic curve $C$ given by a Weierstrass equation of the form \eqref{Equation: Simplified Weierstrass equation} is an elliptic curve defined over $K$ with distinguished point $\infty = (0 : 1 : 0)$. Nonsingularity is equivalent to the discriminant
		\begin{equation}\label{Equation: Discriminant of Weierstrass equation}
		\Delta(E) = \Delta(A, B) = -16(4A^3 + 27B^2)
		\end{equation}
		 of the model being nonzero.
	\end{enumalph}
\end{theorem}

\begin{remark}
	\cite[Proposition 3.1]{Silverman} only affirms a $K$-isomorphism from $E$ to a curve $C$ in $\bbP^2$ of the form
	\begin{equation}\label{Equation: Unsimplified Weierstrass equation}
	y^2 + a_1 xy + a_3 y = x^3 + a_2 x^2 + a_4 x + a_6;
	\end{equation}
	however, because we have assumed $K$ is not of characteristic 2 or 3, a linear change of variables reduces \eqref{Equation: Unsimplified Weierstrass equation} to \eqref{Equation: Simplified Weierstrass equation}.
\end{remark}

\Cref{Theorem: Weierstrass equation for an elliptic curve}(a) is a consequence of the Riemann--Roch theorem; like many results about elliptic curves, it is proven using tools from algebraic geometry. Here and later in the thesis, we will have no direct need of divisors, Jacobians, the Riemann--Roch theorem, or anything of the sort: while we will freely reference results proven using such techniques, we elide the techniques themselves as a technical distraction.

If $E$ has \eqref{Equation: Simplified Weierstrass equation} as its Weierstrass equation, we define the $j$-invariant $j(E)$ of $E$ as follows:
\begin{equation}\label{Equation: j-invariant for E}
j(E) = j(A, B) \colonequals -1728 \frac{(4 A)^3}{\Delta(A, B)} = \frac{2^8 \cdot 3^3 \cdot A^3}{4 A^3+27 B^2}.
\end{equation}
Two elliptic curves have the same $j$-invariant if and only if they are $\bbQalg$-isomorphic \cite[Proposition 1.4(b)]{Silverman}.

\Cref{Theorem: Weierstrass equation for an elliptic curve} gives us an infinite family of models for every elliptic curve, but each such elliptic curve has a canonical \defi{minimal} model over $\bbQ$. Indeed, as noted in \cref{Section: Results}, every elliptic curve $E/ \bbQ$ has a unique Weierstrass model of the form 
\begin{equation}
E \colon y^2 = x^3 + Ax + B,
\end{equation}
where $A, B \in \bbZ$, $\Delta(A, B) \neq 0$, and for every prime $\ell$ we have $\ell^4 \nmid A$ or $\ell^6 \nmid B$. Recall that $\scrE$ is the set of Weierstrass models of this form. With this notation in mind, \Cref{Theorem: Weierstrass equation for an elliptic curve} has the following corollary.

\begin{corollary}\label{Corollary: Isomorphisms are quadratic twists for j not 0 and 1728}
    Let $E, E\prm \in \scrE$, and suppose $E$ is $\bbQalg$-isomorphic to $E\prm$. The following statements hold.
    \begin{itemize}
        \item If $j(E) \neq 0, 1728$, then $E\prm$ is a quadratic twist of $E$.
        \item If $j(E) = 0$, then $E\prm$ is a sextic twist of $E$.
        \item If $j(E) = 1728$, then $E\prm$ is a quartic twist of $E$.
    \end{itemize}
\end{corollary}

\begin{proof}
    Write
    \begin{equation}
    E : y^2 = x^3 + A x + B.
    \end{equation}
    By \Cref{Theorem: Weierstrass equation for an elliptic curve}(b), the $\bbQalg$-isomorphism of $E$ with $E\prm$ implies
    \begin{equation}
    E\prm = E^{(u^2)} : y^2 = x^3 + u^4 A x + u^6 B
    \end{equation}
    for some $u \in (\bbQalg)^\times$. On the other hand, as $E\prm \in \scrE$, we see $u^4 A, u^6 B \in \bbQ$. 
    
    If $j(E) \neq 0, 1728$ then $A, B \neq 0$, so comparing these equations we see $u^4, u^6 \in \bbQ^\times$. Thus $u^2 = u^6/u^4 \in \bbQ^\times$, and $E\prm$ is a quadratic twist of $E$ as desired.

    If $j(E) = 0$, then $u^6 \in \bbQ^\times$, so $E\prm$ is a sextic twist of $E$. If $j(E) = 1728$, then $u^4 \in \bbQ^\times$, so $E\prm$ is a quartic twist of $E$.
\end{proof}

\subsection*{The group law for an elliptic curve}\label{Subsection: Group law}

In this subsection, we give a concrete description of the group law for an elliptic curve, and recall Mordell's Theorem (\Cref{Theorem: Mordell's Theorem}) and Mazur's theorem on torsion (\Cref{Theorem: Mazur's Theorem}). The material outlined here is classical.

Let $E/\bbQ$ be an elliptic curve with Weierstrass model \eqref{Equation: Simplified Weierstrass equation}, and let $P, Q \in E(\bbQ)$. Let $L$ denote the line in $\bbP^2$ passing through both $P$ and $Q$ (if $P = Q$, we take $L$ to be tangent to $P$), and let $R$ be the third point of intersection of $L$ with $E$. Let $L\prm$ be the line passing through $R$ and $\infty$, and let $R\prm$ be the third point of intersection of $L\prm$ with $E$. We define $P + Q \colonequals R\prm$. For $P = (x_0, y_0) \in E$, we define $- P = (x_0, -y_0)$; in particular, the unique line $L$ passing through both $P$ and $-P$ also passes through $\infty$. By definition, we have $P + Q = -R$.

\begin{example}
    Let $E : y^2 = x^3 - 3 x + 62$, let $P = (-1, 8)$, and let $Q = (2, -8)$. The line $L$ is given in blue in Figure \ref{Figure: grouplaw}, and passes through the points $P$, $Q$, and $R = (247/9, -3880/27)$. The line $L\prm$ is given in green in Figure \ref{Figure: grouplaw}, and passes through the points $\infty$, $R$, and $R\prm = (247/9, 3880/27)$. Thus
 \begin{equation}
    (-1, 8) + (2, -8) = (247/9, 3880/27)
 \end{equation}
 in $E$.
\end{example}

\jvfigure{Figure: grouplaw}{\includegraphics[scale=0.3]{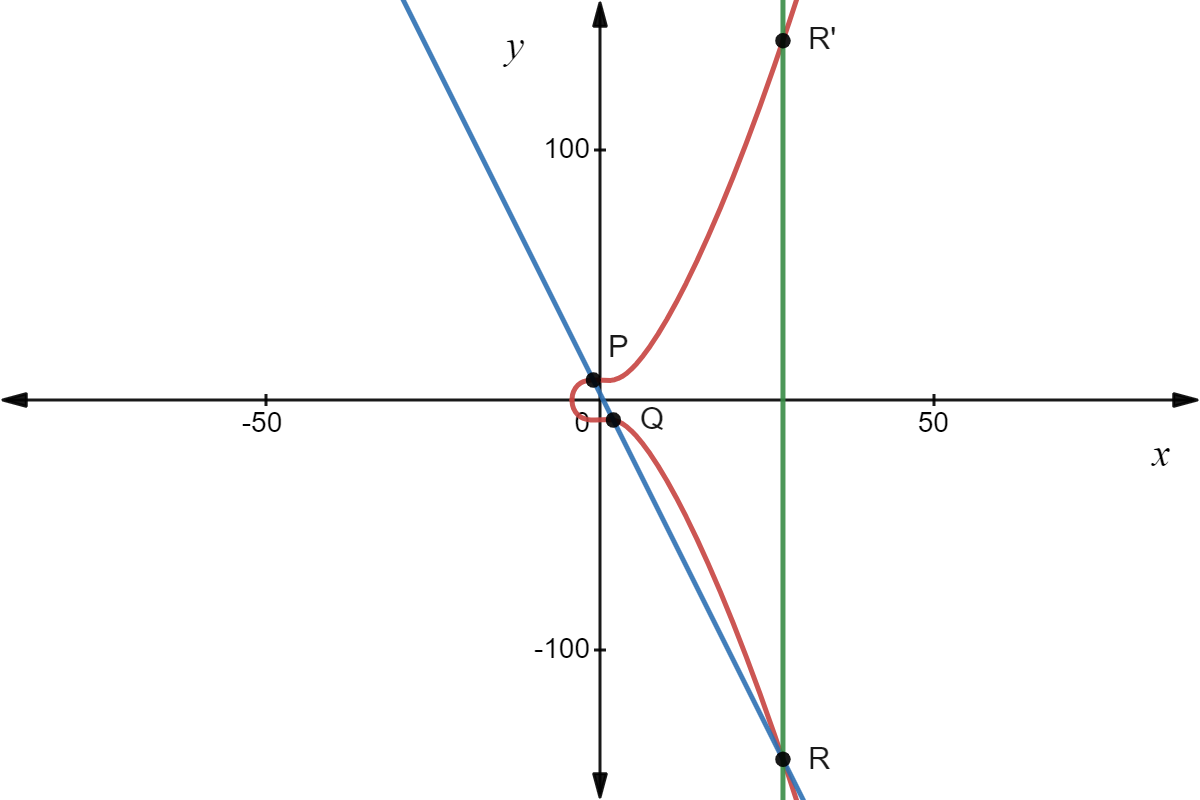}}{The elliptic curve $E : y^2 = x^3 - 3 x + 62$, with some additions}

As our notation suggests, we have the following theorem.

\begin{theorem}[{\cite[Proposition 2.2]{Silverman}}]\label{Theorem: group law}
	Let $E$ be an elliptic curve with Weierstrass equation \eqref{Equation: Simplified Weierstrass equation}. The elliptic curve $E$ is an abelian group under the binary operation $+$ defined above, with identity element $\infty$ and with inversion given by $P \mapsto - P$.
\end{theorem}

\begin{remark}
	Every elliptic curve $E$ is isomorphic to its own Jacobian as an algebraic curve. The Jacobian of an elliptic curve is an abelian group by definition. Although we have given the group law of the elliptic curve geometrically, this group law is also induced on $E$ via its isomorphism with its Jacobian.  
 
    We can derive an explicit formula for calculating the sum of points on the elliptic curve $E$ by examining its Weierstrass equation \cite[Group Law Algorithm 2.3]{Silverman}.
\end{remark}

Note that if $K/\bbQ$ is a field extension of $\bbQ$ and $P, Q \in E(K)$, then $P + Q \in E(K)$ as well. Thus, $E$ defines not a single group, but a family of groups, one for every field extension of $\bbQ$. This is the essential reason that $E$ defines a group scheme. However, in this thesis we will never have occasion to consider groups other than $E(\bbQ)$ and $E(\bbQalg)$, so this rich structure will be largely invisible to us.

\Cref{Theorem: group law} tells us that $E(\bbQ)$ is an abelian group, so it is natural to inquire about and classify its structure. This classification was done almost completely in the following two theorems, due to Mordell and Mazur.

\begin{theorem}[Mordell's Theorem]\label{Theorem: Mordell's Theorem}
	Let $E/\bbQ$ be an elliptic curve. Then the abelian group $E(\bbQ)$ is finitely generated.
\end{theorem}

\begin{proof}
	Mordell \cite{Mordell}. Weil extended Mordell's theorem to number fields \cite{Weil}.
\end{proof}

Due to \Cref{Theorem: Mordell's Theorem} and Weil's extension thereof to number fields \cite{Weil}, we now call $E(\bbQ)$ the \defi{Mordell--Weil group} of $E$. \Cref{Theorem: Mordell's Theorem} shows that we may write 
\begin{equation}
E(\bbQ) \simeq \bbZ^r \times E(\bbQ)_{\textup{tors}},
\end{equation}
where $r \geq 0$ is an integer and $E(\bbQ)_{\textup{tors}}$ is a finite abelian group. The quest to understand the ranks of elliptic curves over $\bbQ$ is over a hundred years old \cite[page 173]{Poincare} and remains an active area of study to this day (see \cite{Claymath, Bhargava-Shankar, Park-Poonen-Voight-Wood}); by contrast, Mazur completely determined the possible torsion of an elliptic curve over $\bbQ$.

\begin{theorem}[Mazur's Theorem on torsion]\label{Theorem: Mazur's Theorem}
	Let $E/\bbQ$ be an elliptic curve. The group $E(\bbQ)_{\textup{tors}}$ is isomorphic to one of the following fifteen groups:
        \begin{equation}
	\begin{aligned}
	\bbZ / N \bbZ \ &\text{with} \ 1 \leq N \leq 10 \ \text{or} \ N = 12, \\
	\bbZ / 2 \bbZ \times \bbZ / 2 N \bbZ \ &\text{with} \ 1 \leq N \leq 4.
	\end{aligned}\end{equation} 
    Moreover, each of these fifteen groups occurs as the torsion group of infinitely many elliptic curves over $\bbQ$.
\end{theorem}

\begin{proof}
    Mazur \cite{Mazur1}.
\end{proof}

\subsection*{Isogenies}\label{Subsection: isogenies}

In this subsection, we define isogenies between elliptic curves and outline several of their remarkable properties. We finish the section by recalling Mazur's theorem on isogenies (\Cref{Theorem: Kenku-Mazer Theorem}). The material outlined here is classical.

If we agree that elliptic curves are interesting objects, it is natural to ask what constitutes a map between elliptic curves. If $E_1, E_2 \in \scrE$ are elliptic curves, such a map should at least give a morphism of curves $E_1 \to E_2$, and should also preserve the distinguished point $\infty$, so $\infty \mapsto \infty$. In the previous subsection, we also showed that the points of an elliptic curve have a natural group structure. Should we demand that our map of elliptic curves preserve this group structure? As it turns out, we get this property for free.

\begin{theorem}[{\cite[Theorem 4.8]{Silverman}}]\label{Theorem: isogenies preserve groups}
	Let $E_1, E_2 \in \scrE$, and let $\phi : E_1 \to E_2$ be a morphism of curves (over $\bbQ$) with $\phi(\infty) = \infty$. Then for all points $P, Q$ of $E_1$, we have
	\begin{equation}
	\phi(P + Q) = \phi(P) + \phi(Q).
	\end{equation}
\end{theorem}

We refer to a morphism of curves $\phi : E_1 \to E_2$ for which $\phi(\infty) = \infty$ as an \defi{isogeny}. Unless otherwise specified, we restrict our attention to nonconstant isogenies, which are necessarily surjective as maps over $\bbQalg$. \Cref{Theorem: Weierstrass equation for an elliptic curve} and the definition of $\scrE$ already give us the following result.

\begin{proposition}\label{Proposition: only isomorphisms are signs}
	Let $E_1, E_2 \in \scrE$ be elliptic curves over $\bbQ$, and let $\phi : E_1 \to E_2$ be an isomorphism of elliptic curves over $\bbQ$. Then $E_1 = E_2$, and $\phi$ is either the identity map $P \mapsto P$ or the negative map $P \mapsto - P$.
\end{proposition}

\begin{proof}
	By definition, $\scrE$ only contains one elliptic curve from each $\bbQ$-isomorphism class, so $E_1 = E_2$. Now by \Cref{Theorem: Weierstrass equation for an elliptic curve}(b), $\phi$ is of the form 
	\begin{equation}
	\phi : (x, y) \mapsto (u^2 x, u^3 y)
	\end{equation}
	for some $u \in \bbQ^\times$, and we have $A = u^4 A$ and $B = u^6 B$. As a consequence, $u$ is a fourth root of unity or a sixth root of unity, but the only roots of unity in $\bbQ$ are $\pm 1$, so $u = \pm 1$. The claim follows.
\end{proof}

To avoid excessive clutter in our writing, we assume all isogenies given are isogenies over $\bbQ$ unless otherwise specified. 

Somewhat remarkably, isogenies induce an equivalence relation on elliptic curves $E \in \scrE$. More precisely, we have the following proposition.

\begin{proposition}\label{Proposition: dual isogenies}
	Let $E_1, E_2$ be elliptic curves, and let $\phi : E_1 \to E_2$ be a nonconstant degree $m$ isogeny. There is a unique isogeny $\widehat{\phi} : E_2 \to E_1$ with 
		\begin{equation}
		\widehat{\phi} \circ \phi : P \mapsto m P \ \text{and} \ \phi \circ \widehat{\phi} : P \mapsto m P.
		\end{equation}
\end{proposition}

\begin{proof}
	Silverman \cite[Theorem 6.1]{Silverman}.
\end{proof}

We say $E, E\prm \in \scrE$ are \defi{isogenous} if there is a nonconstant isogeny $\phi : E \to E\prm$. By \Cref{Proposition: dual isogenies}, this is an equivalence relation, which is coarser than $\bbQ$-isomorphism.

If $\phi, \widehat{\phi}$ are as in \Cref{Proposition: dual isogenies}, we refer to $\widehat{\phi}$ as the \defi{dual isogeny} of $\phi$. We define the dual of the constant isogeny $0 : E_1 \to E_2$ to be the constant isogeny $0 : E_2 \to E_1$.

A morphism of curves is of \defi{degree $m$} if it is generically $m$-to-one. If $\phi : E_1 \to E_2$ is a degree $m$ isogeny, we say $\phi$ is an \defi{$m$-isogeny}. Note that if $\phi : E_1 \to E_2$ is an $m$-isogeny and $\psi : E_2 \to E_3$ is an $n$-isogeny, then $\psi \circ \phi$ is an $mn$-isogeny. We now report some charming additional properties of dual isogenies.

\begin{theorem}\label{Theorem: properties of the dual isogeny}
	Let $E_1, E_2 \in \scrE$, and let $\phi : E_1 \to E_2$ be a (nonconstant) isogeny. The following statements hold.
	\begin{enumalph}
		\item If $\phi$ is an $m$-isogeny, then $\widehat{\phi}$ is an $m$-isogeny.
		\item Let $\psi : E_2 \to E_3$ be another isogeny. We have 
		\begin{equation}
		\widehat{\psi \circ \phi} = \widehat{\phi} \circ \widehat{\psi}.
		\end{equation}
		\item Let $\phi\prm : E_1 \to E_2$ be another isogeny. We have
		\begin{equation}
		\widehat{\phi + \phi\prm} = \widehat{\phi} + \widehat{\phi\prm}.
		\end{equation}
		\item We have $\widehat{\widehat{\phi}} = \phi$.
		\item If $[m]$ is defined by
		\begin{equation}
		[m] : P \mapsto m P,
		\end{equation}
		then $\deg [m] = m^2$, and $\widehat{[m]}$ is also defined by
		\begin{equation}
		\widehat{[m]} : P \mapsto m P.
		\end{equation}
	\end{enumalph}
\end{theorem}

\begin{proof}
	Silverman \cite[Theorem 6.2]{Silverman} proves a more general result.
\end{proof}

If $\phi : E_1 \to E_2$ is an isogeny, we write $\ker \phi \subseteq E_1(\bbQalg)$ for the kernel of $\phi$. It turns out that $\ker \phi$ determines almost the entire behavior of $\phi$.

\begin{theorem}\label{Theorem: Isogenies factor}
	Let $E_1, E_2 \in \scrE$, and let $\phi : E_1 \to E_2$ be an isogeny. The following statements hold.
	\begin{enumalph}
		\item If $\phi$ is an $m$-isogeny, then $\# \ker \phi = m$.  
		\item Suppose $E_2\prm \in \scrE$ and $\phi\prm : E_1 \to E_2\prm$ is another isogeny, and suppose that $\ker \phi \subseteq \ker \phi\prm$. There is a unique isogeny $\psi : E_2 \to E_2\prm$ such that $\psi \circ \phi = \phi\prm$.
	\end{enumalph}
\end{theorem}

\begin{proof}
	Silverman \cite[Theorem 4.10, Corollary 4.11]{Silverman} proves a more general result.
\end{proof}

We define an \defi{unsigned isogeny} $\phi : E_1 \to E_2$ to be an isogeny up to postcomposition by $\pm 1$. By \Cref{Proposition: only isomorphisms are signs}, an isogeny over $\bbQ$ up to postcomposition by the $\bbQ$-automorphisms of $E_2$ is the same as an unsigned isogeny over $\bbQ$.

\begin{theorem}\label{Theorem: every kernel gives rise to an essentially unique isogeny}
	Let $E \in \scrE$, and let $\Phi \subseteq E(\bbQalg)$ be a finite group which is stabilized by the absolute Galois group $\Gal(\bbQalg\,|\,\bbQ)$ of $\bbQ$. Then there is a unique elliptic curve $E\prm \in \scrE$, and a unique unsigned isogeny $\phi : E \to E\prm$ with $\ker \phi = \Phi$.
\end{theorem}

\begin{proof}
	Silverman \cite[Proposition 4.12]{Silverman} (see also V\'{e}lu \cite{Velu} for a more direct construction).
\end{proof}

As a consequence of \Cref{Theorem: every kernel gives rise to an essentially unique isogeny}, for each elliptic curve $E \in \scrE$, there is a bijection between (nonconstant) unsigned isogenies with domain $E$ and finite subgroups of $E(\bbQalg)$ which are stabilized by the absolute Galois group $\Gal(\bbQalg\,|\,\bbQ)$, which we henceforth abbreviate as $\Gal_\bbQ$.

\begin{corollary}\label{Corollary: isogenies are stablized by quadratic twist}
    Let $E \in \scrE$ be an elliptic curve, and let $E\prm$ be a quadratic twist of $E$. If $\Phi \subseteq E(\bbQalg)$ is a finite abelian group stabilized by $\Gal_{\bbQ}$, then the image $\Phi\prm$ of $\Phi$ in $E\prm$ under that quadratic twist is also an abelian group stabilized by $\Gal_\bbQ$. Moreover, $\Phi \simeq \Phi\prm$ as abelian groups.
\end{corollary}

\begin{proof}
    Let $\Phi \subseteq E(\bbQalg)$ be a finite subgroup, and let
    \begin{equation}
    f(t) = \prod_{\substack{\pm P \in \Phi \\ P \neq \infty}} (t - x(P)) \in \bbQalg[t]
    \end{equation}
    be the product of all linear terms $t - x(P)$ as $x$ varies over the $x$-coordinates of the affine points in $\Phi$. We see that $f(t) \in \bbQ[t]$ if and only if $\Phi$ is stable under $\Gal_\bbQ$.

    Now write $E\prm = E^{(c)}$ for appropriately chosen $c \in \bbQ^\times$, and let $\iota^{(c)} : E \to E\prm$ be the $\bbQalg$-isomorphism $\iota^{(c)} : (x, y) \mapsto (c^2 x, c^3 y)$. By definition, we have $\Phi\prm = \iota^{(c)}(\Phi)$, so $\Phi\prm \simeq \Phi$ as abelian groups. On the other hand, the polynomial
    \begin{equation}
    f^{(c)}(t) = \prod_{\substack{\pm P \in \Phi\prm \\ P \neq \infty}} (t - x(P)) = \prod_{\substack{\pm P \in \Phi \\ P \neq \infty}} (t - c^2 x(P)) \in \bbQalg[t]
    \end{equation}
    is in $\bbQ[t]$ if and only if $\Phi\prm$ is stable under $\Gal_\bbQ$. But 
    \begin{equation}
    f^{(c)}(t) = c^{2 \deg f} f(t/c^2) \in \bbQ[t].
    \end{equation}
    Thus $\Phi\prm$ is stabilized under $\Gal_\bbQ$ as desired. 
\end{proof}

By \Cref{Corollary: isogenies are stablized by quadratic twist}, an $m$-isogeny $\phi$ of an elliptic curve $E \in \scrE$ induces an $m$-isogeny $\phi\prm$ with isomorphic kernel for every quadratic twist $E\prm$ of $E$. By contrast, if $j(E) = 0$ or $j(E) = 1728$, quartic or sextic twists $E\prm$ of $E$ might or might not preserve a given isogeny $\phi$: in other words, the induced isogeny $\phi\prm$ on $E\prm$ might or might not be defined over $\bbQ$. 

\begin{example}
The elliptic curve
\begin{equation}
E : y^2=x^3+1
\end{equation}
admits both a $2$-isogeny and a $3$-isogeny. Every sextic twist of $E$ preserves its $3$-isogeny (see \cite[Lemma 2.7]{Pizzo-Pomerance-Voight}), but only quadratic twists of $E$ preserve its $2$-isogeny (see the proof of \Cref{Lemma: parameterizing m-isogenies}).
\end{example}

We now turn our attention to a special class of isogenies, which will occupy us for the remainder of the manuscript. We say that an isogeny $\phi$ is \defi{cyclic} if $\ker \phi$ is a cyclic group. 

\begin{proposition}
    Let $E, E\prm \in \scrE$, and let $\phi : E \to E\prm$ be a cyclic $m$-isogeny. Then $\widehat{\phi} : E\prm \to E$ is also a cyclic $m$-isogeny.
\end{proposition}

\begin{proof}
    We prove that if $\phi$ is not cyclic then $\widehat{\phi}$ is also not cyclic. Indeed, if $\phi$ is not cyclic then for some prime $p$, the $p$-torsion subgroup $E[p]$ of $E(\bbQalg)$ is contained in $\ker \phi$. But $E[p]$ is the kernel of the $p$-multiplication map $[p] : P \mapsto p P$, so by \Cref{Theorem: Isogenies factor} we may write $\phi = \phi\prm \circ [p]$ for some isogeny $\phi\prm : E \to E\prm$. Now by \Cref{Theorem: properties of the dual isogeny}(b) and \Cref{Theorem: properties of the dual isogeny}(e), we see $\widehat{\phi} = [p] \circ \widehat{\phi\prm}$, where we now interpret $[p] : P \mapsto p P$ as an endomorphism of $E\prm$. As $[p]$ is not cyclic, $\widehat{\phi}$ cannot be cyclic, and our claim follows.
\end{proof}

\begin{proposition}\label{Proposition: factoring cyclic m-isogenies}
	Let $E, E\prm \in \scrE$, and let $\phi : E \to E\prm$ be a cyclic $m$-isogeny with $m > 1$. Write $m = p_1 \dots p_r$, where $p_1, \dots, p_r$ are (not necessarily distinct) primes. Then we may write
	\begin{equation}
	\phi = \phi_r \circ \phi_{r-1} \circ \dots \circ \phi_1,
	\end{equation}
	where each $\phi_i$ is a cyclic $p_i$-isogeny.
\end{proposition}

\begin{proof}
	We proceed by induction on $m$. If $m = p$ is prime, we may take $\phi_1 = \phi$. Otherwise, write $m = p m\prm$, with $p$ prime. As the $p$-torsion elements of $\ker \phi$ form a group of order $p$ that is stabilized by $\Gal_\bbQ$, we may apply \Cref{Theorem: Isogenies factor} and \Cref{Theorem: every kernel gives rise to an essentially unique isogeny} to write $\phi = \phi\prm \circ \phi_1$, where $\phi_1$ is a cyclic $p$-isogeny and $\phi\prm$ is a cyclic $m\prm$-isogeny. By induction hypothesis, our claim follows.
\end{proof}

\begin{proposition}\label{Proposition: Compositions of cyclic p-isogenies}
	Let $E, E\prm, E^{\prime \prime} \in \scrE$, and let $\phi : E \to E\prm$ and $\psi : E\prm \to E^{\prime \prime}$ be cyclic $p$-isogenies (over $\bbQ$) for a given prime $p$. Then exactly one of the following conditions holds:
	\begin{itemize}
		\item the composition $\psi \circ \phi$ is a cyclic $p^2$-isogeny;
		\item $E = E^{\prime \prime}$, $\psi$ and $\phi$ are dual isogenies up to sign, and $\psi \circ \phi$ is multiplication by $\pm p$ in $E$.
	\end{itemize}
\end{proposition}

\begin{proof}
	The kernel of $\psi \circ \phi$ in $E(\bbQalg)$ is an abelian group of order $p^2$, and thus either isomorphic to $\bbZ /p^2 \bbZ$ or $\bbZ / p \bbZ \times \bbZ / p \bbZ$. In the former case, $\psi \circ \phi$ is cyclic by definition. In the latter case, $\psi \circ \phi$ has the same kernel as multiplication by $p$, and so is the same map up to sign and isomorphism. But $\scrE$ contains only one elliptic curve from each $\bbQ$-isomorphism class, so $E = E^{\prime \prime}$, and $\psi$ and $\phi$ are dual isogenies up to sign.
\end{proof}

\begin{proposition}\label{Proposition: Compositions of cyclic m and n isogenies}
	Let $E, E\prm, E^{\prime \prime} \in \scrE$, let $\phi : E \to E\prm$ be a cyclic $m$-isogeny (over $\bbQ$), and let $\psi : E\prm \to E^{\prime \prime}$ be a cyclic $n$-isogeny (over $\bbQ$). If $\gcd(m, n) = 1$, then $\psi \circ \phi$ is a cyclic $mn$-isogeny.
\end{proposition}

\begin{proof}
	The kernel of $\psi \circ \phi$ in $E(\bbQalg)$ is an abelian group of order $mn$ which contains both $\bbZ / m \bbZ$ and $\bbZ / n \bbZ$. If $\gcd(m, n) = 1$, then this kernel is necessarily $\bbZ / m \bbZ \times \bbZ / n \bbZ \simeq \bbZ / mn \bbZ$, which is cyclic, and our claim follows.
\end{proof}

We have proven some interesting trivia about cyclic isogenies, but a larger question looms: for what $m \in \bbZ_{>0}$ do cyclic $m$-isogenies exist? Building on earlier work of Ogg \cite{Ogg1, Ogg2}, Mazur \cite{Mazur2} reduced this problem to the assertion that $Y_0(m)(\bbQ) = \emptyset$ for $m \in \set{39, 65, 91, 125, 169}$, and Kenku verified this assertion \cite{Kenku1, Kenku2, Kenku3, Kenku4}.

\begin{theorem}[{\cite[Theorem 1]{Kenku5}, Mazur's theorem on isogenies}]\label{Theorem: Kenku-Mazer Theorem}
	Let $m \in \bbZ_{>0}$. The following are equivalent:
	\begin{itemize}
		\item There exist elliptic curves $E, E\prm / \bbQ$ with a cyclic $m$-isogeny $\phi : E \to E\prm$;
		\item $m \in \set{1, \ldots, 19, 21, 25, 27, 37, 43, 67, 163}$.
	\end{itemize}
\end{theorem}

In \cite{Chiloyan-Lozano-Robledo}, Chiloyan and Lozano-Robledo give a precise description of all possible interrelationships of cyclic isogenies and underlying Mordell--Weil groups for elliptic curves $E/\bbQ$. Using \Cref{Proposition: factoring cyclic m-isogenies}, we extract the following result from \cite{Chiloyan-Lozano-Robledo}.

\begin{theorem}\label{Theorem: elliptic curves with multiple cyclic m-isogenies}
	Let $m \in \bbZ_{>0}$. The following are equivalent:
	\begin{enumerate}
		\item There exist elliptic curves $E, E_1, E_2 \in \scrE$ and cyclic $m$-isogenies $\phi_1 : E \to E_1$ and $\phi_2 : E \to E_2$ with distinct kernels;
		\item $m \in \set{2, 3, 4, 5, 6, 8, 12, 16}$.
	\end{enumerate}
\end{theorem}

\begin{proof}
	By \Cref{Proposition: factoring cyclic m-isogenies} and \Cref{Theorem: Kenku-Mazer Theorem}, it suffices to consider the following levels. Namely, $m \in \set{2,3,4,5,6,8,9,10,12,15,16,18,25,27}$. The cases $m\in\set{2,4,8,16}$, $m=3$, $m=5$, and $m\in\set{6,12}$ follow by inspection of the graphs $T8$, $L3(9)$, $L3(25)$, and $S$, respectively \cite[Tables 1,2,4]{Chiloyan-Lozano-Robledo}. For $m=9,27$, the classified graphs contain at most one cyclic subgroup of the relevant order from any vertex. For $m=10,15$, the rectangular graph $R_4(m)$ has unique rational primary components, hence at most one cyclic $m$-kernel. For $m=18$, the middle vertices of $R_6(18)$ have no cyclic $9$-isogeny, while at a corner the unique $9$- and $2$-primary kernels determine the $18$-kernel. Finally, $L_3(25)$ gives at most one cyclic $25$-kernel \cite[Tables 1,3]{Chiloyan-Lozano-Robledo}.
\end{proof}

\begin{corollary}\label{Corollary: Elliptic curves for which twNeq(X) = twNad(X)}
	Let 
\[
 m \in \set{1, 7, 9, 10, 11, 13, 14, 15, 17, 18, 19, 21, 25, 27, 37, 43, 67, 163}.
 \]
 For all $X > 0$, we have
	\begin{equation}
	\twNeq m (X) = \twNad m (X), \ \text{and} \ \NQeq m (X) = \NQad m (X).
	\end{equation}
\end{corollary}

\begin{proof}
    Immediate from \Cref{Theorem: elliptic curves with multiple cyclic m-isogenies}.
\end{proof}

For $m \in \set{2, 3, 4, 5, 6, 8, 12, 16}$, see \Cref{Theorem: Parameterizing elliptic curves equipped with pairs of isogenies} and \Cref{Lemma: Difference between twNeq and twNad}.

\section{Counting elliptic curves: a brief history}\label{Section: History}

In this section, we briefly recount the history of arithmetic statistics of elliptic curves%
. We place special emphasis on how one can define the size of an elliptic curve, and on discussing the asymptotics of elliptic curves with level structure (specifically, the asymptotics of elliptic curves with a given torsion group or a cyclic $m$-isogeny for particular $m$).

\subsection*{What is the size of an elliptic curve?}\label{Subsection: What is the size of an elliptic curve?}

Mathematicians have long been curious about the statistical behavior of families of elliptic curves. Intuitively natural questions abound. What is the average rank of an elliptic curve over $\bbQ$ \cite{Bektemirov-Mazur-Stein}? What proportion of elliptic curves over $\bbQ$ have trivial torsion \cite{Harron-Snowden}? In order to make sense of these and similar questions, mathematicians needed to endow elliptic curves with a notion of size.

One natural notion is the (absolute value of the) discriminant $\Delta(E)$ of an elliptic curve $E$ \eqref{Equation: Discriminant of Weierstrass equation}. The discriminant $\Delta(E)$ is a nonzero integer, and it is divisible precisely by the primes for which $E(\bbF_p)$ is singular. Guided by the intuition that for reasonably chosen regions $\calR {} \subseteq \bbR^2$, we should expect $\# \parent{\calR {} \cap \bbZ^2} \sim \text{Area}(\calR {})$ (see \Cref{Theorem: Principle of Lipschitz}), Brumer and McGuinness conjectured \cite[section 5]{Brumer-McGuinness} that 
\begin{equation}\label{Equation: count of elliptic curves by discriminant}
\# \set{E \in \scrE : \abs{\Delta(E)} \leq X} \sim \frac{\alpha_\Delta}{\zeta(10)} X^{5/6},
\end{equation}
where
\begin{equation}
\alpha_\Delta \colonequals \frac{3 + \sqrt{3}}{10} \int_1^\infty \frac{\mathrm{d}u}{\sqrt{u^3 - 1}} = \frac{2 \sqrt{\pi } \parent{3 + \sqrt{3}} \Gamma \left(7/6\right)}{10 \Gamma \left(2/3\right)} = 1.149\,249\,826\ldots.
\end{equation}

Elliptic curves have a second invariant, the conductor $\mathrm{cond}(E)$ of $E$. Like the discriminant $\Delta(E)$, the conductor $\mathrm{cond}(E)$ is a nonzero integer, and it is divisible precisely by the primes for which $E(\bbF_p)$ is singular. For such $p$, the conductor conveys more about the singularities of $E(\bbF_p)$ than the discriminant does. In \cite{Brumer-Silverman}, Brumer and Silverman proved coarse bounds for the number of elliptic curves with given (and hence with bounded) conductor. In \cite{Watkins}, Watkins built on \cite{Brumer-McGuinness}, re-establishing the heuristic \eqref{Equation: count of elliptic curves by discriminant}, and making several other heuristic claims. Notably, Watkins conjectured \cite[Heuristic 4.1]{Watkins} that for an explicit $\alpha_{\textrm{cond}} > 0$ we have
\begin{equation}\label{Equation: count of elliptic curves by conductor}
\# \set{E \in \scrE : \textrm{cond}(E) \leq X} \sim \frac{\alpha_{\textrm{cond}}}{\zeta(10)} X^{5/6}.
\end{equation}

In \cite{Faltings1, Faltings2}, Faltings introduced a third invariant of an elliptic curve $E$ (or, more generally, an abelian variety $A$), which we now refer to as the Faltings height $\hht_F(E)$ of $E$. The Faltings height of an elliptic curve also relates to its arithmetic properties, albeit in a subtler way than the discriminant or the conductor do. In 2016, Hortsch \cite[Theorem 1.2]{Hortsch} used a reformulation of the Faltings height for elliptic curves due to Silverman \cite[Proposition 1.1]{Silverman3} to prove that there is an explicit constant $\alpha_{\hht_F} > 0$ such that
\begin{equation}\label{Equation: count of elliptic curves by Faltings height}
\# \set{E \in \scrE : \hht_F(E) \leq X} = \frac{\alpha_{\hht_F}}{\zeta(10)} X^{5/6} + O(X^{1/2} \log^3 X),
\end{equation}
where
\begin{equation}
\alpha_{\hht_F} \approx 348\,716
\end{equation}
is given by an explicit integral.

The na\"ive height is another natural notion of size for an elliptic curve. For $E \in \scrE$ with Weierstrass equation as in \eqref{Equation: Weierstrass equation} or \eqref{Equation: Simplified Weierstrass equation}, we define
\begin{equation}\label{Equation: naive height according to Bhargava}
\hht(E) \colonequals \max(4 \abs{A}^3, 27 \abs{B}^2).
\end{equation}
This notion of size is much easier to study and bound than the discriminant, conductor, or Faltings height. Unlike the discriminant, the na\"ive height cannot fall victim to drastic cancellation between $4 A^3$ and $27 B^2$. Unlike the conductor, the na\"ive height does not require an understanding of the behavior of the elliptic curve at primes of bad reduction. Unlike the Faltings height, the na\"ive height is an integer. One can easily derive the power-saving asymptotic
\begin{equation}
\# \set{E \in \scrE : \hht(E) \leq X} = \frac{\alpha_{\hht}}{\zeta(10)} X^{5/6} + O(X^{1/2})
\end{equation}
(see \Cref{Theorem: Count of elliptic curves} below), where
\begin{equation}
\alpha_{\hht} \colonequals \frac{2^{4/3}}{3^{3/2}} = 0.484\,943\,838\ldots.
\end{equation}
At least conjecturally, then, elliptic curves counted by na\"ive height grow at the same rate as elliptic curves counted by discriminant, conductor, and Faltings height. The na\"ive height naturally arises from the embedding of the coordinates $(A : B)$ from \eqref{Equation: Weierstrass equation} in $\bbP(4, 6)(\bbQ)$ (see \Cref{rmk:moduli0}). It is also natural to view the na\"ive height as arising from the ``two components'' of the discriminant. Throughout the remainder of this thesis, when we write about the ``height'' of an elliptic curve, we refer \emph{only} to its na\"ive height, \emph{not} its Faltings height.

In 1992, Brumer showed that the generalized Birch–Swinnerton-Dyer conjecture and Riemann hypothesis together imply that the average rank of the elliptic curves in $\scrE$, ordered by height, is at most $2.3$ \cite{Brumer}. Under the same hypotheses, Heath-Brown improved this bound from $2.3$ to $2$ \cite{Heath-Brown}, and Young improved it to $25/14 = 1.785\,714\ldots$ \cite{Young}. In 2010, Bhargava and Shankar proved unconditionally that the average rank of elliptic curves over $\bbQ$, ordered by height, is at most $3/2 = 1.5$ \cite{Bhargava-Shankar}. 

\begin{remark}
	Some authors (notably, Duke \cite{Duke} and Harron and Snowden \cite{Harron-Snowden}) define the na\"ive height of an elliptic curve $E \in \scrE$ to be
	\begin{equation}\label{Equation: naive height according to Harron and Snowden}
	\max(\abs{A}^3, \abs{B}^2).
	\end{equation}
	We follow \cite{Bhargava-Shankar} in using \eqref{Equation: naive height according to Bhargava}, but our results would all go through without issue if we used \eqref{Equation: naive height according to Harron and Snowden}. Only the coefficients, \emph{not the exponents}, in our asymptotics would be changed.
\end{remark}

Since then, numerous authors have studied elliptic curves ordered by the na\"ive height \cite{Boggess-Sankar, Bruin-Najman, Harron-Snowden, Phillips1, Pizzo-Pomerance-Voight, Pomerance-Schaefer}. Throughout the remainder of this thesis, we concentrate on elliptic curves counted by height, and neglect the interesting parallel questions about how elliptic curves with a cyclic $m$-isogeny would look counted by discriminant, or conductor, or Faltings height.

\subsection*{Counting elliptic curves with torsion by height} 

In 2013, Harron and Snowden \cite{Harron-Snowden} studied the arithmetic statistics of elliptic curves with a given torsion group (see also previous work of Duke \cite{Duke} and Grant \cite{Grant}). They showed that for each finite abelian group $T$ given by Mazur's theorem on torsion (\Cref{Theorem: Mazur's Theorem}), there is an explicit rational number $d(T)$ such that
\begin{equation}\label{Equation: Count of elliptic curves with given torsion}
\#\set{E \in \scrE : E(\bbQ)_{\textup{tors}} \simeq T, \ \hht(E) \leq X} \asymp X^{1/d(T)}.
\end{equation}
Moreover, for $T \in \set{0, \bbZ / 2 \bbZ, \bbZ / 3\bbZ}$, they gave an asymptotic for the left-hand side of \eqref{Equation: Count of elliptic curves with given torsion} with a power-saving error term \cite[Theorem 5.5]{Harron-Snowden}. For $T \in \set{0, \bbZ / 2 \bbZ, \bbZ / 3\bbZ}$, their argument proceeds in three steps: first, they establish parameterizations for those elliptic curves
\begin{equation}
E : y^2 = x^3 + A x + B
\end{equation}
equipped with each of these torsion groups, then they use the Principle of Lipschitz (\Cref{Theorem: Principle of Lipschitz}) and some elementary sieving to estimate the number of elliptic curves in these families, and finally they show that the discrepancy between elliptic curves with torsion containing $T$ and torsion exactly $T$ is relatively small. For larger $T \neq \bbZ / 2 \bbZ \times \bbZ / 2 \bbZ$, they use the theory of modular curves to establish a universal elliptic curve over an open subset of the affine line $\mathbb{A}^1$ equipped with a subgroup isomorphic to $T$; essentially equivalently, they found an elliptic surface
\begin{equation}
E_T : y^2 = x^3 + f_T(t) x + g_T(t)
\end{equation}
over $\bbP^1$ for which evaluation at $t \in \bbP^1(\bbQ)$ (away from a finite set) yields all elliptic curves with torsion containing $T$. Once they had this universal elliptic curve, they used elementary analytic and algebraic arguments to obtain upper and lower bounds for the number of elliptic curves with bounded height occurring in this family. The case $T = \bbZ / 2 \bbZ \times \bbZ / 2 \bbZ$ requires its own argument.

In 2022, Cullinan, Kenney, and Voight \cite{Cullinan-Kenney-Voight} improved on Harron and Snowden's work by attaining asymptotics for the left-hand side of \eqref{Equation: Count of elliptic curves with given torsion} with power-saving error for all torsion groups $T$ given by Mazur's theorem on torsion. Moreover, they provided satisfactory interpretations of the exponent of $X$ and the constants appearing in these asymptotics. Their work was in fact more general, and gave an asymptotic with power-saving error for counts of elliptic curves for which the Galois representation on the $N$-torsion subgroup $E[N]$ of $E(\bbQalg)$ is isomorphic to a fixed group $G$, so long as the associated modular group $\Gamma_G$ is torsion-free, and the associated modular curve $X_G$ has genus $0$ and no irregular cusps \cite[Theorem 1.3.3]{Cullinan-Kenney-Voight}. %
Their proof extends the arguments given in Harron and Snowden: after establishing a universal elliptic curve, applying the Principle of Lipschitz, and doing some elementary sieving, they also needed to address the discrepancy between counting elliptic curves equipped with level structure, and counting the elliptic curves themselves. There is a negligible contribution here from curves with Galois representation properly contained within $G$, and an integer factor arising from the index of $G$ inside its normalizer.

\subsection*{Counting elliptic curves with a cyclic isogeny by height}

In this subsection, we highlight several known results regarding counts of elliptic curves with a cyclic isogeny.

Just as it is natural to study asymptotics for elliptic curves with a given torsion group, as in \Cref{Theorem: Mazur's Theorem}, it is natural to study asymptotics for elliptic curves admitting a cyclic $m$-isogeny, as in \Cref{Theorem: Kenku-Mazer Theorem}. Equivalently, we wish to study elliptic curves $E$ with a cyclic subgroup of $E(\bbQalg)$ with size $m$ that is stable under $\Gal_\bbQ$  (\Cref{Theorem: every kernel gives rise to an essentially unique isogeny}). 

Recall the definitions of $\NQeq m (X)$ and $\NQad m (X)$ from \eqref{Equation: counts of m-isogenies}. We note in passing that the map $(E, \phi) \mapsto E$ defines a surjection
\begin{equation*}
\begin{tikzcd}
\set{(E, \phi) : E \in \scrEX \ \text{and} \ \phi : E \to E\prm \ \textup{an unsigned cyclic $m$-isogeny}} \arrow[d, twoheadrightarrow] \\
\set{E \in \scrEX : \textup{$E$ admits a cyclic $m$-isogeny}}
\end{tikzcd}
\end{equation*}
from the set that $\NQeq m (X)$ counts to the set that $\NQad m (X)$ counts.

We are interested in the asymptotics of $\NQeq m (X)$ and $\NQad m (X)$.

The case $m = 1$ is well-known (\Cref{Theorem: Count of elliptic curves}). The case $m = 2$ is handled by Harron and Snowden's work (\Cref{Theorem: asymptotics for N(X) for m = 2}); indeed, if $E(\bbQalg)$ has a cyclic subgroup of size $2$ which is stabilized by $\Gal_\bbQ$, that subgroup is necessarily defined over $\bbQ$. 

Counting elliptic curves with given torsion and counting elliptic curves with a cyclic $m$-isogeny are both examples of counting elliptic curves with level structure. However, elliptic curves with a cyclic $m$-isogeny have a feature that elliptic curves with torsion do not share. If two elliptic curves are twist equivalent, then one admits a cyclic $m$-isogeny if and only if the other does (\Cref{Corollary: isogenies are stablized by quadratic twist}). Therefore, we must count both the (noncuspidal) points on the modular curve $X_0(m)$ and their quadratic twists.

In 2019, Pizzo, Pomerance, and Voight \cite{Pizzo-Pomerance-Voight} gave a power-saving asymptotic for $\NQeq 3 (X)$ and $\NQad 3(X)$ (see \Cref{Theorem: asymptotics for N(X) = 3} below). Pizzo, Pomerance, and Voight characterized elliptic curves with a cyclic 3-isogeny in terms of the factorization of the 3-division polynomial, rather than directly producing a universal elliptic curve with a cyclic 3-isogeny. They parameterized the elliptic curves (with $j(E) \neq 0$) that have a cyclic 3-isogeny in terms of a triple $(u, v, w)$ of integers satisfying certain arithmetic conditions, then carefully summed over such triples using techniques from analytic number theory. In 2020, Pomerance and Schaefer \cite{Pomerance-Schaefer} provided asymptotics for $\NQad 4 (X)$ and $\NQeq 4(X)$, among other results. Both \cite{Pizzo-Pomerance-Voight} and \cite{Pomerance-Schaefer} used intricate analytic number-theoretic arguments (including a refinement of the Principle of Lipschitz by Huxley \cite{Huxley1}) to improve their errors as far as feasible.

In 2020, building on the arguments outlined above, Boggess and Sankar provided at least an order of growth for $\NQad m (X)$ for $m \in \set{2, 3, 4, 5, 6, 8, 9, 12, 16, 18}$: these $m$, together with $m = 7, 10, 13, 25$, are precisely the integers $m$ from \Cref{Theorem: Kenku-Mazer Theorem} for which the associated modular curve $X_0(m)$ is of genus $0$. Although they still use the Principle of Lipschitz, their approach is qualitatively different from those taken in \cite{Cullinan-Kenney-Voight, Harron-Snowden, Pomerance-Schaefer, Pizzo-Pomerance-Voight}. Boggess and Sankar leverage a modular curve $X_{1/2}(m)$, which is a degree 2 cover of $X_0(m)$; this modular curve lets them track the twists of the elliptic curves arising from $X_0(m)$, and thereby establish orders of growth for $m \in \set{2, 3, 4, 5, 6, 8, 9, 12, 16, 18}$. The ring of modular forms for $X_0(5)$ is especially difficult to handle, and they do so separately from the other $m$ in their list. In 2022, the author and John Voight obtained asymptotics for $\NQeq 7(X) = \NQad 7 (X)$ using the methods which we develop in this thesis \cite{Molnar-Voight}.

Mathematicians are also interested in counting elliptic curves with a cyclic $m$-isogeny over number fields \cite{Bruin-Najman, Phillips1} and even global fields \cite{Banerjee-Park-Schmitt, Phillips2, PhillipsThesis}. In 2020, Bruin and Najman \cite{Bruin-Najman} gave an order of growth estimate for the number of elliptic curves over a number field $K$ admitting a $G$-level structure whenever the associated modular curve $X_G$ over $K$ is isomorphic to a weighted projective line. In 2022, Phillips \cite{Phillips1} gave an asymptotic for the number of elliptic curves over a number field $K$ admitting a $G$-level structure under similar hypotheses. He recovered the asymptotics for $\NQeq m (X)$ for $m \in \set{2, 4}$, and gave the leading term for the asymptotics of $\NQeq m (X)$ for $m \in \set{6, 8, 9, 12, 16, 18}$, albeit without an explicit error term (power-saving or otherwise). In \cite{PhillipsThesis}, he gave similar asymptotics for these $m$ over well-behaved function fields of characteristic greater than 3.

\section{Previous results for \texorpdfstring{$m \leq 5$}{m <= 5}}\label{Section: Previous results}

In this section, we state the asymptotics of $\NQeq m (X)$ and $\NQad m (X)$ for $1 \leq m < 5$, and the order of growth for $\NQeq 5 (X)$ and $\NQad 5 (X)$. Combined with the findings in \cref{Section: Results}, these furnish estimates for the number of elliptic curves with a cyclic $m$-isogeny for all $m \in \bbZ_{>0}$.

Here and throughout the remainder of the thesis, the notation $c \approx \ldots$ indicates that the estimate given for $c$ is numerically supported, but does not have clear and well-bounded error.

The first theorem is folklore (but see \cite[Theorem 2.1]{Hortsch} for a nice treatment).

\begin{theorem}\label{Theorem: Count of elliptic curves}
	We have
	\begin{equation}
	\NQeq 1 (X) = \NQad 1 (X) = \frac{2^{4/3}}{3^{3/2}\zeta(10)} X^{5/6} + O\parent{X^{1/2}}
	\end{equation}
	for $X \geq 1$.
\end{theorem}

Essentially the same argument yields the asymptotic
\begin{equation}\label{Equation: Count of elliptic curves up to twist}
\twistNeq 1 (X) = \twistNad 1 (X) = \frac{2^{4/3}}{3^{3/2}\zeta(5)} X^{5/6} + O\parent{X^{1/2}}.
\end{equation}
The equalities in \Cref{Theorem: Count of elliptic curves} and \eqref{Equation: Count of elliptic curves up to twist} are exact, because no elliptic curve has more than one unsigned isomorphism over $\bbQ$. 

\begin{theorem}\label{Theorem: asymptotics for N(X) for m = 2}
    There exists an effectively computable constant $\constad 2$ such that
    \begin{equation}
    \NQeq 2 (X), \NQad 2 (X) = \constad 2 X^{1/2} + O(X^{1/3})
    \end{equation}
    for $X \geq 1$.
\end{theorem}

\begin{proof}
    Harron and Snowden prove the claim for $\NQad 2 (X)$ in \cite[Theorem 5.5]{Harron-Snowden}, and \cite[Table 1]{Harron-Snowden} assures us that counting elliptic curves with distinct 2-isogenies separately can contribute at most $O(X^{1/3})$ more to this sum.
\end{proof}

\begin{theorem}[{\cite[Theorem 1.3, Proposition 2.9]{Pizzo-Pomerance-Voight}}]\label{Theorem: asymptotics for N(X) = 3}
	There exist effectively computable constants 
    \begin{equation} 	\begin{aligned}
        \constad 3 &= 0.107\,437\,255\,02\ldots \ \text{and} \\ \constadprm 3 &\approx 0.16
    \end{aligned}\end{equation} 
    such that
	\begin{equation}
	\NQeq 3 (X), \NQad 3 (X) = \frac{2}{3^{3/2} \zeta(6)} X^{1/2} + \constad 3 X^{1/3} \log X + \constadprm 3 X^{1/3} + O(X^{7/24})
	\end{equation}
 for $X \geq 1$.
\end{theorem}

\begin{theorem}[{\cite[Theorem 4.2, Theorem 5.11]{Pomerance-Schaefer}}]\label{Intro Theorem: asymptotics for N(X) = 4}
	There exist effectively computable constants 
 \begin{equation} 	\begin{aligned}
     \constad 4 &= 0.957\,400\,377\,047\ldots, \\
     \consteqprm 4 &= - 1.671\,470\,808\,106\ldots, \ \text{and} \\ \constadprm 4 &= -0.871\,250\,852\,030\ldots
 \end{aligned}\end{equation} 
 such that
	\begin{equation}
	\NQeq 4 (X) = 2 \constad 4 X^{1/3} + \consteqprm 4 X^{1/6} + O(X^{21/200})
	\end{equation}
	and
	\begin{equation}
	\NQad 4 (X) = \constad 4 X^{1/3} + \constadprm 4 X^{1/6} + O(X^{21/200})
	\end{equation}
	for $X \geq 1$.
\end{theorem}

\begin{theorem}[{\cite[Proposition 5.9]{Boggess-Sankar}}]\label{Intro Theorem: asymptotics for N(X) = 5}
	We have
	\begin{equation}
	\NQeq 5 (X) \sim \NQad 5 (X) \asymp X^{1/6} \log^2 X
	\end{equation}
	as $X \to \infty$.
\end{theorem}

For completeness, we also record the following folklore theorem in our notation (this result was known to Kenku \cite{Kenku5} and Mazur \cite{Mazur2}).

\begin{theorem}[\Cref{Theorem: asymptotic for twN(X) for m of nonzero genus}]\label{Intro Theorem: asymptotic for twN(X) for m of nonzero genus}
	For $X$ sufficiently large, we have the following identities:
		\begin{equation} 
  \begin{aligned}
		\twistNad {11} (X) =& 3, \ \twistNad {14} (X) = 2, \ \twistNad {15}(X) = 4, \ \twistNad {17} (X) = 2, \\
    \twistNad {19} (X) =& 1, \twistNad {21}(X) = 4, \twistNad {27} (X) = 1, \ \twistNad {37} (X) = 2, \\ 
    \twistNad {43} (X) =& 1, \ \twistNad {67} (X) = 1, \ \twistNad {163} (X) = 1.
		\end{aligned}
  \end{equation} 
\end{theorem}

For $m \not\in \set{1, \ldots, 19, 21, 25, 27, 37, 43, 67, 163}$, we have
\begin{equation}
\twistNeq m (X) = \twistNad m (X) = 0 %
\end{equation}
identically by Mazur's theorem on isogenies (\Cref{Theorem: Kenku-Mazer Theorem}).

%% file: Ch3_Preliminaries.tex
\chapter{Technical preliminaries}\label{Chapter: Preliminaries}

In this chapter, we recall and build upon a number of results from the literature \cite{Bingham-Goldie-Teugels, Cullinan-Kenney-Voight, Davenport, Faltings1, Faltings2, Horn-Johnson, Ivic, Landau1915, Liu, Molnar-Voight, Phillips1, Pomerance-Schaefer, Roux, Tenenbaum, Walfisz, Widder} which will be important in proving our main theorems (\Cref{Intro Theorem: asymptotic for NQ(X) for 5 < m <= 9}, \Cref{Intro Theorem: asymptotic for NQ(X) for m > 9}, \Cref{Intro Theorem: asymptotic for twN(X) for m of genus $0$}, \Cref{Intro Theorem: asymptotic for twN(X) for m of nonzero genus}). 

In \cref{Section: Height and defect}, we examine the na\"ive height and the twist height of an elliptic curve more carefully, and we establish additional notation for later use in this thesis. In \cref{Section: Parameterizing elliptic curves with a cyclic m-isogeny}, for 
\begin{equation}
m \in \set{4, 5, 6, 7, 8, 9, 10, 12, 13, 16, 18, 25},
\end{equation}
we give a universal elliptic curve
\begin{equation}\label{Equation: universal model up to quadratic twist}
y^2 = x^3 + \f m (t) x + \g m (t)
\end{equation}
over an open subset of the affine line $\mathbb{A}^1$ equipped with a cyclic $m$-isogeny. For each proper divisor $n$ of $m$, we also report the modular curve parameterizing elliptic curves equipped with a pair of cyclic $m$-isogenies whose kernels have intersection of order $n$, as long as this modular curve is defined over $\bbQ$; when this modular curve is of genus $0$ but is not $X_0(m)$, we give a universal elliptic curve
\begin{equation}
y^2 = x^3 + \f {m, n} (t) x + \g {m, n} (t)
\end{equation}
over an open subset of the affine line $\mathbb{A}^1$ equipped with such pairs of isogenies.
In \cref{Section: Lattices and the principle of Lipschitz} we recall the Principle of Lipschitz, and use it to obtain counts of Weierstrass models arising from the universal elliptic curves given by \eqref{Equation: universal model up to quadratic twist}. In \cref{Section: Some analytic trivia} we recall a number of results from analytic number theory which we will require later in the thesis. In \cref{Section: Our approach revisited}, we outline the approach we take in resolving this problem. The material in this section mirrors and abstracts material from \cite[\S 2.1, \S 4.1, \S 5.1]{Molnar-Voight}. 

Impatient readers may restrict their attention to \cref{Section: Height and defect}, \cref{Section: Parameterizing elliptic curves with a cyclic m-isogeny}, and \cref{Section: Our approach revisited}, and refer back to the intervening sections of this chapter as needed.

\section{Height, minimality, and defect} \label{Section: Height and defect}

In this section, we define the height and twist height of an elliptic curve over $\bbQ$ without reference to minimal models. We also define the minimality defect and twist minimality defect of a Weierstrass equation, which measure the discrepancy between the size of the coefficients of a given Weierstrass model and the height or twist height of the associated elliptic curve. 

Let $E$ be an elliptic curve over $\bbQ$, and let
\begin{equation}\label{eqn:yaxb}
y^2 = x^3 + Ax + B
\end{equation}
be a Weierstrass model for $E$ with $A,B \in \bbZ$. Unlike in \cref{Section: Results} and \cref{Section: A primer on elliptic curves}, we do not assume $E \in \scrE$; that is, we do not assume that this model is minimal. Define 
\begin{equation} \label{eqn:HAB} 
H(A,B) \colonequals \max(\abs{4A^3},\abs{27B^2}).
\end{equation}
Note that $H(A, B)$ depends on our choice of model, and not only on our choice of elliptic curve $E$.

\begin{example}\label{Example: H for equivalent elliptic curves}
    The Weierstrass equations
    \begin{equation}\label{Equation: minimal model example}
    y^2 = x^3 + 4 x + 8 
    \end{equation}
    and
    \begin{equation}\label{Equation: nonminimal model example}
    y^2 = x^3 + 324 x + 5832
    \end{equation}
    are models for the same elliptic curve over $\bbQ$. However, $\rawheight(4, 8) = 1728$ while on the other hand $\rawheight(324, 5832) = 918\,330\,048$.
\end{example}

The largest $d \in \Z_{>0}$ such that $d^4 \mid A$ and $d^6 \mid B$ is called the \defi{minimality defect} $\md(A,B)$ of the model \eqref{eqn:yaxb}. Explicitly, we have
\begin{equation}\label{Equation: defect powers}
\md(A,B) \colonequals \prod_{\ell} \ell^{v_\ell}, \quad \text{where} \ v_\ell \colonequals \lfloor \min(\ord_\ell(A)/4,\ord_\ell(B)/6) \rfloor,
\end{equation}
with the product over all primes $\ell$.
We now redefine the \defi{(na\"ive) height} of $E$ to be
\begin{equation}
\hht(E)=\hht(A,B) \colonequals \frac{H(A,B)}{\md(A,B)^{12}}, \label{eqn:justheight}
\end{equation}
which is well-defined up to $\bbQ$-isomorphism. Indeed, if we write $d = \md (A, B)$, the integral Weierstrass equation
\begin{equation} 
y^2=x^3+(A/d^4)x+(B/d^6)
\end{equation}
has minimality defect 1, and is the unique model for $E$ such that for every prime $\ell$ we have $\ell^4 \nmid A$ or $\ell^6 \nmid B$. Our new definition \eqref{eqn:justheight} for the na\"ive height thus agrees with the definition \eqref{Equation: Intro height} we gave in the introduction.

We may similarly consider all integral Weierstrass equations which define an elliptic curve that is twist equivalent to $E$---these are the quadratic twists of $E$ (defined over $\bbQ$). We call the largest $e\in \Z_{>0}$ such that $e^2 \mid A$ and $e^3 \mid B$ the \defi{twist minimality defect} of the model \eqref{eqn:yaxb}, denoted $\tmd(A,B)$. Explicitly, we have
\begin{equation}\label{Equation: twist defect powers}
\tmd(A,B) \colonequals \prod_{\ell} \ell^{v_\ell}, \quad \text{where} \ v_\ell \colonequals \lfloor \min(\ord_\ell(A)/2,\ord_\ell(B)/3) \rfloor,
\end{equation}
with the product over all primes $\ell$. As above, we then define the \defi{twist height} of $E$ to be
\begin{equation} \label{eqn:htdef}
\twht(E)=\twht(A,B) \colonequals \frac{H(A,B)}{\tmd(A,B)^{6}},
\end{equation}
which is well-defined up to $\bbQ$-isomorphism, and even up to $\bbQalg$-isomorphism when $j(E) \neq 0, 1728$. Indeed, if we write $e = \tmd (A, B)$, the integral Weierstrass equations
\begin{equation} 
y^2=x^3+(A/e^2)x\pm(B/e^3)
\end{equation}
have twist minimality defect 1, and are the only twist minimal models in the quadratic twist class of $E$ such that for every prime $\ell$ we have $\ell^2 \nmid A$ or $\ell^3 \nmid B$. We refer to Weierstrass models with twist minimality defect 1 as \defi{twist minimal}: just as the height of $E/\bbQ$ is the height of its minimal model, the twist height of $E/\bbQ$ is the twist height of its twist minimal model.

\begin{example}
    The Weierstrass equations \eqref{Equation: minimal model example} and \eqref{Equation: nonminimal model example} given in \Cref{Example: H for equivalent elliptic curves} have na\"ive height $1728$ and twist height $27$. The equation \eqref{Equation: minimal model example} is a minimal model for the elliptic curve described by these equations, and the two twist minimal models for this elliptic curve are as follows:
    \begin{equation}
    y^2 = x^3 + x \pm 1.
    \end{equation}
\end{example}

Recall from \cref{Section: Motivation and setup} that $\twistE$ consists of the elliptic curves in $\scrE$ up to twist equivalence. By the remarks in the previous paragraph, the twist minimal models
\begin{equation}
E : y^2 = x^3 + A x + B
\end{equation}
for which $B > 0$ or $(A, B) = (1, 0)$ provide a collection of representatives for $\twistE$, and we identify $\twistE$ with this set of Weierstrass models for the remainder of the thesis.

With this convention established, we make an elementary observation. We have
\begin{equation}\label{Equation: scrE as twists of twistE}
\scrE = \set{E^{(c)} : E \in \twistE \ \text{and} \ c \in \bbZ \ \text{squarefree}}.
\end{equation}
Here, as in \eqref{Equation: Defining E^(c)}, $E^{(c)}$ is the quadratic twist of $E$ by $c$. Note that for $E^{(c)}$ as in \eqref{Equation: scrE as twists of twistE}, we have
\begin{equation}\label{Equation: quadratic twists multiply height by c^6}
\hht(E^{(c)}) = c^6 \hht(E) = c^6 \twistheight(E).
\end{equation}

\begin{remark} \label{rmk:moduli0}
This setup records in a direct manner the more intrinsic notions of height coming from moduli stacks. The moduli stack $Y(1)_\bbQ$ of elliptic curves admits an open immersion into a weighted projective line $Y(1) \hookrightarrow \PP(4,6)_\bbQ$ by $E \mapsto (A:B)$ for any choice of Weierstrass model \eqref{eqn:yaxb}, and the height of $E$ is the height of the point $(A:B) \in \PP(4,6)_\bbQ$ associated to $\scrO_{\PP(4,6)}(12)$ (with coordinates harmlessly scaled by $4,27$): see Bruin--Najman \cite[\S 2, \S 7]{Bruin-Najman} and Phillips \cite[\S 2.2]{Phillips1}. Similarly, the height of the twist minimal model is given by the height of the point $(A:B) \in \PP(2,3)_\bbQ$ associated to $\scrO_{\PP(2,3)}(6)$, which is almost but not quite the height of the $j$-invariant (in the usual sense). 
\end{remark}

We remark in passing that for $E/\bbQ$ and $c, c\prm \in \bbQ^\times$, we have
\begin{equation}
\parent{E^{(c)}}^{(c\prm)} = E^{(c c\prm)},
\end{equation}
and thus
\begin{equation}
\parent{E^{(c)}}^{(c)} = E^{(c^2)} = E,
\end{equation}
up to $\bbQ$-isomorphism. 

\section{Parameterizing elliptic curves equipped with a cyclic \texorpdfstring{$m$}{m}-isogeny}\label{Section: Parameterizing elliptic curves with a cyclic m-isogeny}

In this section, we use the theory of modular curves to parameterize all elliptic curves with a cyclic $m$-isogeny when $X_0(m)$ is of genus $0$ and $m > 3$. We also produce modular curves parameterizing elliptic curves equipped with pairs of distinct cyclic $m$-isogenies, whenever such pairs exist over $\bbQ$.

We gather the necessary input from the theory of modular curves. The modular curve $Y_0(m) \subseteq X_0(m)$, defined over $\bbQ$, parameterizes pairs $(E,\phi)$ of elliptic curves $E$ equipped with a cyclic $m$-isogeny $\phi$ up to isomorphism, or equivalently, a cyclic subgroup of order $m$ stable under the absolute Galois group $\Gal_\bbQ \colonequals \Gal(\bbQalg\,|\,\bbQ)$ (see \Cref{Theorem: every kernel gives rise to an essentially unique isogeny}).  For $m \in \set{4, 5, 6, 7, 8, 9, 10, 12, 13, 16, 18, 25},$ we observe that $Y_0(m) \subseteq X_0(m)$ is affine open in $\PP^1$. In these cases, the objects of interest are parameterized by a coordinate $t$ in this affine open subset.

\begin{lemma}\label{Lemma: parameterizing m-isogenies}
	Let $m \in \set{4, 5, 6, 7, 8, 9, 10, 12, 13, 16, 18, 25}$. Then the set of elliptic curves $E$ over $\bbQ$ that admit a cyclic $m$-isogeny (defined over $\bbQ$) is precisely the set of those of the form 
	\begin{equation}\label{Equation: parameterizing elliptic curves with m isogeny}
	E_m^{(c)}(t) \colon y^2 = x^3 + c^2\f m(t)x + c^3\g m(t)
	\end{equation}
	for some $c \in \bbQ^\times$ and $t \in \bbQ$ with $t \not\in \cusps m$, where $\cusps m \subseteq \bbP^1 (\bbQ)$ consists of those elements $t \in \bbP^1(\bbQ)$ for which the Weierstrass equation \eqref{Equation: parameterizing elliptic curves with m isogeny} is singular. Moreover, the set
 \begin{equation}
\set{(c, t) \in \bbZ \times \bbQ : c \ \text{squarefree}, \ t \not\in \cusps m}
 \end{equation}
 is in bijection with the set of elliptic curves equipped with an unsigned cyclic $m$-isogeny via the map $(c, t) \mapsto E_m^{(c)}(t)$.
 
 The polynomials $\f m (t)$ are given in Table \ref{table:fm}, the polynomials $\g m (t)$ are given in Table \ref{table:gm}, and $\cusps m$ is given in Table \ref{table:auxiliaries}.
\end{lemma}

\begin{proof}
    The proof follows for elliptic curves $E$ with $j(E) \neq 0, 1728$ by routine calculations with $q$-expansions for modular forms on the group $\Gamma_0(m) \subseteq \mathrm{SL}_2(\bbZ)$, with the cusps at $t \in \cusps m \subseteq \bbP^1(\bbQ)$. For instance, writing the $j$-function in terms of the Hauptmodul for $X_0(m)$ yields a parameterization
    \begin{equation}\label{Equation: parameterizing elliptic curves with m isogeny up to quadratic twist}
        E_m(t) : y^2 = x^3 + \f m (t) x + \g m (t),
    \end{equation}
    for elliptic curves equipped with an unsigned cyclic $m$-isogeny, up to quadratic twist (see also \cite[Proposition 3.3.16]{Cullinan-Kenney-Voight}). Taking quadratic twists of \eqref{Equation: parameterizing elliptic curves with m isogeny up to quadratic twist}, we obtain \eqref{Equation: parameterizing elliptic curves with m isogeny}.

    We now turn our attention to the circumstance that $j(E) = 0, 1728$. By inspection, for $E_m(t)$ as in \eqref{Equation: parameterizing elliptic curves with m isogeny up to quadratic twist}, $j(E_m(t)) = 0$ if and only if $(m, t) \in \set{(6, 1), (9, 0)}$, and $j(E_m(t)) = 1728$ if and only if $(m, t) = (4, 0)$. In these cases, we must investigate the sextic or quartic twists of $E_m(t)$, respectively (\Cref{Corollary: Isomorphisms are quadratic twists for j not 0 and 1728}).

    First let $m = 6$ and $t = 1$. We compute $\g 6 (1) = -64$, which after a quadratic twist yields the elliptic curve
    \begin{equation}
    E : y^2 = x^3 + 1.
    \end{equation}
    Now let
    \begin{equation}
    E\prm : y^2 = x^3 \pm c
    \end{equation}
    be a sextic twist of $E$. If $E\prm$ has a rational 6-isogeny, then \textit{a fortiori} it also has a rational 2-torsion point, so $c^{1/3} \in \bbQ$. But this implies $E\prm$ is a quadratic twist of $E$, as desired.

    Now let $m = 9$ and $t = 0$. We compute $\g 9 (0) = -16$, which after a quadratic twist yields the elliptic curve
    \begin{equation}
    E : y^2 = x^3 - 2.
    \end{equation}
    Let $\Phi \subseteq E(\bbQalg)$ be the kernel of the 9-isogeny on $E$. Inspecting the division polynomial for $E$, we readily compute the polynomial 
    \begin{equation}
    \prod_{(x, \pm y) \in \Phi} (t - x) = t(t^3 + 6 t^2 - 8).
    \end{equation}
    For $u \in (\bbQalg)^\times$, if $(x, y) \mapsto (u^2 x, u^3 y)$ is a (not necessarily quadratic) twist of $E$ which sends $\Phi$ to another group stabilized by $\Gal_\bbQ$, then
    \begin{equation}
    t (t^3 + 6 u^2 t^2 - 8 u^6)
    \end{equation}
    is a rational polynomial. In particular, $u^2 \in \bbQ^\times$, so $(x, y) \mapsto (u^2 x, u^3 y)$ is a quadratic twist.

    Finally, let $m = 4$ and $t = 0$. In this case, we find $\f 4 (0) = 9$, which after a quadratic twist yields the elliptic curve
    \begin{equation}
    E : y^2 = x^3 + x.
    \end{equation}
    Note that this elliptic curve has two distinct 4-isogenies over $\bbQ$, generated by the points $(1, 2^{1/2})$ and $(-1, (-2)^{1/2})$ respectively. If $(x, y) \mapsto (u^2 x, u^3 y)$ is a twist of $E$ which preserves either of these 4-isogenies, then $\pm u^2 \in \bbQ^\times$, so again we have a quadratic twist, as desired.
\end{proof}

We highlight that the point $t = 0$ encodes \emph{two distinct} cyclic 4-isogenies over $\bbQ$. It turns out that elliptic curves equipped with a cyclic 4-isogeny always come in pairs (see \Cref{Corollary: discrepancy between twNeq m and twNad m when 4 | m} and the paragraphs immediately preceding it), but this is the only case where this pair is indexed by one argument $t$ instead of two arguments $t_1$ and $t_2$.

Of course we can ignore the factor $c$ in \Cref{Lemma: parameterizing m-isogenies} for elliptic curves over $\bbQ$ up to quadratic twist.

For $m \in \set{4, 5, 6, 7, 8, 9, 10, 12, 13, 16, 18, 25}$, we define $\degB m \colonequals \deg \g m (t)$, and note that $\deg \f m (t) = 2 \degB m / 3$. 

\jvtable{table:fm}{
\rowcolors{2}{white}{gray!10}
\begin{tabular}{c | c}
$m$ & $\f m (t)$ \\
\hline\hline
$4$ & $-3 (t^2 + 6 t - 3)$ \\
$5$ & $-3 (t^2 + 1) \cdot (t^2 + 114 t + 124)$ \\
$6$ & $-3 (t - 1) \cdot (t^3 + 33 t^2 - 117 t + 99)$ \\
$7$ & $-3 (t^2 + t + 7) \cdot (t^2 - 231 t + 735)$ \\
$8$ & $-3 (t^4 + 24 t^3 - 88 t^2 + 96 t - 32)$ \\
$9$ & $-3 t (t^3 - 24 t^2 - 24 t - 8)$ \\
$10$ & $-3 (t^2 + 1) \cdot (t^6 - 118 t^5 + 360 t^4 - 240 t^3 + 240 t^2 - 8 t + 4)$ \\
$12$ & $-3 (t^2 - 3) \cdot (t^6 - 108 t^5 - 657 t^4 - 1512 t^3 - 1701 t^2 - 972 t - 243)$ \\
$13$ & $-3 (t^2 + t + 7) \cdot (t^2 + 4) \cdot (t^4 - 235 t^3 + 1211 t^2 - 1660 t + 6256)$ \\
$16$ & $-3 (t^8 + 48 t^7 - 432 t^6 + 1536 t^5 - 2896 t^4 + 3072 t^3 - 1728 t^2 + 384 t + 16)$ \\
$18$ & \begin{tabular}{@{}c@{}}
$-3 (t^3 + 6 t^2 + 4) \cdot (t^9 + 234 t^8 + 756 t^7 + 2172 t^6$ \\
$+ 1872 t^5 + 3024 t^4 + 48 t^3 + 3744 t^2 + 64)$
\end{tabular} \\
$25$ & \begin{tabular}{@{}c@{}}
$-3 (t^2 + 4) \cdot (t^{10} + 240 t^9 + 2170 t^8 + 8880 t^7 + 34835 t^6 + 83748 t^5$ \\
$+ 206210 t^4 + 313380 t^3 + 503545 t^2 + 424740 t + 375376)$ \\
\end{tabular}
\end{tabular} \\
}{$\f m (t)$ for $m$ with $X_0(m)$ of genus $0$}

\jvtable{table:gm}{
\rowcolors{2}{white}{gray!10}
\begin{tabular}{c | c}
$m$ & $\g m (t)$ \\
\hline\hline
$4$ & $2 t \cdot (t^2 - 18 t + 9)$ \\
$5$ & $2 (t^2 + 1)^2 (t^2 - 261 t - 2501)$ \\
$6$ & $2 (t^2 - 3) (t^4 - 96 t^3 + 426 t^2 - 648 t + 333)$  \\
$7$ & $2 (t^2 + t + 7) (t^4 + 518 t^3 - 11025 t^2 + 6174 t - 64827)$ \\
$8$ & $2 (t^2 - 2) (t^4 - 72 t^3 + 248 t^2 - 288 t + 112)$ \\
$9$ & $2 (t^6 + 60 t^5 - 12 t^4 - 124 t^3 - 120 t^2 - 48 t - 8)$ \\
$10$ & \begin{tabular}{@{}c@{}} 
$2 (t^2 + 1)^2 (t^2 - 2 t + 2) (t^2 - 11 t - 1)$ \\
$(t^4 + 268 t^3 - 66 t^2 + 52 t - 4)$ 
\end{tabular} \\
$12$ & \begin{tabular}{@{}c@{}} 
$2 (t^4 - 6 t^3 - 36 t^2 - 54 t - 27) (t^8 + 276 t^7 + 1836 t^6$ \\
$+ 4860 t^5 + 6750 t^4 + 6156 t^3 + 4860 t^2 + 2916 t + 729)$ \\
\end{tabular} \\
$13$ & \begin{tabular}{@{}c@{}}
$2 (t^2 + t + 7) (t^2 + 4)^2 (t^6 + 512 t^5 - 13073 t^4$ \\
$+ 34860 t^3 - 157099 t^2 + 211330 t - 655108)$ \\
\end{tabular} \\
$16$ & 
\begin{tabular}{@{}c@{}} 
$2 (t^4 - 12 t^2 + 24 t - 14) (t^8 - 144 t^7 + 1200 t^6 - 4416 t^5 $ \\
$+ 9152 t^4 - 11520 t^3 + 8832 t^2 - 3840 t + 736)$  \\
\end{tabular} \\
$18$ & \begin{tabular}{@{}c@{}}
$2 (t^6 + 24 t^5 + 24 t^4 + 92 t^3 - 48 t^2 + 96 t - 8) (t^{12} - 528 t^{11}$ \\
$ - 3984 t^{10} - 14792 t^9 - 27936 t^8 - 42624 t^7 - 37632 t^6$ \\ 
$- 52992 t^5 - 25344 t^4 - 43520 t^3 - 6144 t^2 - 6144 t - 512)$ \\
\end{tabular} \\
$25$ & \begin{tabular}{@{}c@{}}
$2 (t^2 + 4)^2 \cdot (t^4 + 6 t^3 + 21 t^2 + 36 t + 61) (t^{10} - 510 t^9$ \\
$ - 13580 t^8 - 36870 t^7 - 190915 t^6 - 393252 t^5 - 1068040 t^4$ \\ 
$- 1508370 t^3 - 2581955 t^2 - 2087010 t - 1885124)$ \\
\end{tabular}
\end{tabular}
}{$\g m (t)$ for $m$ with $X_0(m)$ of genus $0$}

For $m \in \set{4, 5, 6, 7, 8, 9, 10, 12, 13, 16, 18, 25}$, we record $\degB m$ and $\cusps m$ in the following table. For later use, we also record the resultant $\Res(\f m (t), \g m (t))$ of $f_m(t)$ and $g_m(t)$.

\jvtable{table:auxiliaries}{
\rowcolors{2}{white}{gray!10}
\begin{tabular}{c | c c c}
$m$ & $\degB m$ & $\cusps m$ & $\Res(\f m (t), \g m (t))$ \\
\hline\hline
$4$ & $3$ & $\set{1/2, 1, \infty}$ & $2^6 \cdot 3^6$ \\
$5$ & $6$ & $\set{11/2, \infty}$ & $0$ \\
$6$ & $6$ & $\set{3/2, 5/3, 3, \infty}$ & $-2^{24} \cdot 3^{15}$ \\
$7$ & $6$ & $\set{-7, \infty}$ & $0$ \\
$8$ & $6$ & $\set{1, 3/2, 2, \infty}$ & $2^{24} \cdot 3^{12}$ \\
$9$ & $6$ & $\set{-1, \infty}$ & $- 2^{16} \cdot 3^{15}$ \\
$10$ & $12$ & $\set{-2, 0, 1/2, \infty}$ & $0$ \\
$12$ & $12$ & $\set{-3, -2, -3/2, -1, 0, \infty}$ & $2^{48} \cdot 3^{78}$ \\
$13$ & $12$ & $\set{-3, \infty}$ & $0$ \\
$16$ & $12$ & $\set{1, 3/2, 2, \infty}$ & $2^{48} \cdot 3^{24}$ \\
$18$ & $18$ & $\set{-1, 0, 2, \infty}$ & $- 2^{144} \cdot 3^{189}$ \\
$25$ & $18$ & $\set{1, \infty}$ & $0$ \\
\end{tabular}
}{Miscellaneous data about the model $E_m(t) : y^2 = x^3 + \f m (t) x + \g m (t)$}

We emphasize that $\f m (t)$ and $\g m (t)$ are unique only up to a fractional linear transformation, and that such fractional linear transformations will also act on $\cusps m$. Of course, $\# \cusps m$ is independent of our choice of model 
\begin{equation}
E_m(t) : y^2 = x^3 + \f m (t) x + \g m (t).
\end{equation}

To work with integral models, we take $t=a/b$ (in lowest terms) and homogenize, giving the following polynomials in $\Z[a,b]$:
\begin{equation} \label{Equation: defining ABCab for m = 7}
\begin{aligned}
	\A m(a,b) &\colonequals b^{2 \degB m / 3} \f m(a/b), \\
	\B m(a,b) &\colonequals b^{\degB m} \g m(a/b). \\ 
\end{aligned}
\end{equation}

For 
\[
m \in \set{4, 5, 6, 7, 8, 9, 10, 12, 13, 16, 18, 25},
\]
the set $\cusps m \subseteq \bbP^1 (\bbQ)$ described in \Cref{Lemma: parameterizing m-isogenies} consists of those elements $(a : b) \in \bbP^1(\bbQ)$ for which
\begin{equation}
E_m(a, b) : y^2 = x^3 + \A m (a, b) x + \B m (a, b)\label{Equation: elliptic curve defined in terms of fm and gm}
\end{equation}
is singular. Thus
\begin{equation}
	\cusps m = \set{(a : b) \in \bbP^1 (\bbQ) : 4 \A m(a, b)^3 + 27 \B m(a, b)^2 = 0}.
\end{equation}

\begin{definition}\label{Definition: m-groomed pairs}
	We say that a pair $(a, b) \in \bbZ^2$ is \defi{$m$-groomed} if $\gcd(a, b) = 1$, $b > 0$, and $a/b \not\in \cusps m$.
\end{definition}

\begin{remark}
	It would be more technically accurate to define a pair $(a, b) \in \bbZ^2$ to be $m$-groomed if $\gcd(a, b) = 1$, $b > 0$ or $(a, b) = (1, 0)$, and $a/b \not\in \cusps m$. However, \Cref{Definition: m-groomed pairs} is harmless because $\infty \in \cusps m$ for all 
 \[
 m \in \set{4, 5, 6, 7, 8, 9, 10, 12, 13, 16, 18, 25}.
 \]
\end{remark}

Thus \Cref{Lemma: parameterizing m-isogenies} asserts that elliptic curves $E \in \scrE$ that admit a cyclic $m$-isogeny are precisely those with a model 
\begin{equation}
y^2 = x^3 + \frac{c^2 \A m(a, b)}{d^4} x + \frac{c^3 \B m(a, b)}{d^6} \label{Equation: Weierstrass equation for elliptic surface, homogenized}
\end{equation}
where $(a, b)$ is $m$-groomed, $c \in \Z$ is squarefree, and $d=\md(c^2\A m(a,b),c^3\B m(a,b))$.
For $m > 4$, the count $\NQeq m(X)$ can be computed as
\begin{equation} \label{Equation: NQx in terms of abc}
\begin{aligned}
\NQeq m(X) =& \#\left\{(a, b, c) \in \bbZ^3 : 
\begin{minipage}{36ex} $(a, b)$ $m$-groomed, $c$ squarefree, and \\
$\hht(c^2 \A m(a,b),c^3 \B m(a,b)) \leq X$
\end{minipage}
\right\}
\end{aligned}
\end{equation}
with the height defined as in \eqref{eqn:justheight}. 

\begin{remark}
    It is not hard to adapt this asymptotic to estimate $\NQeq 4 (X)$. For $X > 0$ a real number, we let
    \begin{equation}\label{Equation: Definition of squarefree count}
	\kfree 2(X) \colonequals \# \set{n \in \bbZ_{>0} : n \leq X, \ n \ \textup{a squarefree integer}},
    \end{equation}
so $\kfree 2 (X)$ is the number of squarefree positive integers less than or equal to $X$. When $m = 4$, we need only add $2 \kfree 2 ((X/4)^{1/6})$ to the right-hand side of \eqref{Equation: NQx in terms of abc} to account for the fact that there are two cyclic 4-isogenies associated to each triple of the form $(0, 1, c)$. However, \Cref{Intro Theorem: asymptotics for N(X) = 4} addresses the asymptotics of $\NQeq 4 (X)$ and $\NQad 4 (X)$ to our satisfaction, so we need not pursue this case further.
\end{remark}

Similarly, but more simply, for $m \in \set{4, 5, 6, 7, 8, 9, 10, 12, 13, 16, 18, 25}$, the subset of $E \in \twistE$ that admit a cyclic $m$-isogeny consists of models
\begin{equation}
y^2 = x^3 + \frac{\A m(a, b)}{e^2} x + \frac{\abs{\B m(a, b)}}{e^3} \label{Equation: Weierstrass equation for elliptic surface, homogenized, twist}
\end{equation}
with $(a, b)$ $m$-groomed and $e=\twistdefect(\A m(a,b),\B m(a,b))$ the twist minimality defect defined in \eqref{Equation: twist defect powers}. Accordingly, for $m > 4$ we have
\begin{equation} \label{eqn:twistheightcalc}
\begin{aligned}
\twistNeq m(X) =& \# \set{(a, b) \in \bbZ^2 : \begin{minipage}{33ex} $(a, b)$ $m$-groomed, and \\
$\twht(\A m(a,b),\B m(a,b)) \leq X$
\end{minipage}}.
\end{aligned}
\end{equation}
When $m = 4$, the count on the left and the count on the right can differ by at most 1.

\begin{remark} \label{rmk:moduli}
Returning to \Cref{rmk:moduli0}, we conclude that counting elliptic curves over $\bbQ$ equipped with an $m$-isogeny is the same as counting points on $\PP(4,6)_\bbQ$ in the image of the natural map $Y_0(m) \to Y(1) \subseteq \PP(4,6)_\bbQ$.  Counting them up to twist replaces this with the further natural quotient by $\mu_2$, giving $\PP(2,3)_\bbQ$.  
\end{remark}

We now turn our attention to those elliptic curves which admit pairs of cyclic $m$-isogenies. These elliptic curves are not parameterized by a single modular curve; however, they are parameterized by a family of modular curves, one for each proper divisor $n$ of $m$. \Cref{Theorem: Parameterizing elliptic curves equipped with pairs of isogenies} and the comments that follow it may be extracted from \cite{Chiloyan-Lozano-Robledo}, which proves a great deal more; however, we opt to give a self-contained argument.

\begin{theorem}\label{Theorem: Parameterizing elliptic curves equipped with pairs of isogenies}
    Let
    \[
    m \in \set{2, 3, 4, 5, 6, 7, 8, 9, 10, 12, 13, 16, 18, 25}
    \]
    and let $n$ be a proper divisor of $m$. Then Table \ref{table:modularcurves} records each modular curve defined over $\bbQ$ which parameterizes elliptic curves equipped with pairs of unsigned cyclic $m$-isogenies whose kernels have intersection of order $n$.
\end{theorem}

The RSZB labels used in Table \ref{table:modularcurves} are defined in \cite{Rouse-Sutherland-Zureick-Brown}. Let $\Gamma \leq \GL_2(\widehat{\bbZ})$ be an open subgroup, and let $X_\Gamma$ be the associated modular curve. The first three numbers that occur in the RSZB label of $X_\Gamma$ are its level, the index of $\Gamma$ in $\GL_2(\widehat{\bbZ})$, and the genus of $X_\Gamma$. For instance, when $m = 6$ and $n = 2$, we obtain the modular curve $X_0(2) \times_{X(1)} X_{\textup{sp}}(3)$, which has RSZB label 6.36.0.1. Thus $X_0(2) \times_{X(1)} X_{\textup{sp}}(3)$ is level 6, has index $36$, and has genus $0$.

\jvtable{table:modularcurves}{
\rowcolors{2}{white}{gray!10}
\begin{tabular}{c c | c}
$m$ & $n$ & RSZB label \\
\hline
\hline
$2$ & $1$ & 2.6.0.1 \\
\hline
$3$ & $1$ & 3.12.0.1 \\
\hline
$4$ & $1$ & 4.24.0.8 \\
$4$ & $2$ & 4.6.0.1 \\
\hline
$5$ & $1$ & 5.30.0.1 \\
\hline
$6$ & $1$ & 6.72.1.1 \\
$6$ & $2$ & 6.36.0.1 \\
$6$ & $3$ & 6.24.0.1 \\
\hline
$7$ & $1$ & 7.56.1.1 \\
\hline
$8$ & $1$ & 8.96.3.16 \\
$8$ & $2$ & 8.24.0.67 \\
$8$ & $4$ & 8.12.0.5 \\
\hline
$9$ & $1$ & 9.108.4.4 \\
$9$ & $3$ & Not defined over $\bbQ$ \\
\hline
$10$ & $1$ & 10.180.7.1 \\
$10$ & $2$ & 10.90.2.1 \\
$10$ & $5$ & 10.36.1.1 \\
\end{tabular}
\quad \quad \quad
\rowcolors{2}{white}{gray!10}
\begin{tabular}{c c | c}
$m$ & $n$ & RSZB label \\
\hline
\hline
$12$ & $1$ & 12.288.13.3 \\ %
$12$ & $2$ & 12.72.1.1 \\ %
$12$ & $3$ & 12.96.3.4 \\ %
$12$ & $4$ & 12.72.1.1 \\ %
$12$ & $6$ & 12.24.0.3 \\ %
\hline
$13$ & $1$ & 13.182.8.1 \\
\hline
$16$ & $1$ & 16.384.21.47 \\ %
$16$ & $2$ & 16.96.3.226 \\ %
$16$ & $4$ & Not defined over $\bbQ$ \\
$16$ & $8$ & 16.24.0.2 \\ %
\hline
$18$ & $1$ & 18.648.37.11 \\ %
$18$ & $2$ & 18.324.16.5 \\ %
$18$ & $3$ & \multicolumn{1}{c}{Not defined over $\bbQ$} \\
$18$ & $6$ & \multicolumn{1}{c}{Not defined over $\bbQ$} \\
$18$ & $9$ & 18.72.1.1 \\ %
\hline
$25$ & $1$ & 25.750.48.1  \\
$25$ & $5$ &  Not defined over $\bbQ$ \\
\end{tabular}
}{Modular curves parameterizing elliptic curves with pairs of unsigned cyclic $m$-isogenies whose kernels have intersection of order $n$}

\begin{proof}[Proof of \Cref{Theorem: Parameterizing elliptic curves equipped with pairs of isogenies}]
    
    Let $m$ be as above, and suppose $E \in \scrE$ has two distinct unsigned cyclic $m$-isogenies with kernels $\Phi_1, \Phi_2 \subseteq E(\bbQalg)$. If we let $n \colonequals \#(\Phi_1 \cap \Phi_2)$, we must have $n \mid m$ by elementary group theory, and $n < m$ as otherwise $\Phi_1 = \Phi_2$ and our $m$-isogenies would not be distinct up to sign. We first handle the case where $m$ is a power of a prime, and afterwards we will address the case where $m$ is a product of powers of distinct primes.

    Suppose first that $m$ is a power of a prime. There are two subcases: $n = 1$ and $n > 1$. 
    
    In the first subcase, $\Phi_1 + \Phi_2 = E[m]$. As the absolute Galois group $\Gal_\bbQ$ stabilizes both $\Phi_1$ and $\Phi_2$, the modular curve $X_{\textup{sp}}(m)$ associated to the split Cartan group precisely parameterizes elliptic curves equipped with such a pair of unsigned cyclic isogenies. %

    On the other hand, if $n > 1$, write $n\prm = m/n$. Choose $P$ and $Q$ to be generators for $\Phi_1$ and $\Phi_2$ as abelian groups, so that $n\prm P = n\prm Q$ is a generator for $\Phi_1 \cap \Phi_2$. Choose $P_2 \in E[m](\bbQalg)$ so that $\set{P_1 = P, P_2}$ is a basis for $E[m](\bbQalg)$. We therefore have
    \begin{equation}
    Q = a P_1 + b P_2
    \end{equation}
    for some $a, b \in \bbZ$, and
    \begin{equation}
    n\prm P = n\prm Q =  a n\prm P_1 + b n\prm P_2. 
    \end{equation}
    As $m$ is a prime power, $a$ is a unit modulo $m$, so we may replace $Q$ with $a^{-1} Q$ and take $a = 1$. As $n\prm b \equiv 0 \pmod{m}$, $b$ must be a multiple of $n$ in $\bbZ / m \bbZ$. But if $b \bbZ / m \bbZ \subsetneq n \bbZ / m \bbZ$, then $\# \parent{\Phi_1 \cap \Phi_2} > n$, a contradiction. So $b= u n$ for some $u \in (\bbZ / m \bbZ)^\times$, and replacing $P_2$ with $u P_2$ if necessary, we may assume $b = n$.

    Now in the coordinates $\set{P_1, P_2}$, we have
    \begin{equation}
    P = \begin{pmatrix} 1 \\ 0 \end{pmatrix} \ \text{and} \ Q = \begin{pmatrix} 1 \\ n \end{pmatrix}.
    \end{equation}
    The image of $\Gal_\bbQ$ in $\GL_2(\bbZ / m \bbZ)$ preserves each of the cyclic subgroups generated by these two vectors, and hence is contained in the subgroup
    \begin{equation}\label{Equation: Galois image in GLm}
    H_{m,n} \colonequals \set{\begin{pmatrix}
        a & b \\ 0 & d
    \end{pmatrix} \in \GL_2(\bbZ / m \bbZ) : d \equiv a + n b \pmod{n\prm}}.
    \end{equation}
    Conversely, if the image of $\Gal_\bbQ$ on $E[m](\bbQalg)$ is contained in $H_{m,n}$ for some choice of basis $P_1, P_2$, then the cyclic subgroups generated by $P_1$ and $P_1+nP_2$ are Galois-stable and have intersection of order $n$. Thus $E$ has two unsigned cyclic $m$-isogenies with kernels having intersection of size $n$.

    Let us consider a product $m$ of distinct prime powers $p_1^{v_1},\dots,p_r^{v_r}$. To obtain modular curves that parameterize all elliptic curves with repeated cyclic $p_j^{v_j}$-isogenies, we apply the arguments above to each prime power $p_j^{v_j}$. We then take the fiber products of these modular curves over $X(1)$. However, we must also include the modular curve $X_0(p_j^{v_j})$ for each prime power $p_j^{v_j}$ in the fiber product (corresponding to the case where the $p_j$-part of the two cyclic isogenies is identical). Of course, the fiber product of all $X_0(p_j^{v_j})$ is $X_0(m)$, and we discard this case. 
\end{proof}

    A few remarks about Table \ref{table:modularcurves} are in order. 
    
    First, the modular curves for the cases
    \begin{equation}
    (m, n) \in \set{(9, 3), (16, 4), (18, 3), (18, 6), (25, 5)}.
    \end{equation}
    are not defined over $\bbQ$. Indeed, the Weil pairing tells us that the determinant map from the image of $\Gal_\bbQ$ in $\GL(E[m])$ to $(\bbZ / m \bbZ)^\times$ must be surjective (see \cite[Proposition 8.1 and Corollary 8.1.1]{Silverman}), but in each of these cases the determinant image of $H_{m,n}$ from \eqref{Equation: Galois image in GLm} is a proper subgroup of $(\bbZ/m\bbZ)^\times$. In particular, there are no elliptic curves $E/\bbQ$ possessing such a pair. Thus, for instance, there are no elliptic curves $E/\bbQ$ possessing two cyclic $9$-isogenies with kernels having an intersection of size $3$ (this possibility was also precluded by the isogeny graphs of \cite{Chiloyan-Lozano-Robledo}).
    
    Second, we note that when $m = 2^v > 2$ and $m/n = n\prm = 2$, we obtain the modular curve $X_0(m)$. In these cases, cyclic $m$-isogenies invariably come in pairs (see \cite[\S 3]{Pomerance-Schaefer} for a careful treatment of this observation in the case $m = 4$). This property is preserved by fiber products, so cyclic $12$-isogenies also come in pairs. This observation gives rise to the following corollary.

\begin{corollary}\label{Corollary: discrepancy between twNeq m and twNad m when 4 | m}
    Let $m \in \set{12, 16}$. Then we have
    \begin{equation}
    \twistNeq m (X) = 2 \twistNad m (X) + O(1)
    \end{equation}
    for $X \geq 1$.
\end{corollary}

\begin{proof}
    In these cases $Y_0(m)$ parameterizes both elliptic curves equipped with an unsigned cyclic $m$-isogeny, and elliptic curves equipped with a pair of unsigned cyclic $m$-isogenies whose kernels have intersection of order $m/2$. Thus the count of elliptic curves we obtain by substituting $t \in \bbQ$ into the Weierstrass equation
    \begin{equation}
    y^2 = x^3 + \f m (t) x + \g m (t)
    \end{equation}
    is off by a factor of $2$. For all divisors $n$ of $m$ with $1 \leq n < m/2$, the modular curve we obtain has nonzero genus. When the genus is greater than 1, Faltings's theorem (\Cref{Theorem: Faltings Theorem}) gives us a contribution of $O(1)$ to the difference $\twistNeq m (X) - 2 \twistNad m (X)$. Otherwise, when $m = 12$ and $n \in \set{2, 4}$, the associated modular curve $X_{\textup{sp}}(3) \times_{X(1)} X_0(4)$ has genus 1: by inspection, this modular curve is $\bbQ$-isomorphic to the elliptic curve
    \begin{equation}
    y^2 = x^3 + 1,
    \end{equation}
    which has Mordell-Weil group $\bbZ/6 \bbZ$, and again we get a contribution of $O(1)$ to $\twistNeq {12} (X) - 2 \twistNad {12} (X)$.
\end{proof}

    If
    \begin{equation}
    (m, n) \in \set{(2, 1), (3, 1), (4, 1), (5, 1), (6, 2), (6, 3), (8, 2)},
    \end{equation}
    the associated modular curve has genus $0$, but is not the moduli space of elliptic curves equipped with a cyclic $m$-isogeny. In these cases, we can compute universal families
    \begin{equation}
    E_{m, n}(t) : y^2 = x^3 + \f {m, n} (t) x + \g {m, n} (t)
    \end{equation}
    for elliptic curves (over $\bbQ$, up to quadratic twist) with this level structure, in the manner of \Cref{Lemma: parameterizing m-isogenies}. For later use, the polynomials $\f {m, n} (t)$ and $\g {m, n} (t)$ are recorded in Tables \ref{table:fmn} and \ref{table:gmn} below.

\jvtable{table:fmn}{
\rowcolors{2}{white}{gray!10}
\begin{tabular}{c c | c}
$m$ & $n$ & $\f {m, n} (t)$ \\
\hline\hline
$2$ & $1$ & $-(3 t^2 + 1)$ \\
$3$ & $1$ & $-3 t (2 + t) (4 - 2 t + t^2)$ \\
$4$ & $1$ & $-3 \left(t^4-2 t^3+2 t^2+2 t+1\right) \left(t^4+2 t^3+2 t^2-2 t+1\right)$ \\
$5$ & $1$ &  \begin{tabular}{@{}c@{}}$-3 \left(t^2+4\right) \left(t^2-3 t+1\right) \left(t^4-t^3+11 t^2+4 t+16\right)$ \\
$\left(t^4+4 t^3+11 t^2+14 t+31\right)$
\end{tabular} \\
$6$ & $2$ & $-3 \left(t^3-2\right) \left(t^3+6 t-2\right) \left(t^6-6 t^4-4 t^3+36 t^2+12 t+4\right)$ \\
$6$ & $3$ & $-3 \left(t^2+3\right) \left(t^6-15 t^4+75 t^2+3\right)$ \\
$8$ & $2$ & $-3 \left(t^4-8 t^3+2 t^2+8 t+1\right) \left(t^4+8 t^3+2 t^2-8 t+1\right)$ \\
\end{tabular} \\
}{$\f {m, n} (t)$ for $(m, n) \in \set{(2, 1), (3, 1), (4, 1), (5, 1), (6, 2), (6, 3), (8, 2)}$}

\jvtable{table:gmn}{
\rowcolors{2}{white}{gray!10}
\begin{tabular}{c c | c}
$m$ & $n$ & $\g {m, n} (t)$ \\
\hline\hline
$2$ & $1$ & $2 t (t - 1) (t + 1)$ \\
$3$ & $1$ & $2 (t^2 - 2 t - 2) (t^4 + 2 t^3 + 6 t^2 - 4 t + 4)$ \\
$4$ & $1$ & $2 \left(t^2-2 t-1\right) \left(t^2+2 t-1\right) \left(t^4+1\right) \left(t^4+6 t^2+1\right)$ \\
$5$ & $1$ & \begin{tabular}{@{}c@{}}
$2 \left(t^2+4\right)^2 \left(t^2+2 t-4\right) \left(t^4+3 t^2+1\right) \left(t^4-6 t^3+21 t^2-36 t+61\right)$ \\
$ \left(t^4+4 t^3+21 t^2+34 t+41\right)$
\end{tabular}\\
$6$ & $2$ & \begin{tabular}{@{}c@{}} $2 \left(t^2+2 t-2\right) \left(t^4-2 t^3-8 t-2\right) \left(t^4-2 t^3+6 t^2+4 t+4\right) $ \\
$\left(t^8+2 t^7+4 t^6-16 t^5-14 t^4+8 t^3+64 t^2-16 t+4\right)$ 
\end{tabular} \\
$6$ & $3$ & $2 \left(t^4-6 t^2-3\right) \left(t^4-6 t^2-24 t-3\right) \left(t^4-6 t^2+24 t-3\right)$ \\
$8$ & $2$ & $2 \left(t^2-2 t-1\right) \left(t^2+2 t-1\right) \left(t^8+132 t^6-250 t^4+132 t^2+1\right)$ \\
\end{tabular}
}{$\g {m, n} (t)$ for $(m, n) \in \set{(2, 1), (3, 1), (4, 1), (5, 1), (6, 2), (6, 3), (8, 2)}$}

\section{Lattices and the principle of Lipschitz}\label{Section: Lattices and the principle of Lipschitz}

In this section, we recall (a special case of) the Principle of Lipschitz, also known as Davenport's Lemma. We apply this result to estimate the number of Weierstrass equations of the form
\begin{equation}
E : y^2 = x^3 + \A m (a, b) x + \B m (a, b) 
\end{equation}
which satisfy $\rawheight(A, B) \leq X$. This count differs substantially from $\twistNeq m (X)$ because the pairs $(a, b)$ and $(d a, d b)$ give rise to the same elliptic curve, just with different models, and (not unrelatedly) $\rawheight(A, B)$ need not be the twist height of $E$. To obtain $\NQeq m (X)$ we will also need to sum over quadratic twists of these elliptic curves. Nevertheless, the estimates given in this section are essential building blocks in what follows.

\begin{theorem}[Principle of Lipschitz]\label{Theorem: Principle of Lipschitz}
	Let $\mathcal{R} \subseteq \bbR^2$ be a closed and bounded region, with rectifiable boundary $\partial \mathcal{R}$. Then we have
	\begin{equation}
	\#(\mathcal{R} \cap \bbZ^2) = \Area(\mathcal{R}) + O(\len(\partial \mathcal{R})).
	\end{equation}
	The implicit constant depends on the similarity class of $\mathcal{R}$, but not on its size, orientation, or position in the plane $\bbR^2$.
\end{theorem}

\begin{proof}
	See Davenport \cite{Davenport}.
\end{proof}

\begin{remark}
	Davenport's formulation of \Cref{Theorem: Principle of Lipschitz} was substantially stronger than what we have recorded. More precisely, he allowed for $\mathcal{R}$ to be a subset of $\bbR^n$, not just $\bbR^2$, he imposed weaker conditions on $\mathcal{R}$ than we have, and he made his error term explicit.
\end{remark}

\subsection*{Applying the Principle of Lipschitz}

Specializing to the case of interest, for $m \in \set{4, 5, 6, 7, 8, 9, 10, 12, 13, 16, 18, 25}$ and for $X > 0$, let
\begin{equation} \label{eqn: R(X)}
\calR m(X) \colonequals \set{(a, b) \in \bbR^2 :\rawheight(\A m(a, b), \B m(a, b)) \leq X, \ b \geq 0},
\end{equation}
and let 
\begin{equation}\label{Equation: Rm}
\R m \colonequals \Area(\calR m(1)).
\end{equation}

\begin{lemma}\label{Lemma: formula for R(X)}
	Let $m \in \set{4, 5, 6, 7, 8, 9, 10, 12, 13, 16, 18, 25}$. For $X > 0$, we have 
    \begin{equation}
    \Area (\calR m(X)) = \R m X^{1/\degB m}.
    \end{equation}
\end{lemma}

\begin{proof}
	Since $\f m(t)=\A m(t,1)$ and $\g m(t)=\B m(t,1)$ have no common real root, the region $\calR m(X)$ is compact \cite[Proof of Theorem 3.3.1, Step 2]{Cullinan-Kenney-Voight}. The homogeneity
	\begin{equation}
 	\rawheight(\A m(u a, ub), \B m(ua, ub)) = u^{2 \degB m} \rawheight(\A m(a, b), \B m(a, b))
	\end{equation}
	implies
	\begin{equation} 	\begin{aligned}
	\Area(\calR m(X)) &= \Area(\{(X^{1/{2 \degB m}} a, X^{1/{2 \degB m}} b) :(a, b) \in \calR m(1)\}) \\
	&= X^{1/\degB m} \Area(\calR m(1)) \\
	&= \R m X^{1/\degB m}
	\end{aligned}\end{equation} 
	as desired.
\end{proof}

We obtain the following corollary of \Cref{Theorem: Principle of Lipschitz}.

\begin{corollary}\label{Corollary: Estimates for lattice counts}
	Let $m \in \set{4, 5, 6, 7, 8, 9, 10, 12, 13, 16, 18, 25}$, and let $a_0, b_0, d \in \bbZ$ with $d \geq 1$. For $X > 0$, we have
	\begin{equation}\label{Equation: principle of Lipschitz congruent to a0 and b0}
	\#\{(a, b) \in \calR m(X) \cap \bbZ^2 : (a, b) \equiv (a_0, b_0) \psmod d\} = \frac{\R m X^{1/\degB m}}{d^2} + O\parent{\frac{X^{1/{2 \degB m}}}{d}}
	\end{equation}
	for $X, d \geq 1$. The implied constants are independent of $X$, $d$, $a_0,$ and $b_0$. In particular,
	\begin{equation}\label{Equation: calR m X principle of Lipschitz}
	\#(\calR m(X) \cap \bbZ^2) = \R m X^{1/\degB m} + O(X^{1/{2 \degB m}})
	\end{equation}
    for $X \geq 1$.
\end{corollary}

\begin{proof}
    We combine \Cref{Lemma: formula for R(X)} and \Cref{Theorem: Principle of Lipschitz} to obtain \eqref{Equation: principle of Lipschitz congruent to a0 and b0} for $X$ sufficiently large, say for $X \geq X_0$. But as the left-hand side of \eqref{Equation: principle of Lipschitz congruent to a0 and b0} is locally bounded, and $X^{1/{2 \degB m}}/d > 0$ for $X, d \geq 1$ we can choose $C$ large enough 
    that
    \begin{equation}
    \abs{\#\{(a, b) \in \calR m(X) \cap \bbZ^2 : (a, b) \equiv (a_0, b_0) \psmod d\} - \frac{\R m X^{1/\degB m}}{d^2}} \leq C \frac{X^{1/{2 \degB m}}}{d}
    \end{equation}
    for $1 \leq X \leq X_0$, and \eqref{Equation: principle of Lipschitz congruent to a0 and b0} holds for $X, d \geq 1$.
\end{proof}

In the proof of \Cref{Corollary: Estimates for lattice counts}, we turned an assertion of the form
\begin{equation}\label{Equation: f = g + O h for X large}
f(X) = g(X) + O(h(X)) \ \text{for} \ X \ \text{sufficiently large}
\end{equation}
into an assertion of the form
\begin{equation}\label{Equation: f = g + O h for X >= 1}
f(X) = g(X) + O(h(X)) \ \text{for} \ X \geq 1.
\end{equation}
This technique readily generalizes: if $f(X)$ is locally bounded and $h(X) = X^\alpha$, then \eqref{Equation: f = g + O h for X large} implies \eqref{Equation: f = g + O h for X >= 1} essentially by the argument given in the proof of \Cref{Corollary: Estimates for lattice counts}. Similarly, if $f(X)$ is locally bounded and $h(X) = X^\alpha \log^\beta X$, then \eqref{Equation: f = g + O h for X large} implies
\begin{equation}\label{Equation: f = g + O h for X >= 2}
f(X) = g(X) + O(h(X)) \ \text{for} \ X \geq 2.
\end{equation}
We shall use these observations without further comment throughout the remainder of the thesis.

\begin{remark}
    Huxley's work \cite{Huxley1} implies that if the boundary of $\mathcal{R}$ is defined by nonlinear polynomials, then the error term in \eqref{Equation: calR m X principle of Lipschitz} can be improved to $O(X^{1/{2 \degB m} - \delta})$ for some $\delta > 0$ (see \cite[page 7]{Pomerance-Schaefer}). Of course, this does not hold for $\calR m(X)$ as we have defined it, because the line $b = 0$ gives a lower boundary for $\calR m (X)$. However, if we drop the condition $b \geq 0$ from the definition of $\calR m (X)$, Huxley's result will apply.

    Under this convention, the points $(a, b)$ and $(-a, -b)$ correspond to the same elliptic curve, and some additional technical care is needed to attend to those points with $b = 0$. However, such an argument enables us to modestly improve error terms for $\cM m (X; e)$ in \Cref{Lemma: asymptotic for M(X; e) for m = 7}, \Cref{Lemma: asymptotic for M(X; e) for m = 5 and 10 and 25}, and \Cref{Lemma: asymptotic for M(X; e) for m = 13}, and consequently improve the asymptotic error in \Cref{Theorem: asymptotic for twN(X) for m = 7}, \Cref{Theorem: asymptotic for twN(X) for m = 10 and 25}, and \Cref{Theorem: asymptotic for twN(X) for m = 13}. For want of time and ease of exposition, we decline to pursue this insight further in this thesis.
\end{remark}

We now record graphs of the region $\R m (1)$ for 
\[
m \in \set{4, 5, 6, 7, 8, 9, 10, 12, 13, 16, 18, 25}.
\]
For each figure, the region graphed in blue is the set
\begin{equation}
\set{(a, b) \in \bbR^2 : 4 \abs{\A m (a, b)}^3 \leq 1, b \geq 0};
\end{equation}
similarly, the region graphed in red is the set
\begin{equation}
\set{(a, b) \in \bbR^2 : 27 \abs{\B m (a, b)}^2 \leq 1, b \geq 0}.
\end{equation}
By definition, $\calR m (1)$ is the intersection of these two regions.

\jvfigure{Figure: calR4(1)}{\includegraphics[scale=0.3]{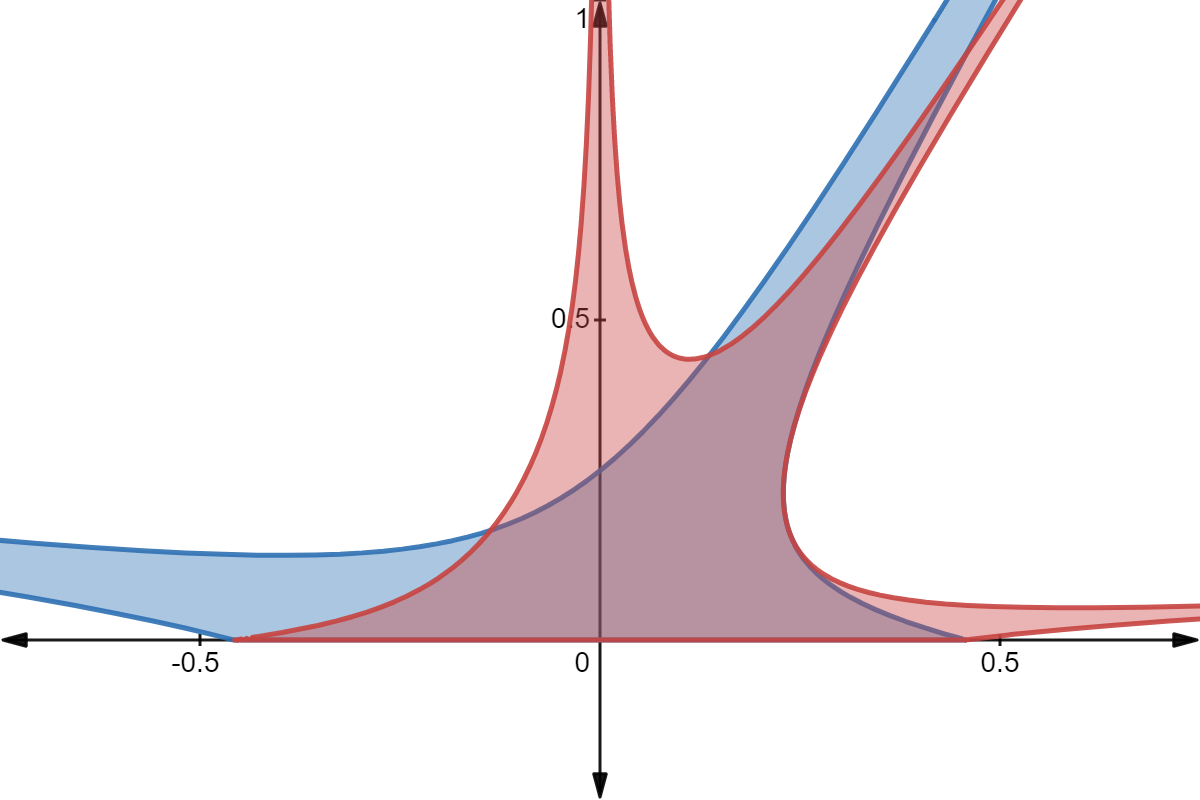}}{The region $\calR 4(1)$}

\jvfigure{Figure: calR5(1)}{\includegraphics[scale=0.3]{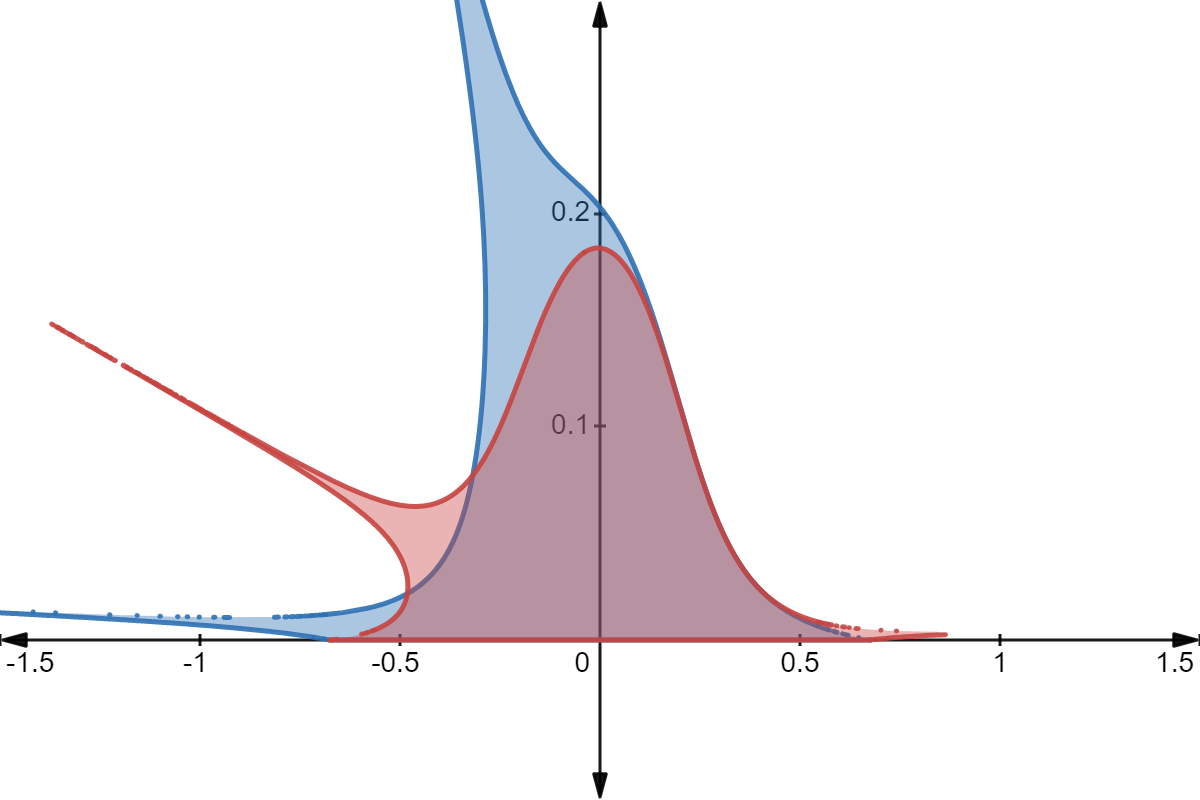}}{The region $\calR 5(1)$}

\jvfigure{Figure: calR6(1)}{\includegraphics[scale=0.3]{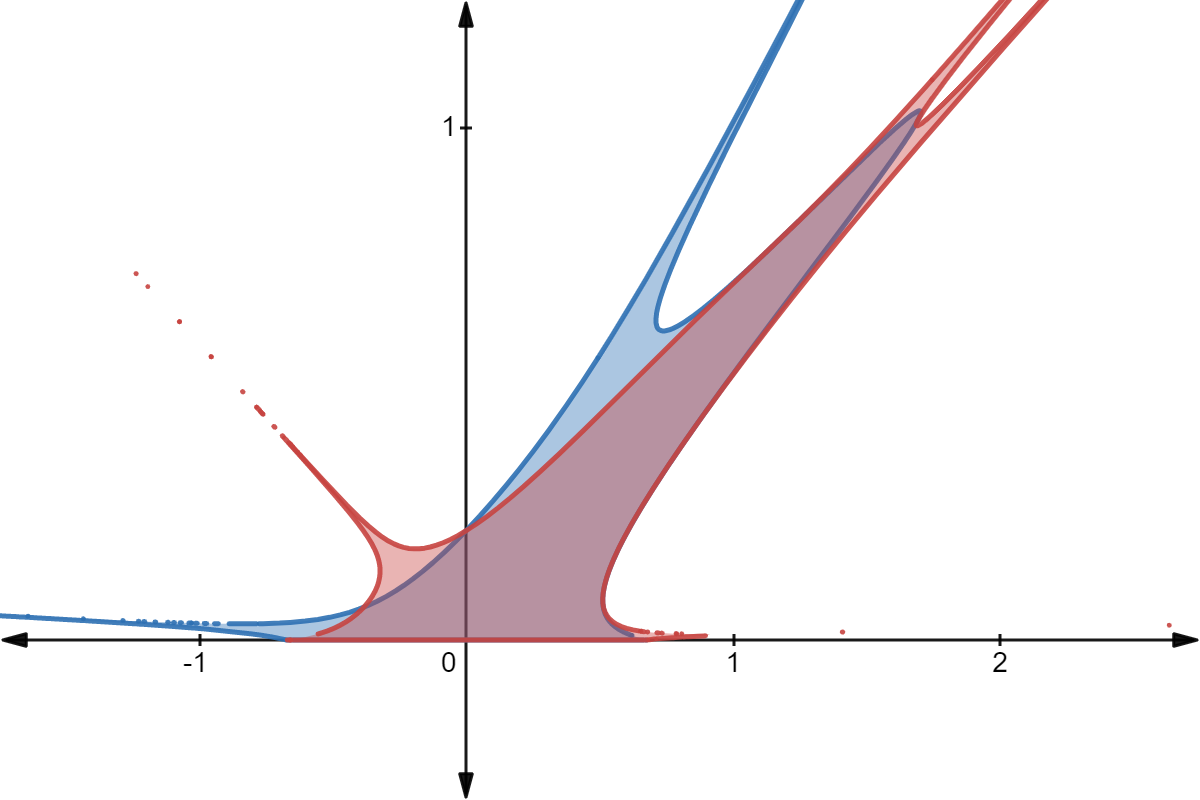}}{The region $\calR 6(1)$}

\jvfigure{Figure: calR7(1)}{\includegraphics[scale=0.3]{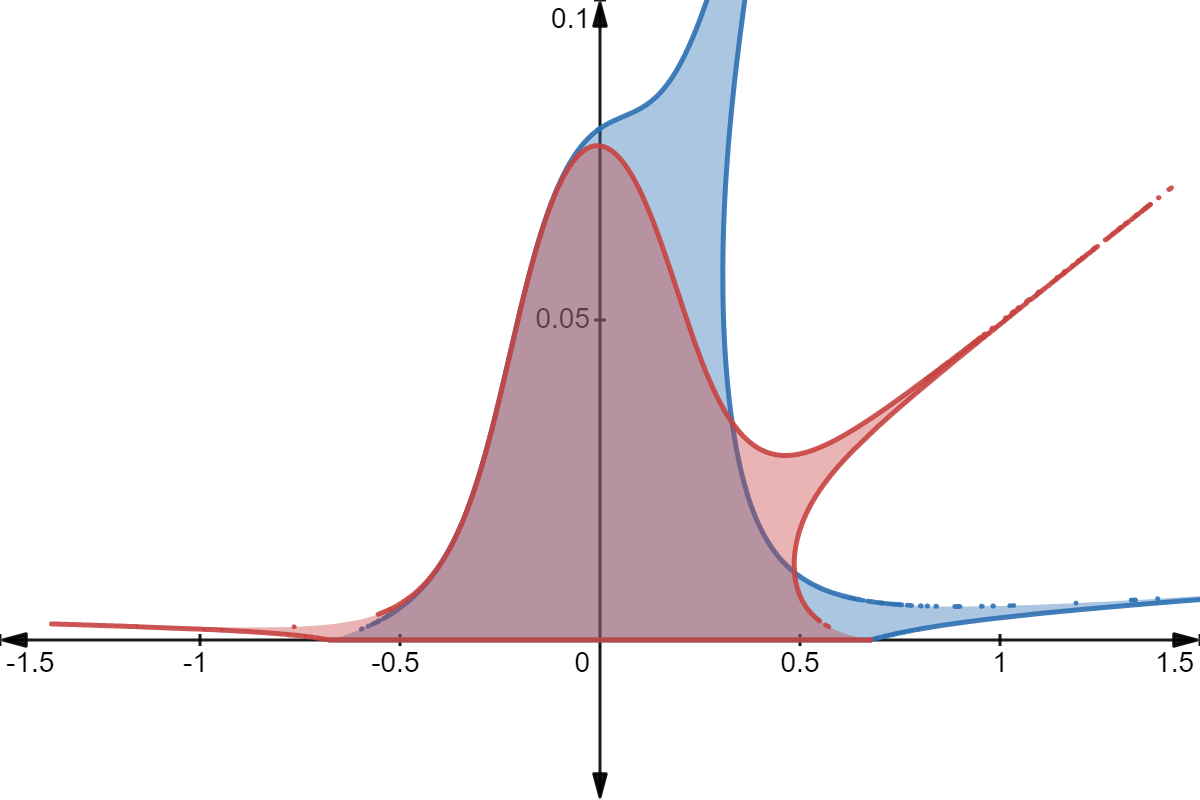}}{The region $\calR 7(1)$}

\jvfigure{Figure: calR8(1)}{\includegraphics[scale=0.3]{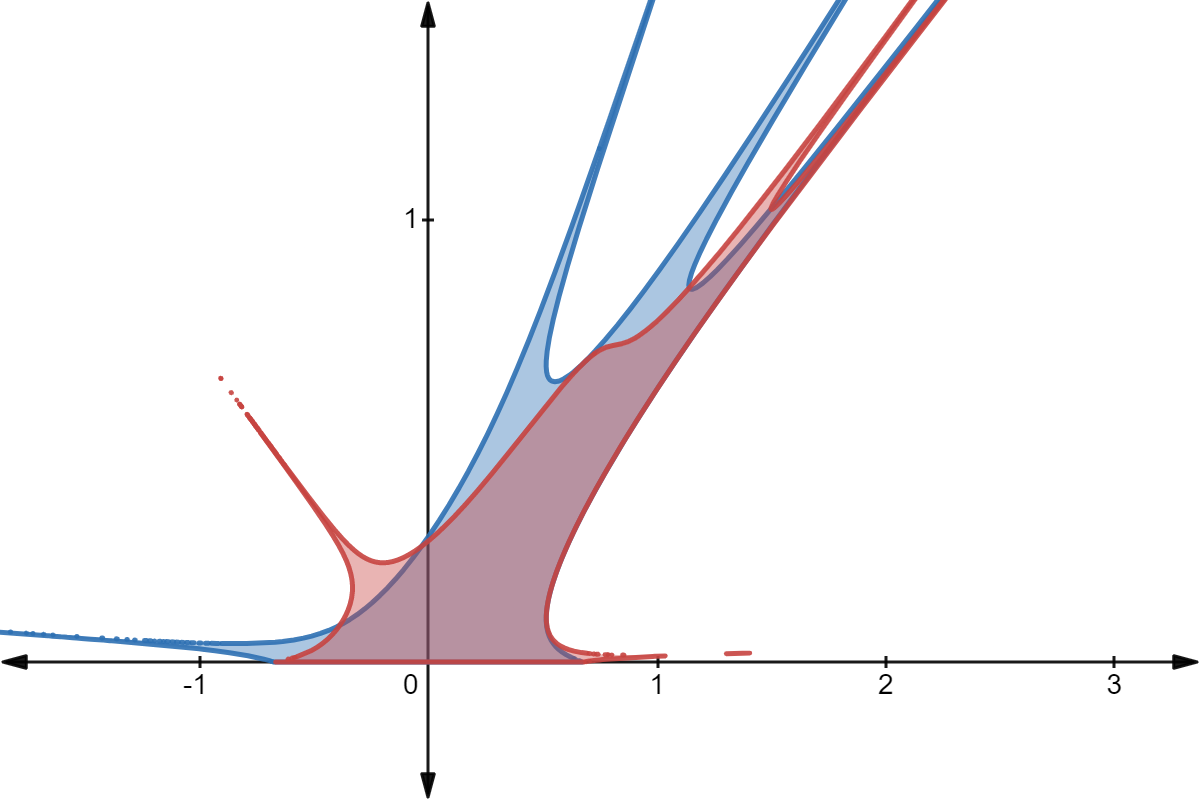}}{The region $\calR 8(1)$}

\jvfigure{Figure: calR9(1)}{\includegraphics[scale=0.3]{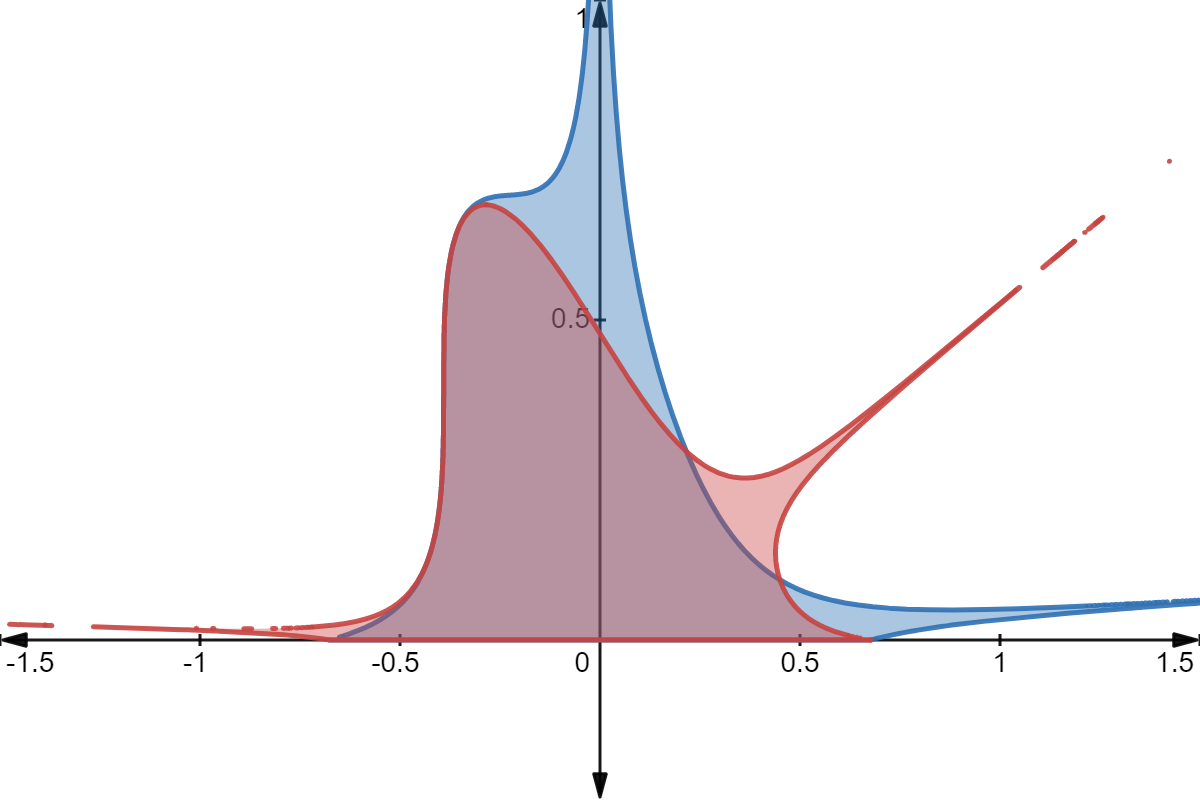}}{The region $\calR 9(1)$}

\jvfigure{Figure: calR10(1)}{\includegraphics[scale=0.3]{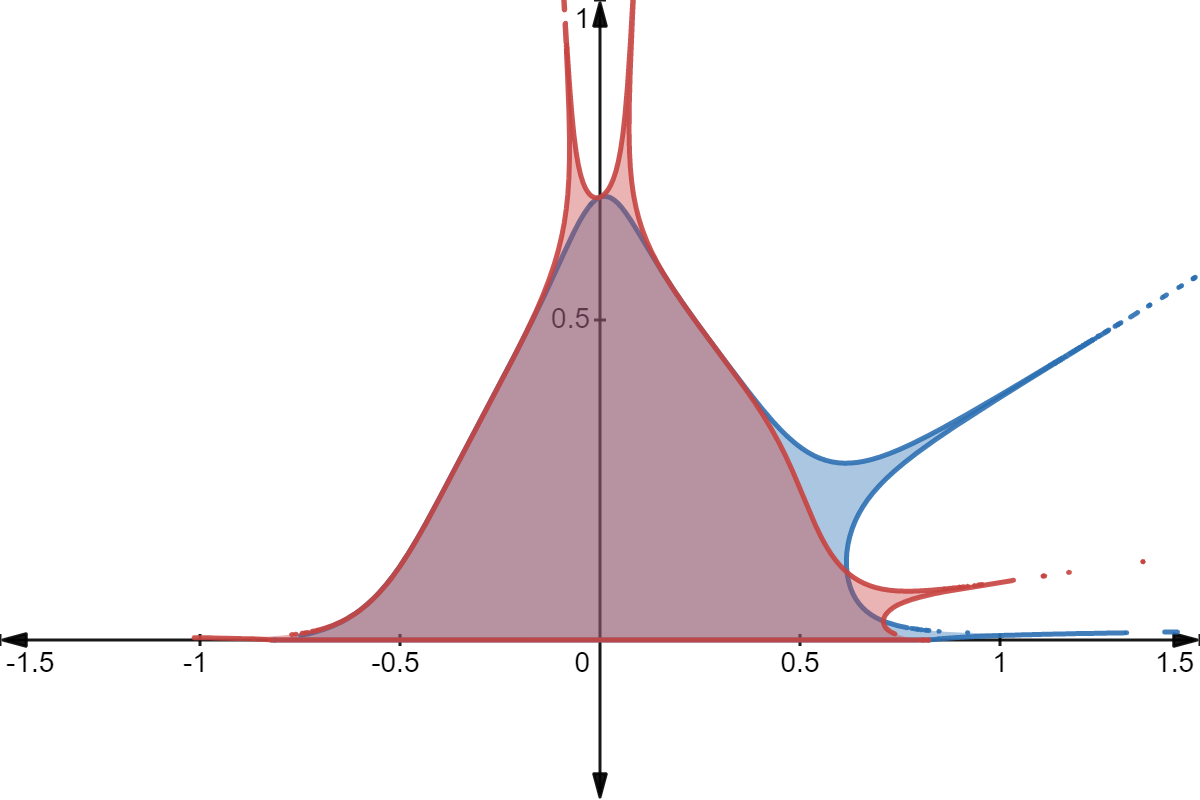}}{The region $\calR {10}(1)$}

\jvfigure{Figure: calR12(1)}{\includegraphics[scale=0.3]{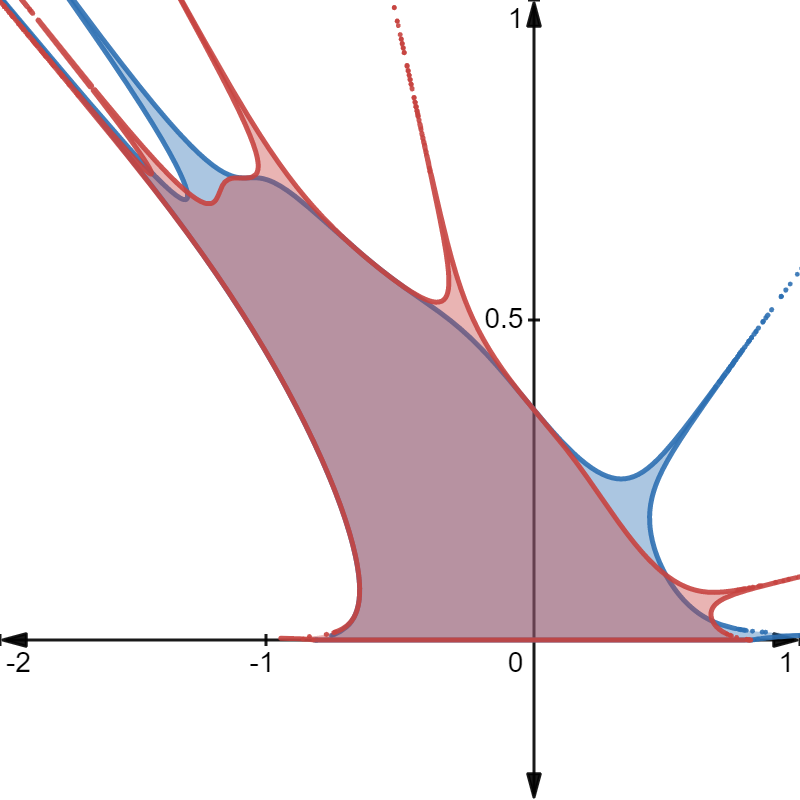}}{The region $\calR {12}(1)$}

\jvfigure{Figure: calR13(1)}{\includegraphics[scale=0.3]{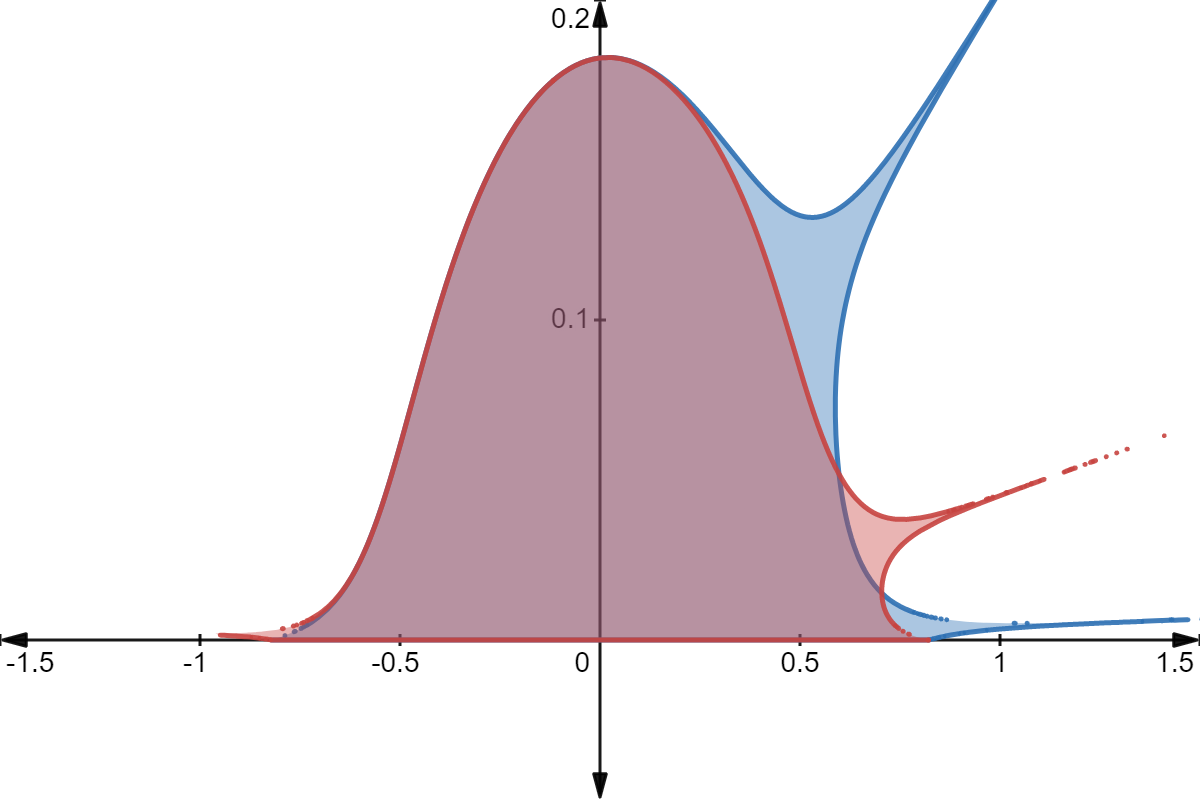}}{The region $\calR {13}(1)$}

\jvfigure{Figure: calR16(1)}{\includegraphics[scale=0.3]{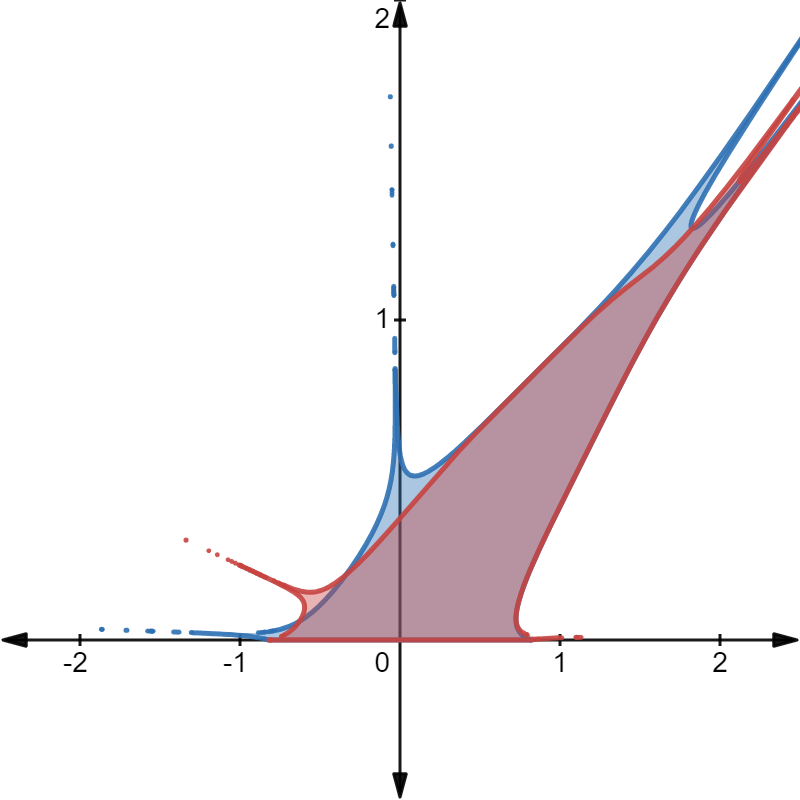}}{The region $\calR {16}(1)$}

\jvfigure{Figure: calR18(1)}{\includegraphics[scale=0.3]{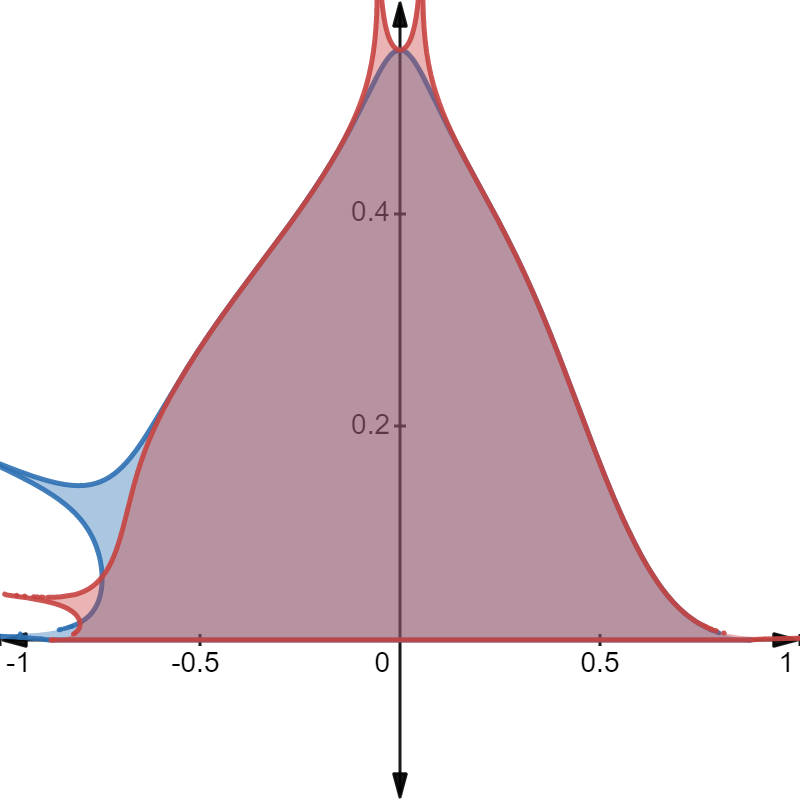}}{The region $\calR {18}(1)$}

\jvfigure{Figure: calR25(1)}{\includegraphics[scale=0.3]{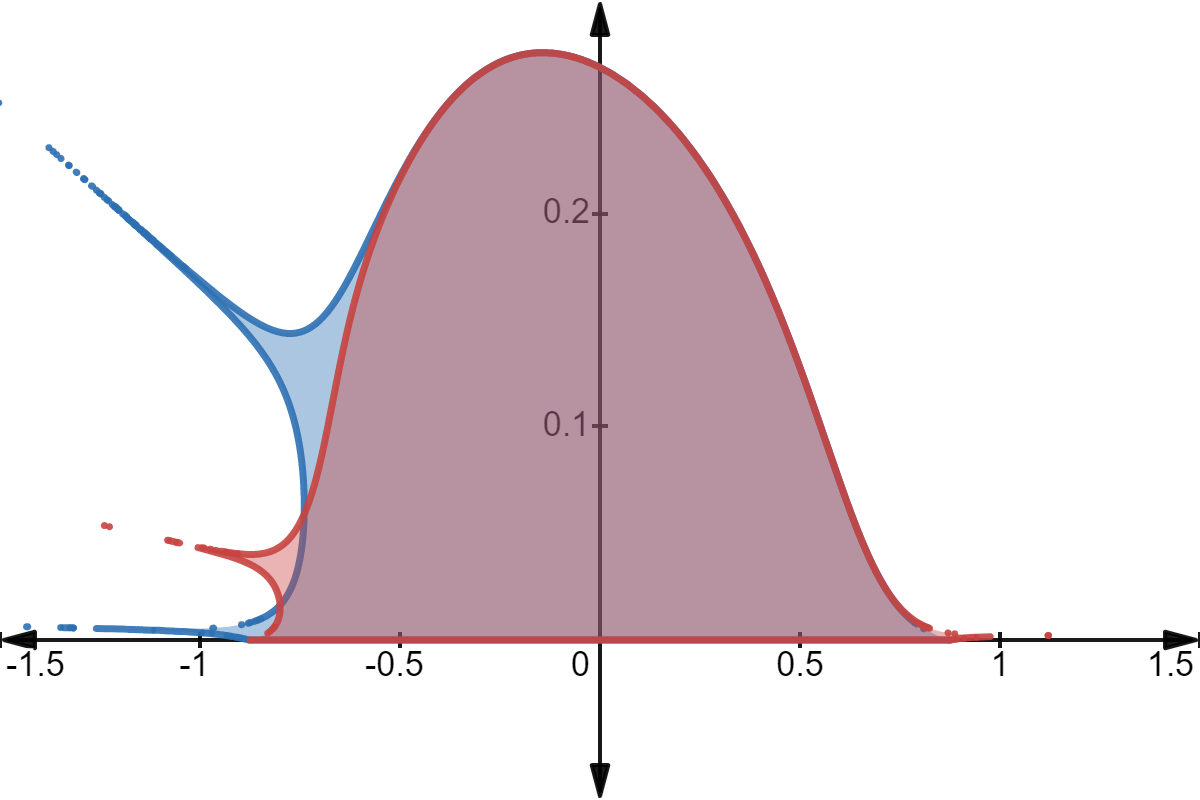}}{The region $\calR {25}(1)$}

\begin{remark}
    The region $\calR 4 (X)$, which we use to parameterize elliptic curves equipped with a cyclic $4$-isogeny, differs substantially from the region $\calR 1 \prm (X)$ from \cite[page 6]{Pomerance-Schaefer}, which Pomerance and Schaefer use to parameterize elliptic curves equipped with a pair of cyclic $4$-isogenies. 
    
    This discrepancy reflects a difference in derivation. Pomerance and Schaefer parameterize elliptic curves $E$ equipped with a pair of cyclic $4$-isogenies essentially by investigating the 2-division polynomial and the $4$-division polynomial for $E$ \cite[Proposition 3.2]{Pomerance-Schaefer}, and thus obtain the parameterization 
    \begin{equation}
    (A_4^{\textup{P-S}}(a, b), B_4^{\textup{P-S}}(a, b)) \colonequals (a^2 - 3 b^2, (a^2 - 2 b^2) b).
    \end{equation}
    We instead parameterize elliptic curves by using the Riemann--Roch theorem to establish an isomorphism from $X_0(4)$ to $\bbP^1$ (\Cref{Lemma: parameterizing m-isogenies}), and thus obtain the parameterization
    \begin{equation}
    (\A 4(a, b), \B 4(a, b)) = (-3(a^2 + 6 a b - 3 b^2), 2 a (a ^2 - 18 a b + 9 b^2)).
    \end{equation}
    This also means that for each elliptic curve $E$ (with $j(E) \neq 0$) equipped with a pair of associated cyclic 4-isogenies, we have two ($4$-groomed) ordered pairs $(a, b) \in \bbZ^2$ (one for each isogeny) where Pomerance and Schaefer only have one. For example, the elliptic curve
    \begin{equation}
    E: y^2 = x^3 - 2 x + 1
    \end{equation}
    arises via our method from the ordered pairs $(3, 1)$ and $(3, 5)$; on the other hand, for Pomerance and Schaefer, this elliptic curve arises from the ordered pair $(1, -1)$.
\end{remark}

\subsection*{Generalizing the Principle of Lipschitz}

As stated, \Cref{Theorem: Principle of Lipschitz} counts points in the lattice $\bbZ^2$, but we may easily extend this result to counts of points in more general lattices.

\begin{corollary}\label{Corollary: Principle of Lipschitz for general lattices}
	Let $\mathcal{R} \subseteq \bbR^2$ be a closed and bounded region, with rectifiable boundary $\partial \mathcal{R}$, let $\mathcal{L} \subseteq \bbR^2$ be a two-dimensional lattice, let $M : \bbR^2 \to \bbR^2$ be a linear map which induces a bijection between $\bbZ^2$ and $\mathcal{L}$, and let $\sigma(M)$ denote the smallest singular value of $M$. We have
	\begin{equation}
	\#(\mathcal{R} \cap \mathcal{L}) = \frac{\Area(\mathcal{R})}{\abs{\det M}} + O\parent{\frac{\len(\partial \mathcal{R})}{\sigma(M)}},
	\end{equation}
	where the implicit constant depends on the similarity class of $\mathcal{R}$, but not on its size, orientation, or position in the plane $\bbR^2$.
\end{corollary}

\begin{proof}
	Let $M : \bbR^2 \to \bbR^2$ be a linear map which induces a bijection from $\bbZ^2$ to $\mathcal{L}$; if $v_1, v_2$ are a $\bbZ$-basis for $\mathcal{L}$, we can take $M = (v_1 \ v_2 )$. By assumption, we have
	\begin{equation}
	\#(\mathcal{R} \cap \mathcal{L}) = \#(\mathcal{R} \cap M \bbZ^2) = \#(M\inv \mathcal{R} \cap \bbZ^2).
	\end{equation}
	By the variational characterization of the singular values of $M$ (\cite[Theorem 7.3.8]{Horn-Johnson}), $M\inv$ stretches $\partial \mathcal{R}$ by at most $1/\sigma(M)$. We now apply \Cref{Theorem: Principle of Lipschitz} to obtain our desired result. 
\end{proof}

\begin{remark}
	Note that $\abs{\det M}$ is the covolume of $\mathcal{L}$, and is independent of our choice of $M$. On the other hand, $\sigma(M)$ depends heavily on our choice of linear transform $M$. This is because shears preserve area but stretch lengths.
 
    For instance, if $\mathcal{L} = \bbZ^2$, then letting $M$ be the identity map recovers \Cref{Theorem: Principle of Lipschitz}. On the other hand, if we foolishly let
	\begin{equation}
	M = \begin{pmatrix} 1 & n \\ 0 & 1 \end{pmatrix},
	\end{equation}
	the error degrades by a factor of more than $n$. In Figure \ref{Figure: shear} below, both the square and the parallelogram are fundamental regions for $\bbZ^2$, but the square has a perimeter of $4$ whereas the parallelogram has a perimeter of $2 + 2 \sqrt{5} = 6.472\ldots > 4$.

 \jvfigure{Figure: shear}{\includegraphics[scale=0.3]{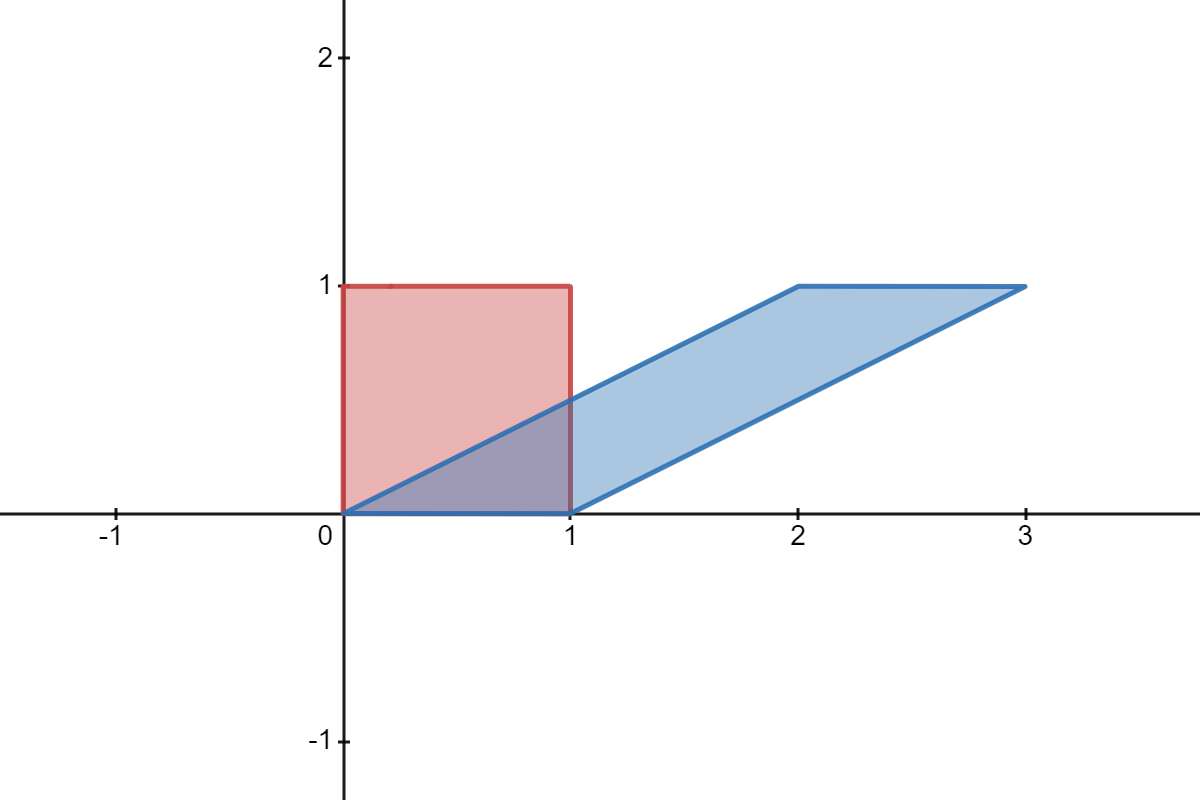}}{The unit square under the transformation $\begin{psmallmatrix} 1 & 2 \\ 0 & 1\end{psmallmatrix}$}
\end{remark}

\section{Analytic ingredients}\label{Section: Some analytic trivia}

In this section, we record several results from analytic number theory for later use. We begin by recording one half of Karamata's integral theorem for regularly varying functions. We then report estimates for counts of squarefree integers. We proceed to give several results from complex analysis about the convergence and growth rate of Dirichlet series. Finally, we record a Tauberian theorem which we will use to translate our estimates for $\twistNeq m (X)$ and $\twistNad m (X)$ to estimates for $\NQeq m (X)$ and $\NQad m (X)$.

\subsection*{Regularly varying functions}

We require a fragment of Karamata's integral theorem for regularly varying functions.

\begin{definition}
	Let $F \colon \bbR_{> 0} \to \bbR$ be measurable and eventually positive. We say that $F$ is \defi{regularly varying of index $\rho \in \bbR$} if for each $\lambda > 0$ we have
	\begin{equation}
	\lim_{y \to \infty} \frac{F(\lambda y)}{F(y)} = \lambda^\rho.
	\end{equation}
\end{definition}

\begin{theorem}[Karamata's integral theorem]\label{Theorem: Karamata's integral theorem}
	Let $F \colon \bbR_{> 0} \to \bbR$ be locally bounded and regularly varying of index $\rho$. Let $\sigma \in \bbR$. Then the following statements hold.
	\begin{enumalph}
		\item If $\sigma > \rho + 1$, then
		\begin{equation}
		\int_y^\infty u^{-\sigma} F(u) \,\mathrm{d}u \sim \frac{y^{1 - \sigma} F(y)}{\abs{\sigma - \rho - 1}}
		\end{equation}
		as $y \to \infty$.
		\item If $\sigma < \rho + 1$, then
		\begin{equation}
		\int_0^y u^{-\sigma} F(u) \,\mathrm{d}u \sim \frac{y^{1 - \sigma} F(y)}{\abs{\sigma - \rho - 1}}
		\end{equation}
	as $y \to \infty$.
	\end{enumalph}
\end{theorem}

\begin{proof}
	See Bingham--Goldie--Teugels \cite[Theorem 1.5.11]{Bingham-Goldie-Teugels}. (Karamata's integral theorem also includes a converse.)
\end{proof}

\begin{corollary}\label{Corollary: tail of sum of f/n^sigma}
	Let $\alpha \colon \bbZ_{> 0} \to \bbR$ be an arithmetic function, and suppose that for some $\kappa, \rho, \tau \in \bbR$ with $\kappa \neq 0$, we have
	\begin{equation}
	F(y) \colonequals \sum_{n \leq y} \alpha(n) \sim \kappa y^\rho \log^\tau y\label{Equation: asymptotic for F in corollary of Karamata's Theorem}
	\end{equation}
	as $y \to \infty$. Let $\sigma, \rho > 0$.
	Then the following statements hold, as $y \to \infty$.
	\begin{enumalph}
	\item If $\sigma > \rho > 0$, then
	\begin{equation}
	\sum_{n > y} n^{-\sigma} \alpha(n) \sim \frac{\rho y^{-\sigma} F(y)}{\abs{\sigma - \rho}} \sim \frac{\kappa \rho y^{\rho - \sigma} \log^\tau y}{\abs{\sigma - \rho}}.
	\end{equation}
	\item If $\rho > \sigma > 0$, then \begin{equation}
	\sum_{n \leq y} n^{-\sigma} \alpha(n) \sim \frac{\rho y^{-\sigma} F(y)}{\abs{\sigma - \rho}} \sim \frac{\kappa \rho y^{\rho - \sigma} \log^\tau y}{\abs{\sigma - \rho}}.
	\end{equation}
	\end{enumalph}
\end{corollary}

\begin{proof}
	Replacing $\alpha$ and $F$ with $-\alpha$ and $-F$ if necessary, we may assume $\kappa > 0$. As a partial sum of an arithmetic function, $F(y)$ is measurable and locally bounded; by \eqref{Equation: asymptotic for F in corollary of Karamata's Theorem}, $F(y)$ is eventually positive. Now for any $\lambda > 0$, we compute
	\begin{equation}
	\lim_{y \to \infty} \frac{F(\lambda y)}{F(y)} = \lim_{y \to \infty} \frac{\kappa (\lambda y)^\rho \log^\tau (\lambda y)}{\kappa y^\rho \log^\tau y} = \lambda^\rho,
	\end{equation}
	so $F$ is regularly varying of index $\rho$.

	Suppose first $\sigma > \rho > 0$. Since
	\begin{equation}
	y^{-\sigma} F(y) \sim \kappa y^{\rho - \sigma} \log^\tau y \to 0
	\end{equation}
	as $y \to \infty$, Abel summation yields
	\begin{equation}
		\sum_{n > y} n^{-\sigma} \alpha(n) = - y^{-\sigma} F(y) + \sigma \int_y^\infty u^{-\sigma - 1} F(u) \,\mathrm{d}u.
	\end{equation}
	Clearly $\sigma + 1 > \rho + 1$, so \Cref{Theorem: Karamata's integral theorem}(a) tells us 
	\begin{equation}
	\int_y^\infty u^{-\sigma - 1} F(u) \,\mathrm{d}u \sim \frac{y^{-\sigma} F(y)}{\abs{\sigma - \rho}} \sim \frac{\kappa y^{\rho - \sigma} \log^\tau y}{\abs{\sigma - \rho}}
	\end{equation}
	and thus
	\begin{equation}
		\sum_{n > y} n^{-\sigma} \alpha(n) \sim \frac{\rho y^{-\sigma} F(y)}{\abs{\sigma - \rho}}
	\end{equation}
	as $y \to \infty$.
	
	The case $\rho > \sigma > 0$ is similar.
\end{proof}

\subsection*{Counting squarefree integers}\label{Subsection: counting squarefree integers}

    Recall from \eqref{Equation: Definition of squarefree count} that for $X > 0$ a real number, we have
    \begin{equation}
	\kfree 2(X) \colonequals \# \set{n \in \bbZ_{>0} : n \leq X, \ n \ \text{a squarefree integer}}.
    \end{equation}
    In this subsection, we record estimates for $\kfree 2 (X)$ due to Walfisz \cite{Walfisz} and Liu \cite{Liu}. These estimates are obtained by examining the zero-free region of the Riemann zeta function.

\begin{theorem}\label{Theorem: Estimating squarefree(X)}
	Let $X > 0$ be a real number. Then for some constant $\kappa > 0$, we have
	\begin{equation}
	\kfree 2(X) = \frac{X}{\zeta(2)} + O\parent{X^{1/2} e^{-\kappa \frac{\log^{3/5} X}{\log^{1/5} \log X}}}
	\end{equation}
	for $X \geq 1$.
\end{theorem}

\begin{proof}
	Walfisz \cite[Satz V.6.1]{Walfisz} proves a stronger result.
\end{proof}

If the Riemann hypothesis holds, we can say substantially more about $\kfree 2 (X)$.

\begin{theorem}\label{Theorem: Estimating squarefree(X) given RH}
	Let $X > 0$ be a real number. If the Riemann hypothesis holds, then for any $\epsilon > 0$, we have
	\begin{equation}
	\kfree 2(X) = \frac{X}{\zeta(2)} + O\parent{X^{11/35 + \epsilon}}
	\end{equation}
	for $X \geq 1$. The implicit constant depends on $\epsilon$.
\end{theorem}

\begin{proof}
	Liu \cite[Theorem 1]{Liu}.
\end{proof}

\subsection*{Dirichlet series}

In this subsection, we record several analytic results about Dirichlet series.

The following theorem is attributed to Stieltjes.

\begin{theorem}\label{Theorem: product of Dirichlet series converges}
	Let $\alpha, \beta : \bbZ_{> 0} \to \bbR$ be arithmetic functions. If the Dirichlet series $L_\alpha(s) \colonequals \sum_{n \geq 1} \alpha(n) n^{-s}$ and $L_\beta(s) \colonequals \sum_{n \geq 1} \beta(n) n^{-s}$ both converge for $s = \sigma + i t$ with $\sigma > \sigma_0$, and one of these two series converges absolutely, then
	\begin{equation}
	L_{\alpha \ast \beta}(s) \colonequals \sum_{n \geq 1} \parent{\sum_{d \mid n} \alpha(d) \beta\parent{\frac nd}} n^{-s}
	\end{equation}
	converges for $s = \sigma + i t$ with $\sigma > \sigma_0$. If both $L_\alpha(s)$ and $L_\beta(s)$ converge absolutely when $\sigma > \sigma_0$, then so does $L_{\alpha \ast \beta}(s)$.
\end{theorem}

\begin{proof}
	Widder \cite[Theorems 11.5 and 11.6b]{Widder} proves a more general result, or see Tenenbaum \cite[proof of Theorem II.1.2, Notes on p.~204]{Tenenbaum}.
\end{proof}

Let 
\begin{equation}\label{Equation: Euler-Mascheroni constant}
\gamma \colonequals \lim_{y \to \infty} \bigl(\sum_{n \leq y} 1/n - \log y \bigr)
\end{equation}
be the Euler--Mascheroni constant. 

\begin{theorem}\label{Theorem: Laurent series expansion for zeta(s)}
	The difference
	\begin{equation}
	\zeta(s) - \parent{\frac{1}{s - 1} + \gamma}
	\end{equation}
	is entire on $\bbC$ and vanishes at $s = 1$. 
\end{theorem}

\begin{proof}
	Ivi\'c \cite[page 4]{Ivic} proves a more general result.
\end{proof}

Recall that a complex function $F(s)$ has \defi{finite order} on a domain $D \subseteq \bbC$ if there exists $\xi \in \bbR_{>0}$ such that 
	\begin{equation}
	F(\sigma \pm i t) = O(1 + \abs{t}^\xi)
	\end{equation}
	whenever $\sigma \pm i t \in D$. We emphasize that this is a statement about $F(\sigma \pm i t)$ for $t$ large, not about $F$ on or close to the real line. If $F$ is of finite order on a right half-plane, we define
	\begin{equation}
	\mu_F(\sigma) \colonequals \inf\{\xi \in \bbR_{\geq 0} :F(\sigma + i t) = O(1 + \abs{t}^\xi)\}
	\end{equation}
	as $t \to \infty$, where the implicit constant depends on $\sigma$ and $\xi$.

\begin{proposition}\label{Proposition: mu acts like a negative valuation}
	Let $F$ and $G$ be complex functions of finite order on a right half-plane. We have
	\begin{equation}
	\mu_{F + G}(\sigma) \leq \max(\mu_F(\sigma), \mu_G(\sigma)),
	\end{equation}
	and
	\begin{equation}
	\mu_{FG}(\sigma) \leq \mu_F(\sigma) + \mu_G(\sigma).
	\end{equation}
\end{proposition}

\begin{proof}
    Immediate from the definition of $\mu_F(\sigma)$.
\end{proof}

\begin{theorem}\label{Theorem: absolutely convergent Dirichlet series have mu = 0}
    Let $L(s)$ be a Dirichlet series with abscissa of absolute convergence $\sigma_a$. Then we have $\mu_L(\sigma)=0$ for all $\sigma > \sigma_a$, and $\mu_L(\sigma)$ is nonincreasing (as a function of $\sigma$) on any region where $L$ has finite order.
\end{theorem}

\begin{proof}
	Tenenbaum \cite[Theorem II.1.21]{Tenenbaum}.
\end{proof}

\begin{theorem}\label{Theorem: muzeta(sigma)}
	Let $\zeta(s)$ be the Riemann zeta function, and let $\sigma \in \bbR$. Then we have
	\begin{equation}
	\mu_\zeta(\sigma) \leq \begin{cases}
	\frac 12 - \sigma, & \text{if} \ \sigma \leq 0; \\
	\frac 12 - \frac{29}{45} \sigma, & \text{if} \ 0 \leq \sigma \leq \frac 12; \\
	\frac{13}{42}(1 - \sigma), & \text{if} \ \frac 12 \leq \sigma \leq 1; \\
	0 & \text{if} \ \sigma \geq 1.
	\end{cases}
	\end{equation}
	Moreover, equality holds if $\sigma<0$ or $\sigma>1$.
\end{theorem}

\begin{proof}
	Tenenbaum \cite[page 235]{Tenenbaum} proves the claim when $\sigma<0$ or $\sigma>1$. Now 
	\begin{equation}
	\mu_\zeta(1/2) \leq 13/84 \label{Equation: Bourgain's bound on muzeta}
	\end{equation}
	by Bourgain \cite[Theorem 5]{Bourgain}, and our result follows from the subconvexity of $\mu_\zeta$ \cite[Theorem II.1.20]{Tenenbaum}.
\end{proof}

\begin{remark}
	We use \Cref{Theorem: muzeta(sigma)} as an input to \Cref{Theorem: Landau's Tauberian theorem}. Although we have stated it in the strongest form we know, for our applications, we could replace \eqref{Equation: Bourgain's bound on muzeta} with the much weaker statement
	\begin{equation}
	\mu_\zeta(1/2) < 1/2.
	\end{equation}
	However, appropriate refinements to \Cref{Theorem: Landau's Tauberian theorem} might enable us to leverage the full strength of \Cref{Theorem: muzeta(sigma)} (see \Cref{Remark: true bound for twistN(X) for m = 7} below).
\end{remark}

\subsection*{A Tauberian theorem}\label{Subsection: Perron's formula}

We now present a Tauberian theorem, due in essence to Landau \cite{Landau1915}, and in this formulation to Roux \cite{Roux}.

\begin{definition}\label{Definition: admissible sequences}
	Let $\parent{\alpha(n)}_{n \geq 1}$ be a sequence with $\alpha(n) \in \bbR_{\geq 0}$ for all $n$, and let $L_\alpha(s) \colonequals \sum_{n \geq 1} \alpha(n) n^{-s}$.
	We say the sequence $\parent{\alpha(n)}_{n \geq 1}$ is \defi{admissible} with (real) parameters $\parent{\sigma_a, \delta, \xi}$ if the following hypotheses hold:
	\begin{enumroman}
		\item $L_\alpha(s)$ has abscissa of absolute convergence $\sigma_a$.\label{Condition: abscissa of absolute convergence sigma}
		\item The function $L_\alpha(s)/s$ has meromorphic continuation to 
  \begin{equation}
      \set{s = \sigma + i t \in \bbC : \sigma > \sigma_a - \delta},
  \end{equation}
  and only finitely many poles in this region. \label{Condition: meromorphic extension}
		\item For $\sigma > \sigma_a - \delta$, we have $\mu_{L_\alpha}(\sigma) \leq \xi$.
		\label{Condition: mu is bounded}
	\end{enumroman}
\end{definition}

	If $\parent{\alpha(n)}_n$ is admissible, let $s_1, \dots, s_r$ denote the poles of $L_\alpha(s)/s$ with real part greater than $\sigma_a - \delta/(\floor{\xi} + 2)$.

The following theorem is essentially an application of Perron's formula \cite[\S II.2.1]{Tenenbaum}, which is itself an inverse Mellin transform.

\begin{theorem}[Landau's Tauberian Theorem]\label{Theorem: Landau's Tauberian theorem}
	Let $\parent{\alpha(n)}_{n \geq 1}$ be an admissible sequence (\textup{\Cref{Definition: admissible sequences}}), and write $N_\alpha(X) \colonequals \sum_{n \leq X} \alpha(n)$. Then for all $\epsilon>0$,
		\begin{equation}
		N_\alpha(X) = \sum_{j = 1}^r \res_{s=s_j}\parent{\frac{L_\alpha(s) X^{s}}{s}} + O\!\parent{X^{\sigma_a - \frac{\delta}{\floor{\xi} + 2} + \epsilon}},
		\end{equation}
		as $X \to \infty$, where the main term is a sum of residues. The implicit constant depends on $\epsilon$.
\end{theorem}

\begin{proof}
	See Roux \cite[Theorem 13.3, Remark 13.4]{Roux}.
\end{proof}

\begin{remark}
	Landau's original theorem \cite{Landau1915} was fitted to a more general context, and allowed sums of the form
	\begin{equation}
	\sum_{n \geq 1} \alpha(n) \ell(n)^{-s}
	\end{equation}
	as long as $\parent{\ell(n)}_{n \geq 1}$ was increasing and tended to $\infty$. However, Landau also required that $L_\alpha(s)$ has a meromorphic continuation to all of $\bbC$, and Roux \cite[Theorem 13.3, Remark 13.4]{Roux} relaxes this assumption. Both Landau and Roux wrote out the expression
	\begin{equation}
	\res_{s=s_j}\parent{\frac{L_\alpha(s) X^{s}}{s}}
	\end{equation}
	in terms of the Laurent series expansion for $L_\alpha(s)$ around $s = s_j$, but we believe expressing the result in terms of residues is more transparent.
\end{remark}

We now illustrate the applicability and power of \Cref{Theorem: Landau's Tauberian theorem}. Let $\omega(n)$ denote the number of distinct prime divisors of $n$, and recall the definition of the Euler--Mascheroni constant $\gamma$ \eqref{Equation: Euler-Mascheroni constant}.

\begin{corollary}\label{Corollary: sum of 2^omega(n)}
	For any $\epsilon > 0$, we have
	\begin{equation} 	\begin{aligned}
	\sum_{n \leq y} 2^{\omega(n)} &= \frac{1}{\zeta(2)}y \log y + \frac{(2 \gamma - 1) \zeta(2) - 2 \zeta\prm(2)}{\zeta(2)^2} y + O\parent{y^{3/4 + \epsilon}} \\
	&\approx 0.607927 y \log y + 0.786872 y + O\parent{y^{3/4 + \epsilon}}
	\end{aligned}\end{equation} 
	for $y \geq 1$. The implicit constant depends only on $\epsilon$.
\end{corollary}

\begin{proof}
Recall that
\begin{equation}
L_{2^\omega}(s) \colonequals \frac{\zeta(s)^2}{\zeta(2s)} = \sum_{n \geq 1} \frac{2^{\omega(n)}}{n^s}.
\end{equation}
For $\sigma > 1/2$, this function is holomorphic except at $s = 1$, where it has a pole of order 2. Write $\zeta_a(s) \colonequals \zeta(a s)$. By \Cref{Theorem: absolutely convergent Dirichlet series have mu = 0}, $\mu_{1/\zeta_2}(\sigma) = 0$ for $\sigma > \frac 12$, and by \Cref{Theorem: muzeta(sigma)}, $\mu_{\zeta}(\sigma) \leq \frac{13}{42}\parent{1 - \sigma}$ for $\frac 12 \leq \sigma \leq 1$. Thus for any $\epsilon > 0$ and all $\sigma > \frac 12 + \epsilon$, \Cref{Proposition: mu acts like a negative valuation} tells us 
\begin{equation}
\mu_{L_{2^\omega}}(\sigma) \leq \frac{2 \cdot 13}{84} + 0 = \frac {13}{42}.
\end{equation}
Consequently, the sequence $\parent{2^{\omega(n)}}_{n \geq 1}$ is admissible with parameters $\parent{1, 1/2, 13/42}$.

We compute the residue of $L_{2^\omega}(s) \cdot \frac{y^s}{s}$ at $s = 1$, which is
\begin{equation}
	\frac{1}{\zeta(2)}y \log y + \frac{(2 \gamma - 1) \zeta(2) - 2 \zeta\prm(2)}{\zeta(2)^2} y%
\end{equation}
Applying Theorem \ref{Theorem: Landau's Tauberian theorem}, we conclude
\begin{equation}
\sum_{n \leq y} 2^{\omega(n)} = \frac{1}{\zeta(2)}y \log y + \frac{(2 \gamma - 1) \zeta(2) - 2 \zeta\prm(2)}{\zeta(2)^2} y + O\parent{y^{3/4 + \epsilon}}
\end{equation}
for any $\epsilon > 0$.
\end{proof}

\begin{remark}
	Better estimates for $\sum_{n \leq y} 2^{\omega(n)}$ are possible (see \cite[Exercise I.3.54]{Tenenbaum}), but we shall not require them.
\end{remark}

Although we do not require the strength of the asymptotics of \Cref{Corollary: sum of 2^omega(n)}, it will be useful to have an order of growth for $\kappa^{\omega(n)}$ in greater generality.

\begin{theorem}\label{Theorem: Asymptotics of k^omega(n)}
	Let $\kappa > 0$ be a real number. Then we have
	\begin{equation}
	\sum_{n \leq y} \kappa^{\omega(n)} \asymp y \log^{\kappa - 1} y.
	\end{equation}
\end{theorem}

\begin{proof}
	Ivi\'c proves a stronger result \cite[Theorem 14.10]{Ivic}.
\end{proof}

\section{Our approach revisited}\label{Section: Our approach revisited}

In this section, we elaborate on the intuitions laid out in \cref{Section: Our approach} and set up a general framework that will enable us to determine the asymptotic number of elliptic curves with a cyclic $m$-isogeny over $\bbQ$, up to both quadratic twist and $\bbQ$-isomorphism, when $m \in \set{7, 10, 13, 25}$. Our method also applies even more easily when $m \in \set{4,6,8,9,12,16,18}$. When $X_0(m)$ is of nonzero genus, i.e., when $m \in \set{11,14,15,17,19,21,27,37,43,67,163}$, there are four or fewer elliptic curves in $\twistE$ with a cyclic $m$-isogeny, so there is no need to sieve to obtain $\twistNeq m (X) = \twistNad m (X)$ (see \cref{Chapter: m of nonzero genus} for the asymptotics of $\twistNad m (X)$ and $\NQad m (X)$ for these $m$).

Choose $m$ so that $X_0(m)$ is of genus $0$ and $m > 5$ or $m = 4$. We now elaborate on the steps described in \cref{Section: Our approach}. The first step, obtaining a parameterization for the family of elliptic curves with a cyclic $m$-isogeny up to quadratic twist, was completed in \cref{Section: Parameterizing elliptic curves with a cyclic m-isogeny}. The second step, estimating
\begin{equation}
    \#\{(a, b) \in \calR m(X) \cap \bbZ^2 : (a, b) \equiv (a_0, b_0) \psmod d\}
\end{equation}
for $a_0, b_0, d \in \bbZ$, was completed in \cref{Section: Lattices and the principle of Lipschitz}.

We erect a more abstract framework for addressing the third step and the fifth step of our approach (the fourth step is addressed by \Cref{Corollary: Elliptic curves for which twNeq(X) = twNad(X)} and \Cref{Lemma: Difference between twNeq and twNad}).

\subsection*{Estimating $\twNad \calE (X)$}

In this subsection, we abstract the strategy developed in \cite[section 4.1]{Molnar-Voight} and formulate it in a way that will enable us to estimate $\twNeq m (X)$ and $\twNad m (X)$ for 
\begin{equation}
m \in \set{4, 6, 7, 8, 9, 10, 12, 13, 16, 18, 25}.
\end{equation}

We establish some notation for brevity and ease of exposition. Suppose we have a sequence of real-valued functions $\parent{\alpha(X; n)}_{n \geq 1}$, and an eventually positive function $\phi : \bbR_{>0} \to \bbR$. We write
\begin{equation}
\sum_{n \geq 1} \alpha(X; n) = \sum_{n \ll \phi(X)} \alpha(X; n)
\end{equation}
if there is a positive constant $\kappa$ such that for all $X \in \bbR_{> 0}$ and all $n > \kappa \phi(X)$, we have $\alpha(X; n) = 0$.

Let $\calE$ denote a multiset of Weierstrass models for elliptic curves over $\bbQ$ with integral coefficients, as in \eqref{eqn:yaxb}. Thus $E \in \calE$ is given by a Weierstrass equation
\begin{equation}\label{Equation: Unreduced Weierstrass equation}
E : y^2 = x^3 + A x + B
\end{equation}
with $A, B \in \bbZ$ and $4 A^3 + 27 B^2 \neq 0$, and we assume nothing further about $A$ and $B$. For $E$ as in \eqref{Equation: Unreduced Weierstrass equation}, we abuse notation and write $\rawheight(E)$ for $\rawheight(A, B)$ and $\tmd(E)$ for $\tmd(A, B)$. For the definition of $\rawheight$ and $\tmd$, see \eqref{eqn:HAB} and \eqref{Equation: twist defect powers}.

Let 
\begin{equation}
\twNad {\calE} (X) \colonequals \# \set{E \in \calE : \twistheight(E) \leq X}.
\end{equation}
We assume that $\twNad {\calE} (X)$ is finite for all $X > 0$. Note that this framework encompasses both the counts $\twistNeq m (X)$ and the counts $\twistNad m (X)$: the former by letting $\calE$ denote the multiset of elliptic curves in $\twistE$ that admit a cyclic $m$-isogeny, counted with repetition if the same elliptic curve admits multiple unsigned cyclic $m$-isogenies, and the latter by letting $\calE$ denote simply the set of elliptic curves in $\twistE$ that possess a cyclic $m$-isogeny. %

Let 
\begin{equation}\label{Equation: defining cM}
	\cM \calE (X; e) \colonequals \# \set{E \in \calE : H(E) \leq X, e \mid \tmd(E)}.
\end{equation}
Note that the points counted by $\twNad \calE (X)$ have \emph{twist height} bounded by $X$, but the points counted by $\cM \calE(X; e)$ have only the function $H$ bounded by $X$.

For our applications, we estimate $\cM \calE (X; e)$ using \Cref{Corollary: Estimates for lattice counts} and sometimes \Cref{Corollary: Principle of Lipschitz for general lattices}, together with some sieving. There are important differences between the sieving necessary for $m = 7$, $m = 10$ and $m = 25$, and $m = 13$ (see \Cref{Lemma: asymptotic for M(X; e) for m = 7}, \Cref{Lemma: asymptotic for M(X; e) for m = 5 and 10 and 25}, and \Cref{Lemma: asymptotic for M(X; e) for m = 13}). These differences stem from certain geometric differences between $X_0(7)$, $X_0(10)$ and $X_0(25)$, and $X_0(13)$: to wit, the elliptic surface parameterizing elliptic curves equipped with a cyclic $7$-isogeny has a place of type II additive reduction, the elliptic surfaces parameterizing elliptic curves equipped with a cyclic $10$-isogeny and cyclic $25$-isogeny have places of type III additive reduction, and the elliptic surface parameterizing elliptic curves equipped with a cyclic $13$-isogeny has places of both type II and type III additive reduction.

Once we have estimates for $\cM \calE (X; e)$, however, these geometric disparities play no further role in our technique. Our results in this thesis depend upon the following application of the M\"{o}bius sieve.

\begin{lemma}\label{Lemma: fundamental sieve}
	We have
	\begin{equation}\label{Equation: twN calE (X) as a sum}
		\twNad {\calE} (X) = \sum_{n \geq 1} \sum_{e \mid n} \mu(n/e) \cM \calE (e^6 X; n).
	\end{equation}
\end{lemma}

\begin{proof}
	As an intermediate step, we define
	\begin{equation}\label{Equation: Defining twN(X e)}
	\twNad \calE (X; e) = \# \set{E \in \calE : \twistheight(E) \leq X, \tmd (E) = e}.
	\end{equation}
	By \eqref{eqn:htdef}, 
	\begin{equation}
	\twNad {\calE} (X) = \sum_{e \geq 1} \twNad \calE (X; e).\label{Equation: twN(X) in terms of twN(X e)}
	\end{equation}
	Moreover, \eqref{eqn:htdef} yields
	\begin{equation}
	\cM \calE (X; e) = \sum_{ef \ll X^{1/6}} \twNad {\calE} (X/(ef)^6; e f).\label{Equation: cM(X e) in terms of twN(X e)}
	\end{equation}

	For the inversion, we must keep the raw-height cutoff fixed.  A model counted by $\twNad\calE(X;e)$ has $\tmd(E)=e$ and $\rawheight(E)\leq e^6X$.  Inclusion--exclusion among the multiples of $e$ therefore gives
	\begin{equation}
		\twNad {\calE} (X; e) = \sum_{f \geq 1} \mu(f) \cM \calE (e^6 X; e f), \label{Equation: twN(X e) in terms of cM(X e)}
	\end{equation}
	where the sum is finite.  Substituting \eqref{Equation: twN(X e) in terms of cM(X e)} into \eqref{Equation: twN(X) in terms of twN(X e)} and letting $n = ef$, we obtain our desired result.
\end{proof}

Before proceeding, we record an absolute-value bound on the summands of \eqref{Equation: twN calE (X) as a sum}, which will be of use in the proofs of \Cref{Lemma: bound on twN>y(X) for m = 7}, \Cref{Lemma: bound on twN>y(X) for m = 10 and 25}, and \Cref{Lemma: bound on twN>y(X) for m = 13}.

\begin{proposition}\label{Proposition: bounding the summands of twN calE X}
For all $n\geq1$ and $X>0$,
\begin{equation}
\abs{\sum_{e\mid n}\mu(n/e)\cM\calE(e^6X;n)}\leq\sum_{e\mid n}\abs{\mu(n/e)}\cM\calE(e^6X;n).
\end{equation}
Every model on the right has twist height at most $X$, and after summing over $n$ a model of twist minimality defect $D$ occurs with multiplicity at most $\tau(D)^2$.
\end{proposition}

\begin{proof}
The displayed inequality is the triangle inequality. If $E$ is counted with $e\mid n\mid D\colonequals\tmd(E)$, then $\rawheight(E)\leq e^6X\leq D^6X$; its total multiplicity is at most $\sum_{n\mid D}2^{\omega(n)}\leq\tau(D)^2$.
\end{proof}

\begin{lemma}\label{Lemma: weighted 2omega divisor factor}
For fixed $\alpha>0$, $\beta>1$, and $y\geq2$,
\begin{equation}
\sum_{n>y}\frac{2^{\omega(n)}\prod_{p\mid n}(1+p^{-\alpha})}{n^\beta}=O_{\alpha,\beta}\parent{\frac{\log y}{y^{\beta-1}}}.
\end{equation}
\end{lemma}

\begin{proof}
Since $\prod_{p\mid n}(1+p^{-\alpha})\leq\sum_{d\mid n}d^{-\alpha}$ and $2^{\omega(dm)}\leq2^{\omega(d)}2^{\omega(m)}$, \Cref{Corollary: sum of 2^omega(n)} gives partial sums $O_\alpha(x\log x)$; partial summation proves the claim.
\end{proof}

\begin{definition}\label{Definition: twNadly and twNadgy}
Define
\begin{equation}
\twNadly {\calE} (X) \colonequals \sum_{n \leq y} \sum_{e \mid n} \mu(n/e) \cM \calE (e^6 X; n)\label{Equation: defining twNadly X}
\end{equation}
and
\begin{equation}
\twNadgy {\calE} (X) \colonequals \sum_{n > y} \sum_{e \mid n} \mu(n/e) \cM \calE (e^6 X; n).\label{Equation: defining twNadgy X}
\end{equation}
\end{definition}

By \Cref{Lemma: fundamental sieve}, we have
\begin{equation}
\twNad {\calE} (X) = \twNadly {\calE} (X) + \twNadgy {\calE} (X).\label{Equation: twN X in terms of twNly X and twNgy X}
\end{equation}

Utilizing our estimates for $\cM \calE (X; e)$ (\Cref{Lemma: asymptotic for M(X; e) for m = 7}, \Cref{Lemma: asymptotic for M(X; e) for m = 5 and 10 and 25}, \Cref{Lemma: asymptotic for M(X; e) for m = 13}, and \Cref{Lemma: asymptotic for M(X; e) for m of genus $0$}), we can obtain asymptotics for $\twNadly \calE (X)$ and $\twNadgy \calE (X)$. Choosing $y$ to minimize the error in \eqref{Equation: twN X in terms of twNly X and twNgy X}, we hope to obtain a good asymptotic for $\twNad \calE (X)$ as $X \to \infty$.

This concludes our treatment of the third step of our approach.

\subsection*{Estimating $\NQad \calE (X)$}

In this subsection, we adapt and abstract the notation and arguments in \cite[\S 5.1]{Molnar-Voight}, enabling us to derive the asymptotics for $\NQeq m (X)$ and $\NQad m (X)$ from those of $\twistNeq m (X)$ and $\twistNad m (X)$. This will clarify the final step described in \cref{Section: Our approach}. We continue to use the assumptions and notation set out in \cref{Section: Our approach revisited}. %

For $X > 0$, we define
\begin{equation}
\NQad \calE (X) \colonequals \# \set{E^{(c)} : E \in \calE, \ c \in \bbZ \ \text{squarefree,} \ \text{and} \ \hht(E^{(c)}) \leq X},
\end{equation}
where $E^{(c)}$ is as in \eqref{Equation: Defining E^(c)}. Here we consider the right-hand side to be the count of a multiset.

The function $\NQad \calE (X)$ counts (with repetition) all elliptic curves in $\scrE$ with rational height less than or equal to $X$ which are quadratic twists of elliptic curves in $\calE$.

 We wish to obtain asymptotics for $\NQad \calE (X)$ as $X \to \infty$. By \Cref{Lemma: parameterizing m-isogenies}, this framework encompasses both the counts $\NQeq m (X)$ and the counts $\NQad m (X)$: the former by letting $\calE$ denote the multiset of elliptic curves in $\twistE$ that possess a cyclic $m$-isogeny, counted with repetition if the same elliptic curve admits multiple unsigned cyclic $m$-isogenies, and the latter by letting $\calE$ denote simply the set of elliptic curves in $\twistE$ that admit a cyclic $m$-isogeny.

We adopt the convention $\twistNad \calE (0) = \NQad \calE (0) = 0$. For $n \geq 1$, write
\begin{equation}\label{Equation: twisth calE}
	\twisthad \calE (n) \colonequals \twistNad \calE (n) - \twistNad \calE (n-1);
\end{equation}
likewise, write
\begin{equation}\label{Equation: hQ calE}
	\hQad \calE (n) \colonequals \NQad \calE (n) - \NQad \calE (n-1).
\end{equation}
Then $\twisthad \calE(n)$ counts the number of elliptic curves $E \in \calE$ of twist height $n$, and $\hQad \calE (n)$ counts the number of elliptic curves $E \in \scrE$ of height $n$ that are quadratic twists of elements of $\calE$.

We define
    \begin{equation}\label{Equation: defining twistL calE X}
    \twistLad \calE (s) \colonequals \sum_{n \geq 1} \twisthad \calE (n) n^{-s}
    \end{equation}
    and
    \begin{equation}\label{Equation: defining LQ calE X}
	\LQad \calE (s) \colonequals \sum_{n \geq 1} \hQad \calE (n) n^{-s}
    \end{equation}
wherever these Dirichlet series converge.

Note that we have $\twistNad \calE (X) = \sum_{n \leq X} \twisthad \calE (n)$, and conversely we have 
$\twistLad \calE(s) = \int_0^\infty u^{-s} \,\mathrm{d}\twistNad \calE(u)$. A good asymptotic understanding of $\twistNad \calE (X)$ therefore enables us to develop a good analytic understanding of $\twistLad \calE (s)$ (for instance, see \Cref{Corollary: twL(s) has a meromorphic continuation for m = 7}). Similarly, we have $\NQad \calE (X) = \sum_{n \leq X} \hQad \calE (n)$, and conversely we have 
$\LQad \calE(s) = \int_0^\infty u^{-s} \,\mathrm{d}\NQad \calE(u)$.

\begin{theorem}\label{Theorem: relationship between twistL(s) and L(s)}
The following statements hold.
\begin{enumerate}
	\item We have
	\begin{equation}
	\hQad \calE(n) = 2 \sum_{c^6 \mid n} \abs{\mu(c)} \twisthad \calE \parent{n/c^6}
	\end{equation}
	\item We have
	\begin{equation}\label{Equation: Relationship between LcalE and twLcalE}
	\LQad \calE (s) = \frac{2 \zeta(6s) \twistLad \calE(s)}{\zeta(12s)}
	\end{equation}
	wherever both sides converge.
	\end{enumerate}
\end{theorem}

\begin{proof}
For (a), we first collect the terms that contribute to $\hQad \calE(n)$ by the quadratic twist factor $c$: 
	\begin{equation} 
	\hQad \calE^{\parent{c}}(n) \colonequals \#\set{E \in \calE : \hht(E^{(c)})= n}
	\end{equation}
By \eqref{Equation: quadratic twists multiply height by c^6} we have $\height(E^{(c)})=c^6 \height(E)$, so
\begin{equation}
	\hQad \calE^{\parent{c}}(n) = 
	\begin{cases}
	\twisthad \calE(n/c^6), & \text{if $c^6 \mid n$;} \\
	0, & \text{otherwise.}
	\end{cases}
	\end{equation}
	 Therefore
	\begin{equation}
	\hQad \calE (n) = \sum_{c \ \text{squarefree}} \hQad \calE^{\parent{c}}(n) \\
	= 2 \sum_{c \geq 1} \abs{\mu(c)} \hQad \calE^{\parent{c}}(n) \\
= 2 \sum_{c^6 \mid n} \abs{\mu(c)} \twisthad \calE\parent{n/c^6} \end{equation}
	proving (a). 
	
	Part (b) simply reframes (a) in the language of Dirichlet series rather than Dirichlet convolutions.
\end{proof}

Thus, we can leverage our understanding of $\twistLad \calE (s)$ to obtain information about $\LQad \calE (s)$. Finally, using Landau's Tauberian theorem (\Cref{Theorem: Landau's Tauberian theorem}), we can transform analytic information about $\LQad \calE (s)$ into asymptotic information about $\NQad \calE (X)$.

\begin{remark}
	In this thesis, we apply a Tauberian theorem (\Cref{Theorem: Landau's Tauberian theorem}) to \eqref{Equation: Relationship between LcalE and twLcalE} in order to obtain asymptotics for $\NQad \calE (X)$. In doing so, we implicitly invoke the apparatus of complex analysis, which is used in the proof of Perron's formula and of Landau's Tauberian theorem. However, we believe an elementary argument applying Dirichlet's hyperbola method \cite[Theorem I.3.1]{Tenenbaum} to \Cref{Theorem: relationship between twistL(s) and L(s)}(a) could achieve similar asymptotics, and perhaps even modestly improve on the error term. 
\end{remark}

We now specialize to the notation we have developed in this section to our problems of interest. 

When $\calE$ is the multiset of elliptic curves in $\twistE$ that possess a cyclic $m$-isogeny, counted with repetition if the same elliptic curve possesses multiple cyclic $m$-isogenies, we have $\twistNad \calE (X) = \twistNeq m (X)$. For this choice of $\calE$, we write
\begin{equation}\label{Equation: Applying N calE to Neq m}
\begin{aligned}
\twistNeq m (X; e) &\colonequals \twistNad \calE (X; e), \\
\twistNeqly m (X) &\colonequals \twistNadly {\calE} (X), \\
\twistNeqgy m (X) &\colonequals \twistNadgy {\calE} (X), \\
\twistheq m (n) &\colonequals \twisthad \calE (n), \\
\twistLeq m (s) &\colonequals \twistLad \calE (s), \\
\hQeq m (n) &\colonequals \hQad \calE (n), \ \text{and} \\
\LQeq m (s) &\colonequals \LQad \calE (s).
\end{aligned}
\end{equation}

Thus for instance $\twistLeq m (s)$ is the height zeta function whose $n$th coefficient $\twistheq m (n)$ is the number of pairs $(E, \phi)$ of elliptic curves $E$ up to quadratic twist with an unsigned cyclic $m$-isogeny $\phi$.

Similarly, when $\calE$ is the set of elliptic curves in $\twistE$ admitting a cyclic $m$-isogeny, we have $\twistNad \calE (X) = \twistNad m (X)$. In analogy with \eqref{Equation: Applying N calE to Neq m}, for this $\calE$ we write
\begin{equation}\label{Equation: Applying N calE to Nad m}
\begin{aligned}
\twistNad m (X; e) &\colonequals \twistNad \calE (X; e), \\
\twistNadly m (X) &\colonequals \twistNadly {\calE} (X), \\
\twistNadgy m (X) &\colonequals \twistNadgy {\calE} (X), \\
\twisthad m (n) &\colonequals \twisthad \calE (n), \\
\twistLad m (s) &\colonequals \twistLad \calE (s), \\
\hQad m (n) &\colonequals \hQad \calE (n), \ \text{and} \\
\LQad m (s) &\colonequals \LQad \calE (s).
\end{aligned}
\end{equation}

Thus for instance $\twistLad m (s)$ is the height zeta function whose $n$th coefficient $\twisthad m (n)$ is the number of elliptic curves $E$ up to quadratic twist admitting an unsigned cyclic $m$-isogeny.

%% file: Ch4_m7.tex
\chapter{Counting elliptic curves with a cyclic \texorpdfstring{$m$}{m}-isogeny for \texorpdfstring{$m = 7$}{m = 7}}\label{Chapter: m = 7}

In this chapter, we prove \Cref{Intro Theorem: asymptotic for NQ(X) for 5 < m <= 9} (\Cref{Theorem: asymptotic for NQ(X) for m = 7}), and \Cref{Intro Theorem: asymptotic for twN(X) for m of genus $0$} (\Cref{Theorem: asymptotic for twN(X) for m = 7}) when $m = 7$. The arguments in this chapter are taken (in many cases verbatim) from \cite{Molnar-Voight}. However, we make several refinements and adjustments to the arguments of \cite{Molnar-Voight}, which account for our much improved error term. Notably, we have furnished a second proof for \Cref{Lemma: asymptotic for M(X; e) for m = 7} which improves its error term. We also establish improved estimates for the order of growth for $\twistLad 7 (s)$. Both improvements propagate through the rest of the chapter.

In \cref{Section: Establishing notation when m = 7}, we establish notations pertaining to $\f 7 (t)$ and $\g 7 (t)$ which will be used throughout the remainder of the chapter. In \cref{Section: twist minimality defect for m = 7}, we develop bounds relating the twist minimality defect to the greatest common divisor of $\f 7 (t)$ and $\g 7 (t)$. In \cref{Section: Estimating twN(X) for m = 7}, we apply the framework developed in \cref{Section: Our approach revisited} to prove \Cref{Intro Theorem: asymptotic for twN(X) for m of genus $0$} when $m = 7$, with a slightly improved error term. In \cref{Section: Working over the rationals for m = 7}, we prove \Cref{Intro Theorem: asymptotic for NQ(X) for 5 < m <= 9} when $m = 7$. In \cref{Section: Computations for m = 7}, we enumerate the elliptic curves with a cyclic $7$-isogeny and twist height at most $10^{42}$, estimate the constants appearing in \Cref{Theorem: asymptotic for twN(X) for m = 7} and \Cref{Theorem: asymptotic for NQ(X) for m = 7}, and empirically confirm that the count of elliptic curves with a cyclic $7$-isogeny aligns with our theoretical estimate.

\section{Establishing notation for \texorpdfstring{$m = 7$}{m = 7}}\label{Section: Establishing notation when m = 7}

By \Cref{Corollary: Elliptic curves for which twNeq(X) = twNad(X)}, we have
\begin{equation}
\twNeq 7 (X) = \twNad 7 (X) \ \text{and} \ \NQeq 7 (X) = \NQad 7 (X)
\end{equation}
for all $X > 0$, so we may use either notation interchangeably. We opt to work with $\twNad 7 (X)$ and related functions.

Pursuant to the notation established in \cref{Section: Parameterizing elliptic curves with a cyclic m-isogeny}, we define 
\begin{equation}
\h 7(t) = \gcd(\f 7 (t), \g 7 (t)),    
\end{equation}
and we define $\tf 7 (t)$ and $\tg 7 (t)$ so that
\begin{equation} \label{equation: Definition of tf and tg for 7-isogenies}
\f 7 (t) = \tf 7 (t) \h 7 (t) \ \text{and} \ \g 7 (t) = \tg 7 (t) \h 7 (t).
\end{equation}
We emphasize that $\tf 7$ and $\tg 7$ are \emph{not} the derivatives of $\f 7$ and $\g 7$. We have
\begin{equation}\label{Equation: Explicit expressions for h7 and tf7 and tg7}
\begin{aligned}
	\h 7(t) &= t^2 + t + 7, \\
	\tf 7 (t) &= -3 (t^2 - 231 t + 735), \ \text{and} \\
	\tg 7 (t) &= 2 (t^4 + 518 t^3 - 11025 t^2 + 6174 t - 64827).
\end{aligned}
 \end{equation}
To work with integral models, as discussed in \cref{Section: Our approach}, we take $t=a/b$ (in lowest terms) and homogenize, obtaining
\begin{equation} \label{equation: defining Cab for m = 7}
\begin{aligned}
	{\C 7}(a, b) &\colonequals b^2 \h 7(a/b)=a^2 + ab + 7 b^2, \\
	\tA 7 (a, b) &\colonequals b^2 \tf 7 (a/b) = -3 (a^2 - 231 a b + 735 b^2), \ \text{and} \\
	{\tB 7} (a, b) &\colonequals b^4 \tg 7 (a/b) = 2 (a^4 + 518 a^3 b - 11025 a^2 b^2 + 6174 a b^3 - 64827 b^4)
\end{aligned}
\end{equation}
In analogy with \eqref{Equation: Explicit expressions for h7 and tf7 and tg7}, we have
\begin{equation}
\begin{aligned}
{\C 7}(a,b) &= \gcd(\A 7(a,b),\B 7(a,b)) \in \Z[a,b], \\
\A 7 (a, b) &= \tA 7 (a, b) {\C 7} (a, b), \ \text{and} \\
\B 7 (a, b) &= \tB 7 (a, b) {\C 7} (a, b).
\end{aligned}
\end{equation}

\section{The twist minimality defect for \texorpdfstring{$m = 7$}{m = 7}}\label{Section: twist minimality defect for m = 7}

In this section, we study the twist minimality defect for
\begin{equation}
E_7(a, b) : y^2 = x^3 + \A 7 (a, b) x + \B 7(a, b)
\end{equation}
using the polynomials $\tA 7 (a, b)$, $\tB 7 (a, b)$, and $\C 7 (a, b)$.

We begin with the following lemma, which shows that when $\gcd(a, b) = 1$, the largest cube dividing $\C 7(a, b)$ almost determines the twist minimality defect.

\begin{lemma}\label{Lemma: 3 and 7 are the ungroomed primes for m = 7}
	Let $(a,b) \in \Z^2$ be $7$-groomed, let $\ell$ be prime, and let $v \in \Z_{>0}$. Then the following statements hold.
	\begin{enumalph}
	\item If $\ell \neq 3, 7$, then $\ell^v \mid \twistdefect(\A 7(a,b),\B 7(a,b))$ if and only if $\ell^{3v} \mid {\C 7}(a, b)$.
	\item $\ell^{3v} \mid {\C 7}(a,b)$ if and only if $\ell \nmid b$ and $\h 7(a/b) \equiv 0 \psmod{\ell^{3v}}$.
	\item If $\ell \neq 3$, then $\ell \mid {\C 7}(a,b)$ implies $\ell \nmid (\partial {\C 7}/\partial a)(a,b) = 2a+b$.
	\end{enumalph}
\end{lemma}

We give two proofs of \Cref{Lemma: 3 and 7 are the ungroomed primes for m = 7}, one using resultants and the other using the polynomial division algorithm.

\begin{proof}[First proof of \Cref{Lemma: 3 and 7 are the ungroomed primes for m = 7}]
We argue as in Cullinan--Kenney--Voight \cite[Proof of Theorem 3.3.1, Step 3]{Cullinan-Kenney-Voight}. 
	For part (a), we compute the resultants
\begin{equation}
\Res(\tA 7(t,1),{\tB 7}(t,1))=\Res(\tf 7(t),\tg 7(t))=-2^8\cdot 3^7 \cdot 7^{14} = \Res(\tA 7(1,u),{\tB 7}(1,u)).
\end{equation}
If $\ell \neq 2,3,7$, then $\ell \nmid \gcd(\tA 7(a,b),\tB 7(a,b))$; so by \eqref{Equation: twist defect powers}, if $\ell^v \mid \tmd(\A 7(a,b),\B 7(a,b))$ then $\ell^{2v} \mid {\C 7}(a,b)$. But also
\begin{equation}
\Res(\tB 7(t,1),{\C 7}(t,1))=\Res(\tg 7(t),\h 7(t)) = 2^8 \cdot 3^3 \cdot 7^7 = \Res(\tB 7(1,u),{\C 7}(1,u)),
\end{equation}
so $\ell \nmid \gcd(\tB 7(a,b), {\C 7}(a, b))$ and thus $\ell^v \mid \tmd(\A 7(a,b),\B 7(a,b))$ if and only if $\ell^{3v} \mid {\C 7}(a,b)$. If $\ell = 2$, then $\C 7(a,b)$ and $\tA 7(a,b)$ are both odd whenever $(a,b)$ is $7$-groomed, so $\A 7(a,b)$ is odd and our claim also holds in this case. 

For (b), by homogeneity it suffices to show that $\ell \nmid b$: this holds since if $\ell \mid b$ then $\A 7(a,0) \equiv -3a^4 \equiv 0 \psmod{\ell}$ and $\B 7(a,0) \equiv 2a^6 \equiv 0 \psmod{\ell}$ so $\ell \mid a$, a contradiction.
	
Part (c) follows from (b) and the fact that $\h 7(t)$ has discriminant $\disc(\h 7(t))=-3^3$. 
\end{proof}

\begin{proof}[Second proof of \Cref{Lemma: 3 and 7 are the ungroomed primes for m = 7}]
	Let $F \neq G$ be homogeneous polynomials in $a$ and $b$ whose homogenizations are coprime over $\bbQ$. Applying the polynomial division algorithm to $F$ and $G$ to cancel all occurrences of the variable $a$, we obtain coprime homogeneous polynomials $\gcdP FG 1(a, b) \in \bbZ[a, b]$ and $\gcdQ FG 1(a, b) \in \bbZ[a, b]$, and a positive integer $\gcdm FG1$, such that
	\begin{equation}
	\gcdP FG 1 (a, b) F(a, b) + \gcdQ FG 1 (a, b) G(a, b)
	\end{equation}
    equals $\gcdm FG1$ times a power of $b$. Likewise, applying the polynomial division algorithm to $F$ and $G$ to cancel all occurrences of the variable $b$, we obtain coprime homogeneous polynomials $\gcdP FG 2(a, b) \in \bbZ[a, b]$ and $\gcdQ FG 2(a, b) \in \bbZ[a, b]$, and a positive integer $\gcdm FG2$, such that
	\begin{equation}
	\gcdP FG 2 (a, b) F(a, b) + \gcdQ FG 2 (a, b) G(a, b)
	\end{equation}
	equals $\gcdm FG2$ times a power of $a$.

	Explicitly, when $(F, G) \in \set{({\tA 7}, {\tB 7}), ({\tB 7}, {\C 7}), ({\C 7}, {\frac{\partial}{\partial a} \C 7})}$, we have %
	\begin{equation} 	\begin{aligned}
		\gcdP {\tA 7} {\tB 7} 1 (a, b) &= 2\parent{13 a^3 + 6776 a^2 b - 121422 a b^2 - 303555 b^3}, \\
		\gcdQ {\tA 7} {\tB 7} 1 (a, b) &= 3\parent{13 a - 2961 b}, \\
		\gcdP {\tA 7} {\tB 7} 2 (a, b) &= -2\parent{1835 a^3 - 44037 a^2 b + 23373 a b^2 - 259308 b^3}, \\
		\gcdQ {\tA 7} {\tB 7} 2 (a, b) &= 9\parent{303 a - 980 b}, \\
		\gcdP {\tB 7} {\C 7} 1 (a, b) &= -1763 a + 239 b, \\
		\gcdQ {\tB 7} {\C 7} 1 (a, b) &= 2 \parent{1763 a^3 + 911232 a^2 b - 20484450 a b^2 + 27625563 b^3}, \\
		\gcdP {\tB 7} {\C 7} 2 (a, b) &= 13 a + 20 b, \ \text{and} \\
		\gcdQ {\tB 7} {\C 7} 2 (a, b) &= 2 \parent{10571 a^3 - 17325 a^2 b + 76293 a b^2 + 185220 b^3}, \\
		\gcdP {\C 7} {{\frac{\partial}{\partial a} \C 7}} 1 (a, b) &= 2^2, \\
		\gcdQ {\C 7} {{\frac{\partial}{\partial a} \C 7}} 1 (a, b) &= - (2 a + b), \\
		\gcdP {\C 7} {{\frac{\partial}{\partial a} \C 7}} 2 (a, b) &= 1, \\
		\gcdQ {\C 7} {{\frac{\partial}{\partial a} \C 7}} 2 (a, b) &= 13 a - 7 b,
	\end{aligned}\end{equation} 
	so
	\begin{equation} 	\begin{aligned}
		\gcdm {\tA 7} {\tB 7} 1 = 2^4 \cdot 3^3 \cdot 7^8 \ &\text{and} \ \gcdm {\tA 7} {\tB 7} 2 = 2^4 \cdot 3 \cdot 7^3, \\ 
		\gcdm {\tB 7} {\C 7} 1 = 2^4 \cdot 3^3 \cdot 7^7 \ &\text{and} \ \gcdm {\tB 7} {\C 7} 2 = 2^4 \cdot 3^3 \cdot 7^2, \ \text{and} \\ 
		\gcdm {\C 7} {{\frac{\partial}{\partial a} \C 7}} 1 = 3^3 \ &\text{and} \ \gcdm {\C 7} {{\frac{\partial}{\partial a} \C 7}} 2 = 3^3. 
	\end{aligned}\end{equation} 

	Let $p$ be a prime not dividing $\gcdm FGi$ for each $F \neq G$ chosen from $\set{{\tA 7}, {\tB 7}, {\C 7}}$ and $i \in \set{1, 2}$. In other words, $p \neq 2, 3, 7$. Let $a$ and $b$ be coprime integers, and suppose $v$ is a positive integer such that $p^{v} \mid \twistdefect(\A 7(a, b), \B 7(a, b))$. 
	
	Suppose by way of contradiction that $p^{2 v} \nmid {\C 7}(a, b)$. As $p^{2v} \mid {\A 7}(a, b), {\B 7} (a, b)$, this implies that $p \mid {\tA 7}(a, b), {\tB 7} (a, b)$, and therefore that
	\begin{equation}
	p \mid \gcdP {{\tA 7}}{{\tB 7}} 1 (a, b) {\tA 7}(a, b) + \gcdQ {{\tA 7}}{{\tB 7}} 1 (a, b) {\tB 7}(a, b) = 2^4 \cdot 3^3 \cdot 7^8 \cdot b^5,
	\end{equation}
	so $p \mid b^5$. Similarly, $p \mid a^5$, but $\gcd(a, b) = 1$, so we have obtained a contradiction, and we conclude $p^{2v} \mid {\C 7}(a, b)$. Now if $p^{3v} \nmid {\C 7}(a, b)$, then $p \mid {\tB 7}(a, b), {\C 7}(a, b)$, and an argument of exactly the same style using $\gcdP {{\tB 7}}{{\C 7}} 1, \gcdQ {{\tB 7}}{{\C 7}} 1, \gcdP {{\tB 7}}{{\C 7}} 2, \gcdQ {{\tB 7}}{{\C 7}} 2$ gives us another contradiction. We conclude that $p^{3 v} \mid {\C 7}(a, b)$ as desired. Finally, comparing ${\C 7}(a, b)$ against ${\frac{\partial}{\partial a} \C 7}(a, b) = 2a + b$, we find that if $p$ divides the former integer it cannot divide the latter, and we have proven our lemma in the case $p \neq 2, 3, 7$.

Now let $p = 2$. For $a$ and $b$ coprime, ${\tA 7}(a,b) \not\equiv 0 \psmod 2$ and ${\C 7}(a,b) \not\equiv 0 \psmod 2$, so ${\A 7}(a,b)$ is odd and our claim holds vacuously in this case. 
\end{proof}

This proof gives bounds on the discrepancy $e\prm$ between the largest cube dividing ${\C 7} (a, b)$ and the cube of the twist minimality defect. Indeed, if $\C 7 (a, b) = e_0^3 n_0$ with $n_0$ cubefree, and $\tmd(\A 7(a, b), \B 7 (a, b)) = e_0 e\prm$, then
\begin{equation}
(e\prm)^3 \mid \lcm(\gcdm {\tA 7} {\tB 7} 1, \gcdm {\tA 7} {\tB 7} 2) \cdot \lcm(\gcdm {\tB 7} {\C 7} 1, \gcdm {\tB 7} {\C 7} 2) = 2^8 \cdot 3^6 \cdot 7^{15}.
\end{equation}
In conjunction with the last paragraph of our second proof, we conclude
\begin{equation}
e\prm \mid 3^2 \cdot 7^5.\label{Equation: bad bound on e'}
\end{equation}
This bound is not sharp, as we shall see.

\begin{definition}\label{Definition: T(e) for m = 7}
For $e \in \bbZ_{>0}$, let $\tcalT 7(e)$ denote the image of
\begin{equation}
\set{(a, b) \in \bbZ^2 : (a, b) \ \text{$7$-groomed}, \ e \mid \twistdefect(\A 7 (a, b), \B 7 (a, b))}
\end{equation}
under the projection
\begin{equation}
\bbZ^2 \to (\bbZ / e^3 \bbZ)^2,
\end{equation}
and let $\tcT 7(e) \colonequals \# \tcalT 7(e)$. Similarly, let $\calT 7(e)$ denote the image of
\begin{equation}
\set{t \in \bbZ : e^2 \mid \f 7(t) \ \text{and} \ e^3 \mid \g 7(t)}
\end{equation}
under the projection
\begin{equation}
\bbZ \to \bbZ / e^3 \bbZ,
\end{equation}
and let $\cT 7(e) \colonequals \#\calT 7(e)$.
\end{definition}

\begin{lemma}\label{Lemma: bound on T(e) for m = 7}
The following statements hold.
\begin{enumalph}
\item $\tcalT 7(e)$ consists of those pairs $(a, b) \in (\bbZ / e^3 \bbZ)^2$ which satisfy the following conditions:%
\begin{itemize}
	\item $\A 7 (a, b) \equiv 0 \psmod {e^2}$ and $\B 7 (a, b) \equiv 0 \psmod {e^3}$, and
	\item $\ell \nmid \gcd(a,b)$ for all primes $\ell \mid e$. %
\end{itemize}
\item Let $(a, b) \in \bbZ^2$. If $(a, b) \psmod{e^3} \in \tcalT {7}(e)$ then $e \mid \tmd(\A {7} (a, b), \B {7} (a, b))$. 
\item The functions $\tcT 7(e)$ and $\cT 7(e)$ are multiplicative, and $\tcT 7(e) = \varphi(e^3) \cT 7(e)$.
\item For every prime $\ell \neq 3,7$ and all $v \geq 1$, we have
	\begin{equation}
	\cT 7(\ell^v) = \cT 7(\ell) = 1 + \left(\frac{\ell}{3}\right).
	\end{equation}
\item The nonzero values of $\cT {7}(3^v)$ are given in Table \ref{table:T7(3)} below. We have $\cT 7 (3^v) = 0$ for $v \geq 3$. The values $\cT {7}(7)$ and $\cT {7} (7^2)$ are given in Table \ref{table:T7(7)} below. We have $\cT 7 (7^v) = 1 + 7^7$ for $v \geq 3$.
\item 
	We have $\cT 7(e) =O(2^{\omega(e)})$ for $e \geq 1$, where $\omega(e)$ is the number of distinct prime divisors of $e$.
	\end{enumalph}
\end{lemma}

\begin{proof}
	Parts (a) and (b) are immediate from \Cref{Definition: T(e) for m = 7}.

	For part (c), multiplicativity follows from the CRT (Sun Zi theorem). For the second statement, let $\ell$ be a prime, and let $e = \ell^v$ for some $v \geq 1$. Consider the injective map
\begin{equation}
\begin{aligned}
\calT 7(\ell^v) \times (\bbZ/ \ell^{3 v})^\times &\to \tcalT 7(\ell^v) \\
(t,u) &\mapsto (tu,u)
\end{aligned}
\end{equation}
We observe $A(1, 0) = -3$ and $B(1, 0) = 2$ are coprime, so no pair $(a, b)$ with $b \equiv 0 \psmod \ell$ can be a member of $\tcalT 7(\ell^v)$. Surjectivity of the given map follows, and counting both sides gives the result.
	
Now part (d). For $\ell \neq 3, 7$, \Cref{Lemma: 3 and 7 are the ungroomed primes for m = 7}(a)--(b) yield
	\begin{equation}
	\calT 7(\ell^v) = \set{t \in \bbZ/\ell^{3v} \bbZ : \h 7(t) \equiv 0 \psmod{\ell^{3v}}}.
	\end{equation}
By \Cref{Lemma: 3 and 7 are the ungroomed primes for m = 7}(c), $\h 7(t) \equiv 0 \psmod \ell$ implies $\frac{\textup{d}}{\textup{d}t}\h 7(t) \not\equiv 0 \psmod \ell$, so Hensel's lemma applies and we need only count roots of $\h 7(t)$ modulo $\ell$. By quadratic reciprocity, this count is 
\begin{equation}
1 + \parent{\frac{-3}{\ell}} = 1 + \parent{\frac{\ell}{3}} = \begin{cases}
2, & \ \text{if} \ \ell \equiv 1 \psmod 3; \\
0, & \ \text{else.}
\end{cases}
\end{equation}

Next, part (e). For $\ell = 3$, we just compute 
$\cT 7(3) = 18$, $\cT 7(3^2) = 27$, and $\cT 7(3^3) = 0$; the observation $\cT 7(3^3) = 0$ implies $\cT 7(3^v) = 0$ for all $v \geq 3$. 
For $\ell = 7$, we compute
\begin{equation}
\cT 7(7) = 1+7^2, \ \cT 7(7^2) = 1 + 7^4, \ \cT 7(7^3) = \dots = \cT 7(7^6) = 1 + 7^7.
\end{equation}
Hensel's lemma still applies to $\h 7(t)$: let $t_0,t_1$ be the roots of $\h 7(t)$ in $\Z_7$ with $t_0 \colonequals 248\,044 \psmod{7^7}$ (so that $t_1=-1-t_0$). The branch $t_0 \equiv -1 \psmod 7$ has $\tg 7(t_0) \not\equiv 0 \psmod 7$, while on the branch $t_1 \equiv 0 \psmod 7$ one checks that $7^4 \mid \tf 7(t)$ and $7^7 \Vert \tg 7(t)$ whenever $t \equiv 42 \psmod{49}$. Consequently, for $v \geq 3$ we have
\begin{equation}
\calT 7(7^v) = \set{t_0} \sqcup \set{t_1 + 7^{3v - 7} u \in \bbZ / 7^{3v} \bbZ :u \in \bbZ / 7^7 \bbZ},\label{equation: decomposition of T(7^v) for m = 7}
\end{equation}
because the first branch requires $7^{3v} \mid \h 7(t)$, whereas on the second branch $7^{3v-7} \mid \h 7(t)$ is necessary and sufficient for the two divisibility conditions defining $\calT 7(7^v)$. Thus $\cT 7(7^v) = 1 + 7^7 = 823544$.

Finally, part (f). From (c)--(e) we conclude
\begin{equation}
\cT 7(e) \leq \frac{27 \cdot 823\,544}{4} \cdot \prod_{\substack{\ell \mid e \\ \ell \neq 3,7}} \parent{1 + \parent{\frac{\ell}{3}}} \leq 5\,558\,922 \cdot 2^{\omega(e)}
\end{equation}
so $\cT 7(e) = O(2^{\omega(e)})$ as claimed.
\end{proof}

\jvtable{table:T7(3)}{
\rowcolors{2}{white}{gray!10}
\begin{tabular}{c c}
$\cT {7} (3^1)$ & $\cT {7} (3^2)$ \\
\hline\hline
$2 \cdot 3^2$ & $3^3$ \\
\end{tabular}
}{All nonzero $\cT {7} (3^v)$}

\jvtable{table:T7(7)}{
\rowcolors{2}{white}{gray!10}
\begin{tabular}{c c}
$\cT {7} (7^1)$ & $\cT {7} (7^2)$ \\
\hline\hline
$1 + 7^2$ & $1 + 7^4$ \\
\end{tabular}
}{$\cT {7} (7^v)$ for $v = 1,2$}

\subsection*{The common factor \texorpdfstring{$\C 7(a, b)$}{C7ab}}

In view of \Cref{Lemma: 3 and 7 are the ungroomed primes for m = 7}, the twist minimality defect away from the primes $2,3,7$ is determined by the quadratic form $\C 7(a,b)=a^2+ab+7b^2=b^2 \h 7(a/b)$. We define
\begin{equation}\label{Equation: Definition of calcc 7}
    \calcc 7 (e) \colonequals \set{(a, b) \in \bbZ^2 : \C 7 (a, b) = e \ \text{and} \ \gcd(a, b) = 1},
\end{equation}
and note $\# \calcc 7 (e) \leq 2^{\omega(e) + 1}$.

Fortunately, $\C 7 (a, b)$ is the norm form of a quadratic order of class number $1$, namely $\bbZ[3 \zeta_6]$, where $\zeta_6$ is a primitive $6$th root of unity. We record some elementary algebraic observations about $\C 7 (a, b)$ and the order $\bbZ[3 \zeta_6]$.

\begin{lemma}\label{Lemma: Algebraic structure of Z[3 zeta6] and C}
	The following statements hold.
	\begin{enumalph}	
		\item The right regular representation of $\bbZ[3 \zeta_6]$ in the basis $\set{1, -1 + 3 \zeta_6}$ induces the map $\repr 7 : \bbZ^2 \to \textup{M}_2(\bbZ)$ given by
		\begin{equation}\label{Equation: Defining repr for m = 7}
		\repr 7 : (a, b) \mapsto \begin{pmatrix}
		a & b \\ -7 b & a + b
		\end{pmatrix}.
		\end{equation}
		\item For all $a, b, c, d, e \in \bbZ$, we have the following implication:
		\begin{equation}
		\C 7 (a, b) = e \implies e \mid \C 7 ((c, d) \cdot \repr 7 (a, b)).
		\end{equation}
		\item Conversely, if $c\prm, d\prm, e, k$ are integers such that $k \geq 1$, $e^k \mid \C 7 (c\prm, d\prm)$, and 
    \begin{equation}
    \gcd(c\prm, d\prm, e) = \gcd(3, e) = 1,
    \end{equation}
    then there are integers $a, b, c, d \in \bbZ$ with $(a, b) \in \calcc 7 (e)$ and
		\begin{equation}
		(c\prm, d\prm) = (c, d) \cdot \repr 7 (a, b)^k.
		\end{equation}
	\end{enumalph}
\end{lemma}

\begin{proof}
	Part (a) is a short computation.
	
	Part (b) follows from the observation that $\C 7 (a, b)$ is the norm on $\bbZ[3 \zeta_6]$ in the basis $\set{1, -1 + 3 \zeta_6}$.
	
	Let $e, k \in \bbZ_{>0}$ and $\alpha\prm \in \bbZ[3 \zeta_6]$. Part (c) will follow if we can prove the following statement. If no inert prime divides both $e$ and $\alpha\prm \in \bbZ[3 \zeta_6]$, $e^k \mid \Nm(\alpha\prm)$, and $\gcd(3, e) = 1$, then there are algebraic integers $\alpha, \beta \in \bbZ[3 \zeta_6]$ such that $\alpha\prm = \alpha \beta^k$ and $\Nm(\beta) = e$. 
 
        We now prove this assertion. The order $\bbZ[3 \zeta_6]$ is a suborder of the Euclidean domain $\bbZ[\zeta_6]$ of conductor 3, and it inherits the following almost unique factorization: up to sign, every nonzero $\alpha \in \Z[3\zeta_6]$ can be written uniquely %
	as 
	\begin{equation}
	\alpha = \beta \pi_1^{e_1} \cdots \pi_r^{e_r},
	\end{equation}
	where $\Nm(\beta)$ is a power of $3$, $\pi_1, \dots, \pi_r$ are distinct irreducibles coprime to $3$, and $e_1, \dots, e_r$ are positive integers.

    Write
    \begin{equation}
    \alpha\prm = \beta\prm \pi_1^{e_1} \cdots \pi_r^{e_r}.
    \end{equation}
    As every prime dividing both $\alpha\prm$ and $e$ splits, $\Nm(\pi_j) = p_j$ is a rational prime whenever $\pi_j \mid e$. Moreover, because $\gcd(c\prm, d\prm, e) = 1$, if a prime $p$ divides $e$ then $\Nm(\pi_i) = \Nm(\pi_j) = p$ implies $\pi_i = \pi_j$. Write $e = \prod_{j = 1}^r p_j^{f_j}$, and let $\beta = \prod_{j = 1}^r \pi_j^{f_j}$. Necessarily, $\beta^k \mid \alpha\prm$; letting $\alpha = \beta\prm \prod_{j = 1}^r \pi_j^{e_j - k f_j}$, our claim follows.
\end{proof}

The twist minimality defect measures the disparity between $\rawheight(A, B)$, which is easy to compute, and $\twistheight(A, B)$, which is of arithmetic interest: this disparity cannot be too large compared to $\C 7(a,b)$, as the following theorem shows.

\begin{theorem}\label{Theorem: Controlling size of twist minimality defect for m = 7}
	The following statements hold.
	\begin{enumalph}
	\item For all $(a, b) \in \bbR^2$, we have
	\begin{equation}\label{Equation: upper and lower bounds for H for m = 7}
	108 \C 7(a, b)^6 \leq \rawheight(\A 7(a, b), \B 7(a, b)) \leq \upperratio 7 \C 7(a, b)^6,
	\end{equation}
	where the constant $\upperratio 7 = 311\,406\,871.990\,204\ldots$ is an algebraic number given by evaluating the function $H(\A 7(a, b), \B 7(a, b))$ at appropriate roots of \eqref{equation: Roots defining upperratio for m = 7}.
	\item If $\C 7(a, b) = e_0^3 n_0$, with $n_0$ cubefree, then $\twistdefect(\A 7(a, b), \B 7(a, b)) = e_0 e\prm$ for some $e\prm \mid 3 \cdot 7^3$, and
	\begin{equation}
	\frac{2^2}{3^3 \cdot 7^{18}} e_0^{12} n_0^6 \leq \twistheight(\A 7(a, b), \B 7(a, b)) \leq \upperratio 7 e_0^{12} n_0^6.
	\end{equation}
	\end{enumalph}
\end{theorem}

\begin{proof}
	We first prove (a). We wish to find the extrema of the ratio 
 \begin{equation}\label{Equation: ratio for m = 7}
 \rawheight(\A 7(a, b), \B 7(a, b))/\C 7(a, b)^6.
 \end{equation}
 As \eqref{Equation: ratio for m = 7} is homogeneous of degree 0, and $\C 7(a, b)$ is positive definite, we may assume without loss of generality that $\C 7(a, b) = 1$. Using the theory of Lagrange multipliers, and examining the critical points of $\rawheight(\A 7(a, b), \B 7(a, b))$ subject to $\C 7(a, b) = 1$, we verify that \eqref{Equation: upper and lower bounds for H for m = 7} holds. Moreover, the lower bound is attained at $(1, 0)$, and the upper bound is attained when $a$ and $b$ are appropriately chosen roots of
	\begin{equation}\label{equation: Roots defining upperratio for m = 7}
	\begin{aligned}
		 &1296 a^8 - 2016 a^6 + 2107 a^4 - 1596 a^2 + 252  \\
		 =& 2^4 \cdot 3^4 \cdot a^8 - 2^5 \cdot 3^2 \cdot 7 \cdot a^6 + 7^2 \cdot 43 \cdot a^4 - 2^2 \cdot 3 \cdot 7 \cdot 19 \cdot a^2 + 2^2 \cdot 3^2 \cdot 7, \ \text{and} \\
		 &1\,067\,311\,728 b^8 - 275\,298\,660 b^6 + 43\,883\,077 b^4 - 3\,623\,648 b^2 + 1849 \\
		 =& 2^4 \cdot 3^4 \cdot 7^7 \cdot b^8 - 2^2 \cdot 3^2 \cdot 5 \cdot 7^6 \cdot 13 \cdot b^6 + 7^6 \cdot 373 \cdot b^4 - 2^5 \cdot 7^2 \cdot 2311 \cdot b^2 + 43^2
	\end{aligned}
	\end{equation}
	respectively. For $(a, b) = (0.450\,760\,996\,604\,693\,04\ldots, -0.371\,118\,011\,382\,744\,86\ldots)$, the arguments that maximize the ratio 
 \begin{equation}
 \rawheight(\A 7(a, b), \B 7(a, b))/\C 7(a, b)^6,
 \end{equation}
 we have $27 \abs{\B 7(a, b)}^2 > 4 \abs{\A 7(a, b)}^3$. %
	
	We now prove (b). Write $\C 7(a, b) = e_0^3 n_0$ with $n_0$ cubefree, and write 
 \begin{equation}
 \twistdefect(\A 7(a, b), \B 7(a, b)) = e_0 e\prm.
 \end{equation}
 By \Cref{Lemma: 3 and 7 are the ungroomed primes for m = 7}, $e\prm = 3^v \cdot 7^w$ for some $v, w \geq 0$; a short computation shows $v \in \set{0, 1}$, and \eqref{equation: decomposition of T(7^v) for m = 7} shows $w \leq \ceil{7/3} = 3$. 
	
	As 
	\begin{equation}
	\rawheight(\A 7(a, b), \B 7(a, b)) = e_0^6 \parent{e\prm}^6 \twistheight(\A 7(a, b), \B 7(a, b)),
	\end{equation}
	we see
	\begin{equation}
	\frac{108}{(e\prm)^6} e_0^{12} n_0^6 \leq \twistheight(\A 7(a, b), \B 7(a, b)) \leq \frac{\upperratio 7}{(e\prm)^6} e_0^{12} n_0^6.
	\end{equation}
	Rounding $e\prm$ up to $3 \cdot 7^3$ on the left, and rounding down to $1$ on the right gives the desired result.
\end{proof}

Unsurprisingly, \Cref{Theorem: Controlling size of twist minimality defect for m = 7} shows that the bound on $e\prm$ given by \eqref{Equation: bad bound on e'} is not sharp.

\begin{corollary}\label{Corollary: bound twist defect in terms of twist height for m = 7}
	Let $(a, b)$ be a $7$-groomed pair. We have
\begin{equation}
	\twistdefect(\A 7(a, b), \B 7(a, b)) \leq \frac{3^{5/4} \cdot 7^{9/2}}{2^{1/6}} \twistheight(\A 7(a, b), \B 7(a, b))^{1/12} \end{equation}
where $3^{5/4} \cdot 7^{9/2} / 2^{1/6} = 22\,344.227\,186\ldots$ 
\end{corollary}

\begin{proof}
	In the notation of \Cref{Theorem: Controlling size of twist minimality defect for m = 7}(b),
	\begin{equation}
	e_0^{12} n_0^6 \leq \frac{3^3 \cdot 7^{18}}{2^2} \twistheight(\A 7(a, b), \B 7(a, b)).
	\end{equation}
	Multiplying through by $(e\prm)^{12}$, rounding $n_0$ down to $1$ on the left, rounding $e\prm$ up to $3 \cdot 7^3$ on the right, and taking $12$th roots of both sides, we obtain the desired result.
\end{proof}

\begin{remark}
	 We could instead prove \Cref{Theorem: Controlling size of twist minimality defect for m = 7}(a) as follows. We assume without loss of generality that $\C 7(a, b) = 1$. Now the level set
	\begin{equation}
	\set{(a, b) \in \bbR^2 : \C 7(a, b) = a^2 + ab + 7 b^2 = 1}
	\end{equation}
	is parameterized by the function 
	\begin{equation}
	\phi : \theta \mapsto \parent{\cos \theta - \frac{1}{3 \sqrt{3}} \sin \theta, \frac{2}{3 \sqrt{3}} \sin \theta}.
	\end{equation}
	We can now use a computer to show
	\begin{equation} 	\begin{aligned}
	\min_\theta \rawheight(\A 7(\phi(\theta)), \B 7(\phi(\theta))) &= 108, \text{and} \\
	\max_\theta \rawheight(\A 7(\phi(\theta)), \B 7(\phi(\theta))) &\eqqcolon \upperratio 7 = 311\,406\,871.990\,204\ldots,
	\end{aligned}\end{equation} 
	with the minimum attained when $\theta = 0$ and $(a, b) = (1, 0)$, and the maximum attained when $\theta = 4.980\,802\,4\ldots$ and $(a, b) = (0.450\ldots, -0.371\ldots)$. It is straightforward to show $\rawheight(\A 7(1, 0), \B 7(1, 0)) = 108$ but 
	\begin{equation}
	\frac{\rawheight(\A 7(a, 1), \B 7(a, 1))}{\C 7(a, 1)^6} > 108
	\end{equation}
	for all $a$, so this minimum value is exact. This argument does not express $\upperratio 7$ as an algebraic number, however.
\end{remark}

\section{Estimates for twist classes for \texorpdfstring{$m = 7$}{m = 7}} \label{Section: Estimating twN(X) for m = 7}

In this section, we use \cref{Section: Our approach revisited} to estimate $\twistNad 7(X)$, counting the number of twist minimal elliptic curves over $\bbQ$ admitting a cyclic $7$-isogeny.

Recall \eqref{Equation: defining cM}, as well as \eqref{Equation: Applying N calE to Neq m} and \eqref{Equation: Applying N calE to Nad m}. By \cref{Section: Parameterizing elliptic curves with a cyclic m-isogeny}, $\cM 7(X; e)$ counts pairs $(a, b) \in \bbZ^2$ with
\begin{itemize}
	\item $(a, b)$ $7$-groomed,
	\item $\rawheight(\A 7(a,b),\B 7(a,b)) \leq X$, and
	\item $e \mid \twistdefect(\A 7(a, b), \B 7(a, b))$.
\end{itemize}

The following proposition refines \Cref{Lemma: fundamental sieve}, and specifies both an order of growth and an explicit upper bound past which the summands of \eqref{Equation: twN calE (X) as a sum} vanish when $m = 7$.

\begin{proposition}\label{Proposition: fundamental sieve for for m = 7}
	We have
	\begin{equation}\label{Equation: twN(X) in terms of M(X; e) for m = 7}
	\twNad 7 (X) = \sum_{n \ll X^{1/12}} \sum_{e \mid n} \mu(n/e) \cM 7 (e^6 X; n);
	\end{equation}
	more precisely, we can restrict our sum to
	\begin{equation}
	n \leq \frac{3^{5/4} \cdot 7^{9/2}}{2^{1/6}} \cdot X^{1/12}.
	\end{equation}
\end{proposition}

\begin{proof}
	Let $(a, b) \in \bbZ^2$, and suppose 
 \begin{equation}
 \rawheight(\A 7(a,b),\B 7(a,b)) \leq e^6 X \ \text{and} \ e \mid \twistdefect(\A 7(a, b), \B 7(a, b)).
 \end{equation}
	If we can prove 
	\begin{equation}\label{Equation: bound on e for m = 7}
	e \leq \frac{3^{5/4} \cdot 7^{9/2}}{2^{1/6}} \cdot X^{1/12},
	\end{equation}
	then our claim will follow.
	
	Write $\C 7(a, b) = e_0^3 n_0$, with $n_0$ cube-free. By \Cref{Theorem: Controlling size of twist minimality defect for m = 7}(a), we have
	\begin{equation}
	108 e_0^{18} n_0^6 \leq e^6 X.\label{Equation: First inequality between e0 and X for m = 7}
	\end{equation}
	On the other hand, by \Cref{Theorem: Controlling size of twist minimality defect for m = 7}(b), we have $e \mid 3 \cdot 7^3 \cdot e_0$, and \textit{a fortiori}
	\begin{equation}
	e \leq 3 \cdot 7^3 e_0.\label{Equation: e0 versus e for m = 7}
	\end{equation}
	Multiplying \eqref{Equation: First inequality between e0 and X for m = 7} through by $(3 \cdot 7^3)^{18}$ and utilizing \eqref{Equation: e0 versus e for m = 7}, we conclude
	\begin{equation}
	2^2 \cdot 3^3 e^{18} \cdot n_0^6 \leq 3^{18} \cdot 7^{54} e^6 X.
	\end{equation}
	Rounding $n_0$ down to $1$ and rearranging, we obtain \eqref{Equation: bound on e for m = 7}.
\end{proof}

Recall that a pair $(a, b) \in \bbZ^2$ is $7$-groomed if $\gcd(a, b) = 1$, $b > 0$, and $a/b \not\in \cusps 7 = \set{-7, \infty}$ (see \Cref{Definition: m-groomed pairs} and Table \ref{table:auxiliaries}). In order to estimate $\cM 7(X; e)$, we further unpack the $7$-groomed condition on pairs $(a, b)$. We therefore let $\cM 7(X; d, e)$ denote the number of pairs $(a, b) \in \bbZ^2$ with
	\begin{itemize}
		\item $\gcd(da, db, e) = 1$, $b > 0$, and $a/b \not\in \cusps 7$,
		\item $\rawheight(\A 7(d a, d b), \B 7(d a, db)) \leq X$, and
		\item $e \mid \twistdefect(\A 7(d a, d b), \B 7(da, db))$.
	\end{itemize}
By \Cref{Theorem: Controlling size of twist minimality defect for m = 7}, and because $\rawheight(\A 7(a, b), \B 7(a, b))$ is homogeneous of degree 12, a M\"{o}bius sieve yields
\begin{equation}\label{equation: cM(X;e) in terms of cM(X; d, e) for m = 7}
	\cM 7(X; e) = \sum_{\substack{d \ll X^{1/12} \\ \gcd(d, e) = 1}} \mu(d) \cM 7(X; d, e);
\end{equation}	
more precisely, we can restrict our sum to
	\begin{equation}
	d \leq \frac{1}{2^{1/6} \cdot 3^{1/4}} \cdot X^{1/12}.
	\end{equation} 

Before proceeding, we give an outline of the argument employed in this section. In \Cref{Lemma: asymptotic for M(X; e) for m = 7}, we use the Principle of Lipschitz to estimate $\cM 7(X; d, e)$, then piece these estimates together using \eqref{equation: cM(X;e) in terms of cM(X; d, e) for m = 7} to estimate $\cM 7(X; e)$. Heuristically,
\begin{equation}
\cM 7(X; d, e) \sim \frac{\R 7 \cT 7(e) X^{1/6}}{d^2 e^3} \prod_{\ell \mid e} \parent{1 - \frac{1}{\ell}}
\end{equation}
(where $\R 7$ is the area of \eqref{eqn: R(X)} when $m = 7$ and $\cT 7$ is the arithmetic function investigated in \Cref{Lemma: bound on T(e) for m = 7})
by summing over the congruence classes modulo $e^3$ that satisfy $e \mid \twistdefect(\A 7(d a, d b), \B 7(da, db))$. Then \eqref{equation: cM(X;e) in terms of cM(X; d, e) for m = 7} suggests
\begin{equation}
\cM 7(X; e) \sim \frac{\R 7 \cT 7(e) X^{1/6}}{\zeta(2) e^3 \prod_{\ell \mid e} \parent{1 + \frac{1}{\ell}}}.\label{equation: Main term for cM(X; e) for m = 7}
\end{equation}
Substituting \eqref{equation: Main term for cM(X; e) for m = 7} into \Cref{Lemma: fundamental sieve}, we obtain the heuristic estimate
\begin{equation}
\twistNad 7(X) \sim \frac{\Q 7 \R 7 X^{1/6}}{\zeta(2)},
\end{equation}
where
\begin{equation}\label{Equation: Q for m = 7}
\Q 7 \colonequals \sum_{n \geq 1} \frac{\cT 7(n) \varphi(n) }{n^3 \prod_{\ell \mid n} \parent{1 + \frac{1}{\ell}}}.
\end{equation}

	To make this estimate for $\twistNad 7(X)$ rigorous, and to get a better handle on the order of growth of its error term, we now decompose \eqref{Equation: twN(X) in terms of M(X; e) for m = 7} in accordance with \Cref{Definition: twNadly and twNadgy}, so
	\begin{equation} \twistNad 7(X) = \twistNadly 7(X) + \twistNadgy 7(X). \end{equation} 
	We then estimate $\twistNadly 7(X)$ in \Cref{Lemma: asymptotic for twN<=y(X) for m = 7}, and treat $\twistNadgy 7(X)$ as an error term which we bound in \Cref{Lemma: bound on twN>y(X) for m = 7}. Setting the error from our estimate equal to the error arising from $\twistNadgy 7(X)$, we obtain \Cref{Theorem: asymptotic for twN(X) for m = 7}.

In the remainder of this section, we follow the outline suggested here by successively estimating $\cM 7(X; d, e)$, $\cM 7(X; e)$, $\twistNadly 7(X)$, $\twistNadgy 7(X)$, and finally $\twistNad 7(X)$. 

We first estimate $\cM 7(X; d, e)$ and $\cM 7(X; e)$.

\begin{lemma}\label{Lemma: asymptotic for M(X; e) for m = 7}
The following statements hold.
\begin{enumalph}
\item	If $\gcd(d, e) > 1$, then $\cM 7(X; d, e) = 0$. If $\gcd(d, e) = 1$, we have 
	\begin{equation}
	\cM 7(X; d, e) = \frac{\R 7 \cT 7(e) X^{1/6}}{d^2 e^3} \prod_{\ell \mid e} \parent{1 - \frac{1}{\ell}} + O\parent{\frac{2^{\omega(e)} X^{1/12}}{d e^{3/2}}}
	\end{equation}
        for $X, d, e \geq 1$. Here, $\R 7$ is the area of \eqref{eqn: R(X)} when $m = 7$.
\item We have
	\begin{equation}
	\cM 7(X; e) = \frac{\R 7 \cT 7(e) X^{1/6}}{\zeta(2) e^3 \prod_{\ell \mid e} \parent{1 + \frac{1}{\ell}}} + O\parent{\frac{2^{\omega(e)} X^{1/12} \log X}{e^{3/2}}}
	\end{equation}
        for $X \geq 2$, $e \geq 1$.
\end{enumalph}
In both cases, the implied constants are independent of $d$, $e$, and $X$.
\end{lemma}

We give two partial proofs of \Cref{Lemma: asymptotic for M(X; e) for m = 7}. The first proof gives an intuitive interpretation of the coefficient of $X^{1/6}$, and generalizes readily to other elliptic surfaces with type II additive reduction, but yields only a degraded error term. The second proof leverages the observation that $\C 7 (a, b)$ is the norm of the order $\bbZ[3 \zeta_6]$ to give the full error term, but does not make the leading coefficient as explicit.

\begin{proof}[First proof of \Cref{Lemma: asymptotic for M(X; e) for m = 7}]
	We begin with (a) and examine $\cM 7(X; d, e)$. If $d$ and $e$ are not coprime, then $\cM 7(X; d, e) = 0$ because $\gcd(da, db, e) \geq \gcd(d, e) > 1$. On the other hand, if $\gcd(d, e) = 1$, we have a bijection from the pairs counted by $\cM 7(X; 1, e)$ to the pairs counted by $\cM 7(d^{12} X; d, e)$ given by $(a, b) \mapsto (d a, d b)$.

For $X\geq 1$ and $e, a_0, b_0 \in \bbZ$, we write
\begin{equation}
    L_7(X; e, a_0, b_0) \colonequals \#\{(a, b) \in \calR 7(X) \cap \bbZ^2 : (a, b) \equiv (a_0, b_0) \psmod {e^3}, (a:b) \not\in \cusps 7 \}
\end{equation}
(this notation will not be used outside of this proof). By \Cref{Corollary: Estimates for lattice counts}, we have
\begin{equation}
L_{7}(X; e, a_0, b_0) = \frac{\R {7} X^{1/6}}{e^6} + O \parent{\frac{X^{1/12}}{e^3}}.
\end{equation}
Now by \Cref{Lemma: bound on T(e) for m = 7}(c), we have
\begin{equation}\label{Equation: cM (X 1 e) for m = 7}
\begin{aligned}
	\cM 7(X; 1, e) &= \sum_{(a_0, b_0) \in \tcalT 7(e)} L_7(X; e, a_0, b_0) \\
	&= \varphi(e^3) \cT 7(e) \parent{\frac{\R 7 X^{1/6}}{e^6} + O\parent{\frac{X^{1/12}}{e^3}}} \\
	&= \frac{\R 7 \cT 7(e) X^{1/6}}{e^3} \prod_{\ell \mid e} \parent{1 - \frac{1}{\ell}} + O(\cT 7(e) X^{1/12}).
\end{aligned}
\end{equation}
Scaling by $d$ and invoking \Cref{Lemma: bound on T(e) for m = 7}(f), we obtain
\begin{equation}
	\cM 7(X; d, e) = \frac{\R 7 \cT 7(e) X^{1/6}}{d^2 e^3} \prod_{\ell \mid e} \parent{1 - \frac{1}{\ell}} + O\parent{\frac{2^{\omega(e)} X^{1/12}}{d}}.
\end{equation}

For part (b), we compute
\begin{equation} \label{equation: cmXe for m = 7}
\begin{aligned}
	\cM 7(X; e) &= \sum_{\substack{d \ll X^{1/12} \\ \gcd(d, e) = 1}} \mu(d) \cM 7(X; d, e) \\
	&= \sum_{\substack{d \ll X^{1/12} \\ \gcd(d, e) = 1}} \mu(d) \parent{\frac{\cT 7(e) \R 7 X^{1/6}}{d^2 e^3} \prod_{\ell \mid e} \parent{1 - \frac{1}{\ell}} + O\parent{\frac{2^{\omega(e)} X^{1/12}}{d}}} \\
	&= \frac{\R 7 \cT 7(e) X^{1/6}}{e^3} \prod_{\ell \mid e} \parent{1 - \frac{1}{\ell}} \sum_{\substack{d \ll X^{1/12} \\ \gcd(d, e) = 1}} \frac{ \mu(d)}{d^2} + O\parent{2^{\omega(e)} X^{1/12} \sum_{\substack{d \ll X^{1/12} \\ \gcd(d, e) = 1}} \frac 1d}.
	\end{aligned}
	\end{equation}
We plug the straightforward estimates
\begin{equation}
\sum_{\substack{d \ll X^{1/12} \\ \gcd(d, e) = 1}} \frac{ \mu(d)}{d^2} = \frac{1}{\zeta(2)} \prod_{\ell \mid e} \parent{1 - \frac{1}{\ell^2}}\inv + O(X^{-1/12})
\end{equation}
and
\begin{equation} 
\sum_{\substack{d \leq X^{1/12}}} \frac 1d = \frac{1}{12}\log X+ O(1) \end{equation}
into \eqref{equation: cmXe for m = 7}, along with \Cref{Lemma: bound on T(e) for m = 7}(f). Simplifying now gives
\begin{equation}
\begin{aligned}
\cM 7(X;e)	
	&= \frac{\R 7 \cT 7(e) X^{1/6}}{\zeta(2) e^3 \prod_{\ell \mid e} \parent{1 + \frac{1}{\ell}}} + O(2^{\omega(e)} X^{1/12} \log X)
\end{aligned}
\end{equation}
proving (b) with a degraded error term.
\end{proof}

We now give our second proof of \Cref{Lemma: asymptotic for M(X; e) for m = 7}.

\begin{proof}[Second proof of \Cref{Lemma: asymptotic for M(X; e) for m = 7}]
	As in the previous proof, we may restrict our attention to $\cM 7 (X; 1, e)$. Throughout this proof, $d$ will not refer to the second argument of $\cM 7 (X; d, e)$.
 
	Let $e \in \bbZ_{>0}$, and let $e_0$ be the smallest integer for which $e \mid 3 \cdot 7^3 e_0$. By \Cref{Theorem: Controlling size of twist minimality defect for m = 7}(b), if $e \mid \tmd(\A 7 (a, b), \B 7 (a, b))$, then $e_0^3 \mid \C 7 (a, b)$. By \Cref{Lemma: Algebraic structure of Z[3 zeta6] and C}(b), if $\gcd(3, e_0) = 1$, we have the following implications for all $a_0, b_0, c, d \in \bbZ$:
	\begin{equation}\label{Equation: C = e implies e^3 divides C}
	\C 7 (a_0, b_0) = e_0 \implies e_0^3 \mid \C 7 ((c, d) \cdot \repr 7(a_0,b_0)^3),
	\end{equation}
	where $\repr 7 : \bbZ^2 \to \textup{M}_2(\bbZ)$ is defined in \eqref{Equation: Defining repr for m = 7}. By \Cref{Lemma: Algebraic structure of Z[3 zeta6] and C}(c), if $\gcd(c\prm, d\prm, e_0) = 1$, we also have the converse implication
	\begin{equation}\label{Equation: e^3 divides C implies C = e}
	e_0^3 \mid \C 7 (c\prm, d\prm) \implies (c\prm, d\prm) = (c, d) \cdot \repr 7(a_0, b_0)^3
	\end{equation}
	for some $(a_0, b_0) \in \calcc 7 (e_0)$ and $(c, d) \in \bbZ$. Our aim is to use \eqref{Equation: C = e implies e^3 divides C} and \eqref{Equation: e^3 divides C implies C = e}, in tandem with \Cref{Corollary: Principle of Lipschitz for general lattices}, to improve on the error term given in the last proof.
		
	For $e \geq 1$, let $\tcalT 7(a_0, b_0, e)$ denote the image of
	\begin{equation}
	\set{(c\prm, d\prm) \in \bbZ^2 \cdot \repr 7(a_0, b_0)^3 : e \mid \twistdefect(\A 7 (c\prm, d\prm), \B 7 (c\prm, d\prm)) \ \text{and} \ (c\prm, d\prm) \ \text{$7$-groomed}}
	\end{equation}
	under the projection
	\begin{equation}
	\bbZ^2 \to (\bbZ / e^3 \bbZ)^2.
	\end{equation}
	We also let $\tcT 7(a_0, b_0, e) \colonequals \# \tcalT 7(a_0, b_0, e)$. If $e_0 \mid e$, we have the straightforward equality
	\begin{equation}
	\# \parent{e_0^3 \bbZ / e^3 \bbZ} = e^3/e_0^3.
	\end{equation}
	Thus if $\C 7 (a_0, b_0) = e_0$ and $e/e_0$ is an integer dividing $3 \cdot 7^3$, we have the bound $\tcT 7(a_0, b_0, e) \leq 3^6 \cdot 7^{18}$.
	
	If $\gcd(3^2, e) > 3$, then $\cM 7(X; 1, e) = 0$ by \Cref{Lemma: bound on T(e) for m = 7}. Otherwise, we let $e_0$ be the smallest integer for which $e \mid 3 \cdot 7^3 \cdot e_0$. In this case, $\cM 7(X; 1, e)$ is the sum over $(a_0, b_0) \in \calcc 7(e_0)$ and $(c_0, d_0) \in \tcalT 7(a_0, b_0, e)$ of
	\begin{equation}
         \# \set{(c\prm, d\prm) \in \calR 7(X) \cap (\bbZ^2 \cdot \repr 7(a_0, b_0)^3): (c\prm, d\prm) \equiv (c_0, d_0) \psmod {e^3}, \ c\prm/d\prm \not\in \cusps 7}.
	\end{equation}
    By \Cref{Corollary: Principle of Lipschitz for general lattices}, we therefore have
	\begin{equation}
	\cM 7(X; 1, e) = \sum_{(a_0, b_0) \in \calcc 7(e_0)} \sum_{(c_0, d_0) \in \tcalT 7(a_0, b_0, e)} \parent{\frac{\R 7 X^{1/6}}{(\det \repr 7(a_0, b_0))^3} + O\parent{\frac{X^{1/12}}{\sigma(\repr 7 (a_0, b_0))^3}}}.
	\end{equation}
	But $\det \repr 7 (a_0, b_0) = e_0$ by assumption, and on the other hand the singular values of $\repr 7 (a_0, b_0)$ are
 \begin{equation}
\sigma_{\pm}(a_0, b_0) =  \parent{\frac{2 a_0^2 + 2 a_0 b_0 + 51 b_0^2 \pm \abs{b_0} \sqrt{148 a_0^2 + 148 a_0 b_0 + 2405 b_0^2}} 2}^{1/2}.
 \end{equation}
 We use Lagrange multipliers to find the extrema of $\sigma_{-}(a_0, b_0)$ subject to the constraint $\C 7(a_0, b_0) = 1$, and thus of $\sigma_{-}(a_0, b_0)/\C 7 (a_0, b_0)^{1/2}$. We thereby obtain
	\begin{equation}
	\parent{\frac{101 - 16 \sqrt{37}}{27}}^{1/2} e_0^{1/2} \leq \sigma_{-}(\repr 7 (a_0, b_0)) \leq e_0^{1/2}.
	\end{equation}
	The lower extremum is attained when $a_0/b_0 = -1/2$, and the upper extremum is attained when $b_0 = 0$.
 
 Now as $\#\calcc 7(e_0) = O(2^{\omega(e)})$ and $\tcT 7(a_0, b_0, e) = O(1)$, we have
	\begin{equation}
	\cM 7(X; 1, e) = \frac{\R 7 X^{1/6}}{e_0^3} \sum_{(a_0, b_0) \in \calcc 7(e_0)} \tcT 7(a_0, b_0, e) + O\parent{2^{\omega(e)} \frac{X^{1/12}}{e^{3/2}}}.\label{Equation: New expression for cM(X e 1) when m = 7}
	\end{equation}
	By considering the limit
 \begin{equation}
 \lim_{X \to \infty} \frac{\cM 7 (X; 1, e)}{X^{1/6}},
 \end{equation}
 we deduce
	\begin{equation}
	\frac{\R 7}{e_0^3} \sum_{(a_0, b_0) \in \calcc 7(e_0)} \tcT 7(a_0, b_0, e) = \frac{\R 7 \cT 7(e)}{e^3}.
	\end{equation}

    We turn our attention at last to $\cM 7 (X; d, e)$: scaling by $d$ as in the previous proof, we conclude
	\begin{equation}
	\cM 7(X; d, e) = \frac{\R 7 \cT 7(e) X^{1/6}}{d^2 e^3} \prod_{\ell \mid e} \parent{1 - \frac{1}{\ell}} + O\parent{2^{\omega(e)} \frac{X^{1/12}}{d e^{3/2}}},
	\end{equation}
	where the implicit constant is independent of $X,$ $d$, and $e$.
	
	Following the proof of part (b) above, we obtain
	\begin{equation}
	\cM 7(X; e) = \frac{\R 7 \cT 7(e) X^{1/6}}{\zeta(2) e^3 \prod_{\ell \mid e} \parent{1 + \frac{1}{\ell}}} + O\parent{\frac{2^{\omega(e)} X^{1/12} \log X}{e^{3/2}}}.
	\end{equation}
\end{proof}

We let
	\begin{equation}
	\Q 7 \colonequals \sum_{n \geq 1} \frac{\varphi(n) \cT 7(n)}{n^3 \prod_{\ell \mid n} \parent{1 + \frac{1}{\ell}}}, %
	\end{equation}
	and we let
\begin{equation}\label{Equation: twconst for m = 7}
	\twconst 7 \colonequals \frac{\Q 7 \R 7}{\zeta(2)}.
\end{equation}
Here, as always, $\R 7$ is the area of the region
	\begin{equation}
	\calR 7(1) = \set{(a, b) \in \bbR^2 : \rawheight(\A 7 (a, b), \B 7 (a, b)) \leq 1, b \geq 0}.
	\end{equation}

We are now in a position to estimate $\twistNadly 7(X)$.

\begin{prop}\label{Lemma: asymptotic for twN<=y(X) for m = 7}
	Suppose $y \ll X^{\frac{1}{12}}$. Then
	\begin{equation}
	\twistNadly 7(X) =  \twconst 7 X^{1/6} + O\parent{\max\parent{\frac{X^{1/6} \log y}{y}, X^{1/12} \log X \log^4 y }}
	\end{equation}
	for $X, y \geq 2$. The constant $\twconst 7$ is given in \eqref{Equation: twconst for m = 7}.
\end{prop}

\begin{proof}
	Substituting the asymptotic for $\cM 7(X; e)$ from \Cref{Lemma: asymptotic for M(X; e) for m = 7}(b) into the defining series \eqref{Equation: defining twNadly X} for $\twistNadly 7(X)$, we have
	\begin{equation}
		\twistNadly 7(X) = \sum_{n \leq y} \sum_{e \mid n} \mu\parent{n/e} \parent{\frac{\R 7 \cT 7(n) e X^{1/6}}{\zeta(2) n^3 \prod_{\ell \mid n} \parent{1 + \frac{1}{\ell}}} + O\parent{\frac{2^{\omega(e)} e^{1/2} X^{1/12} \log (e^6 X)}{n^{3/2}}}}.
	\end{equation}
	
	We handle the main term and the error of this expression separately. For the main term, we have
	\begin{equation}
	\begin{aligned}
	\sum_{n \leq y} \sum_{e \mid n} \mu\parent{n/e} \frac{\R 7 \cT 7(n) e X^{1/6}}{\zeta(2) n^3 \prod_{\ell \mid n} \parent{1 + \frac{1}{\ell}}} &= \frac{\R 7 X^{1/6}}{\zeta(2)} \sum_{n \leq y} \frac{\cT 7(n)}{n^3 \prod_{\ell \mid n} \parent{1 + \frac{1}{\ell}}} \sum_{e \mid n} \mu\parent{n/e} e \\
	&= \frac{\R 7 X^{1/6}}{\zeta(2)} \sum_{n \leq y} \frac{\varphi(n) \cT 7(n)}{n^3 \prod_{\ell \mid n} \parent{1 + \frac{1}{\ell}}}.
	\end{aligned}
	\end{equation}
	By \Cref{Lemma: bound on T(e) for m = 7}(f), we see
	\begin{equation}
	\frac{\varphi(n) \cT 7(n)}{n^3 \prod_{\ell \mid n} \parent{1 + \frac{1}{\ell}}} = O\parent{\frac{2^{\omega(n)}}{n^2}}.
	\end{equation}
	
	By \Cref{Corollary: tail of sum of f/n^sigma} and \Cref{Corollary: sum of 2^omega(n)}, we have
	\begin{equation}
	\sum_{n > y} \frac{2^{\omega(n)}}{n^2} \sim \frac{\log y}{\zeta(2) y}.
	\end{equation}
	\textit{A fortiori,}
	\begin{equation}
	\sum_{n > y} \frac{\varphi(n) \cT 7(n)}{n^3 \prod_{\ell \mid n} \parent{1 + \frac{1}{\ell}}} = O\parent{\sum_{n > y} \frac{2^{\omega(n)}}{n^2}} = O\parent{\frac{\log y}{y}},
	\end{equation}
	so the series
	\begin{equation}
	\sum_{n \geq 1} \frac{\varphi(n) \cT 7(n)}{n^3 \prod_{\ell \mid n} \parent{1 + \frac{1}{\ell}}} = \Q 7 \label{equation: sum for Q for m = 7}
	\end{equation}
	is absolutely convergent, and 
	\begin{equation}
	\begin{aligned}
	\sum_{n \leq y} \sum_{e \mid n} \mu\parent{n/e} \parent{\frac{\R 7 \cT 7(n) e X^{1/6}}{\zeta(2) n^3 \prod_{\ell \mid n} \parent{1 + \frac{1}{\ell}}}} &= \frac{\R 7 X^{1/6}}{\zeta(2)} \parent{\Q 7 + O\parent{\frac{\log y}{y}}} \\
	&= \twconst 7 X^{1/6} + O\parent{\frac{X^{1/6} \log y}{y}}.
	\end{aligned}
	\end{equation}
	
	As the summands of \eqref{equation: sum for Q for m = 7} constitute a nonnegative multiplicative arithmetic function, we can factor $\Q 7$ as an Euler product. For $p$ prime, \Cref{Lemma: bound on T(e) for m = 7} yields
	\begin{equation}\label{Equation: Q 7 terms}
	\Q 7 (p) \colonequals \sum_{a \geq 0} \frac{\varphi(p^a) \cT 7(p^a)}{p^{3a} \parent{1 + \frac{1}{p}}} = \begin{cases}
	1 + \displaystyle{\frac{2}{(p+1)^2}}, & \textup{if $p \equiv 1 \psmod{3}$ and $p \neq 7$;} \\
	13/6, & \text{if $p=3$;} \\
	63/8, & \text{if $p=7$;} \\ 
	1 & \text{else}.
	\end{cases}
	\end{equation}
	Thus
	\begin{equation}
	\Q 7 = \prod_{\textup{$p$ prime}} \Q 7 (p) = \Q 7 (3) \Q 7 (7) \prod_{\substack{p \neq 7 \ \textup{prime} \\ p \equiv 1 \psmod {3}}} \parent{1 + \frac{2}{(p+1)^2}}. \label{equation: product for Q for m = 7}
	\end{equation}
	
	We now turn to the error term. Since $y \ll X^{1/12}$, for $e \leq y$ we have $\log (e^6 X) \ll \log X$. We obtain
	\begin{equation} \label{equation: partial simplification twN<=y(X) for m = 7}
	\begin{aligned}
		&\sum_{n \leq y} \sum_{e \mid n} \mu\parent{n/e} O\parent{\frac{2^{\omega(n)} e^{1/2} X^{1/12} \log \parent{e^6 X}}{n^{3/2}}} \\
        =& O\parent{X^{1/12} \log X \sum_{n \leq y} \frac{2^{\omega(n)}}{n^{3/2}} \sum_{e \mid n} \abs{\mu\parent{n/e}} e^{1/2} }.
	\end{aligned}
	\end{equation}
	The inequality
	\begin{equation}
	\sum_{e \mid n} \abs{\mu\parent{n/e}} e^{1/2} < 2^{\omega(n)} \sqrt{n}
	\end{equation}
	implies
	\begin{equation}\label{Equation: sum of 2^omega/n^3/2 * sum for m = 7}
	\sum_{n \leq y} \frac{2^{\omega(n)}}{n^{3/2}} \sum_{e \mid n} \abs{\mu\parent{n/e}} e^{1/2} = O\parent{\sum_{n \leq y} \frac{4^{\omega(n)}}{n}}.
	\end{equation}
	But \Cref{Theorem: Asymptotics of k^omega(n)} together with Abel summation imply that \eqref{Equation: sum of 2^omega/n^3/2 * sum for m = 7} is $O(\log^4 y)$, yielding our desired result.
\end{proof}

We emphasize that \eqref{Equation: Q 7 terms} and \eqref{equation: product for Q for m = 7} from the proof of \Cref{Lemma: asymptotic for twN<=y(X) for m = 7} have given us the following Euler product expansion for $\Q 7$:
\begin{equation}
\Q 7 = \Q 7(3) \Q 7 (7) \prod_{\substack{p \neq 7 \ \textup{prime} \\ p \equiv 1 \psmod {3}}} \parent{1 + \frac{2}{(p+1)^2}}, 
\end{equation}
where $\Q 7 (3) = 13/6$ and $\Q 7 (7) = 63/8$.

We are now ready to prove \Cref{Intro Theorem: asymptotic for twN(X) for m of genus $0$} when $m = 7$, which we restate here with an improved error term in the notations we have established. We give two proofs of this important statement. The first proof is an easy argument using \Cref{Proposition: fundamental sieve for for m = 7} and \Cref{Lemma: asymptotic for twN<=y(X) for m = 7}. The second proof requires deriving a bound on $\twistNadgy 7 (X)$ which is in some sense superfluous; however, this proof is more typical of our arguments in the remaining chapters of this thesis.

\begin{theorem}\label{Theorem: asymptotic for twN(X) for m = 7}
	We have 
	\begin{equation}
	\twistNad 7(X) = \twconst 7 X^{1/6} + O(X^{1/12} \log^5 X)
	\end{equation}
	for $X \geq 2$. The constant $\twconst 7$ is given in \eqref{Equation: twconst for m = 7}.
\end{theorem}

\begin{proof}[First Proof of \Cref{Theorem: asymptotic for twN(X) for m = 7}]
	Let $X > 0$, and let $y$ be slightly larger than $\frac{3^{5/4} \cdot 7^{9/2}}{2^{1/6}} \cdot X^{1/12}$. By \Cref{Proposition: fundamental sieve for for m = 7}, 
	\begin{equation}
	\twistNadly 7 (X) = \twistNad 7 (X),
	\end{equation}
	and the result is now immediate from \Cref{Lemma: asymptotic for twN<=y(X) for m = 7}.
\end{proof}

We now bound $\twistNadgy 7(X)$ as a step towards our alternate proof. The proof below is somewhat cleaner than that given in \cite{Molnar-Voight}.

\begin{lemma}\label{Lemma: bound on twN>y(X) for m = 7}
	We have
	\begin{equation}
	\twistNadgy 7(X) = O\parent{\frac{X^{1/6} \log y}{y}}
	\end{equation}
        for $X, y \geq 2$.
\end{lemma}

\begin{proof}
	By \Cref{Lemma: bound on T(e) for m = 7}(f), \Cref{Lemma: asymptotic for M(X; e) for m = 7}, and \Cref{Proposition: bounding the summands of twN calE X},
	\begin{equation}
	\abs{\twistNadgy 7(X)} \ll X^{1/6}\sum_{n>y}\frac{2^{\omega(n)}}{n^2}\prod_{p\mid n}\parent{1+\frac1p}.
	\end{equation}
	The result follows from \Cref{Lemma: weighted 2omega divisor factor} with $\alpha=1$ and $\beta=2$.
\end{proof}

\begin{proof}[Second Proof of \Cref{Theorem: asymptotic for twN(X) for m = 7}]
	Let $y$ be a positive quantity with $y \ll X^{1/12}$; in particular, $\log y \ll \log X$. \Cref{Lemma: asymptotic for twN<=y(X) for m = 7} and \Cref{Lemma: bound on twN>y(X) for m = 7} together tell us
	\begin{equation}
	\twistNad 7(X) = \twconst 7 X^{1/6} + O\parent{\max\parent{\frac{X^{1/6} \log y}{y}, X^{1/12} \log X \log^4 y}}.
	\end{equation}
	Now letting $y = X^{1/12}$, our claim follows.
\end{proof}

\subsection*{$L$-series}\label{Subection: L-series for m = 7}\label{Subsection: L-series for m = 7}

To conclude this section, we set up \cref{Section: Working over the rationals for m = 7} by interpreting
\Cref{Theorem: asymptotic for twN(X) for m = 7} in terms of Dirichlet series. Recall \eqref{Equation: twisth calE}, \eqref{Equation: hQ calE}, \eqref{Equation: defining twistL calE X}, and \eqref{Equation: defining LQ calE X}.

\begin{cor}\label{Corollary: twL(s) has a meromorphic continuation for m = 7}
    The following statements hold.
    \begin{enumerate}
        \item The Dirichlet series $\twistLad 7(s)$ has abscissa of (absolute) convergence $\sigma_a=\sigma_c = 1/6$ and has a meromorphic continuation to the region
	\begin{equation}
	\set{s = \sigma + i t \in \bbC : \sigma > 1/12}. \label{equation: domain of twistL(s) for m = 7}
	\end{equation}
        \item The function $\twistLad 7(s)$ has a simple pole at $s = 1/6$ with residue 
	\begin{equation} 
        \res_{s=\frac{1}{6}} \twistLad 7(s) = \frac{\twconst 7}{6} ;
        \end{equation}
        it is holomorphic elsewhere on the region \eqref{equation: domain of twistL(s) for m = 7}.
    \item We have
    \begin{equation}
    \mu_{\twistLad 7}(\sigma) \leq \max\{0, 2-12\sigma\}
     \end{equation}
    for $\sigma > 1/12$.
    \end{enumerate}
\end{cor}

\begin{proof}
	We first prove part (a). Let $s = \sigma + i t \in \bbC$ be given with $\sigma > 1/6$. Abel summation yields
	\begin{equation}
	\begin{aligned}
		\sum_{n \leq X} \twisthad 7(n) n^{-s} &= \twistNad 7(X) X^{-s} + s \int_1^X \twistNad 7(u) u^{-s-1} \,\mathrm{d}u \\
		&= O\parent{X^{1/6 - \sigma} + s \int_1^X u^{- 5/6 - \sigma} \,\mathrm{d}u};
	\end{aligned}
	\end{equation}
	as $X \to \infty$ the first term vanishes and the integral converges. Thus, when $\sigma > 1/6$,
	\begin{equation}
	\sum_{n \geq 1} \twisthad 7(n) n^{-s} = s \int_1^\infty \twistNad 7(u) u^{-1-s}\,\mathrm{d}u
	\end{equation}
	and this integral converges. A similar argument shows that the sum defining $\twistLad 7(s)$ diverges when $\sigma < 1/6$. We have shown $\sigma_c = 1/6$ is the abscissa of convergence for $\twistLad 7(s)$, but as $\twisthad 7(n) \geq 0$ for all $n$, $1/6$ is also the abscissa of \emph{absolute} convergence $\sigma_a=\sigma_c$.
	
	Now define $\twistLadR 7(s)$ so that
	\begin{equation}
		\twistLad 7(s) = \twconst 7 \zeta(6s) + \twistLadR 7(s).\label{equation: Defining twLR for m = 7}
	\end{equation}		
	Abel summation and the substitution $u \mapsto u^{1/6}$ yields the following equality for $\sigma>1/6$:
	\begin{equation}
	\zeta(6s) = s \int_1^\infty \floor{u^{1/6}} u^{- 1 - s}\,\textrm{d}u = s \int_1^\infty \parent{u^{1/6} + O(1)} u^{- 1 - s} \,\textrm{d}u.
	\end{equation}
	Let 
	\begin{equation}
	\chi_6(n) \colonequals \begin{cases} 1, & \text{if} \ n = k^6 \ \text{for some} \ k \in \bbZ; \\
	0, & \text{else.}
	\end{cases}
	\end{equation}
	Then
	\begin{equation}
	\begin{aligned}
	\twistLadR 7(s) &= \sum_{n \geq 1} \parent{\twisthad 7(n) - \twconst 7 \chi_6(n)} n^{-s}\\
	&= s \int_1^\infty \parent{\twistNad 7(u) - \twconst 7 \floor{u^{1/6}}} u^{-1-s} \,\textrm{d}u \label{equation: twistLR(s) as an integral for m = 7}
	\end{aligned}
	\end{equation}
	when $\sigma > 1/6$. But then for any $\epsilon > 0$,
	\begin{equation}
	\twistNad 7(u) - \twconst 7 \floor{u^{1/6}} = O(u^{1/12 + \epsilon}) \label{equation: twistN(u) - ctw7 floor(u^1/6) for m = 7}
	\end{equation}
	by \Cref{Theorem: asymptotic for twN(X) for m = 7}. 
 Substituting \eqref{equation: twistN(u) - ctw7 floor(u^1/6) for m = 7} into \eqref{equation: twistLR(s) as an integral for m = 7}, we obtain
	\begin{align}
	\twistLadR 7(s) = s \int_1^\infty \parent{\twistNad 7(u) - \twconst 7 \floor{u^{1/6}}} u^{-1-s}\,\textrm{d}u &= O\left(s \int_1^\infty u^{-11/12 - \sigma + \epsilon} \,\textrm{d}u\right) \label{equation: bound on twLR for m = 7}
	\end{align}
	where the integral converges whenever $\sigma > 1/12 + \epsilon$. Letting $\epsilon \to 0$, we obtain an analytic continuation of $\twistLadR 7(s)$ to the region \eqref{equation: domain of twistL(s) for m = 7}.
	
    We proceed to part (b). The Dirichlet series $\zeta(6s)$ has meromorphic continuation to $\bbC$ with a simple pole at $s=1/6$ with residue $1/6$. Thus looking back at \eqref{equation: Defining twLR for m = 7}, we find that 
	\begin{equation}
	\twistLad 7(s) = \twconst 7 \zeta(6s) + s \int_1^\infty \parent{\twistNad 7(u) - \twconst 7 \floor{u^{1/6}}} u^{-1-s} \,\textrm{d}u
	\end{equation}
	when $\sigma > 1/6$, but in fact the right-hand side of this equality defines a meromorphic function on the region \eqref{equation: domain of twistL(s) for m = 7} with a simple pole at $s = 1/6$ and no other poles in this region. 
 
    Finally, we prove part (c). For fixed $1/12<\sigma_0<1/6$, \eqref{equation: bound on twLR for m = 7} gives $\mu_{\twistLadR 7}(\sigma_0)\leq 1$, while absolute convergence gives $\mu_{\twistLadR 7}(\sigma)=0$ for $\sigma>1/6$. By subconvexity of $\mu$ \cite[Theorem II.1.20]{Tenenbaum},
    \begin{equation}
        \mu_{\twistLadR 7}(\sigma)\leq 2-12\sigma \qquad (1/12<\sigma<1/6).
    \end{equation}
    Since \Cref{Theorem: muzeta(sigma)} gives $\mu_{\zeta_6}(\sigma)\leq \frac{13}{42}(1-6\sigma)<2-12\sigma$ on the same strip, \eqref{equation: Defining twLR for m = 7} and \Cref{Proposition: mu acts like a negative valuation} yield $\mu_{\twistLad 7}(\sigma)\leq\max\{0,2-12\sigma\}$ for $\sigma>1/12$.
\end{proof}

\section{Estimates for rational isomorphism classes for \texorpdfstring{$m = 7$}{m = 7}}\label{Section: Working over the rationals for m = 7}

In \cref{Section: Estimating twN(X) for m = 7}, we counted the number of elliptic curves over $\bbQ$ with a (cyclic) $7$-isogeny up to quadratic twist (\Cref{Theorem: asymptotic for twN(X) for m = 7}). In this section, we count all isomorphism classes over $\bbQ$ by enumerating over twists using Landau's Tauberian theorem (\Cref{Theorem: Landau's Tauberian theorem}). We first describe the analytic behavior of $\LQad 7 (s)$.

\begin{theorem}\label{Theorem: relationship between twistL(s) and L(s) for m = 7}
The following statements hold.
\begin{enumalph}
	\item The Dirichlet series $\LQad 7(s)$ has a meromorphic continuation to the region \eqref{equation: domain of twistL(s) for m = 7} with a double pole at $s = 1/6$ and no other singularities on this region. 
	\item The principal part of $\LQad 7(s)$ at $s = 1/6$ is
	\begin{equation}
	\frac{1}{3 \zeta(2)} \parent{\frac{\twconst 7}{6} \parent{s - \frac 16}^{-2} + \parent{\ell_{7, 0} + \twconst 7 \parent{\gamma - \frac{2 \zeta\prm(2)}{\zeta(2)}}} \parent{s - \frac{1}{6}}\inv},
	\end{equation}
	where $\twconst 7$ is given in \eqref{Equation: twconst for m = 7}, and
	\begin{equation}
	\ell_{7, 0} \colonequals \twconst 7 \gamma + \frac {1}6 \int_1^\infty \parent{\twistNad 7(u) - \twconst 7 \floor{u^{1/6}}} u^{-7/6} \,\mathrm{d}u\label{Equation: Definition of ell0 for m = 7}
	\end{equation}
	is the constant term of the Laurent expansion for $\twistLad 7(s)$ around $s = 1/6$.
	\end{enumalph}
\end{theorem}

\begin{proof}
	For part (a), since $\zeta(s)$ is nonvanishing when $\sigma > 1$, the ratio $\zeta(6s)/\zeta(12s)$ is a meromorphic function for $\sigma > 1/12$. But \Cref{Corollary: twL(s) has a meromorphic continuation for m = 7} gives a meromorphic continuation of $\twistLad 7(s)$ to the region \eqref{equation: domain of twistL(s) for m = 7}. The function $\LQad 7(s)$ is a product of these two meromorphic functions on \eqref{equation: domain of twistL(s) for m = 7}, and so it is a meromorphic function on this region. The holomorphy and singularity for $\LQad 7(s)$ then follow from those of $\twistLad 7(s)$ and $\zeta(s)$. 
	
	We deduce part (b) by computing Laurent expansions. We readily verify
\begin{equation}\label{Equation: Laurent expansion for zeta(6s)/zeta(12s)}
\frac{\zeta(6s)}{\zeta(12s)} = \frac{1}{\zeta(2)}\parent{\frac{1}{6}\parent{s - \frac 16}\inv + \parent{\gamma - \frac{2 \zeta\prm(2)}{\zeta(2)}} + \ldots},
\end{equation}
whereas the Laurent expansion for $\twistLad 7(s)$ at $s = 1/6$ begins
\begin{equation}
\twistLad 7(s) = \frac{\twconst 7}{6}\parent{s - \frac 16}\inv + \ell_{7, 0} + \dots,
\end{equation}
with $\ell_{7, 0}$ given by \eqref{Equation: Definition of ell0 for m = 7}. Multiplying the Laurent series tails gives the desired result.
\end{proof}

Using \Cref{Theorem: relationship between twistL(s) and L(s) for m = 7}, we deduce the following lemma.

\begin{lemma}\label{Lemma: Delta NQ(n) is admissible for m = 7}
	For every $\delta$ with $0<\delta<7/97$, the sequence $\parent{\hQad 7(n)}_{n \geq 1}$ is admissible \textup{(\Cref{Definition: admissible sequences})} with parameters
    \begin{equation}
        \parent{1/6,\delta,97\delta/7}.
    \end{equation}
\end{lemma}

\begin{proof}
	We check each condition in \Cref{Definition: admissible sequences}. Since $\hQad 7(n)$ counts objects, we indeed have $\hQad 7(n) \in \Z_{\geq 0}$. 
	
	For (i), \Cref{Corollary: twL(s) has a meromorphic continuation for m = 7} tells us that $\twistLad 7(s)$ has $1/6$ as its abscissa of absolute convergence. Likewise, $\displaystyle{\frac{\zeta(6s)}{\zeta(12s)}}$ has $1/6$ as its abscissa of absolute convergence. By \Cref{Theorem: relationship between twistL(s) and L(s)}(b), 
	\begin{equation}
	\LQad 7(s) = \frac{2\zeta(6s) \twistLad 7(s)}{\zeta(12s)},
	\end{equation}
	and by \Cref{Theorem: product of Dirichlet series converges} this series converges absolutely for $\sigma > 1/6$, so the abscissa of absolute convergence for $\LQad 7(s)$ is at most $1/6$. But for $\sigma < 1/6$, $\LQad 7(\sigma) > \twistLad 7(\sigma)$ by termwise comparison of coefficients, so the Dirichlet series for $\LQad 7(s)$ diverges when $\sigma < 1/6$, and (i) holds with $\sigma_a = 1/6$.

For (ii), \Cref{Corollary: twL(s) has a meromorphic continuation for m = 7} tells us that $\twistLad 7(s)$ has a meromorphic continuation when $\sigma=\repart(s)>1/12$; on the other hand, as $\zeta(12s)$ is nonvanishing for $\sigma > 1/12$, we see that $\zeta(6s)/\zeta(12s)$ has a meromorphic continuation to $\sigma>1/12$. Since $\delta<7/97<1/12$, condition (ii) holds on the smaller half-plane $\sigma>1/6-\delta$. (The only pole of $\LQad 7(s)/s$ there is the double pole at $s = 1/6$ indicated in \Cref{Theorem: relationship between twistL(s) and L(s) for m = 7}(b).)
	
For (iii), suppose first that $1/6-\delta<\sigma\leq1/6$, and put $d=1/6-\sigma$, so $0\leq d<\delta$. By \Cref{Corollary: twL(s) has a meromorphic continuation for m = 7},
\begin{equation}
    \mu_{\twistLad 7}(\sigma) \leq 12d.
\end{equation}
Writing $\zeta_a(s)=\zeta(as)$, \Cref{Theorem: muzeta(sigma)} gives
\begin{equation}\label{Equation: bound on zeta6(s)}
    \mu_{\zeta_6}(\sigma)=\mu_\zeta(6\sigma)\leq \frac{13}{42}(1-6\sigma)=\frac{13}{7}d,
\end{equation}
and, as $\zeta(12s)^{-1}$ is absolutely convergent for $\sigma>1/12$, we have $\mu_{\zeta_{12}^{-1}}(\sigma)=0$. Therefore \Cref{Theorem: relationship between twistL(s) and L(s)}(b) and \Cref{Proposition: mu acts like a negative valuation} imply
\begin{equation}
    \mu_{\LQad 7}(\sigma) \leq \parent{12+\frac{13}{7}}d=\frac{97}{7}d<\frac{97}{7}\delta.
\end{equation}
For $\sigma>1/6$, all three factors are absolutely convergent, so $\mu_{\LQad 7}(\sigma)=0$. Thus condition (iii) holds with $\xi=97\delta/7$, proving the claim.
\end{proof}

We now prove \Cref{Intro Theorem: asymptotic for NQ(X) for 5 < m <= 9} when $m = 7$, which we restate here in this special case in our established notation.

\begin{theorem}\label{Theorem: asymptotic for NQ(X) for m = 7}
	We define
        \begin{equation}\label{Equation: constants for m = 7}
	\begin{aligned}
	\constad 7 &\colonequals \frac{\Q 7\R 7}{3 \zeta(2)^2}, \\
	\constadprm 7 &\colonequals \frac{2}{\zeta(2)}\parent{\ell_{7, 0} + \twconst 7 \parent{\gamma - 1 - \displaystyle{\frac{2 \zeta\prm(2)}{\zeta(2)}}}},
	\end{aligned}
        \end{equation}
	where $\twconst 7$ is defined in \eqref{Equation: twconst for m = 7}, and $\ell_{7, 0}$ is defined in \eqref{Equation: Definition of ell0 for m = 7}. Then for all $\epsilon > 0$, we have
		\begin{equation}
	\NQad 7(X) = \constad 7 X^{1/6} \log X + \constadprm 7 X^{1/6} + O\parent{X^{38/291 + \epsilon}}
	\end{equation}
	for $X \geq 1$. The implicit constant depends only on $\epsilon$. 
\end{theorem}

\begin{proof}
	Apply \Cref{Theorem: Landau's Tauberian theorem} and let $\delta\to7/97$ in \Cref{Lemma: Delta NQ(n) is admissible for m = 7}; the exponent tends to $1/6-7/194=38/291$.
\end{proof}

\begin{remark}\label{Remark: true bound for twistN(X) for m = 7}
	The loss from the twist-counting error $O(X^{1/12}\log^5 X)$ to the exponent $38/291$ above comes from the vertical growth of the analytically continued remainder in \eqref{equation: twistLR(s) as an integral for m = 7}. We believe an elementary argument applying Dirichlet's hyperbola method to \Cref{Theorem: relationship between twistL(s) and L(s)}(a) could achieve a comparable, and perhaps modestly improved, error term without this analytic-growth loss.
\end{remark}

\section{Computations for \texorpdfstring{$m = 7$}{m = 7}}\label{Section: Computations for m = 7}

In this section, we furnish computations that render \Cref{Theorem: asymptotic for twN(X) for m = 7} and \Cref{Theorem: asymptotic for NQ(X) for m = 7} completely explicit.

\subsection*{Enumerating elliptic curves with a cyclic $7$-isogeny}\label{Subsection: Computing elliptic curves with 7-isogeny}

We begin by outlining an algorithm for computing all elliptic curves (up to quadratic twist) with twist height at most $X$ that admit a cyclic 7-isogeny. In a nutshell, we iterate over possible factorizations $e^3 m$ with $m$ cubefree to find all $7$-groomed pairs $(a, b)$ for which $\C 7(a, b) = e^3 m$, then check if $\twistheight(\A 7(a, b), \B 7(a, b)) \leq X$.

In detail, our algorithm proceeds as follows. 
\begin{enumalg}
	\item We list all primes $p \equiv 1 \psmod 3$ up to $\parent{3^3\cdot 7^{18}X/4}^{1/6}$ (this bound follows from \Cref{Theorem: Controlling size of twist minimality defect for m = 7}(b)).
	\item For each pair $(a, b) \in \bbZ^2$ with $b > 0$, $\gcd(a, b) = 1$, and $\C 7(a, b)$ coprime to $3$ and less than $Y$, we compute $\C 7(a, b)$. We organize the results into a lookup table, so that for each $c$ we can find all pairs $(a, b)$ with $b > 0$, $\gcd(a, b) = 1$, and $\C 7(a, b) = c$. We append $1$ to our table with lookup value $(1, 0)$. For each $c$ in our lookup table, we record whether $c$ is cubefree by sieving against the primes we previously computed.
	\item For positive integer pairs $(e_0, m)$ with $4e_0^{12}m^6 \leq 3^3\cdot 7^{18}X$ and $m$ cubefree, we find all $7$-groomed pairs $(a, b) \in \bbZ^2$ with $\C 7(a, b) = e_0^3 m$. If $\gcd(e_0, 3) = \gcd(m, 3) = 1$, we can do this as follows. If $e_0^3 < Y$, we iterate over $7$-groomed pairs $(a_e, b_e)$ and $(a_m, b_m)$ yielding $\C 7(a_e, b_e) = e_0^3$ and $\C 7(a_m, b_m) = m$ respectively, and taking the product 
	\begin{equation}
	(a_e + b_e\parent{-1+3\zeta_6}) (a_m + b_m\parent{-1+3\zeta_6}) = a + b \parent{-1+3\zeta_6} \in \bbZ[3\zeta_6]
	\end{equation}
	as in the proof of \Cref{Lemma: Algebraic structure of Z[3 zeta6] and C}. If $e_0^3 > Y$, we iterate over $7$-groomed pairs $(a_e\prm, b_e\prm)$ with $\C 7(a_e\prm, b_e\prm) = e_0$ instead of over $7$-groomed pairs $(a_e, b_e)$, and compute 
	\begin{equation}
	(a_e\prm + b_e\prm\parent{-1+3\zeta_6})^3 (a_m + b_m\parent{-1+3\zeta_6}) = a + b \parent{-1+3\zeta_6} \in \bbZ[3\zeta_6].
	\end{equation}
	If $\gcd(e_0, 3) > 1$ or $\gcd(m, 3) > 1$, we perform the steps above for the components of $e_0$ and $m$ coprime to $3$, and then postmultiply by those $7$-groomed pairs $(a_3, b_3) \in \bbZ^2$ with $\C 7(a_3, b_3)$ an appropriate power of $3$ (which is no greater than $27$, by \Cref{Lemma: bound on T(e) for m = 7}).
	\item For each pair $(a, b)$ with $\C 7(a, b) = e_0^3 m$, obtained in the previous step, we compute $\rawheight(\A 7(a, b), \B 7(a, b))$. We compute the $3$-component of the twist minimality defect $e_3$, the $7$-component of the twist minimality defect $e_7$, and thereby compute the twist minimality defect $e = \lcm(e_0, e_3, e_7)$. We compute the twist height using the reduced pairs $(\A 7(a, b) / e^2, \abs{\B 7(a, b)} / e^3)$. If this result is less than or equal to $X$, we report $(a, b)$, together with their twist height and any auxiliary information we care to record.
\end{enumalg}

We list the first few twist minimal elliptic curves admitting a (cyclic) $7$-isogeny in Table \ref{table:firstfew7}.

\jvtable{table:firstfew7}{\small
\rowcolors{2}{white}{gray!10}
\begin{tabular}{c c|c c}
$(A, B)$ & $(a, b)$ & $\twistheight(E)$ & $\twistdefect(E)$ \\
\hline\hline
$(-3, 62)$ & $(14, 5)$ & $103788$ & $1029$ \\
$(13, 78)$ & $(21, 4)$ & $164268$ & $1029$ \\
$(37, 74)$ & $(42, 1)$ & $202612$ & $1029$ \\
$(-35, 98)$ & $(0, 1)$ & $259308$ & $21$ \\
$(45, 18)$ & $(35, 2)$ & $364500$ & $1029$ \\
$(-43, 166)$ & $(7, 13)$ & $744012$ & $3087$ \\
$(-75, 262)$ & $(-7, 8)$ & $1853388$ & $1029$ \\
$(-147, 658)$ & $(-56, 1)$ & $12706092$ & $1029$ \\
$(-147, 1582)$ & $(7, 6)$ & $67573548$ & $343$ \\
$(285, 2014)$ & $(28, 3)$ & $109517292$ & $343$ \\
$(-323, 2242)$ & $(-21, 10)$ & $135717228$ & $1029$ \\
$(-395, 3002)$ & $(-63, 2)$ & $246519500$ & $1029$ \\
$(-155, 3658)$ & $(21, 11)$ & $361286028$ & $1029$ \\
$(357, 5194)$ & $(7, 1)$ & $728396172$ & $21$ \\
$(-595, 5586)$ & $(-14, 1)$ & $842579500$ & $63$ \\
$(285, 5662)$ & $(91, 1)$ & $865572588$ & $1029$ \\
$(-603, 5706)$ & $(-28, 11)$ & $879077772$ & $1029$ \\
\end{tabular}
}{$E \in \twistE$ with a cyclic $7$-isogeny and $\twht E \leq 10^{9}$}

Running this algorithm out to $X = 10^{42}$ in Python took us approximately $34$ CPU hours on a single core, producing 4\,582\,079 
elliptic curves admitting a (cyclic) 7-isogeny in $\twistE_{\leq 10^{42}}$. To check the accuracy of our code, we confirmed that the $j$-invariants of these curves are distinct. We also confirmed that the 7-division polynomial of each curve has a linear or cubic factor over $\bbQ$; this took 3.5 CPU hours. For $X = 10^{42}$, we have
\begin{equation}\label{Equation: Ratio for twN7}
\frac{ \twistNad 7(10^{42})}{\twconst 7 (10^{42})^{1/6}} = 0.99996\ldots,
\end{equation}
which is close to 1. We compute $\twconst 7 = \Q 7 \R 7 /\zeta(2)$ below.

Reorganizing \Cref{Theorem: relationship between twistL(s) and L(s)}(a), we find
\begin{equation}
\NQad 7(X)=2\sum_{n\leq X}\twisthad 7(n)\,\kfree 2\parent{(X/n)^{1/6}}.
\end{equation}
Letting $X = 10^{42}$ and using our list of 4\,582\,079  elliptic curves admitting a (cyclic) 7-isogeny, we compute that there are $88\,157\,174$ elliptic curves admitting a (cyclic) 7-isogeny in $\scrE_{\leq 10^{42}}$.

\subsection*{Computing \texorpdfstring{$\twconst 7$}{ctw7}}

In this subsection, we estimate the constant $\twconst {7}$ appearing in \Cref{Theorem: asymptotic for twN(X) for m = 7} by estimating $\Q {7}$ and $\R {7}$.

We begin with $\Q 7$, given by \eqref{equation: product for Q for m = 7}. Truncating the Euler product as a product over $p \leq Y$ gives us a lower bound
\begin{equation}
\Qly 7 \colonequals \frac{273}{16} \prod_{\substack{7 < p \leq Y \\ p \equiv 1 \psmod 3}} \parent{1 + \frac{2}{(p+1)^2}}
\end{equation}
for $\Q 7$. To obtain an upper bound, we compute
\begin{equation}
	\Q 7 < \Qly 7 \exp\parent{2 \sum_{\substack{p > Y \\ p \equiv 1 \psmod 3}} \frac{1}{p^2+1}}.
\end{equation}
For $e, a_0 \in \bbZ$ and $X \geq 1$, write
\begin{equation}\label{Equation: prime counting function}
    \pi(X; e, a_0) \colonequals \# \set{p \in \bbZ_{>0} : p \ \text{prime}, \ p \leq X, \ p \equiv a_0 \psmod e}. 
\end{equation}
Suppose $Y \geq 8 \cdot 10^9$. Using Abel summation together with explicit estimates for $\pi(Y; 3, 1)$ furnished by Bennett--Martin--O'Bryant--Rechnitzer \cite[Theorem 1.4]{Bennett-Martin-OBryant-Rechnitzer}, we obtain
\begin{equation} 	\begin{aligned}
\sum_{\substack{p > Y \\ p \equiv 1 \psmod 3}} \frac{1}{p^2+1} &= -\frac{\pi(Y;3,1)}{Y^2 + 1} + 2 \int_Y^\infty \frac{\pi(u; 3, 1) u} {(u^2 + 1)^2} \,\mathrm{d}u \\
&< -\frac{Y}{2\parent{Y^2 + 1} \log Y} + \parent{\frac 1{\log Y} + \frac{5}{2 \log^2 Y}} \int_Y^\infty \frac{u^2} {(u^2 + 1)^2} \,\mathrm{d}u \\
&= \frac 12 \parent{\frac{5Y}{2 (Y^2 + 1) \log Y} + \parent{\frac 1{\log Y} + \frac{5}{2 \log^2 Y}} \parent{\frac{\pi}{2} - \tan\inv(Y)}}
\end{aligned}\end{equation} 
so
\begin{equation}
\Q 7 <\Qly 7 \cdot \exp\parent{\frac{5Y}{2 (Y^2 + 1) \log Y} + \parent{\frac 1{\log Y} + \frac{5}{2 \log^2 Y}} \parent{\frac{\pi}{2} - \tan\inv(Y)}}.
\end{equation}
In particular, letting $Y = 8 \cdot 10^9$, we compute
\begin{equation}
17.460\,405\,236\,640\,390 < \Q 7 < 17.460\,405\,236\,985\,885.
\end{equation}

We now turn our attention to $\R 7$, given in \eqref{Equation: Rm}. We observe 
\begin{equation}\label{Equation: R7 inside a rectangle}
\calR 7(1) \subseteq [-0.677, 0.677] \times [0, 0.078],
\end{equation}
so we can estimate $\calR 7(1)$ by performing rejection sampling on the rectangle given in \eqref{Equation: R7 inside a rectangle}, which has area $0.105612$. 

We find $r_7 \colonequals 243\,228\,665\,965$ of our first $s_7 \colonequals 595\,055\,000\,000$ samples lie in $\calR 7(1)$, so
\begin{equation}
\R 7 \approx 0.105612 \cdot\frac{r_7}{s_7}=0.043168893413\ldots
\end{equation}
with standard error $6.73053\ldots\cdot 10^{-8}$. This took 11 CPU weeks to compute. Using the midpoint of our interval for $\Q 7$, we obtain $\twconst 7\approx0.4582228478\ldots$ with Monte Carlo standard error $7.15\cdot10^{-7}$.

\subsection*{Computing \texorpdfstring{$\constad {7}$ and $\constadprm 7$}{and c7 and c7'}}

In this subsection, we estimate the constants $\constad 7$ and $\constadprm 7$, which are defined in \eqref{Equation: constants for m = 7} and used in \Cref{Theorem: asymptotic for NQ(X) for m = 7}.

We have the identity $\constad 7 = \twconstad 7 / 3 \zeta(2)$, so $\constad 7\approx0.092\,855\,36263\ldots$ with Monte Carlo standard error $1.45\cdot10^{-7}$.

We now turn our attention to $\constadprm 7$. Truncating \eqref{Equation: Definition of ell0 for m = 7} at $X=10^{42}$ and using our list of $4\,582\,079$ twist-minimal curves gives
\begin{equation}
\ell_{7,0}\approx-0.463\,534\,394.
\end{equation}
In \Cref{Theorem: asymptotic for twN(X) for m = 7}, we have shown that for some $M>0$ and all $u>X$,
\begin{equation}
\abs{\twistNad 7(u)-\twconst 7\floor{u^{1/6}}}<M u^{1/12}\log^5u.
\end{equation}
Thus the omitted tail in $\ell_{7,0}$ satisfies
\begin{equation}
\begin{aligned}
&\abs{\frac16\int_X^\infty \parent{\twistNad 7(u)-\twconst 7\floor{u^{1/6}}}u^{-7/6}\,\mathrm{d}u}\\
&\qquad<2MX^{-1/12}(\log^5X+60\log^4X+2880\log^3X\\
&\qquad\qquad+103680\log^2X+2488320\log X+29859840).
\end{aligned}
\end{equation}
We do not know the exact value of $M$. For the large values of $u$ examined in our computation up to $10^{42}$, we find
\begin{equation}
-5.11\cdot10^{-6}\leq\frac{\twistNad 7(u)-\twconst 7\floor{u^{1/6}}}{u^{1/12}\log^5u}\leq6.29\cdot10^{-7}.
\end{equation}
If we assume these inequalities continue to hold for $u>10^{42}$, the resulting tail error is less than $56.73$.

We can do better by sidestepping the logarithms. Taking $\epsilon\colonequals10^{-4}$, over the values examined up to $10^{42}$ we find
\begin{equation}
-1.2174\leq\frac{\twistNad 7(u)-\twconst 7\floor{u^{1/6}}}{u^{1/12+\epsilon}}\leq0.52272.
\end{equation}
If we assume these inequalities continue to hold for larger $u$, the resulting tail error is less than $7.79\cdot10^{-4}$. Propagating the Monte Carlo standard error in $\twconstad 7$ through the truncated integral contributes $1.09\cdot10^{-5}$, giving reported uncertainties $56.73$ or $7.90\cdot10^{-4}$, respectively.

Using \eqref{Equation: constants for m = 7}, we obtain
\begin{equation}
\constadprm 7\approx-0.164\,049\,781,
\end{equation}
with reported uncertainties $68.98$ or $9.59\cdot10^{-4}$ under the two extrapolations above. As a sanity check,
\begin{equation}
\frac{\NQad 7(10^{42})}{10^7}-42\constad 7\log10=-0.164\,192\,299\ldots\approx\constadprm 7,
\end{equation}
which agrees to three decimal places with the estimate for $\constadprm 7$ we gave above.

%% file: Ch5_m10_25.tex
\chapter{Counting elliptic curves with a cyclic \texorpdfstring{$m$}{m}-isogeny for \texorpdfstring{$m = 10, 25$}{m = 10, 25}}\label{Chapter: m = 10 and 25}

In this chapter, we prove \Cref{Intro Theorem: asymptotic for NQ(X) for m > 9} (\Cref{Theorem: asymptotic for NQ(X) for m = 10 and 25}), and \Cref{Intro Theorem: asymptotic for twN(X) for m of genus $0$} (\Cref{Theorem: asymptotic for twN(X) for m = 10 and 25}) when $m = 10, 25$. These results are new, but our arguments mirror those in \cref{Chapter: m = 7} (and thus also \cite{Molnar-Voight}), and we encourage anyone reading to skim them on a first perusal of this thesis. There is one major new complication, however: viewed as an elliptic curve over $\bbQ(t)$, the elliptic surfaces describing elliptic curves with a cyclic $5$-isogeny, with a cyclic $10$-isogeny, and with a cyclic $25$-isogeny exhibit potential type III additive reduction rather than potential type II additive reduction. This forces us to change how we define $\cT m (e)$ and related functions (see \Cref{Definition: T(e) for m = 10 and 25}), and changes the details of our sieving somewhat.

Although we are unable to derive asymptotics for $\twNeq 5 (X)$ or $\NQeq 5 (X)$, this case is structurally similar enough to $m = 10, 25$ that we opt to provide some preliminary information about the structure of $\cM 5 (X)$ and related functions.

In \cref{Section: Establishing notation when m = 10 and 25}, for $m \in \set{5, 10, 25}$, we establish notations pertaining to $\f m (t)$ and $\g m (t)$ which will be used throughout the remainder of the chapter. In \cref{Section: Estimating twN(X) for m = 5 and 10 and 25}, we develop bounds relating the twist minimality defect to the greatest common divisor of $\f m (t)$ and $\g m (t)$. In \cref{Section: Estimating twN(X) for m = 5 and 10 and 25}, we apply the framework developed in \cref{Section: Our approach revisited} to prove 
\Cref{Intro Theorem: asymptotic for twN(X) for m of genus $0$} for $m = 10, 25$. In \cref{Section: Working over the rationals for m = 10 and 25}, we prove \Cref{Intro Theorem: asymptotic for NQ(X) for m > 9} for $m = 10, 25$. In \cref{Section: Computations for m = 10 and 25}, we produce supplementary computations to estimate the constants appearing in \Cref{Theorem: asymptotic for twN(X) for m = 10 and 25} and \Cref{Theorem: asymptotic for NQ(X) for m = 10 and 25} and empirically confirm that the count of elliptic curves with a cyclic $m$-isogeny aligns with our theoretical estimates when $m = 10, 25$.

\section{Establishing notation for \texorpdfstring{$m \in \set{5, 10, 25}$}{m in 5, 10 25}}\label{Section: Establishing notation when m = 10 and 25}

By \Cref{Corollary: Elliptic curves for which twNeq(X) = twNad(X)}, for $m = 10, 25$, we have
\begin{equation}
\twNeq m (X) = \twNad m (X) \ \text{and} \ \NQeq m (X) = \NQad m (X)
\end{equation}
for all $X > 0$, so we may use either notation interchangeably. On the other hand, $\twistNeq 5 (X) \neq \twistNad 5 (X)$ in general. We work with $\twistNeq m (X)$ until we can proceed no further on the case $m = 5$, and then transition over to using the notation $\twistNad m (X)$.

Let $m \in \set{5, 10, 25}$. Pursuant to the notation established in \cref{Section: Parameterizing elliptic curves with a cyclic m-isogeny}, we define $\h m(t) = \gcd(\f m (t), \g m (t))$, and we define $\tf m (t)$ and $\tg m (t)$ so that
\begin{equation} \label{equation: Definition of tf and tg for m-isogenies for m = 5 and 10 and 25}
\f m (t) = \tf m (t) \h m (t) \ \text{and} \ \g m (t) = \tg m (t) \h m (t)^2.
\end{equation}
Note that $\h m (t)^2$ divides $\g m (t)$, whereas only $\h 7 (t)$ divides $\g 7 (t)$. From \eqref{equation: Definition of tf and tg for m-isogenies for m = 5 and 10 and 25} we conclude
\begin{equation}
\begin{aligned}
	\h {5} (t) &= \h {10} (t) = t^2 + 1, \ \text{and} \\
	\h {25} (t) &= t^2 + 4.
\end{aligned}
\end{equation}
To work with integral models, we take $t=a/b$ (in lowest terms) and homogenize, obtaining
\begin{equation} \label{equation: defining Cab for m = 5 and 10 and 25}
\begin{aligned}
	{\C m}(a, b) &\colonequals b^2 \h m(a/b), \\
	\tA m (a, b) &\colonequals b^{2 \degB m / 3 - 2} \tf m (a/b), \ \text{and} \\
	{\tB m} (a, b) &\colonequals b^{\degB m - 4} \tg m (a/b).
\end{aligned}
\end{equation}
Of course, we have
\begin{equation}
\begin{aligned}
{\C m}(a,b) &= \gcd(\A m(a,b),\B m(a,b)) \in \Z[a,b], \\
\A m (a, b) &= \tA m (a, b) {\C m} (a, b), \ \text{and} \\
\B m (a, b) &= \tB m (a, b) {\C m} (a, b)^2.
\end{aligned}
\end{equation}

Recall $\degB {5} = 6$, $\degB {10} = 12$, and $\degB {25} = 18$.

\section{The twist minimality defect for \texorpdfstring{$m \in \set{5, 10, 25}$}{m in 5, 10 25}}\label{Section: twist minimality defect when m = 5 and 10 and 25}

In this section, for $m \in \set{5, 10, 25}$, we study the twist minimality defect for
\begin{equation}
y^2 = x^3 + \A m (a, b) x + \B m(a, b)
\end{equation}
using the polynomials $\tA m (a, b)$, $\tB m (a, b)$, and $\C m (a, b)$. This section mirrors \cref{Section: twist minimality defect for m = 7}, but our definition of $\cT m(e)$ is changed.

\begin{lemma}\label{Lemma: 2 and 5 are the ungroomed primes for m = 5 and 10 and 25}
	Let $m \in \set{5, 10, 25}$, let $(a,b) \in \Z^2$ be $m$-groomed, let $\ell$ be prime, and let $v \in \Z_{>0}$.  Then the following statements hold.
	\begin{enumalph}
	\item If $\ell \neq 2, 5$, then $\ell^v \mid \twistdefect(\A m (a,b),\B m (a,b))$ if and only if $\ell^{2v} \mid \C m (a, b)$.
	\item We have $\ell^{2v} \mid \C m (a,b)$ if and only if $\ell \nmid b$ and $\h m (a/b) \equiv 0 \psmod{\ell^{2v}}$.
	\item If $\ell \neq 2$, then $\ell \mid \C m (a,b)$ implies $\ell \nmid (\partial \C m/\partial a)(a,b)$.
	\end{enumalph}
\end{lemma}

\begin{proof}
	The proof of part (a) is essentially the same as the proof for \Cref{Lemma: 3 and 7 are the ungroomed primes for m = 7}, with two modifications. First, as $\C m(a, b)^2 \mid \B m (a, b)$, the condition $\ell^{2v} \mid \C m(a, b)$ implies $\ell^v \mid \twistdefect(\A m (a, b), \B m (a, b))$. Second, we have the resultants:
	\begin{equation}	
	\begin{aligned}
&\Res(\tA {5}(t,1),\tB {5}(t,1))=\Res(\tf {5}(t),\tg {5}(t)) \\
&= - 2^2 \cdot 3^4 \cdot 5^{10} = \Res(\tA {5}(1,u), \tB {5}(1,u)), \\
&\Res(\tA {5}(t,1),\C {5}(t,1))=\Res(\tf {5}(t),\h {5}(t)) \\
&= 3^{4} \cdot 5^{5} = \Res(\tA {5}(1,u), \C {5}(1,u)), \\
&\Res(\tA {10}(t,1),\tB {10}(t,1))=\Res(\tf {10}(t),\tg {10}(t)) \\
&= -2^{18} \cdot 3^{16} \cdot 5^{35} = \Res(\tA {10}(1,u), \tB {10}(1,u)), \\
&\Res(\tA {10}(t,1),\C {10}(t,1))=\Res(\tf {10}(t),\h {10}(t)) \\
&= 3^{4} \cdot 5^{5} = \Res(\tA {10}(1,u), \C {10}(1,u)), \\
&\Res(\tA {25}(t,1),\tB {25}(t,1))=\Res(\tf {25}(t),\tg {25}(t)) \\
&= - 2^{38} \cdot 3^{28} \cdot 5^{118} = \Res(\tA {25}(1,u), \tB {25}(1,u)), \\
&\Res(\tA {25}(t,1),\C {25}(t,1))=\Res(\tf {25}(t),\h {25}(t)) \\
&= 2^4 \cdot 3^{4} \cdot 5^{9} = \Res(\tA {25}(1,u), \C {25}(1,u)).
\end{aligned}
\end{equation}
So $2, 3,$ and $5$ are our badly behaved primes. A short computation shows that $3$ divides neither $\twistdefect(\A m (a,b),\B m (a,b))$ nor $\C m (a, b)$, however, which proves (a).

Part (b) proceeds as in the proof of \Cref{Lemma: 3 and 7 are the ungroomed primes for m = 7}.

Part (c) follows from (b) and the fact that $\h 5 (t) = \h {10} (t)$ has discriminant $-2^2,$ and $\h {25} (t)$ has discriminant $\disc(\h {25} (t))= -2^4$.
\end{proof}

\begin{remark}
	We could have adapted the second proof of \Cref{Lemma: 3 and 7 are the ungroomed primes for m = 7} to give an alternate proof of \Cref{Lemma: 2 and 5 are the ungroomed primes for m = 5 and 10 and 25}.%
\end{remark}

\begin{definition}\label{Definition: T(e) for m = 10 and 25}
For $m \in \set{5, 10, 25}$ and for $e \geq 1$, let $\tcalT m(e)$ denote the image of
\begin{equation}
\set{(a, b) \in \bbZ^2 : (a, b) \ m\text{-groomed}, \ e \mid \twistdefect(\A m (a, b), \B m (a, b))}
\end{equation}
under the projection
\begin{equation}
\bbZ^2 \to (\bbZ / e^2 \bbZ)^2.
\end{equation}
For $m \in \set{5, 10, 25}$, let $\tcT m(e) \colonequals \# \tcalT m(e)$.
Similarly, for $m \in \set{5, 10, 25}$, let $\calT m(e) $ denote the image of
\begin{equation}
\set{t \in \bbZ : e^2 \mid \f m (t) \ \text{and} \ e^3 \mid \g m(t)}
\end{equation}
under the projection
\begin{equation}
\bbZ \to \bbZ / e^2 \bbZ,
\end{equation}
and let $\cT m(e) \colonequals \#\calT m(e)$.
\end{definition}

Note that for $m \in \set{5, 10, 25}$, the set $\tcalT m(e)$ is a subset of $(\bbZ / e^2 \bbZ)^2$, whereas $\tcalT 7(e)$ is a subset of $(\bbZ / e^3 \bbZ)^2$! This reflects the difference between potential type II additive reduction and potential type III additive reduction, which manifests in the discrepancy between \Cref{Lemma: 3 and 7 are the ungroomed primes for m = 7}(a) and \Cref{Lemma: 2 and 5 are the ungroomed primes for m = 5 and 10 and 25}(a).

\begin{lemma}\label{Lemma: bound on T(e) for m = 5 and 10 and 25}
Let $m \in \set{5, 10, 25}$. The following statements hold.
\begin{enumalph}
\item If $2^3 \mid e$, then $\tcalT m(e) = \emptyset$. Otherwise, $\tcalT m(e)$ consists of those pairs $(a, b) \in (\bbZ / e^2 \bbZ)^2$ which satisfy the following conditions:
\begin{itemize}
	\item $\A m (a, b) \equiv 0 \psmod {e^2}$, and
	\item $\ell \nmid \gcd(a, b)$ for all primes $\ell \mid e$.
\end{itemize}
\item Let $(a, b) \in \bbZ^2$. If $(a, b) \psmod {e^2} \in \tcalT m(e)$, then $e \mid \tmd(\A m(a, b), \B m(a, b))$.
\item The functions $\tcT m(e)$ and $\cT m(e)$ are multiplicative, and $\tcT m(e) = \varphi(e^2) \cT m(e)$.
\item For all $\ell \neq 2,5$ and $v \geq 1$, 
	\begin{equation} 	\begin{aligned}
	\cT m(\ell^v) &= \cT m(\ell) = 1 + \left(\frac{-1}{\ell}\right).
	\end{aligned}\end{equation} 
\item For $e \in \set{2, 2^2, 5, 5^2, 5^3, 5^4}$, the nonzero values of $\cT m (e)$ are given in Table \ref{table:Tm(2) for m =51025} and Table \ref{table:Tm(5) for m =51025} below. We have
	\begin{equation}
 \begin{aligned}
	\cT {5}(2^v) &= 0 \ \text{for} \ v \geq 1, \\ 
 \cT {10}(2^v) &= 0 \ \text{for} \ v \geq 2, \ \text{and} \\ 
 \cT {25}(2^v) &= 0 \ \text{for} \ v \geq 3;
\end{aligned}
 \end{equation}
	we also have
	\begin{equation}
        \begin{aligned}
	\cT {5}(5^v) &= 1 + 5^{5} \ \text{for} \ v \geq 4, \\ 
 \cT {10}(5^v) &= 1 + 5^5 \ \text{for} \ v \geq 3, \ \text{and} \\ 
 \cT {25}(5^v) &= 1 + 5^9 \ \text{for} \ v \geq 5.
    \end{aligned}
	\end{equation}
\item 
	We have $\cT m(e) =O(2^{\omega(e)})$, where $\omega(e)$ is the number of distinct prime divisors of $e$.
	\end{enumalph}
\end{lemma}

\begin{proof}
	For parts (a) and (b), by the CRT (Sun Zi theorem), it suffices to consider $e = \ell^v$ a power of a prime. For $\ell \neq 2, 5$, both claims follow from \Cref{Lemma: 2 and 5 are the ungroomed primes for m = 5 and 10 and 25}(a)--(b). But a finite computation verifies our claim when $\ell = 2, 5$ as well (see the proof of (e) below).
	
Parts (c) and (d) follow by essentially the same arguments as in the proof of \Cref{Lemma: bound on T(e) for m = 7}.

Next, part (e).  For $\ell = 2$, the claim is a finite computation. 
For $\ell = 5$, we first certify the assertion computationally for $v < 9$. Hensel's lemma still applies to $\h m (t)$: let $t_0,t_1$ be the roots of $\h m (t)$ in the $5$-adic integers $\Z_5$, with $t_0 \equiv 2 \psmod 5$ and with $t_1 \equiv -2 \psmod 5$ if $m = 5, 10$, and $t_0 \equiv -1 \psmod 5$ and $t_1 \equiv 1 \psmod 5$ if $m = 25$. It is easy to verify
\begin{equation}
\tf m (t_0) \not\equiv 0 \psmod 5, \ \tg m (t_0) \not\equiv 0 \psmod 5,
\end{equation}
and on the other hand that for $m = 5, 10$, we have
\begin{equation}
\tf {m} (t_1) \equiv \tg {m} (t_1) \equiv 0 \psmod {5^5} \ \text{but} \ \tf {m} (t_1) \not\equiv 0 \psmod{5^6},
\end{equation}
and that for $m = 25$, we have
\begin{equation}
\tf {25} (t_1) \equiv \tg {25} (t_1) \equiv 0 \psmod {5^9} \ \text{but} \ \tf {25} (t_1) \not\equiv 0 \psmod{5^{10}}.
\end{equation}
For $m = 5, 10$, we therefore have
\begin{equation}
\calT {m} (5^v) = \set{t_0} \sqcup \set{t_1 + 5^{2v-5} u : u \in \bbZ / 5^5 \bbZ} \ \text{for} \ v \geq 5, \label{Equation: decomposition of T(5^v) for m = 5 and 10}
\end{equation}
and for $m = 25$, we have
\begin{equation}
\calT {25} (5^v) = \set{t_0} \sqcup \set{t_1 + 5^{2v-9} u : u \in \bbZ / 5^9 \bbZ} \ \text{for} \ v \geq 9. \label{Equation: decomposition of T(5^v) for m = 25}
\end{equation}
Part (e) is now clear.

Finally, part (f).  From (d)--(e), for $m = 5, 10$ the local factor at $5$ is at most $3126$, the local factor at $2$ is at most $2$, and every other local factor is at most $2$.  Thus
\begin{equation}
\cT {m}(e) \leq 2 \cdot 3126 \cdot 2^{\omega(e)}.
\end{equation}
Similarly, for $m = 25$, the local factors at $2$ and $5$ are at most $8$ and $1953126$, respectively, so
\begin{equation}
\cT {25}(e) \leq 8 \cdot 1953126 \cdot 2^{\omega(e)}.
\end{equation}
Hence $\cT m(e) = O(2^{\omega(e)})$ as claimed.
\end{proof}

\jvtable{table:Tm(2) for m =51025}{
\rowcolors{2}{white}{gray!10}
\begin{tabular}{c | c c}
$m$ & $\cT m (2^1)$ & $\cT m (2^2)$ \\
\hline\hline
$5$ & -- & -- \\
$10$ & $2$ & -- \\
$25$ & $2$ & $2^3$  \\
\end{tabular}
}{All nonzero $\cT m (2^v)$ for $m \in \set{5, 10, 25}$}

\jvtable{table:Tm(5) for m =51025}{
\rowcolors{2}{white}{gray!10}
\begin{tabular}{c | c c c c}
$m$ & $\cT m (5^1)$ & $\cT m (5^2)$ & $\cT m (5^3)$ & $\cT m (5^4)$ \\
\hline\hline
$5$ & $1 + 5$ & $1 + 5^2$ & $1 + 5^3$ & $1 + 5^5$ \\
$10$ & $1 + 5$ & $1 + 5^3$ & $1 + 5^5$ & $1 + 5^5$ \\
$25$ & $1 + 5$ & $1 + 5^3$ & $1 + 5^5$ & $1 + 5^7$ \\
\end{tabular}
}{All $\cT m (5^v)$ for $m \in \set{5, 10, 25}$ and $v \leq 4$}

\subsection*{The common factor \texorpdfstring{$\C m(a, b)$}{Cab}}

In view of \Cref{Lemma: 2 and 5 are the ungroomed primes for m = 5 and 10 and 25}, when $m \in \set{5, 10, 25}$, the twist minimality defect away from the primes $2,5$ is determined by the quadratic form 
\begin{equation}
    \C 5 (a, b) = \C {10} (a,b)=a^2+b^2= b^2 \h {5} (a/b) = b^2 \h {10} (a/b)
\end{equation}
or the quadratic form 
\begin{equation}
\C {25} (a, b) = a^2 + 4 b^2 = b^2 \h {25} (a/b)    
\end{equation}
respectively. For $m \in \set{5, 10, 25}$, we define
\begin{equation}\label{Equation: Definition of calcc m for m = 5 and 10 and 25}
    \calcc m (e) \colonequals \set{(a, b) \in \bbZ^2 : \C m (a, b) = e \ \text{and} \ \gcd(a, b) = 1},
\end{equation}
and note $\# \calcc m (e) \leq 2^{\omega(e) + 1}$.

Just like $\C 7$, the polynomials $\C 5 = \C {10}$ and $\C {25}$ are both norm forms of quadratic orders with class number $1$, namely $\bbZ[\zeta_4]$ or $\bbZ[2 \zeta_4]$, where $\zeta_4$ is a primitive $4$th root of unity, i.e., a square root of $-1$. We record some elementary algebraic observations about $\C {5} (a, b) = \C {10} (a, b)$ and $\C {25} (a, b)$ and the associated orders $\bbZ[\zeta_4]$ and $\bbZ[2 \zeta_4]$.

\begin{lemma}\label{Lemma: Algebraic structure of Z[(2) zeta4] and C}
	Let $m \in \set{5, 10, 25}$. The following statements hold.
	\begin{enumalph}	
		\item The right regular representation of $\bbZ[\zeta_4]$ in the basis $\set{1, \zeta_4}$ induces the map $\repr 5 = \repr {10} : \bbZ^2 \to \textup{M}_2(\bbZ)$ given by
		\begin{equation}\label{Equation: Defining repr for m = 5 and 10}
		\repr 5 = \repr {10} : (a, b) \mapsto \begin{pmatrix}
		a & b \\ - b & a
		\end{pmatrix},
		\end{equation}
  and the right regular representation of $\bbZ[2 \zeta_4]$ in the basis $\set{1, 2 \zeta_4}$ induces the map $\repr {25} : \bbZ^2 \to \textup{M}_2(\bbZ)$ given by
  \begin{equation}\label{Equation: Defining repr for m = 25}
		\repr {25} : (a, b) \mapsto \begin{pmatrix}
		a & b \\ - 4 b & a
		\end{pmatrix}.
		\end{equation}
		\item For all $a, b, c, d, e \in \bbZ$, we have the following implication:
		\begin{equation}
		\C m (a, b) = e \implies e \mid \C m ((c, d) \cdot \repr m (a, b)).
		\end{equation}
		\item Conversely, if $c\prm, d\prm, e, k$ are integers such that $k \geq 1$, $e^k \mid \C m (c\prm, d\prm)$, and 
    \begin{equation}
    \gcd(c\prm, d\prm, e) = \gcd(2, e) = 1,
    \end{equation}
    then there are integers $a, b, c, d \in \bbZ$ with $(a, b) \in \calcc m (e)$ and
		\begin{equation}
		(c\prm, d\prm) = (c, d) \cdot \repr m (a, b)^k.
		\end{equation}
	\end{enumalph}
\end{lemma}

\begin{proof}
    The proof is essentially identical to that of \Cref{Lemma: Algebraic structure of Z[3 zeta6] and C}.
\end{proof}

The twist minimality defect measures the disparity between $\rawheight(A, B)$, which is easy to compute, and $\twistheight(A, B)$, which is of arithmetic interest: this disparity cannot be too large compared to $\C m(a,b)$, as the following theorem shows.

\begin{theorem}\label{Theorem: Controlling size of twist minimality defect for m = 5 and 10 and 25}
	Let $m \in \set{5, 10, 25}$. The following statements hold.
	\begin{enumalph}
	\item For all $(a, b) \in \bbR^2$, we have
	\begin{equation}\label{Equation: upper and lower bounds for H for m = 5 and 10 and 25}
	108 \C m (a, b)^{\degB m} \leq \rawheight(\A m (a, b), \B m (a, b)) \leq \upperratio m \C m (a, b)^{\degB m},
	\end{equation}
	where the constants 
	\begin{equation} 	\begin{aligned}
	\upperratio {5} &= 679\,212\,199.08278056\ldots, \\
	\upperratio {10} &= 211\,362\,386.0164477\ldots, \ \text{and} \\
	\upperratio {25} &= 26\,367\,187\,500, \\
	\end{aligned}\end{equation} 
	are algebraic numbers given by evaluating the function $H(\A {5}(a, b), \B {5}(a, b))$ at appropriate roots of \eqref{Equation: Roots defining upperratio for m = 5}, evaluating the function $H(\A {10}(a, b), \B {10}(a, b))$ at appropriate roots of \eqref{Equation: Roots defining upperratio for m = 10}, and evaluating the function $H(\A {25}(a, b), \B {25}(a, b))$ at the appropriate roots of \eqref{Equation: Roots defining upperratio for m = 25}, respectively.
	\item If $\C m (a, b) = e_0^2 n_0$, with $n_0$ squarefree, then $\twistdefect(\A m (a, b), \B m (a, b)) = e_0 e\prm$, where $e\prm \mid 5^3$ if $m = 5$, $e\prm \mid 2 \cdot 5^3$ if $m = 10$, and $e\prm \mid 2^2 \cdot 5^5$ if $m = 25$. We have
	\begin{equation}\label{Equaton: upper and lower bounds on twist height for m = 5}
	\frac{2^2 \cdot 3^3}{5^{18}} e_0^{6} n_0^{6} \leq \twistheight(\A {5} (a, b), \B {5} (a, b)) \leq \upperratio {5} e_0^{6} n_0^{6},
	\end{equation}
	\begin{equation}\label{Equaton: upper and lower bounds on twist height for m = 10}
	\frac{3^3}{2^4 \cdot 5^{18}} e_0^{18} n_0^{12} \leq \twistheight(\A {10} (a, b), \B {10} (a, b)) \leq \upperratio {10} e_0^{18} n_0^{12},
	\end{equation}
	and
	\begin{equation}\label{Equaton: upper and lower bounds on twist height for m = 25}
	\frac{3^3}{2^{10} \cdot 5^{30}} e_0^{30} n_0^{18} \leq \twistheight(\A {25} (a, b), \B {25} (a, b)) \leq \upperratio {25} e_0^{30} n_0^{18}.
	\end{equation}
	\end{enumalph}
\end{theorem}

\begin{proof}
    The proof follows the contours of \Cref{Theorem: Controlling size of twist minimality defect for m = 7}.

    Part (a) is proven exactly as in the proof of \Cref{Theorem: Controlling size of twist minimality defect for m = 7}. The lower bound $108$ of \eqref{Equation: upper and lower bounds for H for m = 5 and 10 and 25} is attained at $(1, 0)$, and the upper bound $\upperratio m$ is attained when $a$ and $b$ are appropriately chosen roots of
	\begin{equation}\label{Equation: Roots defining upperratio for m = 5}
	\begin{aligned}
		 &312500 a^4 - 312500 a^2 + 841 \\
		 =& 2^{2} \cdot 5^{7} \cdot a^{4} - 2^{2} \cdot 5^{7} \cdot a^{2} + 29^{2}, \ \text{and} \\
		 &312500 b^4 - 312500 b^2 + 841 \\
		 =& 2^{2} \cdot 5^{7} \cdot b^{4} - 2^{2} \cdot 5^{7} \cdot b^{2} + 29^{2},
	\end{aligned}
	\end{equation}
	if $m = 5$, and of
	\begin{equation}\label{Equation: Roots defining upperratio for m = 10}
	\begin{aligned}
		 &253125 a^{12} - 478125 a^{10} + 480625 a^8 - 455000 a^6 \\
		 &+ 252200 a^4 - 49360 a^2 + 16 \\
		 =& 3^{4} \cdot 5^{5} \cdot a^{12} - 3^{2} \cdot 5^{5} \cdot 17 \cdot a^{10} + 5^{4} \cdot 769 \cdot a^{8} \\
		 &- 2^{3} \cdot 5^{4} \cdot 7 \cdot 13 \cdot a^{6} + 2^{3} \cdot 5^{2} \cdot 13 \cdot 97 \cdot a^{4} - 2^{4} \cdot 5 \cdot 617 \cdot a^{2} + 2^{4}, \ \text{and} \\
		 &253125 b^{12} - 1040625 b^{10} + 1886875 b^8 - 1748750 b^6 \\
		 &+ 786575 b^4 - 140665 b^2 + 3481 \\
		 =& 3^{4} \cdot 5^{5} \cdot b^{12} - 3^{2} \cdot 5^{5} \cdot 37 \cdot b^{10} + 5^{4} \cdot 3019 \cdot b^{8} \\
		 &- 2 \cdot 5^{4} \cdot 1399 \cdot b^{6} + 5^{2} \cdot 73 \cdot 431 \cdot b^{4} - 5 \cdot 7 \cdot 4019 \cdot b^{2} + 59^{2}
	\end{aligned}
	\end{equation}
	if $m = 10$, and of
	\begin{equation}\label{Equation: Roots defining upperratio for m = 25}
	\begin{aligned}
		 &53833007812500 a^{28} - 891577148437500 a^{26} + 7403853759765625 a^{24} \\
		 &- 38650180664062500 a^{22} + 139358151855468750 a^{20} - 361638379062500000 a^{18} \\
		 &+ 690434893630859375 a^{16} - 979823552140625000 a^{14} + 1042891876273125000 a^{12} \\
		 &- 839328158831937500 a^{10} + 509588953407434375 a^8 - 227793758883072500 a^6 \\
		 &+ 70659569038784950 a^4 - 13323520820064520 a^2 + 1058114957485041 \\
		 =& 2^{2} \cdot 3^{2} \cdot 5^{15} \cdot 7^{2} \cdot a^{28} - 2^{2} \cdot 3 \cdot 5^{14} \cdot 7 \cdot 37 \cdot 47 \cdot a^{26} + 5^{13} \cdot 19 \cdot 319223 \cdot a^{24} \\
		 &- 2^{2} \cdot 3 \cdot 5^{13} \cdot 13 \cdot 17 \cdot 11939 \cdot a^{22} + 2 \cdot 3 \cdot 5^{13} \cdot 739 \cdot 25747 \cdot a^{20} \\
	   &- 2^{5} \cdot 5^{10} \cdot 197 \cdot 677 \cdot 8677 \cdot a^{18} + 5^{9} \cdot 149 \cdot 2372501111 \cdot a^{16} \\
		 &- 2^{3} \cdot 5^{9} \cdot 31 \cdot 2022861527 \cdot a^{14} + 2^{3} \cdot 3 \cdot 5^{7} \cdot 47 \cdot 61 \cdot 10321 \cdot 18797 \cdot a^{12} \\
        &- 2^{2} \cdot 5^{6} \cdot 11 \cdot 19 \cdot 11251 \cdot 5711029 \cdot a^{10} + 5^{5} \cdot 163068465090379 \cdot a^{8} \\
	&- 2^{2} \cdot 5^{4} \cdot 91117503553229 \cdot a^{6}  + 2 \cdot 5^{2} \cdot 7 \cdot 43 \cdot 4694987975999 \cdot a^{4} \\
        &- 2^{3} \cdot 5 \cdot 13 \cdot 8221 \cdot 3116671381 \cdot a^{2} + 3^{2} \cdot 10842893^{2}, \ \text{and} \\
		 &861328125000000 b^{28} + 551660156250000 b^{26} + 712154541015625 b^{24} \\
		 &- 62065429687500 b^{22} + 109768066406250 b^{20} - 83235166015625 b^{18} \\
		 &+ 33729880859375 b^{16} - 11477761718750 b^{14} + 4976470781250 b^{12} \\
		 &- 1567605921875 b^{10} + 281949340625 b^8 - 27367628750 b^6 \\
		 &+ 1121883700 b^4 - 21082345 b^2 + 7056 \\
		 =& 2^{6} \cdot 3^{2} \cdot 5^{15} \cdot 7^{2} \cdot b^{28} + 2^{4} \cdot 3 \cdot 5^{14} \cdot 7 \cdot 269 \cdot b^{26} + 5^{13} \cdot 583397 \cdot b^{24} \\
		 &- 2^{2} \cdot 3 \cdot 5^{13} \cdot 19 \cdot 223 \cdot b^{22} + 2 \cdot 3 \cdot 5^{13} \cdot 7 \cdot 2141 \cdot b^{20} - 5^{10} \cdot 13 \cdot 655637 \cdot b^{18} \\
		 &+ 5^{9} \cdot 1231 \cdot 14029 \cdot b^{16} - 2 \cdot 5^{9} \cdot 2938307 \cdot b^{14} + 2 \cdot 3 \cdot 5^{7} \cdot 2777 \cdot 3823 \cdot b^{12} \\
		 &- 5^{6} \cdot 7 \cdot 2011 \cdot 7127 \cdot b^{10} + 5^{5} \cdot 2377 \cdot 37957 \cdot b^{8} - 2 \cdot 5^{4} \cdot 7 \cdot 11^{2} \cdot 25849 \cdot b^{6} \\
		 &+ 2^{2} \cdot 5^{2} \cdot 7 \cdot 1602691 \cdot b^{4} - 5 \cdot 4216469 \cdot b^{2} + 2^{4} \cdot 3^{2} \cdot 7^{2}
	\end{aligned}
	\end{equation}
	if $m = 25$. For $m \in \set{5, 10, 25}$, consider the arguments that maximize the ratio
 \begin{equation}
 \rawheight(\A {m}(a, b), \B {m}(a, b))/\C m(a, b)^{\degB m}.
 \end{equation}
 At these arguments, we have $27 \abs{\B {5}(a, b)}^2 > 4 \abs{\A {5}(a, b)}^3$ when $m = 5$, $4 \abs{\A {10}(a, b)}^3 > 27 \abs{\B {10}(a, b)}^2$ when $m = 10$, and $4 \abs{\A {25}(a, b)}^3 = 27 \abs{\B {25}(a, b)}^2$ when $m = 25$. For the reader's information, 
	\begin{equation}
	(a, b) = (0.051946913\ldots, 0.998649847\ldots)
	\end{equation}
	maximizes this ratio when $m = 5$, 
	\begin{equation}
	(a, b) = (-0.766646866\ldots, 0.642068986\ldots)
	\end{equation}
	maximizes this ratio when $m = 10$, and \begin{equation}
	(a, b) = (0.447213595\ldots, 0.447213595\ldots)
	\end{equation}
	maximizes this ratio when $m = 25$. 
		
	We now prove (b). Write $\C m(a, b) = e_0^2 n_0$ with $n_0$ squarefree, and write 
 \begin{equation}
 \twistdefect(\A m(a, b), \B m(a, b)) = e_0 e\prm.
 \end{equation}
 By \ref{Lemma: 2 and 5 are the ungroomed primes for m = 5 and 10 and 25}, $e\prm = 2^v \cdot 5^w$ for some $v, w \geq 0$. A short computation shows that if $m = 5$ then $v = 0$, if $m = 10$ then $v \in \set{0, 1}$, and that if $m = 25$ then $v \in \set{0, 1, 2}$. On the other hand, \eqref{Equation: decomposition of T(5^v) for m = 5 and 10} shows $w \leq \ceil{5/2} = 3$ if $m = 5, 10$, and \eqref{Equation: decomposition of T(5^v) for m = 25} shows $w \leq \ceil{9/2} = 5$ if $m = 25$. 
	
	As 
	\begin{equation}
	\rawheight(\A m (a, b), \B m(a, b)) = e_0^6 \parent{e\prm}^6 \twistheight(\A m(a, b), \B m(a, b)),
	\end{equation}
	we see
	\begin{equation}
	\frac{108}{(e\prm)^6} e_0^{2 \degB m - 6} n_0^{\degB m} \leq \twistheight(\A m(a, b), \B m(a, b)) \leq \frac{\upperratio m}{(e\prm)^6} e_0^{2 \degB m - 6} n_0^{\degB m}.
	\end{equation}
	Rounding $e\prm$ up to $5^3$ (if $m = 5$), $2 \cdot 5^3$ (if $m = 10$), or $2^2 \cdot 5^5$ (if $m = 25$), on the left, and rounding down to $1$ on the right gives the desired result.
\end{proof}

Note that \eqref{Equaton: upper and lower bounds on twist height for m = 5}, in contrast to \eqref{Equaton: upper and lower bounds on twist height for m = 10} and \eqref{Equaton: upper and lower bounds on twist height for m = 25}, has matching exponents for $e_0$ and $n_0$. From a certain perspective, this is what makes the case $m = 5$ difficult to handle.

\begin{corollary}\label{Corollary: bound twist defect in terms of twist height}
	For $m \in \set{5, 10, 25}$, let $(a, b)$ be an $m$-groomed pair. We have
\begin{equation}
\twistdefect(\A {5}(a, b), \B {5} (a, b)) \leq \frac{5^6}{2^{1/3} \cdot 3^{1/2}} \twistheight(\A {5} (a, b), \B {5} (a, b))^{1/6}, 
\end{equation}
where ${5^6}/{(2^{1/3} \cdot 3^{1/2})} = 7\,160.050\ldots$, and
\begin{equation}
\twistdefect(\A {10}(a, b), \B {10} (a, b)) \leq \frac{2^{11/9} \cdot 5^4}{3^{1/6}} \twistheight(\A {10} (a, b), \B {10} (a, b))^{1/18}, 
\end{equation}
where ${2^{11/9} \cdot 5^4}/{3^{1/6}} = 1\,214.186\ldots$, and
\begin{equation}
\twistdefect(\A {25}(a, b), \B {25} (a, b)) \leq \frac{2^{7/3} \cdot 5^6}{3^{1/10}} \twistheight(\A {25} (a, b), \B {25} (a, b))^{1/30},
\end{equation}
where ${2^{7/3} \cdot 5^6}/{3^{1/10}} = 70\,552.307\ldots$.
\end{corollary}

\begin{proof}
	We prove the case $m = 10$; the cases $m = 5$ and $m = 25$ are entirely similar. In the notation of \Cref{Theorem: Controlling size of twist minimality defect for m = 5 and 10 and 25}(b),
	\begin{equation}
	\frac{3^3}{2^4 \cdot 5^{18}} e_0^{18} n_0^{12} \leq \twistheight(\A {10} (a, b), \B {10} (a, b)) \leq \upperratio {10} e_0^{18} n_0^{12}
	\end{equation}
	Multiplying through by $(e\prm)^{18}$, rounding $n_0$ down to $1$ on the left, rounding $e\prm$ up to $2 \cdot 5^3$ on the right, and taking $18$th roots of both sides, we obtain the desired result.
\end{proof}

\section{Estimates for twist classes for \texorpdfstring{$m = 10, 25$}{m = 10, 25}} \label{Section: Estimating twN(X) for m = 5 and 10 and 25}

In this section, we use \cref{Section: Our approach revisited} to estimate $\twistNad m(X)$ for $m = 10, 25$, counting the number of twist minimal elliptic curves over $\bbQ$ admitting an $m$-isogeny for $m = 10, 25$. We also indicate why our method fails for $m = 5$.

Recall \eqref{Equation: defining cM}, \eqref{Equation: Applying N calE to Neq m}, and \eqref{Equation: Applying N calE to Nad m}. By \cref{Section: Parameterizing elliptic curves with a cyclic m-isogeny}, $\cM m(X; e)$ counts pairs $(a, b) \in \bbZ^2$ with
\begin{itemize}
	\item $(a, b)$ $m$-groomed,
	\item $\rawheight(\A m(a,b),\B m(a,b)) \leq X$, and 
	\item $e \mid \twistdefect(\A m(a, b), \B m(a, b))$.
\end{itemize}

The following proposition refines \Cref{Lemma: fundamental sieve}, and specifies both an order of growth and an explicit upper bound past which the summands of \eqref{Equation: twN calE (X) as a sum} vanish when $m \in \set{5, 10, 25}$.

\begin{proposition}\label{Proposition: fundamental sieve for for m = 5 and 10 and 25}
	For $m \in \set{5, 10, 25}$, we have
	\begin{equation}\label{Equation: twN(X) in terms of M(X; e) for m = 5 and 10 and 25}
	\twNeq m (X) = \sum_{n \ll X^{1/(2 \degB m - 6)}} \sum_{e \mid n} \mu(n/e) \cM m (e^6 X; n);
	\end{equation}
	more precisely, if $m = 5$ we can restrict our sum to
	\begin{equation}
	n \leq \frac{5^6}{2^{1/3} \cdot 3^{1/2}} \cdot X^{1/6},
	\end{equation}
	if $m = 10$ we can restrict our sum to
	\begin{equation}
	n \leq \frac{2^{11/9} \cdot 5^{4}}{3^{1/6}} \cdot X^{1/18},
	\end{equation}
	and if $m = 25$, we can restrict our sum to
	\begin{equation}
	n \leq \frac{2^{7/3} \cdot 5^6}{3^{1/10}} X^{1/30}.
	\end{equation}
\end{proposition}

\begin{proof}
	We prove the case $m = 10$ by way of illustration, although the argument precisely mirrors \Cref{Proposition: fundamental sieve for for m = 7}. Let $(a, b) \in \bbZ^2$, and suppose $\rawheight(\A {10}(a,b),\B {10}(a,b)) \leq e^6 X$ and $e \mid \twistdefect(\A {10}(a, b), \B {10}(a, b))$.
	If we can prove 
	\begin{equation}\label{Equation: bound on e for m = 10}
	e \leq \frac{2^{11/9} \cdot 5^{4}}{3^{1/6}} \cdot X^{1/18},
	\end{equation}
	then our claim will follow.
	
	Write $\C {10}(a, b) = e_0^2 n_0$, with $n_0$ square-free. By \Cref{Theorem: Controlling size of twist minimality defect for m = 5 and 10 and 25}(a), we have
	\begin{equation}
	108 e_0^{24} n_0^{12} \leq e^6 X.\label{Equation: First inequality between e0 and X for m = 10}
	\end{equation}
	On the other hand, by \Cref{Theorem: Controlling size of twist minimality defect for m = 5 and 10 and 25}(b), we have $e \mid 2 \cdot 5^3 \cdot e_0$, and \textit{a fortiori} 
	\begin{equation}
	e \leq 2 \cdot 5^3 e_0.\label{Equation: e0 versus e for m = 10}
	\end{equation}
	Multiplying \eqref{Equation: First inequality between e0 and X for m = 10} through by $(2 \cdot 5^3)^{24}$ and utilizing \eqref{Equation: e0 versus e for m = 10}, we conclude
	\begin{equation}
	2^2 \cdot 3^3 e^{24} n_0^{12} \leq 2^{24} \cdot 5^{72} e^6 X.
	\end{equation}
	Rounding $n_0$ down to $1$ and rearranging, we obtain \eqref{Equation: bound on e for m = 10}.
\end{proof}

Let $m \in \set{5, 10, 25}$. As in \cref{Section: Estimating twN(X) for m = 7}, in order to estimate $\cM m(X; e)$, we further unpack the $m$-groomed condition on pairs $(a, b)$. For the reader's convenience, we recall from Table \ref{table:auxiliaries} that $\cusps 5 = \set{11/2, \infty}$, $\cusps {10} = \set{-2, 0, 1/2, \infty}$, and $\cusps {25} = \set{1, \infty}$. Exactly as in \cref{Section: Estimating twN(X) for m = 7}, for $m \in \set{5, 10, 25}$ we let $\cM m(X; d, e)$ denote the number of pairs $(a, b) \in \bbZ^2$ with 
	\begin{itemize}
		\item $\gcd(da, db, e) = 1$, $b > 0$, and $a/b \not\in \cusps m$;
		\item $\rawheight(\A m(d a, d b), \B m(d a, db)) \leq X$;
		\item $e \mid \twistdefect(\A m(d a, d b), \B m(da, db))$;
	\end{itemize}
By \Cref{Theorem: Controlling size of twist minimality defect for m = 5 and 10 and 25}, and because $\rawheight(\A m(a, b), \B m(a, b))$ is homogeneous of degree $2 \degB m$, another M\"{o}bius sieve yields
\begin{equation}\label{equation: cM(X;e) in terms of cM(X; d, e) for m = 5 and 10 and 25}
	\cM m (X; e) = \sum_{\substack{d \ll X^{1/(2 \degB m)} \\ \gcd(d, e) = 1}} \mu(d) \cM m(X; d, e).
\end{equation}

The following lemma gives asymptotics for $\cM m (X)$, which depend on the observation that the largest square dividing $\C m (a, b)$ is essentially the square of the twist minimality defect. Once we have \Cref{Lemma: asymptotic for M(X; e) for m = 5 and 10 and 25}, the rest of our argument proceeds along the lines given in the paragraph after \eqref{equation: cM(X;e) in terms of cM(X; d, e) for m = 7}, and the type III additive reduction for the elliptic surface
\begin{equation}
y^2 = x^3 + \f m (t) x + \g m (t)
\end{equation}
has no additional relevance.

\begin{lemma}\label{Lemma: asymptotic for M(X; e) for m = 5 and 10 and 25}
Let $m \in \set{5, 10, 25}$. The following statements hold.
\begin{enumalph}
\item	If $\gcd(d, e) > 1$, then $\cM {m}(X; d, e) = 0$. If $\gcd(d, e) = 1$, we have 
	\begin{equation}
	\cM m(X; d, e) = \frac{\R m \cT m(e) X^{1/\degB m}}{d^2 e^2} \prod_{\ell \mid e} \parent{1 - \frac{1}{\ell}} + O\parent{\frac{2^{\omega(e)} X^{1/(2 \degB m)}}{d e}}
	\end{equation}
        for $X, d, e \geq 1$,
	where $\R m$ is the area of \eqref{eqn: R(X)}.
\item We have
	\begin{equation}
	\cM m(X; e) = \frac{\R m \cT m(e) X^{1/\degB m}}{\zeta(2) e^2 \prod_{\ell \mid e} \parent{1 + \frac{1}{\ell}}} +  O\parent{\frac{2^{\omega(e)} X^{1/(2 \degB m)} \log X}{e}}
	\end{equation}
    for $X \geq 2$ and $e \geq 1$.
\end{enumalph}
In both cases, the implied constants are independent of $d$, $e$, and $X$.
\end{lemma}

As with \Cref{Lemma: asymptotic for M(X; e) for m = 7}, we prove \Cref{Lemma: asymptotic for M(X; e) for m = 5 and 10 and 25} by means of two partial proofs. The first proof gives an intuitive interpretation of the coefficient of $X^{1/\degB m}$, and generalizes readily to other elliptic surfaces with type III additive reduction. It differs from the first proof of \Cref{Lemma: asymptotic for M(X; e) for m = 7} only in that we sum over congruence classes modulo $e^2$ rather than $e^3$. The second proof leverages the observation that $\C m (a, b)$ is the norm of the order $\bbZ[\zeta_4]$ or $\bbZ[2 \zeta_4]$ to give an enhanced error term.

We write out the first proof to give a flavor for the differences between the case $m = 7$ and the cases $m \in \set{5, 10, 25}$, but largely omit the write-up of the second proof to spare the readers' time.

\begin{proof}[First proof of \Cref{Lemma: asymptotic for M(X; e) for m = 5 and 10 and 25}]
	Our proof mirrors that of \Cref{Lemma: asymptotic for M(X; e) for m = 7}. Let $m \in \set{5, 10, 25}$. We begin with (a) and examine the summands $\cM m(X; d, e)$. If $d$ and $e$ are not coprime, then $\cM m(X; d, e) = 0$ because $\gcd(da, db, e) \geq \gcd(d, e) > 1$. On the other hand, if $\gcd(d, e) = 1$, we have a bijection from the pairs counted by $\cM m(X; 1, e)$ to the pairs counted by $\cM m(d^{2 \degB m} X; d, e)$ given by $(a, b) \mapsto (d a, d b)$.

For $m \in \set{5, 10, 25}$, $X\geq 1$, and $e, a_0, b_0 \in \bbZ$, we write
\begin{equation}
    L_m(X; e, a_0, b_0) \colonequals \#\{(a, b) \in \calR m(X) \cap \bbZ^2 : (a, b) \equiv (a_0, b_0) \psmod {e^2}, (a:b) \not\in \cusps m \}
\end{equation}
(this notation will not be used outside of this proof). By \Cref{Corollary: Estimates for lattice counts}, we have
\begin{equation}
L_{m}(X; e, a_0, b_0) = \frac{\R {m} X^{1/\degB m}}{e^4} + O \parent{\frac{X^{1/(2 \degB m)}}{e^2}}.
\end{equation}
By \Cref{Lemma: bound on T(e) for m = 5 and 10 and 25}(b), we have
\begin{equation}
\begin{aligned}
	\cM m(X; 1, e) &= \sum_{(a_0, b_0) \in \tcalT m(e)} L_m(X; e, a_0, b_0) \\
	&= \varphi(e^2) \cT m(e) \parent{\frac{\R m X^{1/\degB m}}{e^4} + O\parent{\frac{X^{1/(2 \degB m)}}{e^2}}} \\
	&= \frac{\R m \cT m(e) X^{1/\degB m}}{e^2} \prod_{\ell \mid e} \parent{1 - \frac{1}{\ell}} + O(\cT m(e) X^{1/(2 \degB m)}),
\end{aligned}
\end{equation}
and thus
\begin{equation}
	\cM m(X; d, e) = \frac{\R m \cT m(e) X^{1/\degB m}}{d^2 e^2} \prod_{\ell \mid e} \parent{1 - \frac{1}{\ell}} + O\parent{\frac{\cT m(e) X^{1/(2 \degB m)}}{d}}.
\end{equation}

For part (b), we compute
\begin{equation} \label{Equation: cmXe for m = 5 and 10 and 25}
\begin{aligned}
	\cM m(X; e) =& \sum_{\substack{d \ll X^{1/(2 \degB m)} \\ \gcd(d, e) = 1}} \mu(d) \cM m(X; d, e) \\
	=& \sum_{\substack{d \ll X^{1/(2 \degB m)} \\ \gcd(d, e) = 1}} \mu(d) \parent{\frac{\cT m(e) \R m X^{1/\degB m}}{d^2 e^2} \prod_{\ell \mid e} \parent{1 - \frac{1}{\ell}} + O\parent{\cT m(e) \frac{X^{1/(2 \degB m)}}{d}}} \\
	=& \frac{\R m \cT m(e) X^{1/ \degB m}}{e^2} \prod_{\ell \mid e} \parent{1 - \frac{1}{\ell}} \sum_{\substack{d \ll X^{1/(2 \degB m)} \\ \gcd(d, e) = 1}} \frac{ \mu(d)}{d^2} \\
    &+ O\parent{\cT m(e) X^{1/(2 \degB m)} \sum_{\substack{d \ll X^{1/(2 \degB m)} \\ \gcd(d, e) = 1}} \frac 1d}.
	\end{aligned}
	\end{equation}
We plug the straightforward estimates
\begin{equation}
\sum_{\substack{d \ll X^{1/(2 \degB m)} \\ \gcd(d, e) = 1}} \frac{ \mu(d)}{d^2} = \frac{1}{\zeta(2)} \prod_{\ell \mid e} \parent{1 - \frac{1}{\ell^2}}\inv + O(X^{-1/(2 \degB m)})
\end{equation}
and
\begin{equation} 
\sum_{\substack{d \leq X^{1/(2 \degB m)}}} \frac 1d = \frac{1}{2 \degB m}\log X+ O(1) \end{equation}
into \eqref{Equation: cmXe for m = 5 and 10 and 25}, along with \Cref{Lemma: bound on T(e) for m = 5 and 10 and 25}(f). Simplifying, we now obtain
\begin{equation}
\begin{aligned}
\cM m(X;e)	
	&= \frac{\R m \cT m(e) X^{1/\degB m}}{\zeta(2) e^2 \prod_{\ell \mid e} \parent{1 + \frac{1}{\ell}}} + O(2^{\omega(e)} X^{1/(2 \degB m)} \log X),
\end{aligned}
\end{equation}
which proves (b).
\end{proof}

\begin{proof}[Second proof of \Cref{Lemma: asymptotic for M(X; e) for m = 5 and 10 and 25}]
    This proof is, \textit{mutatis mutandis}, the same as the proof we gave for \Cref{Lemma: asymptotic for M(X; e) for m = 7}. However, when we apply \Cref{Lemma: Algebraic structure of Z[(2) zeta4] and C}(c), we do so with $k = 2$ rather than $k = 3$.
\end{proof}

Let $\kappa \in \bbR$. We write
\begin{equation}\label{Equation: generalized totient function}
\varphi_\kappa(n) \colonequals \sum_{d \mid n} \mu(n/d) d^\kappa = n^\kappa \prod_{\ell \mid n} \parent{1 - \frac {1}{\ell^\kappa}}
\end{equation}
for the generalized Jordan totient function. %

For $m = 10, 25$, we let
	\begin{equation}
	\Q m \colonequals \sum_{n \geq 1} \frac{\varphi_{6 / \degB m}(n) \cT m(n)}{n^2 \prod_{\ell \mid n} \parent{1 + \frac{1}{\ell}}}.
	\end{equation}
		Note that the sum defining $\Q m$ diverges when $m = 5$! We let
\begin{equation}\label{Equation: twconst for m = 10 and 25}
	\twconst m \colonequals \frac{\Q m \R m}{\zeta(2)}.
\end{equation}
Here, as always, $\R m$ is the area of the region
	\begin{equation}
	\calR m(1) = \set{(a, b) \in \bbR^2 : \rawheight(\A m (a, b), \B m (a, b)) \leq 1, b \geq 0}.
	\end{equation}
		
We are now in a position to estimate $\twistNadly m(X)$. Here, at last, we must leave $m = 5$ behind.

\begin{lemma}\label{Lemma: asymptotic for twN<=y(X) for m = 10 and 25}
	Let $m \in \set{10, 25}$. Suppose $2 \leq y \leq X$. Then
	\begin{equation}
	\twistNadly m(X) = \frac{\Q m \R m X^{1/\degB m}}{\zeta(2)} + O\parent{\max\parent{\frac{X^{1/\degB m} \log y}{y^{1 - 6/\degB m}}, X^{1/(2 \degB m)} y^{3/\degB m}\log X \log^2 y  }}
	\end{equation}
	for $X, y \geq 2$. The constant $\twconst m$ is given in \eqref{Equation: twconst for m = 10 and 25}.
\end{lemma}

\begin{proof}
	Substituting the asymptotic for $\cM m(X; e)$ from \Cref{Lemma: asymptotic for M(X; e) for m = 5 and 10 and 25}(b) into the defining series \eqref{Equation: defining twNadly X} for $\twistNadly m(X)$, we have
	\begin{equation}
        \begin{aligned}
  \twistNadly m(X) =& \sum_{n \leq y} \sum_{e \mid n} \mu\parent{n/e} \frac{\R m \cT m(n) e^{6/\degB m} X^{1/\degB m}}{\zeta(2) n^2 \prod_{\ell \mid n} \parent{1 + \frac{1}{\ell}}} \\
  &+ \sum_{n \leq y} \sum_{e \mid n} \mu\parent{n/e} O\parent{\frac{2^{\omega(n)} e^{3/\degB m} X^{1/(2 \degB m)} \log X}{n}}.
        \end{aligned}
	\end{equation}
	
	We handle the main term and the error of this expression separately. For the main term, we have
	\begin{equation}\label{Equation: partial sums mu r T e X for m = 10 and 25}
	\begin{aligned}
	\sum_{n \leq y} \sum_{e \mid n} \mu\parent{n/e} \frac{\R m \cT m(n) e^{6/\degB m} X^{1/\degB m}}{\zeta(2) n^2 \prod_{\ell \mid n} \parent{1 + \frac{1}{\ell}}} &= \frac{\R m X^{1 / \degB m}}{\zeta(2)} \sum_{n \leq y} \frac{\varphi_{6/\degB m}(n) \cT m(n)}{n^2 \prod_{\ell \mid n} \parent{1 + \frac{1}{\ell}}}.
	\end{aligned}
	\end{equation}
	By \Cref{Lemma: bound on T(e) for m = 5 and 10 and 25}(f), we see
	\begin{equation}\label{Equation: bounding phi(n) T(n)/n^2 ...}
	\frac{\varphi_{6/\degB m}(n) \cT m(n)}{n^2 \prod_{\ell \mid n} \parent{1 + \frac{1}{\ell}}} = O\parent{\frac{2^{\omega(n)}}{n^{2-6/\degB m}}}.
	\end{equation}
	
	By \Cref{Corollary: tail of sum of f/n^sigma} and \Cref{Corollary: sum of 2^omega(n)}, we have
	\begin{equation}
	\sum_{n > y} \frac{2^{\omega(n)}}{n^{2 - 6/\degB m}} \sim \frac{\degB m \log y}{(\degB m - 6) \zeta(2) y^{1 - 6/\degB m}}
	\end{equation}
	as $y \to \infty$. \textit{A fortiori,}
	\begin{equation}
	\sum_{n > y} \frac{\varphi_{6/\degB m}(n) \cT m(n)}{n^2 \prod_{\ell \mid n} \parent{1 + \frac{1}{\ell}}} = O\parent{\sum_{n > y} \frac{2^{\omega(n)}}{n^{2 - 6/\degB m}}} = O\parent{\frac{\log y}{y^{1 - 6/\degB m}}},
	\end{equation}
	so the series
	\begin{equation}
	\sum_{n \geq 1} \frac{\varphi_{6/\degB m}(n) \cT m(n)}{n^2 \prod_{\ell \mid n} \parent{1 + \frac{1}{\ell}}} = \Q m \label{equation: sum for Q for m = 10 and 25}
	\end{equation}
	is absolutely convergent, and 
	\begin{equation}
	\begin{aligned}
	\sum_{n \leq y} \sum_{e \mid n} \mu\parent{n/e} \frac{\R m \cT m(n) e^{6/\degB m} X^{1/\degB m}}{\zeta(2) n^2 \prod_{\ell \mid n} \parent{1 + \frac{1}{\ell}}} &= \frac{\R m X^{1/\degB m}}{\zeta(2)} \parent{\Q {m} + O\parent{\frac{\log y}{y^{1 - 6/\degB m}}}} \\
	&= \twconst m X^{1/\degB m} + O\parent{\frac{X^{1/\degB m} \log y}{y^{1 - 6/\degB m}}}.
	\end{aligned}
	\end{equation}
	
	As the summands of \eqref{equation: sum for Q for m = 10 and 25} constitute a nonnegative multiplicative arithmetic function, we can factor $\Q m$ as an Euler product. We have
    \begin{equation}
	\Q m = \Q m (2) \Q m (5) \prod_{\substack{p \neq 5 \ \textup{prime} \\ p \equiv 1 \psmod {4}}} \parent{1 + \frac{2p \parent{p^{6/\degB m} - 1}}{(p+1)\parent{p^2 - p^{6/\degB m}}}}. \label{Equation: product for Q for m = 10 and 25};
	\end{equation}
    by \Cref{Lemma: bound on T(e) for m = 5 and 10 and 25} the terms $\Q m (p)$ can be computed as follows:
        \begin{equation}\label{Equation: Euler factors for m = 10 and 25}
	\begin{aligned}
	\Q m (p) &\colonequals 1 + \sum_{a \geq 1} \frac{\varphi_{6/\degB m}(p^a) \cT m(p^a)}{\parent{1 + 1/p} p^{2a} } \\
	&= \begin{cases}
	1 + \frac{2p \parent{p^{6/\degB m} - 1}}{(p+1)\parent{p^2 - p^{6/\degB m}}}, & \textup{if $p \equiv 1 \psmod{4}$ and $p \neq 5$;} \\
	\frac 13 \parent{2 + \sqrt{2}}, & \text{if $m = 10$ and $p=2$;} \\
	\frac{2}{31}\parent{15 + 13 \sqrt 5}, & \text{if $m = 10$ and $p=5$;} \\
	\frac{1}{3} \left(2+2^{2/3}\right), & \text{if $m = 25$ and $p=2$;} \\ 
	\frac{2}{781} \left(375+15 \cdot 5^{1/3} + 313 \cdot 5^{2/3}\right), & \text{if $m = 25$ and $p=5$;} \\ 
	1 & \text{else}.
	\end{cases}
	\end{aligned}
        \end{equation}
        The square and cubic roots appear in \eqref{Equation: Euler factors for m = 10 and 25} because of the generalized Jordan totient functions $\varphi_{1/2}$ and $\varphi_{1/3}$. For instance, for $m = 25$ and $p = 2$ we have
    \begin{equation} 	\begin{aligned}
    \Q {25} (2) &= 1 + \frac{\varphi_{1/3}(2) \cT {25}(2)}{\parent{1 + 1/2} 2^{2} } + \frac{\varphi_{1/3}(2^2) \cT {25}(2^2)}{\parent{1 + 1/2} 2^{4} } \\
    &= 1 + \frac{(2^{1/3} - 1) 2}{\parent{1 + 1/2} 2^{2} } + \frac{\parent{2^{2/3} - 2^{1/3}} 2^3}{\parent{1 + 1/2} 2^{4} } \\
    &= \frac{1}{3} \left(2+2^{2/3}\right).
    \end{aligned}\end{equation} 
	
	We now turn to the error term. Since $y \leq X$, for $e \leq y$ we have $\log (e^6 X) \ll \log X$. We obtain
	\begin{equation} \label{equation: partial simplification twN<=y(X) for m = 10 and 25}
	\begin{aligned}
		&\sum_{n \leq y} \sum_{e \mid n} \mu\parent{n/e} O\parent{\frac{2^{\omega(n)} e^{3/\degB m} X^{1/(2 \degB m)} \log X}{n}} \\
    =& O\parent{X^{1/(2 \degB m)} \log X \sum_{e \leq y} \frac{2^{\omega(e)}}{e^{1-3/\degB m}} \sum_{f \leq y/e} \frac{2^{\omega(f)}}{f}}.
	\end{aligned}
	\end{equation}
	Using \Cref{Corollary: sum of 2^omega(n)} and \Cref{Corollary: tail of sum of f/n^sigma} in tandem, we obtain
	\begin{equation} 	\begin{aligned}
	&O\parent{X^{1/(2 \degB m)} \log X \sum_{e \leq y} \frac{2^{\omega(e)}}{e^{1-3/\degB m}} \sum_{f \leq y/e} \frac{2^{\omega(f)}}{f}} \\
        =& O\parent{X^{1/(2 \degB m)} \log X \sum_{e \leq y} \frac{2^{\omega(e)}}{e^{1-3/\degB m}} \log(y/e)} \\
		=& O\parent{X^{1/(2 \degB m)} y^{3/\degB m} \log X \log^2 y},
	\end{aligned}\end{equation} 
	which proves our desired result.
\end{proof}

We emphasize that the proof of \Cref{Lemma: asymptotic for twN<=y(X) for m = 10 and 25} has given us the following Euler product expansion for $\Q m$:
\begin{equation}
\Q m = \Q m (2) \Q m (5) \prod_{\substack{p \neq 5 \ \textup{prime} \\ p \equiv 1 \psmod {4}}} \parent{1 + \frac{2p \parent{p^{6/\degB m} - 1}}{(p+1)\parent{p^2 - p^{6/\degB m}}}}, %
\end{equation}
where 
\begin{equation}\label{Q for small p for m = 10 and 25}
\begin{aligned}
	\Q {10} (2) &= \frac 13 \parent{2 + \sqrt{2}}, \\
	\Q {10} (5) &= \frac{2}{31}\parent{15 + 13 \sqrt 5}, \\
	\Q {25} (2) &= \frac{1}{3} \left(2+2^{2/3}\right), \ \text{and} \\
	\Q {25} (5) &= \frac{2}{781} \left(375+15 \cdot 5^{1/3} + 313 \cdot 5^{2/3}\right). %
\end{aligned}
\end{equation}

We now bound $\twistNadgy m(X)$ for $m \in \set{10, 25}$. Our proof here follows the archetype set by \Cref{Lemma: bound on twN>y(X) for m = 7}.

\begin{lemma}\label{Lemma: bound on twN>y(X) for m = 10 and 25}
	Let $m = 10, 25$. We have
	\begin{equation}
	\twistNadgy m(X) = O\parent{\frac{X^{1/\degB m} \log y}{y^{1-6/\degB m}}}
	\end{equation}
        for $X, y \geq 2$.
\end{lemma}

\begin{proof}
	By \Cref{Lemma: bound on T(e) for m = 5 and 10 and 25}, $\cT m (e) = O(2^{\omega(e)})$, so by \Cref{Lemma: asymptotic for M(X; e) for m = 5 and 10 and 25}, we have
	\begin{equation}
	\cM m (X; e) = O\parent{\frac{2^{\omega(e)} X^{1/\degB m}}{e^2}}.
	\end{equation}
	Set $\alpha=6/\degB m$. By \Cref{Proposition: bounding the summands of twN calE X},
	\begin{equation}
\begin{aligned}
	\abs{\twistNadgy m(X)}
	&\ll\sum_{n>y}\sum_{e\mid n}\abs{\mu(n/e)}
	\frac{2^{\omega(n)}(e^6X)^{1/\degB m}}{n^2}\\
	&=X^{1/\degB m}\sum_{n>y}
	\frac{2^{\omega(n)}}{n^{2-\alpha}}
	\prod_{p\mid n}\parent{1+p^{-\alpha}}.
\end{aligned}
	\end{equation}
	Since $m=10,25$ gives $0<\alpha<1$, \Cref{Lemma: weighted 2omega divisor factor} with $\beta=2-\alpha>1$ yields
	\begin{equation}
	\twistNadgy m(X) = O\parent{\frac{X^{1/\degB m} \log y}{y^{1 - 6/\degB m}}}
	\end{equation}
	as desired.
\end{proof}

We are now ready to prove \Cref{Intro Theorem: asymptotic for twN(X) for m of genus $0$} for $m = 10, 25$, which we restate here with a modestly improved error term in the notations we have established.

\begin{theorem}\label{Theorem: asymptotic for twN(X) for m = 10 and 25}
	Let $m = 10, 25$. Then we have
	\begin{equation}
	\twistNad m(X) = \twconst m X^{1/\degB m} + O\parent{X^{1/(2(\degB m - 3))} \log^{3(\degB m - 5)/(\degB m - 3)} X}
	\end{equation}
	for $X \geq 2$. The constant $\twconst m$ is given in \eqref{Equation: twconst for m = 10 and 25}. The implicit constant depends only on $m$.
\end{theorem}

\begin{proof}
	Let $m = 10, 25$. \Cref{Lemma: asymptotic for twN<=y(X) for m = 10 and 25} and \Cref{Lemma: bound on twN>y(X) for m = 10 and 25} together tell us
	\begin{equation}
	\twistNad m(X) = \twconst m X^{1/\degB m} + O\parent{\max\parent{\frac{X^{1/\degB m} \log y}{y^{1 - 6/\degB m}}, X^{1/(2 \degB m)} y^{3/\degB m}\log X \log^2 y  }}.
	\end{equation}

	We let
	\begin{equation}
	 y = X^{1/(2(\degB m - 3))} \log^{-2\degB m/(\degB m - 3)} X.
	\end{equation}
	For sufficiently large $X$, this choice satisfies $2 \leq y \leq X$. Moreover, $\log y \asymp \log X$, and a direct calculation gives
	\begin{equation}
	\frac{X^{1/\degB m} \log y}{y^{1 - 6/\degB m}} \asymp X^{1/(2 \degB m)} y^{3/\degB m}\log X \log^2 y \asymp X^{1/(2(\degB m - 3))} \log^{3(\degB m - 5)/(\degB m - 3)} X.
	\end{equation}
	The finitely many smaller values of $X$ may be absorbed into the implicit constant, proving the claim.
\end{proof}

\begin{remark}\label{Remark: Why not m = 5}
    When $m = 5$, we can follow the proof of \Cref{Lemma: asymptotic for twN<=y(X) for m = 10 and 25} up through \eqref{Equation: bounding phi(n) T(n)/n^2 ...}, but here we are stymied by a lack of understanding of the series
    \begin{equation}\label{Equation: bad partial sum for m = 5}
    \sum_{n \leq y} \frac{\varphi(n) \cT 5(n)}{n^2 \prod_{\ell \mid n} \parent{1 + \frac{1}{\ell}}}
    \end{equation}
    appearing on the right-hand side of \eqref{Equation: partial sums mu r T e X for m = 10 and 25}, which is $O(\log y)$ when $m = 5$. We suspect that \eqref{Equation: bad partial sum for m = 5} behaves similarly to the harmonic sum $\sum_{n \leq y} 1/n$, and that for appropriately chosen constants $\Q 5$ and $\Q 5\prm$ we may write
    \begin{equation}\label{Equation: conjectured bad partial sum for m = 5}
    \sum_{n \leq y} \frac{\varphi(n) \cT 5(n)}{n^2 \prod_{\ell \mid n} \parent{1 + \frac{1}{\ell}}} = \Q 5 \log y + \Q 5\prm + O(1/y).
    \end{equation}
    Even assuming \eqref{Equation: conjectured bad partial sum for m = 5}, however, we find ourselves obstructed by \Cref{Lemma: bound on twN>y(X) for m = 10 and 25}: to handle $m = 5$, we would require not only a bound on $\twistNadgy 5 (X)$, but an asymptotic estimate for $\twistNadgy 5 (X)$ with a power-saving error term. This is far more than \Cref{Lemma: bound on twN>y(X) for m = 10 and 25} aspires to offer. We conjecture that for every $\epsilon > 0$ we have
    \begin{equation}\label{Equation: conjectured asymptotic for twN (X) for m = 5}
        \twistNeq 5 (X), \twistNad 5 (X) = \twconst 5 X^{1/6} \log X + \twconstprm 5 X^{1/6} + O(X^{1/12 + \epsilon}).
    \end{equation}
    If so, the associated Dirichlet series $\twistLeq 5 (s)$ and $\twistLad 5 (s)$ will have a double pole at $s = 1/6$, in contrast to the Dirichlet series $\twistLad 7 (s)$ studied in \Cref{Corollary: twL(s) has a meromorphic continuation for m = 7}, which has a simple pole at $s = 1/6$. Given \eqref{Equation: conjectured asymptotic for twN (X) for m = 5}, it would be straightforward to obtain asymptotics for $\NQeq 5 (X)$ and $\NQad 5 (X)$ with power-saving error terms.

    We do not believe that the case $m = 5$ is intractable, but we have little hope that the sieving methods we employ in this thesis will unlock this case. We suspect other methods, such as Poisson summation, may achieve better results.
\end{remark}

\subsection*{$L$-series}\label{Subsection: L-series for m = 10 and 25}

To conclude this section, we set up \cref{Section: Working over the rationals for m = 10 and 25} by interpreting
\Cref{Theorem: asymptotic for twN(X) for m = 10 and 25} in terms of Dirichlet series. Recall \eqref{Equation: twisth calE}, \eqref{Equation: hQ calE}, \eqref{Equation: defining twistL calE X}, and \eqref{Equation: defining LQ calE X}.

\begin{cor}\label{Corollary: twL(s) has a meromorphic continuation for m = 10 and 25}
	Let $m = 10, 25$, and put
	\begin{equation}
	\beta_m \colonequals \frac{1}{2(\degB m - 3)}.
	\end{equation}
	The following statements hold.
 \begin{enumerate}
    \item The Dirichlet series $\twistLad m(s)$ has abscissa of (absolute) convergence $\sigma_a=\sigma_c = 1/\degB m$ and has a meromorphic continuation to the region
	\begin{equation}
	\set{s = \sigma + i t \in \bbC : \sigma > \beta_m}. \label{equation: domain of twistL(s) for m = 10 and 25}
	\end{equation}
	\item The function $\twistLad m(s)$ has a simple pole at $s = 1/\degB m$ with residue 
	\begin{equation} \res_{s=\frac{1}{\degB m}} \twistLad m(s) = \frac{\twconst m}{\degB m}; \end{equation}
	it is holomorphic elsewhere on the region \eqref{equation: domain of twistL(s) for m = 10 and 25}.
    \item For $\sigma>\beta_m$, we have
	\begin{equation}
    \mu_{\twistLad m}(\sigma) \leq \max\left\{0,\frac{1/\degB m-\sigma}{1/\degB m-\beta_m}\right\}.
    \end{equation}
     \end{enumerate}
\end{cor}

\begin{proof}
	The proof is structurally identical to the one given for \Cref{Corollary: twL(s) has a meromorphic continuation for m = 7}, with $1/6$ and $1/12$ replaced by $1/\degB m$ and $\beta_m$.
\end{proof}

\section{Estimates for rational isomorphism classes for \texorpdfstring{$m = 10, 25$}{m = 10, 25}}\label{Section: Working over the rationals for m = 10 and 25}

In \Cref{Section: Estimating twN(X) for m = 5 and 10 and 25}, we counted the number of elliptic curves over $\bbQ$ with a (cyclic) $m$-isogeny up to quadratic twist (\Cref{Theorem: asymptotic for twN(X) for m = 10 and 25}) for $m = 10, 25$.  In this section, we count all isomorphism classes over $\bbQ$ by enumerating over twists using Landau's Tauberian theorem (\Cref{Theorem: Landau's Tauberian theorem}). We first describe the analytic behavior of $\LQad m (s)$ for $m = 10, 25$.

\begin{theorem}\label{Theorem: relationship between twistL(s) and L(s) for m = 10 and 25}
Let $m = 10, 25$. The following statements hold.
\begin{enumalph}
	\item The Dirichlet series $\LQad m(s)$ has a meromorphic continuation to the region 
	\begin{equation}\label{Equation: domain of LQ(s) for m = 10 and 25}
	\set{s = \sigma + i t \in \bbC : \sigma > 1/12}
	\end{equation}
	with a simple pole at $s = 1/6$ and no other singularities on this region. 
	\item The principal part of $\LQad m(s)$ at $s = 1/6$ is
	\begin{equation}
	\frac{\twistLad m (1/6)}{3 \zeta(2)} \parent{s - \frac 16}^{-1}.
	\end{equation}
	\end{enumalph}
\end{theorem}

\begin{proof}
	We proceed as in the proof of \Cref{Theorem: relationship between twistL(s) and L(s) for m = 7}. For (a), since $\zeta(s)$ is nonvanishing when $\sigma > 1$, the ratio $\zeta(6s)/\zeta(12s)$ is a meromorphic function for $\sigma > 1/12$. But \Cref{Corollary: twL(s) has a meromorphic continuation for m = 10 and 25} gives a meromorphic continuation of $\twistLad m(s)$ to the region \eqref{Equation: domain of LQ(s) for m = 10 and 25}. By \Cref{Theorem: relationship between twistL(s) and L(s)}, the function $\LQad m(s)$ is a product of these two meromorphic functions on \eqref{Equation: domain of LQ(s) for m = 10 and 25}, and so it is a meromorphic function on this region. The holomorphy and singularity for $\LQad m(s)$ then follow from those of $\twistLad m(s)$ and $\zeta(s)$. 
	
	We conclude (b) by computing Laurent expansions. We recall \eqref{Equation: Laurent expansion for zeta(6s)/zeta(12s)}, and of course the Laurent expansion for $\twistLad m(s)$ at $s = 1/6$ begins
\begin{equation}
\twistLad m(s) = \twistLad m (1/6) + \dots.
\end{equation}
Multiplying the Laurent series tails gives the desired result.
\end{proof}

Using \Cref{Theorem: relationship between twistL(s) and L(s) for m = 10 and 25}, we deduce the following lemma.

\begin{lemma}\label{Lemma: Delta NQ(n) is admissible for m = 10 and 25}
	Let $m = 10, 25$. The sequence $\parent{\hQad m(n)}_{n \geq 1}$ is admissible \textup{(\Cref{Definition: admissible sequences})} with parameters $(1/6,1/12,13/84)$.
\end{lemma}

\begin{proof}
	The proof is similar to, but simpler than, the one given for \Cref{Lemma: Delta NQ(n) is admissible for m = 7}. The critical difference is this: by \Cref{Corollary: twL(s) has a meromorphic continuation for m = 10 and 25}, the Dirichlet series defining $\twistLad m (s)$ converges absolutely when $\sigma = \repart(s) > 1/12$.

    Let $m \in \set{10, 25}$. We check each condition in \Cref{Definition: admissible sequences}. Since $\hQad m(n)$ counts objects, 
we indeed have $\hQad m(n) \in \Z_{\geq 0}$. 
	
	For (i), note $\displaystyle{\frac{\zeta(6s)}{\zeta(12s)}}$ has $1/6$ as its abscissa of absolute convergence. Now by \Cref{Theorem: relationship between twistL(s) and L(s)}(b), we have
	\begin{equation}
	\LQad m(s) = \frac{2\zeta(6s) \twistLad m(s)}{\zeta(12s)},
	\end{equation}
	and by \Cref{Theorem: product of Dirichlet series converges} this series converges absolutely for $\sigma > 1/6$, so the abscissa of absolute convergence for $\LQad m(s)$ is at most $1/6$. For the reverse inequality, choose any twist-minimal curve in the family and let $n_0$ be its twist height. By \Cref{Theorem: relationship between twistL(s) and L(s)}(a), the coefficients of $\LQad m(s)$ satisfy
 \begin{equation}
 \hQad m(n_0c^6) \geq 2\abs{\mu(c)}
 \end{equation}
 for every positive integer $c$. Hence the Dirichlet series for $\LQad m(s)$ contains the divergent positive subseries
 \begin{equation}
 2n_0^{-\sigma}\sum_{c\geq1}\frac{\abs{\mu(c)}}{c^{6\sigma}}
 \end{equation}
 whenever $\sigma\leq1/6$. Thus its abscissa of absolute convergence is at least $1/6$, and (i) holds with $\sigma_a=1/6$.

For (ii), as $\zeta(12s)$ is nonvanishing for $\sigma > 1/12$, we see that $\zeta(6s)/\zeta(12s)$ has a meromorphic continuation to $\sigma>1/12$, and so (ii) holds with 
\begin{equation}
\delta = 1/6 - 1/12=1/12.
\end{equation}
(The only pole of $\LQad m(s)/s$ with $\sigma > 1/12$ is the simple pole at $s = 1/6$ indicated in 
\Cref{Theorem: relationship between twistL(s) and L(s) for m = 10 and 25}(b).)
	
For (iii), let $\sigma > 1/12$. By \Cref{Theorem: absolutely convergent Dirichlet series have mu = 0}, $\mu_{\twistLad m}(\sigma) = 0$. Recall the notation $\zeta_a(s)=\zeta(as)$. As in \eqref{Equation: bound on zeta6(s)}, we have
\begin{equation}
\mu_{\zeta_6}(\sigma) < \frac{13}{84} 
\end{equation}
if $\sigma \leq 1/6$, and by \Cref{Theorem: absolutely convergent Dirichlet series have mu = 0}, $\mu_{\zeta_6}(\sigma) = 0$ if $\sigma > 1/6$. Finally, as $\zeta(12s)^{-1}$ is absolutely convergent for $\sigma > 1/12$, \Cref{Theorem: absolutely convergent Dirichlet series have mu = 0} tells us $\mu_{{\zeta_{12}}^{-1}}(\sigma) = 0$. Taken together, we see
	\begin{equation}
	\mu_{\LQad m}(\sigma) < 0 + \frac{13}{84} + 0 = \frac{13}{84},
	\end{equation}
	so the sequence $\parent{\hQad m(n)}_{n \geq 1}$ is admissible with final parameter $\xi = 13/84$.
\end{proof}

We now prove \Cref{Intro Theorem: asymptotic for NQ(X) for m > 9} for $m = 10, 25$, which we restate here for ease of reference in our established notation.

\begin{theorem}\label{Theorem: asymptotic for NQ(X) for m = 10 and 25}
	Let $m = 10, 25$, and define
        \begin{equation}
	\begin{aligned}
	\constad m &\colonequals \frac{2 \twistLad m (1/6)}{\zeta(2)}.
	\end{aligned}
 \end{equation}
Then for all $\epsilon > 0$,  we have
		\begin{equation}
	\NQad m(X) = \constad m X^{1/6} + O\parent{X^{1/8 + \epsilon}}
	\end{equation}
	for $X \geq 1$. The implicit constant depends only on $\epsilon$. 
\end{theorem}

\begin{proof}
	By \Cref{Lemma: Delta NQ(n) is admissible for m = 10 and 25}, the sequence $\parent{\hQad m(n)}_{n \geq 1}$ is admissible with parameters $\parent{1/6, 1/12, 13/84}$. We now apply \Cref{Theorem: Landau's Tauberian theorem} to the Dirichlet series $\LQad m(s)$, and our claim follows. 
\end{proof}

\begin{remark}\label{Remark: true bound for twistN(X) for m = 10 and 25}
	We suspect that the true error on $\NQad m(X)$ is at most $O(X^{1/12 + \epsilon})$, and the true error on $\twistNad m (X)$ is at most $O(X^{1/(2 \degB m) + \epsilon})$, but we have been unable to bound the error terms this far using our techniques. See \Cref{Remark: true bound for twistN(X) for m = 7} for some related thoughts. 
\end{remark}

\section{Computations for \texorpdfstring{$m = 10, 25$}{m = 10, 25}}\label{Section: Computations for m = 10 and 25}

In this section, we furnish computations that render \Cref{Theorem: asymptotic for twN(X) for m = 10 and 25} and \Cref{Theorem: asymptotic for NQ(X) for m = 10 and 25} completely explicit.

\subsection*{Enumerating elliptic curves with \texorpdfstring{$m$}{m}-isogeny for \texorpdfstring{$m = 10, 25$}{m = 10, 25}}

The algorithm described in \cref{Section: Computations for m = 7} can be adapted to enumerate elliptic curves admitting a cyclic $m$-isogeny for $m = 10, 25$. Doing so requires paying special attention to the primes $2$ and $5$, rather than $3$ and $7$, and requires writing $\C m (a, b) = e_0^2 n_0$, with $n_0$ squarefree. In step 3, if $e_0^2$ lies in the lookup range, we multiply representatives of norm $e_0^2$ by representatives of norm $n_0$; otherwise we square representatives of norm $e_0$ first. When $m = 10$, the lookup table in step 2 is restricted to pairs $(a,b)$ with $a \geq b$; in step 3 we also use their conjugates. The bounded components at $2$ and $5$ are handled separately.

For $m = 10$, enumerating out to $X = 10^{96}$ produces 106\,785\,276 elliptic curves admitting a cyclic $10$-isogeny. For $X = 10^{96}$, we have
\begin{equation}
\frac{ \twistNad {10}(10^{96})}{\twconst {10}(10^{96})^{1/12}} = 0.999\,671\ldots,\label{Equation: Ratio for twN10}
\end{equation}
which is close to 1. We compute $\twconst {10}$ below.

For $m = 25$, enumerating out to $X = 10^{138}$ produces 35\,007\,063 elliptic curves admitting a cyclic $25$-isogeny. For $X = 10^{138}$, we have
\begin{equation}
\frac{ \twistNad {25}(10^{138})}{\twconst {25}(10^{138})^{1/18}} = 0.999\,937\ldots,\label{Equation: Ratio for twN25}
\end{equation}
which is close to 1. We compute $\twconst {25}$ below. It is interesting to note that the ratios in \eqref{Equation: Ratio for twN7}, \eqref{Equation: Ratio for twN10}, \eqref{Equation: Ratio for twN25} are all less than 1. We do not know if this bias is systematic or coincidental.

We list some twist minimal elliptic curves admitting a cyclic $10$-isogeny in Table \ref{table:firstfew10}, and the first few twist minimal elliptic curves admitting a cyclic $25$-isogeny in Table \ref{table:firstfew25}.

\jvtable{table:firstfew10}{\small
\rowcolors{2}{white}{gray!10}
\begin{tabular}{c c c c}
$(A, B)$ & $(a, b)$ & $\twistheight(E)$ & $\twistdefect(E)$ \\
\hline\hline
$(6, 88)$ & $(3, 1)$ & $209088$ & $125$ \\
$(-66, 200)$ & $(-1, 3)$ & $1149984$ & $125$ \\
$(-435, 4750)$ & $(4, 3)$ & $609187500$ & $250$ \\
$(6981, 92950)$ & $(8, 1)$ & $1360858296564$ & $250$ \\
$(-7635, 256750)$ & $(-3, 4)$ & $1780275091500$ & $125$ \\
$(-8130, 187000)$ & $(1, 7)$ & $2149471188000$ & $125$ \\
$(-4035, 474050)$ & $(2, 1)$ & $6067531867500$ & $2$ \\
$(-26571, 1570426)$ & $(2, 9)$ & $75038421469644$ & $250$ \\
$(-29370, 1937000)$ & $(-7, 1)$ & $101337883812000$ & $125$ \\
$(-30459, 774358)$ & $(-1, 8)$ & $113033431970316$ & $125$ \\
$(-65091, 6383806)$ & $(-4, 7)$ & $1103120162194284$ & $250$ \\
$(-77979, 8511050)$ & $(6, 7)$ & $1955825246767500$ & $250$ \\
$(-46371, 10131550)$ & $(7, 4)$ & $2771504245867500$ & $125$ \\
$(-119235, 15795650)$ & $(-1, 2)$ & $6780648933211500$ & $1$ \\
$(-280227, 56930654)$ & $(-12, 1)$ & $88021734784228332$ & $250$ \\
$(-405507, 980606)$ & $(1, 12)$ & $266719677879435372$ & $125$ \\
$(-418251, 104112250)$ & $(-9, 2)$ & $292665112764269004$ & $125$ \\
$(-504570, 137620600)$ & $(1, 3)$ & $513835691175972000$ & $1$ \\
\end{tabular}
}{Some $E \in \twistE$ with a cyclic $10$-isogeny and $\twht E \leq 10^{18}$}

\jvtable{table:firstfew25}{\small
\rowcolors{2}{white}{gray!10}
\begin{tabular}{c c c c}
$(A, B)$ & $(a, b)$ & $\twistheight(E)$ & $\twistdefect(E)$ \\
\hline\hline
$(-12, 38)$ & $(-4, 1)$ & $38988$ & $12500$ \\
$(-4035, 98750)$ & $(-3, 2)$ & $263292187500$ & $3125$ \\
$(-8634, 308792)$ & $(6, 1)$ & $2574519136416$ & $12500$ \\
$(-11586, 480008)$ & $(-2, 3)$ & $6221007361728$ & $12500$ \\
$(-281532, 57496282)$ & $(0, 1)$ & $89257205983235148$ & $4$ \\
$(622149, 500328938)$ & $(-9, 1)$ & $6758884247405611788$ & $3125$ \\
$(-1768386, 917586232)$ & $(-2, 1)$ & $22733041315210861248$ & $4$ \\
$(-2010243, 1096965250)$ & $(11, 1)$ & $32494186355919275628$ & $15625$ \\
$(-3333819, 2450621162)$ & $(-7, 3)$ & $162149690150340216588$ & $3125$ \\
$(-4367235, 3512882050)$ & $(-1, 1)$ & $333189188024729467500$ & $1$ \\
$(-5840211, 5432389742)$ & $(7, 2)$ & $796793174499269255724$ & $3125$ \\
$(-6208059, 5953630358)$ & $(-1, 4)$ & $957034289871878620428$ & $3125$ \\
$(-6915540, 6999826250)$ & $(8, 3)$ & $1322934323315597856000$ & $12500$ \\
$(8365830, 6918545000)$ & $(-14, 1)$ & $2342001110975069148000$ & $12500$ \\
$(-23656314, 44286231688)$ & $(2, 1)$ & $52954298562326815548576$ & $4$ \\
$(-149675916, 704409673682)$ & $(16, 1)$ & $13412686238635561555901184$ & $12500$ \\
$(18529341, 811299953342)$ & $(-13, 2)$ & $17771605585903747178162028$ & $3125$ \\
$(-273426411, 1988757501158)$ & $(-11, 4)$ & $106789222757129734026206028$ & $3125$ \\
$(-275757339, 2198498350282)$ & $(-3, 1)$ & $130501664897202240375947148$ & $1$ \\
$(-463781604, 3844524236618)$ & $(-8, 7)$ & $399069898360466822606103948$ & $12500$ \\
$(-593007330, 5558251655000)$ & $(2, 7)$ & $834142359428376453675000000$ & $12500$ \\
\end{tabular}
}{$E \in \twistE$ with a cyclic $25$-isogeny and $\twht E \leq 10^{27}$}

\subsection*{Computing \texorpdfstring{$\twconst {10}$ and $\twconst {25}$}{ctw10 and ctw25}}

In this subsection, for $m = 10, 25$, we estimate the constant $\twconst {m}$ appearing in \Cref{Theorem: asymptotic for twN(X) for m = 10 and 25} by estimating $\Q {m}$ and $\R {m}$.

We begin with $\Q m$. Letting $m = 10, 25$, truncating the Euler product \eqref{Equation: product for Q for m = 10 and 25} as a product over $p \leq Y$ gives the lower bound
\begin{equation}
Q_{m,\leq Y} \colonequals \Q m (2) \Q m(5)
\prod_{\substack{p \leq Y,\ p \neq 5 \ \textup{prime} \\ p \equiv 1 \psmod {4}}} \parent{1 + \frac{2p \parent{p^{6/\degB m} - 1}}{(p+1)\parent{p^2 - p^{6/\degB m}}}}.
\end{equation}
Put $\alpha_m=6/\degB m$ and $s_m=2-\alpha_m$. For $p>Y$, the nontrivial local increment is less than $2p^{-s_m}/(1-Y^{-s_m})$. Thus, by Abel summation and Bennett--Martin--O'Bryant--Rechnitzer \cite[Theorem 1.4]{Bennett-Martin-OBryant-Rechnitzer}, for $Y \geq 8\cdot10^9$,
\begin{equation}
\Q m < Q_{m,\leq Y}\exp\parent{
\frac{Y^{1-s_m}}{(s_m-1)\log Y\parent{1-Y^{-s_m}}}
\parent{1+\frac{5s_m}{2\log Y}}}.
\end{equation}
In particular, taking $Y = 10^{11}$ gives
\begin{equation}
3.636\,493\,079\,001\,437\,6 < \Q {10} < 3.636\,494\,121\,479\,969
\end{equation}
and
\begin{equation}
4.244\,853\,881\,138\,272\,6 < \Q {25} < 4.244\,853\,894\,726\,220.
\end{equation}

We now turn our attention to $\R m$ for $m = 10, 25$, given in \eqref{eqn: R(X)}. We compute $\R {10}$ and $\R {25}$ by performing rejection sampling on the rectangles $[-0.8228, 0.8228] \times [0, 0.6934]$ and $[-0.8781, 0.8781] \times [0, 0.2754]$ respectively. 

We find $r_{10} \colonequals 58\,560\,198\,103$ of our first $s_{10} \colonequals 138\,290\,000\,000$ samples lie in $\calR {10}(1)$, so
\begin{equation}
\R {10} \approx 1.141\,059\,04 \cdot\frac{r_{10}}{s_{10}} = 0.483\,192\,157\,275\,428\,47,
\end{equation}
with standard error
\begin{equation}
1.141\,059\,04 \cdot \sqrt{\frac{r_{10}(s_{10} - r_{10})}{s_{10}^3}}
=1.516\,12\ldots \cdot 10^{-6}.
\end{equation}
This took 2 CPU weeks to compute. %

We find $r_{25} \colonequals 245\,430\,977\,211$ of our first $s_{25} \colonequals 406\,130\,000\,000$ samples lie in $\calR {25}(1)$, so
\begin{equation}
\R {25} \approx 0.483\,657\,48 \cdot\frac{r_{25}}{s_{25}} = 0.292\,282\,096\,746\,878\,3,
\end{equation}
with standard error
\begin{equation}
0.483\,657\,48 \cdot \sqrt{\frac{r_{25}(s_{25} - r_{25})}{s_{25}^3}} < 3.8 \cdot 10^{-7}. %
\end{equation}
This took 13 CPU weeks to compute. %

Using the finite Euler products as point estimates, we therefore have
\begin{equation} \begin{aligned}
	\twconst {10} &\approx 1.068\,203\,87 &&\text{with Monte Carlo standard error }3.36\cdot10^{-6}, \\
	\twconst {25} &\approx 0.754\,251\,990 &&\text{with Monte Carlo standard error }9.58\cdot10^{-7}.
\end{aligned}\end{equation}
The corresponding errors from truncating the Euler products are below $3.07\cdot10^{-7}$ and $2.42\cdot10^{-9}$, respectively.

We have computed the constants which appear in \Cref{Theorem: asymptotic for twN(X) for m = 10 and 25}.

\subsection*{Computing \texorpdfstring{$\constad {10}$ and $\constad {25}$}{c10 and c25}}\label{Subsection: computing c for m = 10, 25}

In this subsection, we estimate $\constad {10} = 2 \twistLad {10}(1/6)/\zeta(2)$ and $\constad {25} = 2 \twistLad {25} (1/6)/\zeta(2)$, the constants which appear in \Cref{Theorem: asymptotic for NQ(X) for m = 10 and 25}, using the partial sums
\begin{equation} 	\begin{aligned}
	\sum_{n \leq 10^{96}} \frac{\twisthad {10}(n)}{n^{1/6}} &= 0.411\,595\,409\,503\,702\,919, \ \text{and} \\
	\sum_{n \leq 10^{138}} \frac{\twisthad {25}(n)}{n^{1/6}} &= 0.206\,338\,924\,692\,701\,754.
\end{aligned}\end{equation} 
We empirically confirm that 
\begin{align}
\twistNad {10} (X) &< 1.095\,26 X^{1/12} \ \text{for} \ X \leq 10^{96}, \ \text{and} \label{Equation: bound on twistN10} \\
\twistNad {25} (X) &< 0.909\,77 X^{1/18} \ \text{for} \ X \leq 10^{138}. \label{Equation: bound on twistN25}
\end{align}
If these bounds continue to hold for larger $X$, partial summation gives
\begin{equation}
\sum_{n>Y}\frac{\twisthad m(n)}{n^{1/6}}
=-Y^{-1/6}\twistNad m(Y)+\frac16\int_Y^\infty \twistNad m(x)x^{-7/6}\,\mathrm{d}x.
\end{equation}
Using $\twistNad {10}(10^{96})=106\,785\,276$ and $\twistNad {25}(10^{138})=35\,007\,063$, the two tails are respectively less than $1.12266724\cdot10^{-8}$ and $2.83347\cdot10^{-16}$. Hence
\begin{equation} \begin{aligned}
\twistLad {10}(1/6) &\approx 0.411595415117039 &&\text{with conditional error }<5.62\cdot10^{-9},\\
\twistLad {25}(1/6) &\approx 0.20633892469270190 &&\text{with conditional error }<1.42\cdot10^{-16},
\end{aligned}\end{equation}
and therefore
\begin{equation} \begin{aligned}
	\constad {10} &\approx 0.500440015697013 &&\text{with conditional error }<6.83\cdot10^{-9}, \\
	\constad {25} &\approx 0.25087804897622103 &&\text{with conditional error }<1.73\cdot10^{-16}.
\end{aligned}\end{equation} 

We emphasize that our estimates for $\constad {10}$ and $\constad {25}$ depend on empirical rather than theoretical estimates for the implicit constant in the error term in the asymptotics of $\twistNad {10}(X)$ and $\twistNad {25} (X)$.

%% file: Ch6_m13.tex
\chapter{Counting elliptic curves with a cyclic \texorpdfstring{$m$}{m}-isogeny for \texorpdfstring{$m = 13$}{m = 13}}\label{Chapter: m = 13}

In this chapter, we prove \Cref{Intro Theorem: asymptotic for NQ(X) for m > 9} (\Cref{Theorem: asymptotic for NQ(X) for m = 13}) and \Cref{Intro Theorem: asymptotic for twN(X) for m of genus $0$} (\Cref{Theorem: asymptotic for twN(X) for m = 13}) when $m = 13$. These results are new, but our arguments mirror those in \cref{Chapter: m = 7} and \cref{Chapter: m = 10 and 25}, and we encourage anyone reading to skim them on a first perusal of this thesis. However, handling $m = 13$ is subtler than handling $m \in \set{7, 10, 25}$ because $X_0(13)$ has elliptic points of both orders 2 and 3, and the elliptic surfaces describing elliptic curves with a cyclic $13$-isogeny thus exhibit \emph{both} potential type II additive reduction and potential type III additive reduction.

The organization of this chapter mirrors that of \cref{Chapter: m = 7} and \cref{Chapter: m = 10 and 25}. In \cref{Section: Establishing notation when m = 13}, we establish notations pertaining to $\f {13} (t)$ and $\g {13} (t)$ which will be used throughout the remainder of the chapter. In \cref{Section: twist minimality defect for m = 13}, we develop bounds relating the twist minimality defect to the two factors of the greatest common divisor of $\f {13} (t)$ and $\g {13} (t)$. In \cref{Section: Estimating twN(X) for m = 13}, we apply the framework developed in \cref{Section: Our approach revisited} to prove \Cref{Intro Theorem: asymptotic for twN(X) for m of genus $0$} for $m = 13$. In \cref{Section: Working over the rationals for m = 13}, we prove \Cref{Intro Theorem: asymptotic for NQ(X) for m > 9} for $m = 13$. In \cref{Section: Computations for m = 13}, we produce supplementary computations to estimate the constants appearing in \Cref{Theorem: asymptotic for twN(X) for m = 13} and \Cref{Theorem: asymptotic for NQ(X) for m = 13} and empirically confirm that the count of elliptic curves with a cyclic $13$-isogeny aligns with our theoretical estimate.

\section{Establishing notation for \texorpdfstring{$m = 13$}{m = 13}}\label{Section: Establishing notation when m = 13}

By \Cref{Corollary: Elliptic curves for which twNeq(X) = twNad(X)}, 
\begin{equation}
\twNeq {13} (X) = \twNad {13} (X) \ \text{and} \ \NQeq {13} (X) = \NQad {13} (X)
\end{equation}
for all $X > 0$, so we may use either notation interchangeably. We opt to work with $\twNad {13} (X)$ and related functions.

Note that 
\begin{equation}
\gcd(\f {13} (t), \g {13} (t)) = (t^2 + t + 7) (t^2 + 4)
\end{equation}
factors over $\bbQ$. We define
\begin{equation}
\begin{aligned}
	\hII {13} (t) &\colonequals t^2 + t + 7, \ \text{and} \\
	\hIII {13} (t) &\colonequals t^2 + 4,
\end{aligned}
\end{equation}
so $\gcd(\f {13} (t), \g {13} (t)) = \hII {13} (t)\hIII {13} (t)$.
We define $\tf {13} (t)$ and $\tg {13} (t)$ so that
\begin{equation} \label{equation: Definition of tf and tg for {13}-isogenies}
\f {13} (t) = \tf {13} (t) \hII {13} (t) \hIII {13}(t) \ \text{and} \ \g {13} (t) = \tg {13} (t) \hII {13} (t) \hIII {13} (t)^2.
\end{equation}
Thus
\begin{equation}
\begin{aligned}
	\tf {13} (t) =& -3 \parent{t^4 - 235 t^3 + 1211 t^2 - 1660 t + 6256}, \ \text{and} \\
	\tg {13} (t) =& 2 (t^6 + 512 t^5 - 13073 t^4 + 34860 t^3 - 157099 t^2 + 211330 t - 655108).
\end{aligned}
\end{equation}
As in the previous chapters, to work with integral models, we take $t=a/b$ (in lowest terms) and homogenize, obtaining
\begin{equation} \label{equation: defining Cab for m = {13}}
\begin{aligned}
	{\CII {13}}(a, b) &\colonequals b^2 \hII {13}(a/b)=a^2 + ab + 7 b^2, \\
	\CIII {13}(a, b) &\colonequals b^2 \hIII {13}(a/b) = a^2 + 4 b^2, \\
	\tA {13} (a, b) &\colonequals b^4 \tf {13} (a/b) \\
        &= -3 \parent{a^4 - 235 a^3 b + 1211 a^2 b^2 - 1660 a b^3 + 6256 b^4}, \ \text{and} \\
	{\tB {13}} (a, b) &\colonequals b^6 \tg {13} (a/b) \\
        &= 2 (a^6 + 512 a^5 b - 13073 a^4 b^2 + 34860 a^3 b^3 - 157099 a^2 b^4 + 211330 a b^5 - 655108 b^6).
\end{aligned}
\end{equation}
Of course, we have
\begin{equation} 	\begin{aligned}
\A {13} (a, b) =& \tA {13} (a, b) {\CII {13}} (a, b) {\CIII {13}} (a, b), \ \text{and} \\
\B {13} (a, b) =& \tB {13} (a, b) {\CII {13}} (a, b) {\CIII {13}} (a, b)^2.
\end{aligned}\end{equation} 

\section{The twist minimality defect for \texorpdfstring{$m = 13$}{m = 13}}\label{Section: twist minimality defect for m = 13}

As with the previous chapters, we begin by studying the twist minimality defect. The situation here is complicated somewhat, however, because the twist minimality defect may receive contributions from both $\CII {13} (a, b)$ and $\CIII {13} (a, b)$.

\begin{lemma}\label{Lemma: 2 and 3 and 13 are the ungroomed primes for m = 13}
	Let $(a,b) \in \Z^2$ be $13$-groomed, let $\ell$ be prime, and let $v \in \Z_{>0}$. Then the following statements hold.
	\begin{enumalph}
	\item If $\ell \neq 2, 3, 13$, then $\ell^v \mid \twistdefect(\A {13}(a,b),\B {13}(a,b))$ if and only if 
 \begin{equation}
 \ell^{3v} \mid {\CII {13}}(a, b) \ \text{or} \ \ell^{2v} \mid {\CIII {13}}(a, b). 
 \end{equation}
 Moreover, for $\ell \neq 13$, we cannot have both $\ell \mid {\CII {13}}(a, b)$ and $\ell \mid {\CIII {13}}(a, b)$.
	\item $\ell^{3v} \mid {\CII {13}}(a,b)$ if and only if $\ell \nmid b$ and $\hII {13}(a/b) \equiv 0 \pmod{\ell^{3v}}$. Likewise, $\ell^{2v} \mid {\CIII {13}}(a,b)$ if and only if $\ell \nmid b$ and $\hIII {13}(a/b) \equiv 0 \pmod{\ell^{2v}}$.
	\item If $\ell \neq 3$, then $\ell \mid {\CII {13}}(a,b)$ implies $\ell \nmid (\partial {\CII {13}}/\partial a)(a,b) = 2a+b$. Likewise, if $\ell \neq 2$, then $\ell \mid {\CIII {13}}(a,b)$ implies $\ell \nmid (\partial {\CIII {13}}/\partial a)(a,b) = 2a$.
	\end{enumalph}
\end{lemma}

\begin{proof}
	We argue as in Cullinan--Kenney--Voight \cite[Proof of Theorem 3.3.1, Step 3]{Cullinan-Kenney-Voight}. Our argument is more involved than the proofs of \Cref{Lemma: 3 and 7 are the ungroomed primes for m = 7} or \Cref{Lemma: 2 and 5 are the ungroomed primes for m = 5 and 10 and 25}, however.
	For part (a), we first compute the resultants
	\begin{equation}\label{Equation: Resultant of CII and CIII}
	\begin{aligned}
        &\Res(\CII {13}(t,1),{\CIII {13}}(t,1))=\Res(\hII {13}(t),\hIII {13}(t)) \\
        &= 13 = \Res(\CII {13}(1,u),{\CIII {13}}(1,u)).
        \end{aligned}
	\end{equation}
	Thus if $\ell \neq 13$ is prime, then $\ell$ can divide at most one of $\CII {13}(a, b)$ and $\CIII {13}(a, b)$. We now compute the resultant
	\begin{equation}
        \begin{aligned}
	&\Res(\tA {13}(t,1),{\tB {13}}(t,1))=\Res(\tf {13}(t),\tg {13}(t)) \\
        &= - 2^{14} \cdot 3^{11} \cdot 13^{24} = \Res(\tA {13}(1,u),{\tB {13}}(1,u)).
        \end{aligned}
	\end{equation}
If $\ell \not\in \set{2,3,13}$, then $\ell \nmid \gcd(\tA {13}(a,b),\tB {13}(a,b))$. Recalling the formula for the twist minimality defect \eqref{Equation: twist defect powers}, if $\ell^v \mid \tmd(\A {13}(a,b),\B {13}(a,b))$ then $\ell \mid {\CII {13}}(a,b) {\CIII {13}}(a,b)$. By \eqref{Equation: Resultant of CII and CIII}, we now have two cases: $\ell \mid {\CII {13}}(a,b)$ or $\ell \mid {\CIII {13}}(a,b)$.

Suppose first that $\ell \mid {\CII {13}}(a,b)$. We compute
\begin{equation} 	\begin{aligned}
\Res(\tB {13}(t,1),{\CII {13}}(t,1))=\Res(\tg {13}(t),\hII {13}(t)) =& 2^8 \cdot 3^3 \cdot 13^6 \\
=& \Res(\tB {13}(1,u),{\CII {13}}(1,u)),
\end{aligned}\end{equation} 
so $\ell \nmid \gcd(\tB {13}(a,b), {\CII {13}}(a, b))$, and thus (under our hypotheses) 
\begin{equation}
\ell^v \mid \tmd(\A {13}(a,b),\B {13}(a,b)) \ \text{if and only if} \ \ell^{3v} \mid {\CII {13}}(a,b).
\end{equation}

Suppose instead that $\ell \mid \CIII {13}(a, b)$. We compute
\begin{equation} 	\begin{aligned}
\Res(\tA {13}(t,1),{\CIII {13}}(t,1))=\Res(\tf {13}(t),\hIII {13}(t)) =& 2^4 \cdot 3^4 \cdot 13^4 \\
=& \Res(\tA {13}(1,u),{\CIII {13}}(1,u)),
\end{aligned}\end{equation} 
so $\ell \nmid \gcd(\tA {13}(a,b), {\CIII {13}}(a, b))$, and thus (under our hypotheses) 
\begin{equation}
\ell^v \mid \tmd(\A {13}(a,b),\B {13}(a,b)) \ \text{if and only if} \ \ell^{2v} \mid {\CIII {13}}(a,b).
\end{equation}
This proves (a).

For (b), by homogeneity it suffices to show that $\ell \nmid b$, and indeed this holds since if $\ell \mid b$ then $\A {13}(a,0) \equiv -3a^8 \equiv 0 \pmod{\ell}$ and $\B {13}(a,0) \equiv 2a^{12} \equiv 0 \pmod{\ell}$ so $\ell \mid a$, a contradiction.
	
Part (c) follows from (b) and the fact that $\hII {13}(t)$ has discriminant $-3^3$ and $\hIII {13}(t)$ has discriminant $-2^4$. 
\end{proof}

We now make our main departure from \cref{Chapter: m = 7} and \cref{Chapter: m = 10 and 25}: in contrast with \Cref{Definition: T(e) for m = 7} and \Cref{Definition: T(e) for m = 10 and 25}, we define $\tcalT {13}$ to be a function with two arguments.

\begin{definition}\label{Definition T(e) for m = 13}
	For $e_1, e_2 \geq 1$, let $\tcalT {13}(e_1, e_2)$ denote the image of
	\begin{equation}\label{Equation: Set mapping to tcalT 13}
\set{(a, b) \in \bbZ^2 : \begin{minipage}{52ex} $(a, b) \ \text{$13$-groomed}, \ e_1 e_2 \mid \twistdefect(\A {13} (a, b), \ \B {13} (a, b)),$ \\ 
$\gcd(e_1^3, \CIII {13} (a, b)) \mid 13, \ \gcd(e_2^2, \CII {13} (a, b)) \mid 13$, \\
$3 \mid e_2 \implies 3 \mid \CIII {13} (a, b),$ \\
$13 \mid e_2 \implies 13 \nmid \CII {13} (a, b) \ \text{or} \ 13^2 \mid \CIII {13} (a, b)$
\end{minipage}}
\end{equation}
under the projection
\begin{equation}
\bbZ^2 \to (\bbZ / e_1^3 e_2^2 \bbZ)^2,
\end{equation}
and let $\tcT {13}(e_1, e_2) \colonequals \# \tcalT {13}(e_1, e_2)$. 

Similarly, we let $\calT {13}(e_1, e_2)$ denote the image of
\begin{equation}\label{Equation: Set mapping to calT 13}
\set{t \in \bbZ : \begin{minipage}{46ex} 
$(e_1 e_2)^2 \mid \f {13}(t), \ (e_1 e_2)^3 \mid \g {13}(t),$ \\
$\gcd(e_1^3, \hIII {13} (t)) \mid 13, \ \gcd(e_2^2, \hII {13} (t)) \mid 13$, \\
$3 \mid e_2 \implies 3 \mid \hIII {13} (t),$ \\ 
$13 \mid e_2 \implies 13 \nmid \hII {13} (t) \ \text{or} \ 13^2 \mid \hIII {13} (t)$
\end{minipage}
}
\end{equation}
under the projection
\begin{equation}\label{Equation: Set mapping to calT 13 lesser}
\bbZ \to \bbZ / e_1^3 e_2^2 \bbZ,
\end{equation}
and let $\cT {13}(e_1, e_2) \colonequals \#\calT {13}(e_1, e_2)$.
\end{definition}

By \Cref{Lemma: 2 and 3 and 13 are the ungroomed primes for m = 13}, for $\ell \neq 2, 3, 13$, $e_1$ is the part of the twist minimality defect arising from $\CII {13} (a, b)$ and $e_2$ is the part of the twist minimality defect arising from $\CIII {13} (a, b)$. The final two conditions of \eqref{Equation: Set mapping to tcalT 13} and \eqref{Equation: Set mapping to calT 13} are necessary to avoid double-counting certain pairs $(a, b)$ for which $3 \mid \tmd(\A {13} (a, b), \B{13} (a, b))$ or $13 \mid \tmd(\A {13} (a, b), \B{13} (a, b))$.

\begin{lemma}\label{Lemma: bound on T(e) for m = 13}
The following statements hold.
\begin{enumalph}
\item Suppose $\gcd(e_1, e_2) = 1$, and write $e_1 = 2^{u_1} 3^{v_1} 13^{w_1} e_1\prm$ and $e_2 = 2^{u_2} 3^{v_2} 13^{w_2} e_2\prm$, where $\gcd(e_1\prm, 2 \cdot 3 \cdot 13) = \gcd(e_2\prm, 2 \cdot 3 \cdot 13) = 1$. The set $\tcalT {13}(e_1, e_2)$ consists of those pairs $(a, b) \in (\bbZ / e_1^3 e_2^2 \bbZ)^2$ which satisfy the following conditions:
\begin{itemize}
	\item $\CII {13}(a, b) \equiv 0 \psmod {(e_1\prm)^3}$,
	\item $\CIII {13} (a, b) \equiv 0 \psmod {(e_2\prm)^2}$,
	\item $\ell \nmid \gcd(a, b)$ for all primes $\ell \mid e_1 e_2$,
 	\item if $u_2 > 0$ then $\A {13} (a, b) \equiv 0 \psmod {2^{2u_2}}$ and $\B {13} (a, b) \equiv 0 \psmod {2^{3 u_2}}$, and if $u_1 > 0$ then no pairs are permitted;
        \item if $v_1 > 0$ then $\A {13} (a, b) \equiv 0 \psmod {3^{2v_1}}$ and $\B {13} (a, b) \equiv 0 \psmod {3^{3 v_1}}$, and if $v_2 > 0$ then no pairs are permitted;
        \item If $w_1 > 0$ then $\A {13} (a, b) \equiv 0 \psmod {13^{2w_1}}$ and $\B {13} (a, b) \equiv 0 \psmod {13^{3 w_1}}$, but $\CIII {13} (a, b) \not\equiv 0 \psmod {13^2}$;
        \item If $w_2 > 0$ then $\A {13} (a, b) \equiv 0 \psmod {13^{2w_2}}$ and $\B {13} (a, b) \equiv 0 \psmod {13^{3 w_2}}$, but $\CII {13} (a, b) \not\equiv 0 \psmod {13}$ or $\CIII {13} (a, b) \equiv 0 \psmod {13^2}$.
\end{itemize}
\item Let $(a, b) \in \bbZ^2$. If $(a, b) \psmod{e_1^3 e_2^2} \in \tcalT {13}(e_1, e_2)$, then 
\[e_1 e_2 \mid \tmd(\A {13} (a, b), \B {13} (a, b)).\] 
\item For all $e_1, e_2, e_1\prm, e_2\prm \in \bbZ_{>0}$ with 
\begin{equation}\label{Equation: T13 gcd conditions}
\gcd(e_1, e_2) = \gcd(e_1\prm, e_2\prm) = \gcd(e_1, e_1\prm) = \gcd(e_2, e_2\prm) = 1,
\end{equation}
we have 
\begin{equation}\label{Equation: multiplicativity of T13}
\begin{aligned}
\tcT {13}(e_1 e_1\prm, e_2 e_2\prm) &= \tcT {13}(e_1, e_2) \tcT {13}(e_1\prm, e_2\prm) \ \text{and} \\ 
\cT {13}(e_1 e_1\prm, e_2 e_2\prm) &= \cT {13}(e_1, e_2) \cT {13}(e_1\prm, e_2\prm),
\end{aligned}
\end{equation}
and 
\begin{equation}\label{Equation: Relating tcT 13 to cT 13}
    \tcT {13}(e_1, e_2) = \varphi(e_1^3) \varphi(e_2^2) \cT {13}(e_1, e_2).
\end{equation}
If $\gcd(e_1, e_2) > 1$, then
\begin{equation}\label{Equation: T13 vanishes when gcd is not 1}
\tcT {13} (e_1, e_2) = \cT {13} (e_1, e_2) = 0.
\end{equation}
\item For every prime $\ell \neq 2, 3, 13$ and all $v \geq 1$, we have
	\begin{equation} 	\begin{aligned}
	\cT {13}(\ell^{v}, 1) =& \cT {13}(\ell, 1) = 1 + \left(\frac{\ell}{3}\right), \\
	\cT {13}(1, \ell^{v}) =& \cT {13}(1, \ell) = 1 + \left(\frac{-1}{\ell}\right).
	\end{aligned}\end{equation} 
\item For $e \in \set{2, 2^2, 3, 3^2}$, the nonzero values of $\cT {13} (e, 1)$ and $\cT {13} (1, e)$ are given in Table \ref{table:T13(1 2)} and Table \ref{table:T13(3 1)} below. We have
	\begin{equation}
	\cT {13}(2^v, 1) = 0 \ \text{for} \ v \geq 1, \ \cT {13}(3^v, 1) = 0 \ \text{for} \ v \geq 3,
	\end{equation}
	and
	\begin{equation}
	\cT {13}(1, 2^v) = 0 \ \text{for} \ v \geq 3, \ \cT {13}(1, 3^v) = 0 \ \text{for} \ v \geq 1.
	\end{equation}
	For $v \geq 1$, the $13$-adic values are
	\begin{equation}\label{Equation: cT13 at 13}
	\cT {13}(13^v,1)=
	\begin{cases}
	1+12\cdot 13^{3v-2}, & 1\leq v\leq 3,\\
	1+13^8, & v\geq 4,
	\end{cases}
	\qquad
	\cT {13}(1,13^v)=
	\begin{cases}
	1+13^{2v-2}, & 1\leq v\leq 3,\\
	1+13^5, & v\geq 4.
	\end{cases}
	\end{equation}
	\item We have $\cT {13}(e_1, e_2) = O(2^{\omega(e_1 e_2)})$, where $\omega(e)$ is the number of distinct prime divisors of $e$.
	\item If $(a, b)$ is a $13$-groomed pair and $e \mid \tmd(\A {13}(a, b), \B {13} (a, b))$, then there is a unique factorization $e = e_1 e_2$ with $\gcd(e_1, e_2) = 1$ and $(a, b) \in \tcalT {13}(e_1, e_2)$.
	\end{enumalph}
\end{lemma}

\begin{proof}
	For parts (a) and (b), by the CRT (Sun Zi theorem), it suffices to consider $e = \ell^v$ a power of a prime. For $\ell \neq 2, 3, 13$, both claims follow from \Cref{Lemma: 2 and 3 and 13 are the ungroomed primes for m = 13}(a)--(b). For the exceptional primes, the stated conditions follow directly from \Cref{Definition T(e) for m = 13} and \eqref{Equation: twist defect powers}, apart from the two assertions that no pairs are permitted when $2 \mid e_1$ or $3 \mid e_2$. The first follows by checking the defining conditions modulo $8$, and the second by checking them modulo $27$. Part (b) then follows from the divisibility conditions in part (a).
	
	We now consider part (c). The assertion \eqref{Equation: T13 gcd conditions} implies \eqref{Equation: multiplicativity of T13} is simply multiplicativity in each argument away from the primes dividing the other argument. This follows from the CRT (Sun Zi theorem). For \eqref{Equation: Relating tcT 13 to cT 13}, let $\ell$ be a prime, and let $e = \ell^v$ for some $v \geq 1$. Consider the injective map
\begin{equation}
\begin{aligned}
\calT {13}(\ell^v, 1) \times (\bbZ/ \ell^{3 v})^\times &\to \tcalT {13}(\ell^v, 1) \\
(t,u) &\mapsto (tu,u)
\end{aligned}
\end{equation}
We observe $\A {13}(1, 0) = -3$ and $\B {13}(1, 0) = 2$ are coprime, so no pair $(a, b)$ with $b \equiv 0 \pmod \ell$ can be a member of $\tcalT {13}(\ell^v, 1)$. Surjectivity of the given map follows, and counting both sides gives the result in this component. On the other hand, we can consider the injective map
\begin{equation}
\begin{aligned}
\calT {13}(1, \ell^v) \times (\bbZ/ \ell^{2 v})^\times &\to \tcalT {13}(1, \ell^v) \\
(t,u) &\mapsto (tu,u)
\end{aligned}
\end{equation}
Again, as $\A {13}(1, 0) = -3$ and $\B {13}(1, 0) = 2$ are coprime, no pair $(a, b)$ with $b \equiv 0 \pmod \ell$ can be a member of $\tcalT {13}(1, \ell^v)$, and the desired implication follows. Finally, \eqref{Equation: T13 vanishes when gcd is not 1} whenever $\gcd(e_1, e_2) > 1$ holds by \Cref{Lemma: 2 and 3 and 13 are the ungroomed primes for m = 13} when $\gcd(e_1, e_2)$ is not a power of 13. The case $\ell = 13$ follows from the last condition of \eqref{Equation: Set mapping to tcalT 13} together with the observation that for coprime $(a, b)$ we have $\gcd(\CII {13} (a, b), \CIII {13} (a, b)) \mid 13$.

Now part (d). For $\ell \neq 2, 3, 13$, \Cref{Lemma: 2 and 3 and 13 are the ungroomed primes for m = 13}(a)--(b) yield
	\begin{equation} 	\begin{aligned}
	\calT {13}(\ell^v, 1) =& \set{t \in \bbZ/\ell^{3v} \bbZ : \hII {13}(t) \equiv 0 \psmod{\ell^{3v}}}, \ \text{and} \\
	\calT {13}(1, \ell^v) =& \set{t \in \bbZ/\ell^{2v} \bbZ : \hIII {13}(t) \equiv 0 \psmod{\ell^{2v}}}.
	\end{aligned}\end{equation} 
By \Cref{Lemma: 2 and 3 and 13 are the ungroomed primes for m = 13}(c), $\hII {13}(t) \equiv 0 \pmod \ell$ implies $\frac{\textup{d}}{\textup{d} t} \hII {13} (t) \not\equiv 0 \pmod \ell$, and likewise for $\hIII {13} (t)$ so Hensel's lemma applies and we need only count roots of the corresponding polynomial modulo $\ell$, and our result follows by quadratic reciprocity.

Next, part (e). For $\ell = 2$, we readily verify $\cT {13} (2, 1) = 0$, and hence $\cT {13} (2^v, 1) = 0$ for $v \geq 1$. On the other hand, $\cT {13} (1, 2) = 2$, $\cT {13} (1, 2^2) = 2^3$, and $\cT {13} (1, 2^3) = 0$, so $\cT {13} (1, 2^v) = 0$ for $v \geq 3$. 

For $\ell = 3$, we just compute 
$\cT {13}(3, 1) = 18$, $\cT {13}(3^2, 1) = 27$, and $\cT {13}(3^3, 1) = 0$; the observation $\cT {13}(3^3, 1) = 0$ implies $\cT {13}(3^v, 1) = 0$ for all $v \geq 3$. Similarly, $\cT {13} (1, 3) = 0$ implies $\cT {13} (1, 3^v) = 0$ for all $v \geq 1$.

It remains to treat $\ell=13$.  Modulo $13$ we have
\begin{equation}
\begin{aligned}
\f {13}(t)&\equiv -3(t-2)(t-3)(t+3)^6,\\
\g {13}(t)&\equiv 2(t-2)(t-3)^2(t+3)^9.
\end{aligned}
\end{equation}
Thus only the residue classes $t\equiv2,3,10\pmod {13}$ can contribute.  At $t\equiv2$, the factor $\hII {13}$ has a simple root while $\hIII {13}$, $\tf {13}$, and $\tg {13}$ are units, so Hensel's lemma gives exactly one class in $\calT {13}(13^v,1)$ for every $v$.  Likewise, the simple root $t\equiv3$ of $\hIII {13}$ gives exactly one class in $\calT {13}(1,13^v)$ for every $v$.

For the common residue class, write $t=10+13u$.  Then
\begin{equation}
\begin{aligned}
\hII {13}(t)&=13H_{\mathrm{II}}(u),&
H_{\mathrm{II}}(u)&=13u^2+21u+9,\\
\hIII {13}(t)&=13H_{\mathrm{III}}(u),&
H_{\mathrm{III}}(u)&=13u^2+20u+8,\\
\tf {13}(t)&=13^4F(u),&
F(u)&=-3(u^4-15u^3-31u^2-20u-4),\\
\tg {13}(t)&=13^6G(u),&
G(u)&=2(u^6+44u^5+83u^4+20u^3-59u^2-50u-12).
\end{aligned}
\end{equation}
Modulo $13$ these satisfy
\begin{equation}
H_{\mathrm{II}}(u)\equiv8(u-7),\quad
H_{\mathrm{III}}(u)\equiv7(u-10),\quad
F(u)\equiv-3(u-7)^4,\quad
G(u)\equiv2(u-10)^6.
\end{equation}
Consequently, for $\calT {13}(13^v,1)$ the defining convention forces $H_{\mathrm{III}}(u)$ to be a unit.  If $v\leq3$, the displayed powers $13^6$ and $13^9$ already supply the required divisibility, giving $(13-1)13^{3v-2}$ classes on this branch.  If $v\geq4$, then $G(u)$ is a unit and the condition on $\g {13}$ is exactly
$13^{3v-9}\mid H_{\mathrm{II}}(u)$.  The latter polynomial has a simple root modulo $13$, so there are exactly $13^8$ such classes modulo $13^{3v}$.  Including the branch through $t\equiv2$ gives the first formula of \eqref{Equation: cT13 at 13}.

For $\calT {13}(1,13^v)$, the common branch must instead satisfy $13\mid H_{\mathrm{III}}(u)$, while $H_{\mathrm{II}}(u)$ and $F(u)$ are units.  A further substitution $u=10+13w$ gives
\begin{equation}
13^{-1}G(10+13w)\equiv-5\pmod {13},
\end{equation}
so $v_{13}(G(u))=1$ throughout this branch.  If $r=v_{13}(H_{\mathrm{III}}(u))\geq1$, the two divisibility conditions are therefore
\begin{equation}
6+r\geq2v,\qquad 10+2r\geq3v.
\end{equation}
For $v\leq3$ these impose only $r\geq1$, giving $13^{2v-2}$ classes modulo $13^{2v}$; for $v\geq4$ they reduce to $r\geq2v-6$, giving $13^5$ classes.  Including the branch through $t\equiv3$ gives the second formula of \eqref{Equation: cT13 at 13}.

Part (f) follows from parts (d) and (e).

Finally, for part (g), part (c) assures us that we can take $(e_1, e_2) = (\ell^v, 1)$ or $(e_1, e_2) = (1, \ell^v)$ without loss of generality. If $\ell \neq 2, 3, 13$, the claim now follows from part (d), and if $\ell = 2, 3, 13$, the claim follows from part (e) and by construction.
\end{proof}

\jvtable{table:T13(1 2)}{
\rowcolors{2}{white}{gray!10}
\begin{tabular}{c | c c}
$m$ & $\cT {13} (1, 2^1)$ & $\cT {13} (1, 2^2)$ \\
\hline\hline
$13$ & 2 & $2^3$ \\
\end{tabular}
}{All nonzero $\cT {13} (1, 2^v)$}

\jvtable{table:T13(3 1)}{
\rowcolors{2}{white}{gray!10}
\begin{tabular}{c | c c}
$m$ & $\cT {13} (3^1, 1)$ & $\cT {13} (3^2, 1)$ \\
\hline\hline
$13$ & $2 \cdot 3^2$ & $3^3$ \\
\end{tabular}
}{All nonzero $\cT {13} (3^v, 1)$}

The following theorem gives us the tools to relate the twist height to the twist minimality defect for $m = 13$, in imperfect analogy with \Cref{Theorem: Controlling size of twist minimality defect for m = 7} and \Cref{Theorem: Controlling size of twist minimality defect for m = 5 and 10 and 25}.

\begin{theorem}\label{Theorem: Controlling size of twist minimality defect for m = 13}
	The following statements hold.
	\begin{enumalph}
	\item For all $(a, b) \in \bbR^2$, we have
	\begin{equation}\label{Equation: upper and lower bounds for H for m = 13}
	\begin{aligned}
	108 \CII {13} (a, b)^{12} &\leq \rawheight(\A {13} (a, b), \B {13} (a, b)) \leq \upperratio {II, 13} \CII {13} (a, b)^{12}, \\
	108 \CIII {13} (a, b)^{12} &\leq \rawheight(\A {13} (a, b), \B {13} (a, b)) \leq \upperratio {III, 13} \CIII {13} (a, b)^{12}, \\
	\end{aligned}
	\end{equation}
	where the constants 
	\begin{equation} 	\begin{aligned}
	\upperratio {II, 13} =& 635\,811\,018.28475061\ldots \ \text{and} \\
	\upperratio {III, 13} =& 35\,492\,073\,075.17456568\ldots \\
	\end{aligned}\end{equation} 
	are algebraic numbers given by evaluating $\rawheight(\A {13}(a, b), \B {13}(a, b))$ at appropriate roots of \eqref{Equation: Roots defining upperratio for m = 13 II} and \eqref{Equation: Roots defining upperratio for m = 13 III} respectively.
	\item If $\CII {13} (a, b) = e_{\textup{II}}^3 n_{\textup{II}}$, with $n_{\textup{II}}$ cube-free, and $\CIII {13} (a, b) = e_{\textup{III}}^2 n_{\textup{III}}$ with $n_{\textup{III}}$ square-free, then $\twistdefect(\A {13} (a, b), \B {13} (a, b)) = e_{\textup{II}} e_{\textup{III}} e\prm$, where $e\prm \mid 2 \cdot 3 \cdot 13$. In addition, for all $(a, b) \in \bbR^2$, we have
	\begin{equation}
	\lambda_{13} \CIII {13} (a, b) \leq \CII {13} (a, b) \leq \mu_{13} \CIII {13} (a, b),
	\end{equation}
	where the constants
	\begin{equation}\label{Equation: lambda 13 and mu 13} 	\begin{aligned}
	\lambda_{13} &= 0.92430609\ldots \ \text{and} \\
	\mu_{13} &= 1.82569390\ldots 
	\end{aligned}\end{equation} 
	are algebraic numbers given by evaluating $\CII {13} (a, b)$ at appropriate roots of \eqref{Equation: Roots defining C bounds for m = 13}.
	\end{enumalph}
\end{theorem}

\begin{proof}
        The proof of this theorem is similar to those of \Cref{Theorem: Controlling size of twist minimality defect for m = 7} and \Cref{Theorem: Controlling size of twist minimality defect for m = 5 and 10 and 25}; however, to aid our reader in parsing the contributions of both factors, especially in part (b), we prove it in its entirety.

	We first prove (a). Let $m = 13$. We wish to find the extrema of 
	\begin{equation}
	\rawheight(\A {13}(a, b), \B {13}(a, b))/\CII {13}(a, b)^{12} \ \text{and} \ \rawheight(\A {13}(a, b), \B {13}(a, b))/\CIII {13}(a, b)^{12}.
	\end{equation}
	As these expressions are homogeneous of degree 0, and $\CII {13}(a, b)$ and $\CIII{13} (a, b)$ are positive definite, we may assume without loss of generality that $\CII {13}(a, b) = 1$ or $\CIII {13} (a, b) = 1$ respectively. Using the theory of Lagrange multipliers, and examining the critical points of $\rawheight(\A {13}(a, b), \B {13}(a, b))$ subject to these respective constraints, we verify that \eqref{Equation: upper and lower bounds for H for m = 13} holds. Moreover, the lower bound is attained in both cases at $(1, 0)$, and the upper bound is attained when $a$ and $b$ are appropriately chosen roots of
	\begin{equation}\label{Equation: Roots defining upperratio for m = 13 II}
	\begin{aligned}
		 &105718701441600 a^{20} + 628890736780800 a^{18} + 6862077189805968 a^{16} \\
		 &- 3737927951730336 a^{14} - 7359872595882599 a^{12} - 1358785779700076 a^{10} \\
		 &+ 7533990802873860 a^8 - 2897948832460864 a^6 + 1787484431772288 a^4 \\
		 &- 2069428838131712 a^2 + 643089024640000 \\
		 =& 2^{6} \cdot 3^{4} \cdot 5^{2} \cdot 13^{8} \cdot a^{20} + 2^{9} \cdot 3^{3} \cdot 5^{2} \cdot 13^{7} \cdot 29 \cdot a^{18} \\
	   &+ 2^{4} \cdot 3^{2} \cdot 13^{6} \cdot 2971 \cdot 3323 \cdot a^{16} - 2^{5} \cdot 3 \cdot 7^{2} \cdot 13^{5} \cdot 2140163 \cdot a^{14} \\
	   &- 13^{4} \cdot 103 \cdot 2309 \cdot 1083517 \cdot a^{12} - 2^{2} \cdot 7 \cdot 13^{4} \cdot 199 \cdot 757 \cdot 11279 \cdot a^{10} \\
    &+ 2^{2} \cdot 3 \cdot 5 \cdot 13^{2} \cdot 742997120599 \cdot a^{8} - 2^{6} \cdot 7 \cdot 13 \cdot 6991 \cdot 71175421 \cdot a^{6}\\
		 &+ 2^{7} \cdot 3^{2} \cdot 29311 \cdot 52936979 \cdot a^{4}  - 2^{11} \cdot 61 \cdot 16564972129 \cdot a^{2} + 2^{10} \cdot 5^{4} \cdot 31699^{2}, \ \text{and} \\
		 &105718701441600 b^{20} - 129031030477440 b^{18} + 226457312671512 b^{16} \\
		 &- 162954491664432 b^{14} + 61397224373329 b^{12} - 13194397029476 b^{10} \\
		 &+ 1681210465311 b^8 - 121030573768 b^6 + 4530949623 b^4 \\
		 &- 78302708 b^2 + 28561 \\
		 =& 2^{6} \cdot 3^{4} \cdot 5^{2} \cdot 13^{8} \cdot b^{20} - 2^{7} \cdot 3^{3} \cdot 5 \cdot 7 \cdot 13^{7} \cdot 17 \cdot b^{18} + 2^{3} \cdot 3^{2} \cdot 13^{6} \cdot 613 \cdot 1063 \cdot b^{16} \\
		 &- 2^{4} \cdot 3 \cdot 13^{5} \cdot 103 \cdot 88771 \cdot b^{14} + 13^{4} \cdot 157 \cdot 13692277 \cdot b^{12} - 2^{2} \cdot 13^{4} \cdot 115493129 \cdot b^{10} \\
		 &+ 3 \cdot 13^{2} \cdot 3315996973 \cdot b^{8} - 2^{3} \cdot 13 \cdot 1091 \cdot 1066687 \cdot b^{6} + 3^{3} \cdot 577 \cdot 290837 \cdot b^{4} \\
		 &- 2^{2} \cdot 11 \cdot 1779607 \cdot b^{2} + 13^{4},
	\end{aligned}
	\end{equation}
	if $\CII {13} (a, b) = 1$, and of
	\begin{equation}\label{Equation: Roots defining upperratio for m = 13 III}
	\begin{aligned}
		 &469860895296 a^{16} - 4490785864656 a^{14} + 18528290390389 a^{12} \\
		 &- 42537089721750 a^{10} + 58527314729975 a^8 - 48232472033876 a^6 \\
		 &+ 22080850389507 a^4 - 4364808534790 a^2 + 18869692689 \\
		 =& 2^{6} \cdot 3^{2} \cdot 13^{8} \cdot a^{16} - 2^{4} \cdot 3^{2} \cdot 7 \cdot 13^{7} \cdot 71 \cdot a^{14} + 13^{6} \cdot 3838621 \cdot a^{12} \\
		 &- 2 \cdot 3 \cdot 5^{3} \cdot 13^{5} \cdot 152753 \cdot a^{10} + 5^{2} \cdot 7 \cdot 13^{5} \cdot 23 \cdot 39163 \cdot a^{8} \\
		 &- 2^{2} \cdot 13^{4} \cdot 157 \cdot 1249 \cdot 2153 \cdot a^{6} + 3 \cdot 7 \cdot 13 \cdot 13907 \cdot 5815937 \cdot a^{4} \\
   &- 2 \cdot 5 \cdot 11 \cdot 39680077589 \cdot a^{2} + 3^{4} \cdot 15263^{2}, \ \text{and} \\
		 &30071097298944 b^{16} + 11710379236608 b^{14} + 995577624340 b^{12} \\
		 &- 638359599384 b^{10} - 130255153795 b^8 - 12900146870 b^6 \\
		 &- 775262241 b^4 - 64120726 b^2 + 29241 \\
		 =& 2^{12} \cdot 3^{2} \cdot 13^{8} \cdot b^{16} + 2^{8} \cdot 3^{6} \cdot 13^{7} \cdot b^{14} + 2^{2} \cdot 5 \cdot 13^{6} \cdot 10313 \cdot b^{12} \\
		 &- 2^{3} \cdot 3^{2} \cdot 13^{5} \cdot 23879 \cdot b^{10} - 5 \cdot 13^{5} \cdot 70163 \cdot b^{8} - 2 \cdot 5 \cdot 13^{4} \cdot 31^{2} \cdot 47 \cdot b^{6} \\
		 &- 3^{2} \cdot 13 \cdot 6626173 \cdot b^{4} - 2 \cdot 557 \cdot 57559 \cdot b^{2} + 3^{4} \cdot 19^{2},
	\end{aligned}
	\end{equation}
	if $\CIII {13} (a, b) = 1$. In both cases, $27 \abs{\B {13}(a, b)}^2 > 4 \abs{\A {13}(a, b)}^3$. For the reader's information, 
	\begin{equation}
	(a, b) = (-0.715678818\ldots, 0.320005592\ldots)
	\end{equation}
	maximizes $\rawheight$ subject to $\CII {13} (a, b) = 1$, and 
	\begin{equation}
	(a, b) = (-0.066491149\ldots, 0.498893507\ldots)
	\end{equation}
	maximizes $\rawheight$ subject to $\CIII {13} (a, b) = 1$.
		
	We now prove (b). Write $\CII {13} (a, b) = e_{\textup{II}}^3 n_{\textup{II}}$, with $n_{\textup{II}}$ cube-free, and $\CIII {13} (a, b) = e_{\textup{III}}^2 n_{\textup{III}}$ with $n_{\textup{III}}$ square-free. By \Cref{Lemma: 2 and 3 and 13 are the ungroomed primes for m = 13}(a), $e\prm = 2^u \cdot 3^v \cdot 13^w$ for some $u, v, w \geq 0$. A short computation shows $u, v, w \leq 1$. 
	
	The remainder of the proof of (b) is similar to the proof of (a), but even easier: we wish to find the extrema of $\CII {13} (a, b)$ subject to the constraint $\CIII {13} (a, b) = 1$. We find that these extrema are attained when $a$ and $b$ are appropriately chosen roots of		
		\begin{equation}\label{Equation: Roots defining C bounds for m = 13}
	\begin{aligned}
		 &13 a^4 - 13 a^2 + 1 \ \text{and} \\
		 &208 b^4 - 52 b^2 + 1 = 2^{4} \cdot 13 \cdot b^{4} - 2^{2} \cdot 13 \cdot b^2 + 1.
	\end{aligned}
	\end{equation}
	We maximize the ratio when $(a, b) = (0.289784148\ldots, 0.478546013\ldots)$, and minimize the ratio when $(a, b) = (0.957092026\ldots, -0.144892074\ldots)$.
\end{proof}

\begin{remark}
	Because $\A {13} (a, b)$ and $\B {13} (a, b)$ have both $\CII {13} (a, b)$ and $\CIII {13} (a, b)$ as common factors, \Cref{Theorem: Controlling size of twist minimality defect for m = 7}(b) and \Cref{Theorem: Controlling size of twist minimality defect for m = 5 and 10 and 25}(b) have no perfect analogues. However, \Cref{Theorem: Controlling size of twist minimality defect for m = 13} enables us to bound $\rawheight(\A {13} (a, b), \B {13} (a, b))$ with respect to $\CII {13} (a, b)^k \cdot \CIII {13} (a, b)^{12 - k}$ for any $k \in \bbR$, and therefore enables us to derive a whole family of analogous bounds.
	
	For example, if $\CII {13} (a, b) = e_{\textup{II}}^3 n_{\textup{II}}$, with $n_{\textup{II}}$ cube-free, and $\CIII {13} (a, b) = e_{\textup{III}}^2 n_{\textup{III}}$ with $n_{\textup{III}}$ square-free, then
	\begin{equation}
	\frac{\lambda_{13}^3}{2^4 \cdot 3^3 \cdot 13^6} e_{II}^{21} n_{II}^9 n_{III}^3 \leq \twistheight (\A {13} (a, b), \B {13} (a, b)) \leq \mu_{13}^3 \upperratio {II, 13} e_{II}^{21} n_{II}^9 n_{III}^3.
	\end{equation}
    The constants $\lambda_{13}$ and $\mu_{13}$ are given in \eqref{Equation: lambda 13 and mu 13}.
\end{remark}

\section{Estimates for twist classes for \texorpdfstring{$m = 13$}{m = 13}} \label{Section: Estimating twN(X) for m = 13}

In this section, we use \cref{Section: Our approach revisited} to estimate $\twistNad {13}(X)$, counting the number of twist minimal elliptic curves over $\bbQ$ admitting a cyclic $13$-isogeny.

Recall \eqref{Equation: defining cM}, \eqref{Equation: Applying N calE to Neq m}, and \eqref{Equation: Applying N calE to Nad m}. By \cref{Section: Parameterizing elliptic curves with a cyclic m-isogeny}, $\cM {13}(X; e)$ counts pairs $(a, b) \in \bbZ^2$ with
\begin{itemize}
	\item $(a, b)$ $13$-groomed,
	\item $\rawheight(\A {13}(a,b),\B {13}(a,b)) \leq X$ and
	\item $e \mid \twistdefect(\A {13}(a, b), \B {13}(a, b))$.
\end{itemize}

To avoid technical inconvenience, and in contrast to \Cref{Proposition: fundamental sieve for for m = 7} and \Cref{Proposition: fundamental sieve for for m = 5 and 10 and 25}, we opt not to refine \Cref{Lemma: fundamental sieve}. Instead, we proceed directly to an analogue of \Cref{Lemma: asymptotic for M(X; e) for m = 7} and \Cref{Lemma: asymptotic for M(X; e) for m = 5 and 10 and 25}.

\begin{lemma}\label{Lemma: asymptotic for M(X; e) for m = 13}
The following statements hold.
\begin{enumalph}
\item	If $\gcd(d, e) > 1$, then $\cM {13}(X; d, e) = 0$. If $\gcd(d, e) = 1$, we have 
	\begin{equation}
	\cM {13}(X; d, e) = \frac{\R {13} X^{1/12}}{d^2} \prod_{\ell \mid e} \parent{1 - \frac{1}{\ell}} \sum_{\substack{e_1 e_2 = e \\ \gcd(e_1, e_2) = 1}} \frac{\cT {13}(e_1, e_2)}{e_1^3 e_2^2} + O\parent{\frac{2^{2 \omega(e)} X^{1/24}}{d}}
	\end{equation}
	for $X, d, e \geq 1$. Here, $\R {13}$ is the area of \eqref{eqn: R(X)} for $m = 13$.
\item We have
	\begin{equation} 
	\cM {13}(X;e)
	= \frac{\R {13} X^{1/12}}{\zeta(2) \prod_{\ell \mid e} \parent{1 + \frac{1}{\ell}}} \sum_{\substack{e_1 e_2 = e \\ \gcd(e_1, e_2) = 1}} \frac{\cT {13}(e_1, e_2)}{e_1^3 e_2^2} + O\parent{2^{2 \omega(e)} X^{1/24} \log X}
	\end{equation}
        for $X \geq 2$ and $e \geq 1$.
\end{enumalph}
In both cases, the implied constants are independent of $d$, $e$, and $X$.
\end{lemma}

\begin{proof}
	We begin with (a) and examine the summands $\cM {13}(X; d, e)$. If $d$ and $e$ are not coprime, then $\cM {13}(X; d, e) = 0$ because $\gcd(da, db, e) \geq \gcd(d, e) > 1$. On the other hand, if $\gcd(d, e) = 1$, we have a bijection from the pairs counted by $\cM {13}(X; 1, e)$ to the pairs counted by $\cM {13}(d^{24} X; d, e)$ given by $(a, b) \mapsto (d a, d b)$.

For $X\geq 1$ and $e_1, e_2, a_0, b_0 \in \bbZ$, we write
\begin{equation}
L_{13}(X; e_1, e_2, a_0, b_0) \colonequals \#\set{(a, b) \in \calR {13}(X) \cap \bbZ^2 : \begin{minipage}{28ex} $(a, b) \equiv (a_0, b_0) \psmod {e_1^3 e_2^2},$ \\ $(a:b) \not\in \cusps{13}$ \end{minipage} }
\end{equation}
(this notation will not be used outside of this proof). 
By \Cref{Corollary: Estimates for lattice counts}, we have
\begin{equation}
L_{13}(X; e_1, e_2, a_0, b_0) = \frac{\R {13} X^{1/12}}{e_1^6 e_2^4} + O \parent{\frac{X^{1/24}}{e_1^3 e_2^2}}.
\end{equation}
Now by \Cref{Lemma: bound on T(e) for m = 13}(e)-(g), we have
\begin{equation}\label{Equation: cM (X 1 e) for m = 13}
\begin{aligned}
	\cM {13}(X; 1, e) =& \sum_{\substack{e_1 e_2 = e \\ \gcd(e_1, e_2) = 1}} \sum_{(a_0, b_0) \in \tcalT {13}(e_1, e_2)} L_{13}(X; e_1, e_2, a_0, b_0) \\
	=& \R {13} X^{1/12} \sum_{\substack{e_1 e_2 = e \\ \gcd(e_1, e_2) = 1}} \frac{\varphi(e_1^3 e_2^2) \cT {13}(e_1, e_2)}{e_1^6 e_2^4} \\
 &+ O\parent{X^{1/24} \sum_{\substack{e_1 e_2 = e \\ \gcd(e_1, e_2) = 1}} \frac{\varphi(e_1^3 e_2^2) \cT {13}(e_1, e_2)}{e_1^3 e_2^2} } \\
	=& \R {13} X^{1/12} \prod_{\ell \mid e} \parent{1 - \frac{1}{\ell}} \sum_{\substack{e_1 e_2 = e \\ \gcd(e_1, e_2) = 1}} \frac{\cT {13}(e_1, e_2)}{e_1^3 e_2^2} \\
    &+ O\parent{  X^{1/24} \sum_{\substack{e_1 e_2 = e \\ \gcd(e_1, e_2) = 1}} \cT {13}(e_1, e_2) }.
\end{aligned}
\end{equation}
Scaling by $d$ and invoking \Cref{Lemma: bound on T(e) for m = 13}(f), we obtain
\begin{equation}
	\cM {13}(X; d, e) = \frac{\R {13} X^{1/12}}{d^2} \prod_{\ell \mid e} \parent{1 - \frac{1}{\ell}} \sum_{\substack{e_1 e_2 = e \\ \gcd(e_1, e_2) = 1}} \frac{\cT {13}(e_1, e_2)}{e_1^3 e_2^2} + O\parent{\frac{2^{2 \omega(e)} X^{1/24}}{d}}.
\end{equation}

For part (b), we compute
\begin{equation} \label{equation: cmXe for m = 13}
\begin{aligned}
	\cM {13}(X; e) =& \sum_{\substack{d \ll X^{1/24} \\ \gcd(d, e) = 1}} \mu(d) \cM {13}(X; d, e) \\
	=& \sum_{\substack{d \ll X^{1/24} \\ \gcd(d, e) = 1}} \mu(d) \frac{\R {13} X^{1/12}}{d^2} \prod_{\ell \mid e} \parent{1 - \frac{1}{\ell}} \sum_{\substack{e_1 e_2 = e \\ \gcd(e_1, e_2) = 1}} \frac{\cT {13}(e_1, e_2)}{e_1^3 e_2^2} \\
    &+ \sum_{\substack{d \ll X^{1/24} \\ \gcd(d, e) = 1}} \mu(d) \cdot 
 O\parent{\frac{2^{2 \omega(e)} X^{1/24}}{d}} \\
	=& \R {13} X^{1/12} \prod_{\ell \mid e} \parent{1 - \frac{1}{\ell}} \sum_{\substack{e_1 e_2 = e \\ \gcd(e_1, e_2) = 1}} \frac{\cT {13}(e_1, e_2)}{e_1^3 e_2^2} \sum_{\substack{d \ll X^{1/24} \\ \gcd(d, e) = 1}} \frac{ \mu(d)}{d^2} \\
	&+ O\parent{2^{2 \omega(e)} X^{1/24} \sum_{\substack{d \ll X^{1/24} \\ \gcd(d, e) = 1}} \frac{1}{d}}.
	\end{aligned}
	\end{equation}
Plugging the straightforward estimates
\begin{equation}
\sum_{\substack{d \ll X^{1/24} \\ \gcd(d, e) = 1}} \frac{ \mu(d)}{d^2} = \frac{1}{\zeta(2)} \prod_{\ell \mid e} \parent{1 - \frac{1}{\ell^2}}\inv + O(X^{-1/24})
\end{equation}
and
\begin{equation} 
\sum_{\substack{d \leq X^{1/24}}} \frac 1d = \frac{1}{24}\log X+ O(1) \end{equation}
into \eqref{equation: cmXe for m = 13} then simplifies to give
\begin{equation}
\begin{aligned}
\cM {13}(X;e)	
	=& \frac{\R {13} X^{1/12}}{\zeta(2) \prod_{\ell \mid e} \parent{1 + \frac{1}{\ell}}} \sum_{\substack{e_1 e_2 = e \\ \gcd(e_1, e_2) = 1}} \frac{\cT {13}(e_1, e_2)}{e_1^3 e_2^2} + O\parent{2^{2 \omega(e)} X^{1/24} \log X}
\end{aligned}
\end{equation}
proving (b).
\end{proof}

\begin{remark}
	The alternate proofs for \Cref{Lemma: asymptotic for M(X; e) for m = 7} and \Cref{Lemma: asymptotic for M(X; e) for m = 5 and 10 and 25} do not carry over directly to \Cref{Lemma: asymptotic for M(X; e) for m = 13}, even though both $\CII {13} (a, b)$ and $\CIII {13} (a, b)$ are norms on orders of class number 1, precisely because we have two such factors (i.e., because the elliptic surface under consideration has places of both potential type II and potential type III additive reduction). However, we believe these arguments can be salvaged in part by applying the improved sieving to the larger of $e_1$ and $e_2$ for each factorization $e = e_1 e_2$ occurring in the outer sum of \eqref{Equation: cM (X 1 e) for m = 13}. This will not yield an error term of the same strength as \Cref{Lemma: asymptotic for M(X; e) for m = 7} and \Cref{Lemma: asymptotic for M(X; e) for m = 5 and 10 and 25}, because the other term $e_i$ will be $O(\sqrt{e})$, and must be approximated by summing over congruence classes in the manner indicated in the proof above. Nevertheless, such an argument ought to be able to improve on the error term of \Cref{Lemma: asymptotic for M(X; e) for m = 13}, and therefore of \Cref{Theorem: asymptotic for twN(X) for m = 13} and \Cref{Theorem: asymptotic for NQ(X) for m = 13} below.
\end{remark}

We let
	\begin{equation}
	\Q {13} \colonequals \sum_{n \geq 1} \frac{ \varphi_{1/2} (n)}{\prod_{\ell \mid n} \parent{1 + \frac{1}{\ell}}} \sum_{\substack{n_1 n_2 = n \\ \gcd(n_1, n_2) = 1}} \frac{\cT {13}(n_1, n_2)}{n_1^3 n_2^2},
	\end{equation}
	and we let
\begin{equation}\label{Equation: twconst for m = 13}
	\twconst {13} \colonequals \frac{\Q {13} \R {13}}{\zeta(2)}.
\end{equation}
Here, as always, $\R {13}$ is the area of the region
	\begin{equation}
	\calR {13}(1) = \set{(a, b) \in \bbR^2 : \rawheight(\A {13} (a, b), \B {13} (a, b)) \leq 1, b \geq 0}.
	\end{equation}

We are now in a position to estimate $\twistNadly {13}(X)$. Our argument is similar to those given in \Cref{Lemma: asymptotic for twN<=y(X) for m = 7} and \Cref{Lemma: asymptotic for twN<=y(X) for m = 10 and 25}, but complicated by the necessity of summing over factorizations for the twist minimality defect $e$.

\begin{lemma}\label{Lemma: asymptotic for twN<=y(X) for m = 13}
	Suppose $y \ll X^{1/24}$. Then
	\begin{equation}
	\twistNadly {13}(X) =  \twconst {13} X^{1/12} + O\parent{\max\parent{\frac{X^{1/12} \log^3 y}{y^{1/2}}, X^{1/24} y^{5/4} \log X \log^7 y }}
	\end{equation}
	for $X, y \geq 2$. The constant $\twconst {13}$ is given in \eqref{Equation: twconst for m = 13}.
\end{lemma}

\begin{proof}
	Substituting the asymptotic for $\cM {13}(X; e)$ from \Cref{Lemma: asymptotic for M(X; e) for m = 13}(b) into the defining series \eqref{Equation: defining twNadly X} for $\twistNadly {13}(X)$, we have
	\begin{equation}\label{Equation: main terms for twNly X for m = 13}\begin{aligned}
		\twistNadly {13}(X) =& \sum_{n \leq y} \sum_{e \mid n} \mu\parent{n/e} \frac{\R {13} e^{1/2} X^{1/12}}{\zeta(2) \prod_{\ell \mid n} \parent{1 + \frac{1}{\ell}}} \sum_{\substack{n_1 n_2 = n \\ \gcd(n_1, n_2) = 1}} \frac{\cT {13}(n_1, n_2)}{n_1^3 n_2^2} \\
        &+ \sum_{n \leq y} \sum_{e \mid n} \mu\parent{n/e} 
 O\parent{2^{2 \omega(n)} e^{1/4} X^{1/24} \log e^6 X}.
	\end{aligned}\end{equation} 
	
	We handle the main term and the error of this expression separately. For the main term, recalling the definition of the generalized totient function \eqref{Equation: generalized totient function}, we have
	\begin{equation}
	\begin{aligned}
	&\sum_{n \leq y} \sum_{e \mid n} \mu\parent{n/e} \frac{\R {13} e^{1/2} X^{1/12}}{\zeta(2) \prod_{\ell \mid n} \parent{1 + \frac{1}{\ell}}} \sum_{\substack{n_1 n_2 = n \\ \gcd(n_1, n_2) = 1}} \frac{\cT {13}(n_1, n_2)}{n_1^3 n_2^2} \\
	=& \frac{\R {13} X^{1/12}}{\zeta(2)} \sum_{n \leq y} \frac{ \varphi_{1/2} (n)}{\prod_{\ell \mid n} \parent{1 + \frac{1}{\ell}}} \sum_{\substack{n_1 n_2 = n \\ \gcd(n_1, n_2) = 1}} \frac{\cT {13}(n_1, n_2)}{n_1^3 n_2^2}.
	\end{aligned}
	\end{equation}
	By \Cref{Lemma: bound on T(e) for m = 13}(f), we see
	\begin{equation}
	\frac{ \varphi_{1/2} (n)}{\prod_{\ell \mid n} \parent{1 + \frac{1}{\ell}}} \sum_{\substack{n_1 n_2 = n \\ \gcd(n_1, n_2) = 1}} \frac{\cT {13}(n_1, n_2)}{n_1^3 n_2^2} = O\parent{\frac{4^{\omega(n)}}{n^{3/2}}}.
	\end{equation}
	
	\Cref{Corollary: tail of sum of f/n^sigma} and \Cref{Theorem: Asymptotics of k^omega(n)} together yield
	\begin{equation}
	\sum_{n > y} \frac{4^{\omega(n)}}{n^{3/2}} = O\parent{\frac{\log^3 y}{y^{1/2}}}.
	\end{equation}
	Thus, the series
	\begin{equation}
	\sum_{n \geq 1} \frac{ \varphi_{1/2} (n)}{\prod_{\ell \mid n} \parent{1 + \frac{1}{\ell}}} \sum_{\substack{n_1 n_2 = n \\ \gcd(n_1, n_2) = 1}} \frac{\cT {13}(n_1, n_2)}{n_1^3 n_2^2}= \Q {13} \label{equation: sum for Q for m = 13}
	\end{equation}
	is absolutely convergent, and 
	\begin{equation}
	\begin{aligned}
	\frac{\R {13} X^{1/12}}{\zeta(2)} \sum_{n \leq y} \frac{ \varphi_{1/2} (n)}{\prod_{\ell \mid n} \parent{1 + \frac{1}{\ell}}} \sum_{\substack{n_1 n_2 = n \\ \gcd(n_1, n_2) = 1}} \frac{\cT {13}(n_1, n_2)}{n_1^3 n_2^2} =& \frac{\R {13} X^{1/12}}{\zeta(2)} \parent{\Q {13} + O\parent{\frac{\log^3 y}{y^{1/2}}}} \\
	=& \twconst {13} X^{1/12} + O\parent{\frac{X^{1/12} \log^3 y}{y^{1/2}}}.
	\end{aligned}
	\end{equation}
	
	As the summands of \eqref{equation: sum for Q for m = 13} constitute a nonnegative multiplicative arithmetic function, we can factor $\Q {13}$ as an Euler product. We have
     \begin{equation}
	\Q {13} = \prod_{\textup{$p$ prime}} \Q {13} (p), \label{equation: product for Q for m = 13};
	\end{equation}
	by \Cref{Lemma: bound on T(e) for m = 13}, the terms $\Q {13}(p)$ are computed as follows:
	\begin{equation}\label{Equation: Euler factors for m = 13}
	\begin{aligned}
	\Q {13} (p) &\colonequals 1 + \sum_{a \geq 1} \frac{\varphi_{1/2}(p^a)}{1 + 1/p} \parent{\frac{\cT {13} (p^a, 1)}{p^{3a}} + \frac{\cT {13} (1, p^a)}{p^{2a}}} \\
	=& \begin{cases}
	1 + \displaystyle{\frac{2 p^{1/2} \left(p^2+p^{3/2}+2 p+2 p^{1/2}+2\right)}{(p+1) \left(p+p^{1/2}+1\right) \left(p^2+p^{3/2}+p+p^{1/2}+1\right)}}, & \begin{tabular}{@{}c@{}} 
    \textup{if $p \equiv 1 \psmod{12}$} \\ \textup{and $p \neq 13$;} \\ %
    \end{tabular} \\
	1 + \displaystyle{\frac{2 p^{1/2}}{(p+1) \left(p+p^{1/2}+1\right)}}, & \textup{if $p \equiv 5 \psmod{12}$;} \\ %
	1 + \displaystyle{\frac{2 p^{1/2}}{(p+1) \left(p^2+p^{3/2}+p+p^{1/2}+1\right)}}, & \textup{if $p \equiv -5 \psmod{12}$;} \\ %
	\frac{4}{3}, & \text{if $p=2$;} \\
	\frac{21 + 17 \sqrt{3}}{36}, & \text{if $p=3$;} \\
	\frac{79\,239\,719+73\,639\,763\sqrt{13}}{79\,270\,842}, & \text{if $p=13$;} \\
	1 & \text{if} \ p \equiv -1 \pmod {12}.
	\end{cases}
	\end{aligned}
	\end{equation}
	For $p=13$, substituting \eqref{Equation: cT13 at 13} into the defining series and summing the geometric tails gives
	\begin{equation}
	\Q {13}(13)=\frac{79\,239\,719+73\,639\,763\sqrt{13}}{79\,270\,842}=4.349\,034\,925\,218\,983\,880\ldots.
	\end{equation}
	The square roots in \eqref{Equation: Euler factors for m = 13} appear because of the generalized Jordan totient function $\varphi_{1/2}$. For instance, for $p = 3$ we have
    \begin{equation} 	\begin{aligned}
    \Q {13} (3) =& 1 + \frac{\varphi_{1/2}(3)}{1 + 1/3} \parent{\frac{\cT {13} (3, 1)}{3^{3}} + \frac{\cT {13} (1, 3)}{3^2}} + \frac{\varphi_{1/2}(3^2)}{1 + 1/3} \parent{\frac{\cT {13} (3^2, 1)}{3^{6}} + \frac{\cT {13} (1, 3^2)}{3^4}} \\
    =& 1 + \frac{\sqrt{3} - 1}{1 + 1/3} \cdot \frac{18}{3^{3}} + \frac{3 - \sqrt{3}}{1 + 1/3} \cdot \frac{27}{3^{6}} \\
    =& \frac{21 + 17 \sqrt{3}}{36}.
    \end{aligned}\end{equation} 
	
	We now turn to the error term. Since $y \ll X^{1/24}$, for $e \leq y$ we have $\log (e^6 X) \ll \log X$. The error term of \eqref{Equation: main terms for twNly X for m = 13} is therefore
	\begin{equation} \label{equation: partial simplification twN<=y(X) for m = 13}
	\begin{aligned}
  O\parent{X^{1/24} \log X \sum_{n \leq y} 4^{\omega(n)} \sum_{e \mid n} \abs{\mu\parent{n/e}} e^{1/4}}
	\end{aligned}
	\end{equation}
	As in the proof of \Cref{Lemma: asymptotic for twN<=y(X) for m = 7}, we note
	\begin{equation}
	\sum_{e \mid n} \abs{\mu\parent{n/e}} e^{1/4} \leq n^{1/4} \prod_{p \mid n} \parent{1 + p^{-1/4}} \leq 2^{\omega(n)} n^{1/4};
	\end{equation}
	\Cref{Theorem: Asymptotics of k^omega(n)} tells us $\sum_{n \leq y} 8^{\omega(n)} = O(y \log^7 y)$, so by \Cref{Corollary: tail of sum of f/n^sigma}, we have
	\begin{equation}
	\sum_{n \leq y} 4^{\omega(n)} \sum_{e \mid n} \abs{\mu\parent{n/e}} e^{1/4} = O\parent{y^{5/4} \log^7 y},
	\end{equation}
	and our desired result follows.
\end{proof}

We emphasize that the proof of \Cref{Lemma: asymptotic for twN<=y(X) for m = 13} has given $\Q {13}$ an Euler product expansion 
\begin{equation}\label{Equation: product for Q for m = 13}
    \Q {13} = \prod_p \Q {13} (p), 
\end{equation}
where $\Q {13} (p)$ is given by \eqref{Equation: Euler factors for m = 13}.

\begin{lemma}\label{Lemma: bound on twN>y(X) for m = 13}
For every $\epsilon>0$, we have
\begin{equation}
\twistNadgy {13}(X)=O_\epsilon\parent{\frac{X^{1/12+\epsilon}}{y^{1/6}}+X^{1/14+\epsilon}}
\end{equation}
for $X,y\geq2$.
\end{lemma}

\begin{proof}
We have $\CII {13}=\C 7$ and $\CIII {13}=\C {25}$. The norm-form lattice arguments of \Cref{Lemma: asymptotic for M(X; e) for m = 7,Lemma: asymptotic for M(X; e) for m = 5 and 10 and 25}, together with the Principle of Lipschitz, give
\begin{equation}
\begin{aligned}
\#\{(a,b): e^3\mid\CII {13}(a,b)\text{ for some }e>E\}
&=O_\epsilon(B^{2+\epsilon}E^{-2}+B^{1+\epsilon}E^{-1/2}+B^{2/3+\epsilon}),\\
\#\{(a,b): e^2\mid\CIII {13}(a,b)\text{ for some }e>E\}
&=O_\epsilon(B^{2+\epsilon}E^{-1}+B^{1+\epsilon}),
\end{aligned}
\end{equation}
for primitive pairs with $\abs a,\abs b\leq B$; here the exceptional $3$- and $2$-parts are bounded, since $v_3(\CII {13})\leq3$ and $v_2(\CIII {13})\leq3$ for primitive pairs.

Let $D=e_{\mathrm{II}}e_{\mathrm{III}}e'$ be the actual twist minimality defect, with $e'\mid78$, as in \Cref{Theorem: Controlling size of twist minimality defect for m = 13}(b). By \Cref{Proposition: bounding the summands of twN calE X}, the large-defect sum is bounded, up to $\tau(D)^2$, by pairs of twist height at most $X$ with $D>y$. The same theorem gives $\abs a,\abs b\ll X^{1/14}$ and $D\ll X^{5/42}$, so $\tau(D)^2\ll_\epsilon X^\epsilon$.

On a dyadic interval $D_0\leq D<2D_0$, we have $\abs a,\abs b\ll X^{1/24}D_0^{1/4}$ and either $e_{\mathrm{II}}\gg D_0^{1/3}$ or $e_{\mathrm{III}}\gg D_0^{2/3}$. The preceding estimates therefore give
\begin{equation}
O_\epsilon\parent{X^{1/12+\epsilon}D_0^{-1/6}+X^{1/24+\epsilon}D_0^{1/4}}.
\end{equation}
Summing dyadically for $D_0>y$ and using $D_0\ll X^{5/42}$ proves the lemma.
\end{proof}

We are now ready to prove \Cref{Intro Theorem: asymptotic for twN(X) for m of genus $0$} for $m = 13$, which we restate here in the notations we have established.

\begin{theorem}\label{Theorem: asymptotic for twN(X) for m = 13}
	For every $\epsilon>0$, we have
	\begin{equation}
	\twistNeq {13}(X) = \twistNad {13}(X) = \twconst {13} X^{1/12} + O_\epsilon\parent{X^{4/51+\epsilon}}
	\end{equation}
	for $X \geq 2$. The constant $\twconst {13}$ is given in \eqref{Equation: twconst for m = 13}.
\end{theorem}

\begin{proof}
Combining \Cref{Lemma: asymptotic for twN<=y(X) for m = 13,Lemma: bound on twN>y(X) for m = 13} and taking $y=X^{1/34}$, we have
\begin{equation}
X^{1/24}y^{5/4}=X^{1/12}y^{-1/6}=X^{4/51}.
\end{equation}
The remaining error terms are smaller, and the logarithmic factors are absorbed into $X^\epsilon$.
\end{proof}

\subsection*{$L$-series}

As in the previous two chapters, we set up the next section by interpreting
\Cref{Theorem: asymptotic for twN(X) for m = 13} in terms of Dirichlet series. 

Recall \eqref{Equation: twisth calE}, \eqref{Equation: hQ calE}, \eqref{Equation: defining twistL calE X}, and \eqref{Equation: defining LQ calE X}.

\begin{cor}\label{Corollary: twL(s) has a meromorphic continuation for m = 13}
The following statements hold.
    \begin{enumerate}
    \item The Dirichlet series $\twistLad {13}(s)$ has abscissa of (absolute) convergence $\sigma_a=\sigma_c = 1/12$ and has a meromorphic continuation to the region
	\begin{equation}
	\set{s = \sigma + i t \in \bbC : \sigma > 4/51}. \label{equation: domain of twistL(s) for m = 13}
	\end{equation}
	\item The function $\twistLad {13}(s)$ has a simple pole at $s = 1/12$ with residue 
	\begin{equation} \res_{s=\frac{1}{12}} \twistLad {13}(s) = \frac{\twconst {13}}{12}; \end{equation}
	it is holomorphic elsewhere on the region \eqref{equation: domain of twistL(s) for m = 13}. \item We have
\begin{equation}
    \mu_{\twistLad {13}}(\sigma) \leq \max\{0, 17-204\sigma\}
     \end{equation}
     for $\sigma > 4/51$.
     \end{enumerate}
\end{cor}

\begin{proof}
Abel summation applied to \Cref{Theorem: asymptotic for twN(X) for m = 13} gives (a)--(b). On fixed lines $4/51<\sigma_0<1/12$, the remainder integral is $O(1+\abs t)$, while absolute convergence gives vertical-growth exponent $0$ for $\sigma>1/12$; subconvexity of the $\mu$-function and \Cref{Theorem: muzeta(sigma)} then give (c).
\end{proof}

\section{Estimates for rational isomorphism classes for \texorpdfstring{$m = 13$}{m = 13}}\label{Section: Working over the rationals for m = 13}

In \cref{Section: Estimating twN(X) for m = 13}, we counted the number of elliptic curves over $\bbQ$ with a cyclic $13$-isogeny up to isomorphism over $\bbQalg$ (\Cref{Theorem: asymptotic for twN(X) for m = 13}) for $m = 13$.  In this section, as in \cref{Section: Working over the rationals for m = 7} and \cref{Section: Working over the rationals for m = 10 and 25}, we count all isomorphism classes over $\bbQ$ by enumerating over twists using a Tauberian theorem (\Cref{Theorem: Landau's Tauberian theorem}). We first describe the analytic behavior of $\LQad {13} (s)$.

\begin{theorem}\label{Theorem: relationship between twistL(s) and L(s) for m = 13}
The following statements hold.
\begin{enumalph}
	\item The Dirichlet series $\LQad {13}(s)$ has a meromorphic continuation to the region 
	\begin{equation}\label{Equation: domain of LQ(s) for m = 13}
	\set{s = \sigma + i t \in \bbC : \sigma > 1/12}
	\end{equation}
	with a simple pole at $s = 1/6$ and no other singularities on this region. 
	\item The principal part of $\LQad {13}(s)$ at $s = 1/6$ is
	\begin{equation}
	\frac{\twistLad {13} (1/6)}{3 \zeta(2)} \parent{s - \frac 16}^{-1}.
	\end{equation}
	\end{enumalph}
\end{theorem}

\begin{proof}
	The proof follows by letting $m = 13$ in the argument of \Cref{Theorem: relationship between twistL(s) and L(s) for m = 10 and 25}.
\end{proof}

Using \Cref{Theorem: relationship between twistL(s) and L(s) for m = 13}, we deduce the following lemma.

\begin{lemma}\label{Lemma: Delta NQ(n) is admissible for m = 13}
	The sequence $\parent{\hQad {13}(n)}_{n \geq 1}$ is admissible \textup{(\Cref{Definition: admissible sequences})} with parameters $(1/6,1/12,13/84)$.
\end{lemma}

\begin{proof}
	The proof is structurally identical to the one given for \Cref{Lemma: Delta NQ(n) is admissible for m = 10 and 25}.
\end{proof}

We now prove \Cref{Intro Theorem: asymptotic for NQ(X) for m > 9} for $m = 13$, which we restate here for ease of reference in our established notation.

\begin{theorem}\label{Theorem: asymptotic for NQ(X) for m = 13}
	Define
	\begin{equation} 	\begin{aligned}
	\constad {13} &\colonequals \frac{2 \twistLad {13} (1/6)}{\zeta(2)}.
	\end{aligned}\end{equation} 
For all $\epsilon > 0$, 
	\begin{equation}
	\NQad {13}(X) = \constad {13} X^{1/6} + O\parent{X^{1/8 + \epsilon}}
	\end{equation}
	for $X \geq 1$. The implicit constant depends only on $\epsilon$. 
\end{theorem}

\begin{proof}
	By \Cref{Lemma: Delta NQ(n) is admissible for m = 13}, the sequence $\parent{\hQad {13}(n)}_{n \geq 1}$ is admissible with parameters $\parent{1/6, 1/12, 13/84}$. We now apply \Cref{Theorem: Landau's Tauberian theorem} to the Dirichlet series $\LQad {13}(s)$, and our claim follows. 
\end{proof}

\begin{remark}\label{Remark: true bound for twistN(X) for m = 13}
	We suspect that the true error on $\NQad {13}(X)$ is at most $O(X^{1/12 + \epsilon})$, and the true error on $\twistNad {13} (X)$ is at most $O(X^{1/24 + \epsilon})$, but we have been unable to bound the error terms this far using our techniques. See \Cref{Remark: true bound for twistN(X) for m = 7} for some related thoughts. 
\end{remark}

\section{Computations for \texorpdfstring{$m = 13$}{m = 13}}\label{Section: Computations for m = 13}

In this section, we describe an algorithm to enumerate the elliptic curves with a cyclic $13$-isogeny and twist height at most $X$. We then numerically evaluate the Euler product $\Q{13}$ and the twist-counting constant $\twconst{13}$, and use the list of elliptic curves admitting a cyclic $13$-isogeny to estimate $\constad 13$.

\subsection*{Enumerating elliptic curves with a cyclic $13$-isogeny}

The algorithms described in \cref{Section: Computations for m = 7} and \cref{Section: Computations for m = 10 and 25} cannot be directly adapted to enumerate elliptic curves admitting a cyclic $13$-isogeny, because both the quadratic form $\CII {13} (a, b)$ and the quadratic form $\CIII {13} (a, b)$ can contribute to the twist minimality defect of the pair $(\A {13} (a, b), \B {13} (a, b))$.

We therefore adopt a more na\"ive algorithm. We first obtain bounds \eqref{Equation: bound on a and b for m = 13} on the magnitudes of $a, b \in \bbZ$ subject to the condition $\twht(\A {13} (a, b), \B {13} (a, b)) \leq X$. Recalling \Cref{Theorem: Controlling size of twist minimality defect for m = 13}(b) and using its notation, we have
\begin{equation}\label{Equation: bound on tmd when m = 13}
    \tmd(\A {13} (a, b), \B {13} (a, b)) \leq 2 \cdot 3 \cdot 13 \cdot \CII {13} (a, b)^{1/3} \cdot \CIII {13} (a, b)^{1/2}.
\end{equation}
By \eqref{eqn:htdef} and \eqref{Equation: bound on tmd when m = 13}, if $\twistheight(\A {13}(a, b), \B {13} (a, b)) \leq X$ then
\begin{equation}\label{Equation: bound on H when m = 13}
\frac{\rawheight(\A {13} (a, b), \B {13} (a, b))}{(2 \cdot 3 \cdot 13)^6 \cdot \CII {13} (a, b)^2 \CIII {13} (a, b)^3} \leq X.
\end{equation}
The left-hand side of \eqref{Equation: bound on H when m = 13} is homogeneous of degree 14 in $a$ and $b$; a short computation shows that if $(a, b) \in \bbR^2$ satisfy \eqref{Equation: bound on H when m = 13}, then 
\begin{equation}\label{Equation: bound on a and b for m = 13}
\abs{a} < 5 X^{1/14} \ \text{and} \ \abs{b} < 0.7 X^{1/14}
\end{equation}
(tighter bounds are possible).

For each coprime pair $(a, b) \in \bbZ^2$ satisfying \eqref{Equation: bound on a and b for m = 13} with $b > 0$ and $4\A {13}(a,b)^3+27\B {13}(a,b)^2 \neq 0$, we determine the largest integer $e$ such that $e^6 \mid \gcd(\A {13} (a, b)^3, \B {13} (a, b)^2)$ by computing the prime factorization of this expression. Necessarily this $e$ is the twist minimality defect of $(\A {13} (a, b), \B {13} (a, b))$, and we now use \eqref{eqn:htdef} to compute $\twht(\A {13} (a, b), \B {13} (a, b)$. If the twist height is bounded by $X$, we report $(a, b)$, together with their twist height and any auxiliary information we care to record.

Running our algorithm out to $X = 10^{48}$ in Sage took us approximately $15$ CPU hours on a single core, producing $9644$ elliptic curves admitting a cyclic $13$-isogeny. We list the first few twist minimal elliptic curves admitting a cyclic $13$-isogeny in Table \ref{table:firstfew13}.

\jvtable{table:firstfew13}{\small
\rowcolors{2}{white}{gray!10}
\begin{tabular}{c c c c}
$(A, B)$ & $(a, b)$ & $\twistheight(E)$ & $\twistdefect(E)$ \\
\hline\hline
$(6, 8)$ & $(10, 1)$ & $1728$ & $26364$ \\
$(-84, 322)$ & $(4, 3)$ & $2799468$ & $8788$ \\
$(-35, 350)$ & $(7, 2)$ & $3307500$ & $6591$ \\
$(-338, 2392)$ & $(-2, 5)$ & $154484928$ & $26364$ \\
$(-380, 2850)$ & $(-16, 1)$ & $219488000$ & $26364$ \\
$(-795, 8650)$ & $(1, 4)$ & $2020207500$ & $6591$ \\
$(-2227, 59534)$ & $(11, 5)$ & $95696023212$ & $19773$ \\
$(-9540, 358650)$ & $(-8, 7)$ & $3473005207500$ & $26364$ \\
$(1581, 403310)$ & $(17, 3)$ & $4391791814700$ & $10985$ \\
$(-12818, 745992)$ & $(2, 1)$ & $15025609729728$ & $12$ \\
$(21012, 672590)$ & $(36, 1)$ & $37107540294912$ & $43940$ \\
$(-24474, 1473688)$ & $(-2, 1)$ & $58637420676288$ & $12$ \\
$(-32844, 2292878)$ & $(0, 1)$ & $141946817117868$ & $4$ \\
$(40549, 144566)$ & $(23, 1)$ & $266686134356596$ & $6591$ \\
$(-49851, 4284054)$ & $(-19, 2)$ & $495543307368204$ & $6591$ \\
$(-82739, 9299442)$ & $(5, 7)$ & $2334949780806828$ & $6591$ \\
$(-83595, 9642950)$ & $(1, 1)$ & $2510635086967500$ & $3$ \\
$(-109235, 13896050)$ & $(-1, 1)$ & $5213705551267500$ & $3$ \\
\end{tabular}
}{$E \in \twistE$ with a cyclic $13$-isogeny and $\twht E \leq 10^{16}$}

\subsection*{Computing \texorpdfstring{$\twconst {13}$}{c13tw}}

We first estimate the Euler product $\Q {13}$.  For $p\notin\{2,3,13\}$, set
\begin{equation}
A_p\colonequals \frac{2p^{1/2}}{(p+1)(p+p^{1/2}+1)},\qquad
B_p\colonequals \frac{2p^{1/2}}{(p+1)(p^2+p^{3/2}+p+p^{1/2}+1)}.
\end{equation}
The first case of \eqref{Equation: Euler factors for m = 13} satisfies
\begin{equation}
\frac{2 p^{1/2} \left(p^2+p^{3/2}+2 p+2 p^{1/2}+2\right)}{(p+1) \left(p+p^{1/2}+1\right) \left(p^2+p^{3/2}+p+p^{1/2}+1\right)}=A_p+B_p.
\end{equation}
Consequently, away from the three exceptional primes, the nontrivial local increment is $A_p$ for $p\equiv1\pmod4$, together with $B_p$ for $p\equiv1\pmod3$.  In particular,
\begin{equation}
0<A_p<2p^{-3/2},\qquad 0<B_p<2p^{-5/2}.
\end{equation}
Let $Q_{13,\leq Y}$ denote the Euler product truncated at $p\leq Y$.  Since $\log(1+x)<x$, for $Y\geq8\cdot10^9$ the same Abel-summation argument using the explicit prime-counting estimates of Bennett--Martin--O'Bryant--Rechnitzer \cite[Theorem 1.4]{Bennett-Martin-OBryant-Rechnitzer} gives
\begin{equation}
\begin{aligned}
\log\frac{\Q {13}}{Q_{13,\leq Y}}
<&\frac{Y^{-1/2}}{(1/2)\log Y}
\parent{1+\frac{15}{4\log Y}}\\
&+\frac{Y^{-3/2}}{(3/2)\log Y}
\parent{1+\frac{25}{4\log Y}}.
\end{aligned}
\end{equation}
At $Y=8\cdot10^9$, direct computation gives
\begin{equation}
Q_{13,\leq Y}=9.775\,286\,929\,620\,052\,484\ldots,
\end{equation}
while the exponent in the preceding upper bound is less than
$1.141\,882\cdot10^{-6}$.  Thus
\begin{equation}
9.775\,286\,929\,620\,052<\Q {13}<9.775\,298\,091\,844\,354.
\end{equation}

Using the rejection-sampling estimate for $\R {13}$ below and the finite Euler product as a point estimate, we obtain
\begin{equation}
\twconst {13}=\frac{\Q {13}\R {13}}{\zeta(2)}\approx0.968\,679\,323.
\end{equation}
The Monte Carlo standard error propagated from $\R {13}$ is approximately
$3.58\cdot10^{-6}$, while the deterministic uncertainty from truncating the Euler product is less than $1.11\cdot10^{-6}$.  We report these separately, since the former is a statistical standard error rather than an absolute error bound.

\subsection*{Computing \texorpdfstring{$\constad {13}$}{c13}}

In this subsection, we estimate $\constad {13} = 2 \twistLad {13}(1/6)/\zeta(2)$ by computing the partial sums of $\twistLad {13} (1/6)$:
\begin{equation} 	\begin{aligned}
	\sum_{n \leq 10^{48}} \frac{\twisthad {13}(n)}{n^{1/6}} &= 0.680\,532\,123\,101\,821\,172.
\end{aligned}\end{equation} 
Multiplying by $2/\zeta(2)$, we estimate
\begin{equation}
    \constad {13} \approx 0.827\,427\,843.
\end{equation}

\pagebreak
We estimate $\R {13}$ by rejection sampling on $[-0.8228, 0.8228] \times [0, 0.1822]$. Of $s_{13} \colonequals 61\,749\,000\,000$ samples, $r_{13} \colonequals 33\,570\,382\,383$ lie in $\R {13}$, giving
\begin{equation}
\R {13} \approx 0.299\,828\,32 \cdot\frac{r_{13}}{s_{13}} = 0.163\,004\,281\,067\,749\,86,
\end{equation}
with standard error
\begin{equation}
0.299\,828\,32 \cdot \sqrt{\frac{r_{13}(s_{13} - r_{13})}{s_{13}^3}} < 6.1 \cdot 10^{-7}.
\end{equation}
This took 1 CPU week to compute.

%% file: Ch7_m4_6_8_9_12_16_18.tex
\chapter{Counting elliptic curves with a cyclic \texorpdfstring{$m$}{m}-isogeny for \texorpdfstring{$m \in \set{6, 8, 9, 12, 16, 18}$}{m in 6, 8, 9, 12, 16, 18}}\label{Chapter: other m of genus $0$}

In this chapter, we prove \Cref{Intro Theorem: asymptotic for NQ(X) for 5 < m <= 9} and \Cref{Intro Theorem: asymptotic for NQ(X) for m > 9} when
\[
m \in \set{6, 8, 9, 12, 16, 18}
\]
(\Cref{Theorem: asymptotic for NQ(X) for 5 < m <= 9} and \Cref{Theorem: asymptotic for NQ(X) for m > 9}), and we prove \Cref{Intro Theorem: asymptotic for twN(X) for m of genus $0$} for these $m$ as well as for $m = 4$ (\Cref{Theorem: asymptotic for twN(X) for m of genus $0$} and \Cref{Corollary: asymptotic for twN(X) for m of genus $0$}). These results are new, but our arguments mirror those given in \cref{Chapter: m = 7}, \cref{Chapter: m = 10 and 25}, and \cref{Chapter: m = 13}, and we encourage our readers to skim them on a first perusal of this thesis. Indeed, the sieving required for \Cref{Theorem: asymptotic for twN(X) for m of genus $0$} is much simpler than what we required in previous chapters.

However, in contrast to \cref{Chapter: m = 7}, \cref{Chapter: m = 10 and 25}, and \cref{Chapter: m = 13}, we may have
\begin{equation}
\twistNeq m (X) \neq \twistNad m (X) \ \text{and} \ \NQeq m (X) \neq \NQad m (X)
\end{equation}
for $m \in \set{4, 6, 8, 12, 16}$ (see \Cref{Theorem: elliptic curves with multiple cyclic m-isogenies} and \Cref{Theorem: Parameterizing elliptic curves equipped with pairs of isogenies}). We therefore first estimate $\twistNeq m (X)$ and $\NQeq m (X)$, and then utilize these asymptotics to find estimates for $\twistNad m (X)$ and $\NQad m (X)$ (\Cref{Corollary: asymptotic for twN(X) for m of genus $0$}). We decline to address the cases $m = 2, 3$ for two reasons: firstly, because the asymptotics of $\NQeq 2(X)$, $\NQad 2 (X)$, $\NQeq 3 (X)$, and $\NQad 3 (X)$ are given by \Cref{Theorem: asymptotics for N(X) for m = 2} and \Cref{Theorem: asymptotics for N(X) = 3}, and secondly, because the modular curves $X_0(2)$ and $X_0(3)$ each have exactly one elliptic point, which complicates the application of our method.

The organization of this chapter approximates that of \cref{Chapter: m = 7} and \cref{Chapter: m = 10 and 25}. In \cref{Section: twist minimality defect for other m of genus $0$}, we determine all possible twist minimality defects for 
\begin{equation}\label{Equation: m = 4 and 6 and 8 and 9 and 12 and 16 and 18}
m \in \set{4, 6, 8, 9, 12, 16, 18}.
\end{equation}
In \cref{Section: Estimating twN(X) for m of genus $0$}, we apply the framework developed in \cref{Section: Our approach revisited} to prove \Cref{Intro Theorem: asymptotic for twN(X) for m of genus $0$} for $m$ as in \eqref{Equation: m = 4 and 6 and 8 and 9 and 12 and 16 and 18}. In \cref{Section: Working over the rationals for m of genus $0$}, we prove \Cref{Intro Theorem: asymptotic for NQ(X) for 5 < m <= 9} and \Cref{Intro Theorem: asymptotic for NQ(X) for m > 9} for $m$ as in \eqref{Equation: m = 4 and 6 and 8 and 9 and 12 and 16 and 18} with $m \neq 4$. In \cref{Section: Computations for m of genus $0$}, we produce supplementary computations to estimate the constants appearing in \Cref{Theorem: asymptotic for twN(X) for m of genus $0$} %
\Cref{Theorem: asymptotic for NQ(X) for m > 9}, 
and empirically confirm that the count of elliptic curves with a cyclic $m$-isogeny aligns with our theoretical estimates.

\section{The twist minimality defect for \texorpdfstring{$m \in \set{4, 6, 8, 9, 12, 16, 18}$}{m in 4, 6, 8, 9, 12, 16, 18}}\label{Section: twist minimality defect for other m of genus $0$}

In this section, we bound the twist minimality defect arising from the parameterization
\begin{equation}
y^2 = x^3 + \A m (a, b) x + \B m (a, b)
\end{equation}
for $m \in \set{4, 6, 8, 9, 12, 16, 18}$.

The polynomials $\f m (t)$ and $\g m (t)$ given in Table \ref{table:fm} and Table \ref{table:gm} are coprime when $m \in \set{4, 6, 8, 9, 12, 16, 18}$, so it is markedly easier to handle the twist minimality defect in these cases.

\begin{lemma}\label{Lemma: 2 and 3 are the ungroomed primes for other m of genus $0$}
	Let $m \in \set{4, 6, 8, 9, 12, 16, 18}$, let $(a,b) \in \Z^2$ be $m$-groomed, and let $\ell$ be a prime. If $\ell \mid \twistdefect(\A {m}(a,b),\B m(a,b))$ then $\ell \in \set{2, 3}$.
\end{lemma}

\begin{proof}
	The resultants of $\f m (t)$ and $\g m (t)$ are given in Table \ref{table:auxiliaries}. These resultants are all of the form $\pm 2^v \cdot 3^w$, and the claim follows.
\end{proof}

We now define $\cT m (e)$; this definition is, \textit{mutatis mutandis}, the same as \Cref{Definition: T(e) for m = 7}.

\begin{definition}\label{Definition: T(e) for m of genus $0$}
	Let $m \in \set{4, 6, 8, 9, 12, 16, 18}$. For $e \in \bbZ_{>0}$, let $\tcalT m(e)$ denote the image of
\begin{equation}
\set{(a, b) \in \bbZ^2 : (a, b) \ \text{$m$-groomed}, \ e \mid \twistdefect(\A m (a, b), \B m (a, b))}
\end{equation}
under the projection
\begin{equation}
\bbZ^2 \to (\bbZ / e^3 \bbZ)^2,
\end{equation}
and let $\tcT m(e) \colonequals \# \tcalT m(e)$. Similarly, let $\calT m(e)$ denote the image of
\begin{equation}
\set{t \in \bbZ : e^2 \mid \f m(t) \ \text{and} \ e^3 \mid \g m(t)}
\end{equation}
under the projection
\begin{equation}
\bbZ \to \bbZ / e^3 \bbZ,
\end{equation}
and let $\cT m(e) \colonequals \#\calT m(e)$.
\end{definition}

Let $m \in \set{4, 6, 8, 9, 12, 16, 18}$. By \Cref{Lemma: 2 and 3 are the ungroomed primes for other m of genus $0$}, $\cT m (e) = 0$ whenever $e$ has a prime divisor other than $2$ and $3$. Of course, much more is true.

\begin{lemma}\label{Lemma: bound on T(e) for m of genus $0$}
Let $m \in \set{4, 6, 8, 9, 12, 16, 18}$. The following statements hold.
\begin{enumalph}
\item $\tcalT m(e)$ consists of those pairs $(a, b) \in (\bbZ / e^3 \bbZ)^2$ which satisfy the following conditions:%
\begin{itemize}
	\item $\A m (a, b) \equiv 0 \psmod {e^2}$ and $\B m (a, b) \equiv 0 \psmod {e^3}$, and
	\item $\ell \nmid \gcd(a,b)$ for all primes $\ell \mid e$. %
\end{itemize}
\item Let $(a, b) \in \bbZ^2$. If $(a, b) \psmod{e^3} \in \tcalT {m}(e)$ then $e \mid \tmd(\A {m} (a, b), \B {m} (a, b))$. 
\item The functions $\tcT m(e)$ and $\cT m(e)$ are multiplicative, and $\tcT m(e) = \varphi(e^3) \cT m(e)$.
\item Let $\ell$ be prime and let $v \geq 1$. The nonzero values of $\cT m (\ell^v)$ are given in Table \ref{table:Tm(2)} and Table \ref{table:Tm(3)} below.
\end{enumalph}
\end{lemma}

\begin{proof}
	The proofs of (a)--(c) are exactly as in the proof of \Cref{Lemma: bound on T(e) for m = 7}. Part (d) is a short computation.
\end{proof}

\jvtable{table:Tm(2)}{
\rowcolors{2}{white}{gray!10}
\begin{tabular}{c | c c c c}
$m$ & $\cT m (2^1)$ & $\cT m (2^2)$ & $\cT m (2^3)$ & $\cT m (2^4)$ \\
\hline\hline
$4$ & $2^2$ & --- & --- & --- \\
$6$ & $2^2$ & $2^5$ & --- & --- \\
$8$ & $2^2$ & $2^5$ & --- & --- \\
$9$ & $2^2$ & --- & --- & --- \\
$12$ & $2^2$ & $2^5$ & --- & --- \\
$16$ & $2^2$ & $2^5$ & --- & --- \\
$18$ & $2^2$ & $2^5$ & $2^8$ & $2^{11}$
\end{tabular}
}{All nonzero $\cT m (2^v)$ for $m \in \set{4, 6, 8, 9, 12, 16, 18}$}

\jvtable{table:Tm(3)}{
\rowcolors{2}{white}{gray!10}
\begin{tabular}{c | c c c c c}
$m$ & $\cT m (3^1)$ & $\cT m (3^2)$ & $\cT m (3^3)$ & $\cT m (3^4)$ & $\cT m (3^5)$ \\
\hline\hline
$4$ & $3^2$ & --- & --- & --- & --- \\
$6$ & $3^2$ & --- & --- & --- & --- \\
$8$ & --- & --- & --- & --- & --- \\
$9$ & $3^2$ & --- & --- & --- & --- \\
$12$ & $3^2$ & $3^5$ & $3^8$ & --- & --- \\
$16$ & --- & --- & --- & --- & --- \\
$18$ & $3^2$ & $3^5$ & $3^8$ & $3^{11}$ & $3^{14}$
\end{tabular}
}{All nonzero $\cT m (3^v)$ for $m \in \set{4, 6, 8, 9, 12, 16, 18}$}

\Cref{Theorem: Controlling size of twist minimality defect for m = 7} can have no direct analogue for $m \in \set{4, 6, 8, 9, 12, 16, 18}$, because $\A m (a, b)$ and $\B m (a, b)$ are coprime polynomials for these $m$. Nevertheless, we have the following proposition.

\begin{proposition}\label{Proposition: bounding twist height for m of genus $0$}
	Let $m \in \set{4, 6, 8, 9, 12, 16, 18}$, and let $(a, b)$ be an $m$-groomed pair. We have
	\begin{equation}
	\frac{1}{2^{24} \cdot 3^{30}} \rawheight(\A m (a, b), \B m (a, b)) \leq \twistheight(\A m (a, b), \B m (a, b)) \leq \rawheight(\A m (a, b), \B m (a, b))
	\end{equation}
\end{proposition}

\begin{proof}
	Examining Table \ref{table:Tm(2)} and Table \ref{table:Tm(3)}, we find $\cT m (2^v) = 0$ for $v > 4$ and $\cT m (3^v) = 0$ for $v > 5$. The claim now follows from \eqref{eqn:htdef}.
\end{proof}

\begin{remark}
	For particular $m$, much sharper bounds on $\twistheight(\A m (a, b), \B m (a, b))$ may be derived than those given in \Cref{Proposition: bounding twist height for m of genus $0$}. For instance, if $m = 4$ and $(a, b)$ is a $4$-groomed pair, then
	\begin{equation}
	\frac{1}{2^{6} \cdot 3^{6}} \rawheight(\A 4 (a, b), \B 4 (a, b)) \leq \twistheight(\A 4 (a, b), \B 4 (a, b)) \leq \rawheight(\A 4 (a, b), \B 4 (a, b)).
	\end{equation}
    We have no need to write down these more careful bounds, however.
\end{remark}

\section{Estimates for twist classes for \texorpdfstring{$m \in \set{4, 6, 8, 9, 12, 16, 18}$}{m in 4, 6, 8, 9, 12, 16, 18}} \label{Section: Estimating twN(X) for m of genus $0$}

In this section, we use \cref{Section: Our approach revisited} to estimate $\twistNeq m(X)$ for $m \in \set{4, 6, 8, 9, 12, 16, 18}$, counting the number of twist minimal elliptic curves $E$ over $\bbQ$ equipped with a cyclic $m$-isogeny.

Recall \eqref{Equation: defining cM}, \eqref{Equation: Applying N calE to Neq m}, and \eqref{Equation: Applying N calE to Nad m}. By \cref{Section: Parameterizing elliptic curves with a cyclic m-isogeny}, for $m \in \set{4, 6, 8, 9, 12, 16, 18}$, $\cM {m}(X; e)$ counts pairs $(a, b) \in \bbZ^2$ with
\begin{itemize}
	\item $(a, b)$ $m$-groomed;
	\item $\rawheight(\A m(a,b),\B m(a,b)) \leq X$; 
	\item $e \mid \twistdefect(\A m(a, b), \B m(a, b))$.
\end{itemize}
If $m = 4$, then by \Cref{Lemma: parameterizing m-isogenies}, we double-count multiples of the pair $(0, 1)$ where it appears.

We have the following refinement of \Cref{Lemma: fundamental sieve}.

\begin{lemma}\label{Lemma: fundamental sieve for m of genus $0$}
	Let $m \in \set{4, 6, 8, 9, 12, 16, 18}$. We have
	\begin{equation}\label{Equation: twN(X) in terms of M(X; e) for m of genus $0$}
	\twNeq m (X) = \sum_{n \mid 2^4 \cdot 3^5} \sum_{e \mid n} \mu(n/e) \cM m (e^6 X; n).
	\end{equation}
\end{lemma}

\begin{proof}
	Recall \eqref{eqn:twistheightcalc}. As in the proof of \Cref{Proposition: bounding twist height for m of genus $0$}, we examine Table \ref{table:Tm(2)} and Table \ref{table:Tm(3)}. We find $\cT m (2^v) = 0$ for $v > 4$ and $\cT m (3^v) = 0$ for $v > 5$. By \eqref{Equation: defining cM}, $\cM m (e^6 X; n) = 0$ when $n \nmid 2^4 \cdot 3^5$.
\end{proof}

\begin{remark}
	As with \Cref{Proposition: bounding twist height for m of genus $0$}, sharper bounds are possible. For instance, when $m = 4$, we have
	\begin{equation}
	\twNeq 4 (X) = \sum_{n \mid 2 \cdot 3} \sum_{e \mid n} \mu(n/e) \cM 4 (e^6 X; n).
	\end{equation}
\end{remark}

Let $m \in \set{4, 6, 8, 9, 12, 16, 18}$. In order to estimate $\cM m(X; e)$, we further unpack the $m$-groomed condition on pairs $(a, b)$, as we did before with \eqref{equation: cM(X;e) in terms of cM(X; d, e) for m = 7} and \eqref{equation: cM(X;e) in terms of cM(X; d, e) for m = 5 and 10 and 25}. We therefore let $\cM m(X; d, e)$ denote the number of pairs $(a, b) \in \bbZ^2$ with
	\begin{itemize}
		\item $\gcd(da, db, e) = 1$, $b > 0$, and $a/b \not\in \cusps m$,
		\item $\rawheight(\A m(d a, d b), \B m(d a, db)) \leq X$, and
		\item $e \mid \twistdefect(\A m(d a, d b), \B m(da, db))$.
	\end{itemize}
Because $\rawheight(\A m(a, b), \B m(a, b))$ is homogeneous of degree $2 \degB m$, another M\"{o}bius sieve yields
\begin{equation}
	\cM m(X; e) = \sum_{\substack{d \ll X^{1/(2 \degB m)} \\ \gcd(d, e) = 1}} \mu(d) \cM m(X; d, e).\label{equation: cM(X;e) in terms of cM(X; d, e) for m of genus $0$}
\end{equation}	

\begin{lemma}\label{Lemma: asymptotic for M(X; e) for m of genus $0$}
Let $m \in \set{4, 6, 8, 9, 12, 16, 18}$. The following statements hold.
\begin{enumalph}
\item	If $\gcd(d, e) > 1$, then $\cM {m}(X; d, e) = 0$. If $\gcd(d, e) = 1$, we have 
	\begin{equation}
	\cM m(X; d, e) = \frac{\R m \cT m(e) X^{1/\degB m}}{d^2 e^2} \prod_{\ell \mid e} \parent{1 - \frac{1}{\ell}} + O\parent{\frac{2^{\omega(e)} X^{1/(2 \degB m)}}{d}}
	\end{equation}
	for $X, d, e \geq 1$, where $\R m$ is the area of \eqref{eqn: R(X)}.
\item We have
	\begin{equation}
	\cM m(X; e) = \frac{\R m \cT m(e) X^{1/\degB m}}{\zeta(2) e^3 \prod_{\ell \mid e} \parent{1 + \frac{1}{\ell}}} +  O\parent{2^{\omega(e)} X^{1/(2 \degB m)} \log X}
	\end{equation}
 for $X \geq 2$ and $e \geq 1$.
\end{enumalph}
In both cases, the implied constants are independent of $d$, $e$, and $X$.
\end{lemma}

Because $\gcd(\A m (a, b), \B m (a, b)) = 1$ for $m \in \set{4, 6, 8, 9, 12, 16, 18}$, we are unable to produce proofs akin to the second proofs of \Cref{Lemma: asymptotic for M(X; e) for m = 7} and \Cref{Lemma: asymptotic for M(X; e) for m = 5 and 10 and 25}. We do not need them, however, because by \Cref{Proposition: bounding twist height for m of genus $0$} there are only finitely many possible values for $\tmd(\A m (a, b), \B m (a, b))$.

\begin{proof}
	The proof is, \textit{mutatis mutandis}, the same as the first proof of \Cref{Lemma: asymptotic for M(X; e) for m = 7}.
\end{proof}

In analogy with \eqref{Equation: Q for m = 7}, for $m \in \set{4, 6, 8, 9, 12, 16, 18}$, we let
	\begin{equation}\label{Equation: Qm for m of genus $0$}
	\Q m \colonequals \sum_{n \mid 2^4 \cdot 3^5} \frac{\varphi_{6 / \degB m}(n) \cT m(n)}{n^3 \prod_{\ell \mid n} \parent{1 + \frac{1}{\ell}}},
	\end{equation}
	and we let
\begin{equation}\label{Equation: twconsteq for m of genus $0$}
	\twconsteq m \colonequals \frac{\Q m \R m}{\zeta(2)}.
\end{equation}
Here, as always, $\R m$ is the area of the region
	\begin{equation}
	\calR m(1) = \set{(a, b) \in \bbR^2 : \rawheight(\A m (a, b), \B m (a, b)) \leq 1, b \geq 0}.
	\end{equation}
    We also let
    \begin{equation}\label{Equation: twconst for m of genus $0$}
    \twconst m \colonequals \twconsteq m / \delta_m,
    \end{equation}
    where
    \begin{equation}\label{Equation: defining delta m}
        \delta_m \colonequals \begin{cases} 
        2 & \text{if} \ 4 \mid m, \\
        1 & \text{else}.
        \end{cases}
    \end{equation}
    The factor $\delta_m$ appears because cyclic $m$-isogenies come in pairs when $4 \mid m$ (\Cref{Theorem: Parameterizing elliptic curves equipped with pairs of isogenies}). 

\jvtable{table:Qm}{
\rowcolors{2}{white}{gray!10}
\begin{tabular}{c | c c c}
$m$ & $\Q m$ & $\Q m (2)$ & $\Q m (3)$ \\
\hline\hline
$4$ & $6$ & $2$ & $3$ \\
$6$ & $3$ & $2$ & $3/2$ \\
$8$ & $2$ & $2$ & --- \\
$9$ & $2$ & $4/3$ & $3/2$ \\
$12$ & $1 + \sqrt 3$ & $4/3$ & $3 (1 + \sqrt 3)/4$ \\
$16$ & $4/3$ & $4/3$ & --- \\
$18$ & $\left(1+2^{1/3}\right) \left(1+3^{2/3}\right)/2$ & $2(1 + 2^{1/3})/3$ & $3(1 + 3^{2/3})/4$
\end{tabular}
}{$\Q m$ and its nontrivial Euler factors for $m \in \set{4, 6, 8, 9, 12, 16, 18}$}

We need not estimate $\twistNeqly m (X)$ and $\twistNeqgy m (X)$ in order to estimate $\twistNeq m (X)$, because by \Cref{Lemma: fundamental sieve for m of genus $0$} we have
\begin{equation}
\twistNeq m (X) = \twistNeqly m (X)
\end{equation}
whenever $y > 2^4 \cdot 3^5$. We are now in a position to estimate $\twistNeq m(X)$.

\begin{theorem}\label{Theorem: asymptotic for twN(X) for m of genus $0$}
	Let $m \in \set{4, 6, 8, 9, 12, 16, 18}$. Then we have
	\begin{equation}
	\twistNeq m(X) = \twconsteq m X^{1/\degB m} + O\parent{X^{1/(2 \degB m)} \log X},
	\end{equation}
	for $X \geq 2$. The constant $\twconsteq m$ is given in \eqref{Equation: twconsteq for m of genus $0$}.
\end{theorem}

\begin{proof}
	Substituting the asymptotic for $\cM m(X; e)$ from \Cref{Lemma: asymptotic for M(X; e) for m of genus $0$}(b) into \eqref{Equation: twN(X) in terms of M(X; e) for m of genus $0$} from \Cref{Lemma: fundamental sieve for m of genus $0$}, we have
	\begin{equation} 	\begin{aligned}
		\twistNeq m(X) =& \sum_{n \mid 2^4 \cdot 3^5} \sum_{e \mid n} \mu\parent{\frac ne} \frac{\R m \cT m(n) e^{6/\degB m} X^{1/\degB m}}{\zeta(2) n^3 \prod_{\ell \mid n} \parent{1 + \frac{1}{\ell}}} \\
        &+ \sum_{n \mid 2^4 \cdot 3^5} \sum_{e \mid n} \mu\parent{\frac ne} O\parent{2^{\omega(n)} e^{3/\degB m} X^{1/(2 \degB m)} \log e^6 X}.
	\end{aligned}\end{equation} 
	
	We handle the main term and the error of this expression separately. For the main term, we have
	\begin{equation}\label{Equation: sum for Q for m of genus $0$}
	\begin{aligned}
	\sum_{n \mid 2^4 \cdot 3^5} \sum_{e \mid n} \mu\parent{n/e} \frac{\R m \cT m(n) e^{6/\degB m} X^{1/\degB m}}{\zeta(2) n^3 \prod_{\ell \mid n} \parent{1 + \frac{1}{\ell}}} &= \frac{\R m X^{1 / \degB m}}{\zeta(2)} \sum_{n \mid 2^4 \cdot 3^5} \frac{\cT m (n) \varphi_{6/\degB m}(n)}{n^3 \prod_{\ell \mid n} \parent{1 + \frac{1}{\ell}}} \\
	&= \twconsteq m X^{1/\degB m},
	\end{aligned}
	\end{equation}
	where the last equality follows from \eqref{Equation: Qm for m of genus $0$} and \eqref{Equation: twconsteq for m of genus $0$}.
	
	The summands of \eqref{Equation: sum for Q for m of genus $0$} constitute a nonnegative multiplicative arithmetic function, so we can factor $\Q m$ as an Euler product \begin{equation}
	\Q m = \prod_p \Q m (p).
	\end{equation}
	Moreover, for $p \neq 2, 3$ prime, \Cref{Lemma: bound on T(e) for m of genus $0$} implies $\Q m (p) = 1$, so in fact $\Q m = \Q m (2) \Q m (3)$. The values of $\Q m$, $\Q m (2)$, and $\Q m(3)$ are recorded in Table \ref{table:Qm} above.
	
	We now turn to the error term. As $e \leq 2^4 \cdot 3^5$, we have $\log (e^6 X) \ll \log X$. We obtain
	\begin{equation}
	\begin{aligned}
		&\sum_{n \mid 2^4 \cdot 3^5} \sum_{e \mid n} \mu\parent{n/e} O\parent{2^{\omega(n)} e^{3/\degB m} X^{1/(2 \degB m)} \log e^6 X} \\
  =& O\parent{X^{1/(2 \degB m)} \log X \sum_{n \mid 2^4 \cdot 3^5}  2^{\omega(n)} \sum_{e \mid n} \abs{\mu\parent{n/e}} e^{3/\degB m}} \\
	=& O\parent{X^{1/(2 \degB m)} \log X},
	\end{aligned}
	\end{equation}
	because
	\begin{equation}
	 \sum_{n \mid 2^4 \cdot 3^5}  2^{\omega(n)} \sum_{e \mid n} \abs{\mu\parent{n/e}} e^{3/\degB m} \leq \sum_{n \mid 2^4 \cdot 3^5}  2^{\omega(n)} \prod_{p \mid n} \parent{1 + p} = 25 \cdot 41 < \infty,
	\end{equation}
	independent of $m$. This proves our desired result.
\end{proof}

To finish our proof of \Cref{Intro Theorem: asymptotic for twN(X) for m of genus $0$}, we bound the difference between $\twNeq m (X)$ and $\delta_m \twNad m (X)$ for $m \in \set{4, 5, 6, 8, 9, 12, 16, 18}$. 

\begin{lemma}\label{Lemma: Difference between twNeq and twNad}
    Let $m \in \set{4, 5, 6, 8, 9, 12, 16, 18}$. We have
    \begin{equation}
\begin{alignedat}{2}
\twNeq 4 (X) &= 2 \twNad 4 (X) &&+ O(X^{1/12}), \\
\twNeq 5 (X) &= \twNad 5 (X) &&+ O(X^{1/18}), \\
\twNeq 6 (X) &= \twNad 6 (X) &&+ O(X^{1/12}), \\
\twNeq 8 (X) &= 2 \twNad 8 (X) &&+ O(X^{1/12}),
\end{alignedat}
\qquad
\begin{alignedat}{2}
\twNeq 9 (X) &= \twNad 9 (X), && \\
\twNeq {12} (X) &= 2 \twNad {12} (X) &&+ O(1), \\
\twNeq {16} (X) &= 2 \twNad {16} (X) &&+ O(1), \ \text{and} \\
\twNeq {18} (X) &= \twNad {18} (X). &&
\end{alignedat}
\end{equation}
    for $X \geq 1$. The implicit constant depends on $m$.
\end{lemma}

\begin{proof}
    When $m \in \set{9, 12, 16, 18}$, \Cref{Lemma: Difference between twNeq and twNad} is a consequence of \Cref{Corollary: Elliptic curves for which twNeq(X) = twNad(X)} and \Cref{Corollary: discrepancy between twNeq m and twNad m when 4 | m}.
    
    We recall \Cref{Theorem: Parameterizing elliptic curves equipped with pairs of isogenies}, along with Table \ref{table:modularcurves}, Table \ref{table:fmn}, and Table \ref{table:gmn} to address the remaining cases. For each proper divisor $n$ of $m$, we consider the contribution to $\twistNeq m (X) - \delta_m \twistNad m (X)$ arising from elliptic curves equipped with a pair of unsigned cyclic $m$-isogenies whose kernels have intersection of order $n$. When $(m, n) = (6, 1)$, the associated modular curve has genus 1, and is $\bbQ$-isomorphic to the elliptic curve $y^2 = x^3 + 1$, which has Mordell-Weil group $\bbZ/6\bbZ$, so this contribution is $O(1)$.

    If
    \begin{equation}
    (m, n) \in \set{(2, 1), (3, 1), (4, 1), (5, 1), (6, 2), (6, 3), (8, 2)},
    \end{equation}
    the associated modular curve has genus $0$, but is not $X_0(m)$. In these cases, the universal families
    \begin{equation}
    y^2 = x^3 + \f {m, n} (t) x + \g {m, n} (t)
    \end{equation}
    for elliptic curves (over $\bbQ$, up to quadratic twist) with this level structure are recorded in Table \ref{table:fmn} and Table \ref{table:gmn}.
    
    For such pairs $(m, n)$, we can repeat the proof of \Cref{Theorem: asymptotic for twN(X) for m of genus $0$} essentially verbatim to obtain asymptotics for the number of elliptic curves equipped with two cyclic $m$-isogenies with kernels having order $n$ intersection, with one caveat. When $m = 5$, we must follow the proof of \Cref{Theorem: asymptotic for twN(X) for m = 10 and 25} rather than of \Cref{Theorem: asymptotic for twN(X) for m of genus $0$}, because in this case $\f {5, 1}$ and $\g {5, 1}$ have a common factor $t^2 + 4$, just like $\f {25} (t)$ and $\g {25} (t)$.

    Except when $m = 5$ and $n = 1$, i.e., when our modular curve has elliptic points, the asymptotics we obtain are the same as those predicted by \cite[Theorem 3.3.1]{Cullinan-Kenney-Voight}, although we cannot apply this theorem directly because the congruence groups inducing our modular curves are not torsion-free. They also accord with \cite[Theorem 1.2.2, Theorem 1.2.3]{Phillips1} wherever these theorems apply.

\end{proof}

\begin{remark}
    It is no coincidence that we can count elliptic curves equipped with a single unsigned cyclic $25$-isogeny, or with two unsigned cyclic $5$-isogenies, by essentially the same methods. Indeed, the modular curve $X_{\textup{sp}}(p)$ parameterizing elliptic curves equipped with two unsigned cyclic $p$-isogenies whose kernels have trivial intersection is isomorphic to the curve $X_0(p^2)$, because taking the dual of one of these two isogenies and composing yields a cyclic $p^2$-isogeny (see \cite[Theorem 3.2.1]{Cullinan-Kenney-Voight}).
\end{remark}

\begin{corollary}\label{Corollary: asymptotic for twN(X) for m of genus $0$}
	Let $m \in \set{4, 6, 8, 9, 12, 16, 18}$. We have
	\begin{equation}
	\twistNad m(X) = \twconstad m X^{1/\degB m} + O\parent{X^{1/(2 \degB m)} \log X},
	\end{equation}
	for $X \geq 2$. We have $\twconstad m = \twconsteq m / 2$ if $m \in \set{4, 8, 12, 16}$, and $\twconstad m = \twconsteq m$ if $m \in \set{6, 9, 18}$.
\end{corollary}

\begin{proof}
    Immediate from \Cref{Theorem: asymptotic for twN(X) for m of genus $0$} and \Cref{Lemma: Difference between twNeq and twNad}.
\end{proof}

\subsection*{$L$-series}

To conclude this section, following \Cref{Corollary: twL(s) has a meromorphic continuation for m = 7}, \Cref{Corollary: twL(s) has a meromorphic continuation for m = 10 and 25}, and \Cref{Corollary: twL(s) has a meromorphic continuation for m = 13}, we set up \cref{Section: Working over the rationals for m of genus $0$} by interpreting the asymptotics given by \Cref{Theorem: asymptotic for twN(X) for m of genus $0$} in terms of Dirichlet series.

Recall \eqref{Equation: twisth calE}, \eqref{Equation: hQ calE}, \eqref{Equation: defining twistL calE X}, and \eqref{Equation: defining LQ calE X}.

\begin{cor}\label{Corollary: twL(s) has a meromorphic continuation for m of genus $0$}
	Let $m \in \set{4, 6, 8, 9, 12, 16, 18}$. The following statements hold.
 \begin{enumerate}
     \item The Dirichlet series $\twistLeq m (s)$ and $\twistLad m (s)$ have abscissa of (absolute) convergence $\sigma_a = \sigma_c = 1/\degB m$, and have a meromorphic continuation to the region
		\begin{equation}\label{equation: domain of twistL(s) for m of genus $0$}
	\set{s = \sigma + i t \in \bbC : \sigma > 1/(2 \degB m)}. 
	\end{equation}
    \item The function $\twistLeq m(s)$ has a simple pole at $s = 1/\degB m$ with residue
	 \begin{equation}
	 \res_{s=\frac{1}{\degB m}} \twistLeq m(s) = \frac{\twconsteq m}{\degB m},
	 \end{equation}
	 and is holomorphic elsewhere on this region. Likewise, $\twistLad m(s)$ has a simple pole at $s = 1/\degB m$ with residue
	 \begin{equation}
	 \res_{s=\frac{1}{\degB m}} \twistLad m(s) = \frac{\twconst m}{\degB m},
	 \end{equation}
	 and is holomorphic elsewhere on this region.
  \item We have
\begin{equation}
    \mu_{\twistLad m}(\sigma) \leq \max\{0,2-2\degB m\,\sigma\}
     \end{equation}
     for $\sigma > 1/(2 \degB m)$. The same bound holds for $\twistLeq m(s)$.
     \end{enumerate}
\end{cor}

\begin{proof}
	The proof is structurally identical to the one given for \Cref{Corollary: twL(s) has a meromorphic continuation for m = 7}.
\end{proof}

\section{Estimates for rational isomorphism classes for \texorpdfstring{$m \in \set{6, 8, 9, 12, 16, 18}$}{m in 6, 8, 9, 12, 16, 18}}\label{Section: Working over the rationals for m of genus $0$}

In \cref{Section: Estimating twN(X) for m of genus $0$}, we counted the number of elliptic curves over $\bbQ$ with a (cyclic) $m$-isogeny up to quadratic twist (\Cref{Theorem: asymptotic for twN(X) for m of genus $0$}) for $m \in \set{4, 6, 8, 9, 12, 16, 18}$.  In this section, we count all isomorphism classes over $\bbQ$ by enumerating over twists using Landau's Tauberian theorem (\Cref{Theorem: Landau's Tauberian theorem}). We will not handle the case $m = 4$ as it has already been addressed previously in \cite{Pomerance-Schaefer}; there is no theoretical obstruction to working this case out using our methods, however. As in previous chapters, we first describe the analytic behavior of $\LQeq m (s)$ and $\LQad m (s)$ for $m \in \set{6, 8, 9, 12, 16, 18}$.

\begin{theorem}\label{Theorem: relationship between twistL(s) and L(s) for 5 < m <= 9}
Let $m \in \set{6, 8, 9}$. The following statements hold.
\begin{enumalph}
	\item The Dirichlet series $\LQeq m (s)$ and $\LQad m(s)$ have a meromorphic continuation to the region \eqref{equation: domain of twistL(s) for m of genus $0$} (i.e., \eqref{equation: domain of twistL(s) for m = 7}) with a double pole at $s = 1/6$ and no other singularities on this region. 
	\item The principal part of $\LQeq m (s)$ at $s = 1/6$ is
	\begin{equation}
	\frac{1}{3 \zeta(2)} \parent{\frac{\twconsteq m}{6} \parent{s - \frac 16}^{-2} + \parent{\widetilde{\ell}_{m,0} + \twconsteq m \parent{\gamma - \frac{2 \zeta\prm(2)}{\zeta(2)}}} \parent{s - \frac{1}{6}}\inv},
	\end{equation}
	where $\twconsteq m$ is given in \eqref{Equation: twconsteq for m of genus $0$}, and
	\begin{equation}
	\widetilde{\ell}_{m, 0} \colonequals \twconsteq m \gamma + \frac {1}6 \int_1^\infty \parent{\twistNeq m(u) - \twconsteq m \floor{u^{1/6}}} u^{-7/6} \,\mathrm{d}u\label{equation: Definition of elleq0 for m of genus $0$}
	\end{equation}
	is the constant term of the Laurent expansion for $\twistLeq m(s)$ around $s = 1/6$. Here, $\gamma$ denotes the Euler-Mascheroni constant. Likewise, the principal part of $\LQad m(s)$ at $s = 1/6$ is
	\begin{equation}
	\frac{1}{3 \zeta(2)} \parent{\frac{\twconst m}{6} \parent{s - \frac 16}^{-2} + \parent{\ell_{m,0} + \twconst m \parent{\gamma - \frac{2 \zeta\prm(2)}{\zeta(2)}}} \parent{s - \frac{1}{6}}\inv},
	\end{equation}
	where $\twconst m$ is given in \eqref{Equation: twconst for m of genus $0$}, and
	\begin{equation}
	\ell_{m, 0} \colonequals \twconst m \gamma + \frac {1}6 \int_1^\infty \parent{\twistNad m(u) - \twconst m \floor{u^{1/6}}} u^{-7/6} \,\mathrm{d}u\label{equation: Definition of ell0 for m of genus $0$}
	\end{equation}
	is the constant term of the Laurent expansion for $\twistLad m(s)$ around $s = 1/6$.
	\end{enumalph}
\end{theorem}

\begin{proof}
    The proof is, \textit{mutatis mutandis}, the same as the proof of \Cref{Theorem: relationship between twistL(s) and L(s) for m = 7}; however, we must run through the argument given there twice, once for $\LQeq m (s)$ and once for $\LQad m (s)$.
\end{proof}

\begin{theorem}\label{Theorem: relationship between twistL(s) and L(s) for m > 9}
Let $m \in \set{12, 16, 18}$. The following statements hold.
\begin{enumalph}
	\item The Dirichlet series $\LQeq m (s)$ and $\LQad m(s)$ have a meromorphic continuation to the region 
	\begin{equation}\label{Equation: domain of LQ(s) for m > 9}
	\set{s = \sigma + i t \in \bbC : \sigma > 1/12}
	\end{equation}
	with a simple pole at $s = 1/6$ and no other singularities on this region. 
	\item The principal part of $\LQeq m(s)$ at $s = 1/6$ is
	\begin{equation}
	\frac{\twistLeq m (1/6)}{3 \zeta(2)} \parent{s - \frac 16}^{-1};
	\end{equation}
 likewise, the principal part of $\LQad m(s)$ at $s = 1/6$ is
	\begin{equation}
	\frac{\twistLad m (1/6)}{3 \zeta(2)} \parent{s - \frac 16}^{-1}.
	\end{equation}
	\end{enumalph}
\end{theorem}

\begin{proof}
    The proof is, \textit{mutatis mutandis}, the same as the proof of \Cref{Theorem: relationship between twistL(s) and L(s) for m = 10 and 25}; however, we must run through the argument given there twice, once for $\LQeq m (s)$ and once for $\LQad m (s)$.
\end{proof}

We prove two analogues to \Cref{Lemma: Delta NQ(n) is admissible for m = 7}: one for $m \in \set{6, 8, 9}$ and the other for $m \in \set{12, 16, 18}$.

\begin{lemma}\label{Lemma: Delta NQ(n) is admissible for 5 < m <= 9}
	Let $m \in \set{6, 8, 9}$ and let $0<\delta<7/97$. The sequences 
 \begin{equation}
 \parent{\hQeq m(n)}_{n \geq 1} \ \text{and} \ \parent{\hQad m(n)}_{n \geq 1}
 \end{equation}
 are admissible \textup{(\Cref{Definition: admissible sequences})} with parameters $\parent{1/6,\delta,97\delta/7}$.
\end{lemma}

\begin{proof}
    The proof is structurally identical to the corrected proof of \Cref{Lemma: Delta NQ(n) is admissible for m = 7}; however, we must run through the argument twice, once for $\LQeq m (s)$ and once for $\LQad m (s)$.
\end{proof}

\begin{lemma}\label{Lemma: Delta NQ(n) is admissible for m > 9}
	Let $m \in \set{12, 16, 18}$. The sequences 
 \begin{equation}
 \parent{\hQeq m(n)}_{n \geq 1} \ \text{and} \ \parent{\hQad m(n)}_{n \geq 1}
 \end{equation}
 are admissible \textup{(\Cref{Definition: admissible sequences})} with parameters $(1/6, 1/12, 13/84)$.
\end{lemma}

\begin{proof}
    The proof is structurally identical to the one given for \Cref{Lemma: Delta NQ(n) is admissible for m = 10 and 25}; however, we must run through the argument given there twice, once for $\LQeq m (s)$ and once for $\LQad m (s)$.
\end{proof}

With \Cref{Lemma: Delta NQ(n) is admissible for 5 < m <= 9} and \Cref{Lemma: Delta NQ(n) is admissible for m > 9} at the ready, the proofs of \Cref{Intro Theorem: asymptotic for NQ(X) for 5 < m <= 9} and \Cref{Intro Theorem: asymptotic for NQ(X) for m > 9} are essentially identical.

\begin{theorem}\label{Theorem: asymptotic for NQ(X) for 5 < m <= 9} 
    Let $m \in \set{6, 8, 9}$. We define
	\begin{equation}\label{Equation: constants for 5 < m <= 9}
        \begin{aligned}
         \consteq m &\colonequals \frac{\twconsteq m}{3 \zeta(2)}, \\
	\consteqprm m &\colonequals \frac{2}{\zeta(2)}\parent{\widetilde{\ell}_{m, 0} + \twconsteq m \parent{\gamma - 1 - \displaystyle{\frac{2 \zeta\prm(2)}{\zeta(2)}}}}, \\
	\constad m &\colonequals \frac{\twconst m}{3 \zeta(2)}, \ \text{and} \\
	\constadprm m &\colonequals \frac{2}{\zeta(2)}\parent{\ell_{m, 0} + \twconst m \parent{\gamma - 1 - \displaystyle{\frac{2 \zeta\prm(2)}{\zeta(2)}}}},
        \end{aligned}
	\end{equation}
	where $\twconsteq m$ is defined in \eqref{Equation: twconsteq for m of genus $0$}, $\twconst m$ is defined in \eqref{Equation: twconst for m of genus $0$}, $\widetilde{\ell}_{m, 0}$ is defined in \eqref{equation: Definition of elleq0 for m of genus $0$}, and $\ell_{m, 0}$ is defined in \eqref{equation: Definition of ell0 for m of genus $0$}. Then for all $\epsilon > 0$, we have
    \begin{equation}
    \NQeq m (X) = \consteq m X^{1/6} \log X + \consteqprm m X^{1/6} + O(X^{38/291 + \epsilon})
    \end{equation}
    and
    \begin{equation}
    \NQad m (X) = \constad m X^{1/6} \log X + \constadprm m X^{1/6} + O(X^{38/291 + \epsilon})
    \end{equation}
    for $X \geq 1$. The implicit constants depend on $m$ and $\epsilon$.
\end{theorem}

\begin{proof}
	Fix $\epsilon>0$. By \Cref{Lemma: Delta NQ(n) is admissible for 5 < m <= 9}, for every $0<\delta<7/97$ both $\parent{\hQeq m(n)}_{n \geq 1}$ and $\parent{\hQad m(n)}_{n \geq 1}$ are admissible with parameters $\parent{1/6,\delta,97\delta/7}$. Since $97\delta/7<1$, applying \Cref{Theorem: Landau's Tauberian theorem} and then taking $\delta$ sufficiently close to $7/97$ gives the error $O(X^{38/291+\epsilon})$ in both cases.
\end{proof}

\begin{theorem}\label{Theorem: asymptotic for NQ(X) for m > 9}
    Let $m \in \set{12, 16, 18}$. Define
        \begin{equation}
	\begin{aligned}
        \consteq m &= \frac{2 \twistLeq m (1/6)}{\zeta(2)}, \ \text{and} \\
	\constad m &= \frac{2 \twistLad m (1/6)}{\zeta(2)}.
	\end{aligned}
 \end{equation}
    Then for all $\epsilon > 0$, we have
    \begin{equation}
    \NQeq m (X) = \consteq m X^{1/6} + O(X^{1/8 + \epsilon})
    \end{equation}
    and
    \begin{equation}
    \NQad m (X) = \constad m X^{1/6} + O(X^{1/8 + \epsilon})
    \end{equation}
    for $X \geq 1$. The implicit constants depend on $m$ and $\epsilon$.
\end{theorem}

\begin{proof}
	By \Cref{Lemma: Delta NQ(n) is admissible for m > 9}, both $\parent{\hQeq m(n)}_{n \geq 1}$ and $\parent{\hQad m(n)}_{n \geq 1}$ are admissible with parameters $\parent{1/6, 1/12, 13/84}$. We now apply \Cref{Theorem: Landau's Tauberian theorem} to the Dirichlet series $\LQeq m(s)$ and $\LQad m(s)$, and our claim follows. 
\end{proof}

\section{Computations for \texorpdfstring{$m \in \set{4, 6, 8, 9, 12, 16, 18}$}{m in 4, 6, 8, 9, 12, 16, 18}}\label{Section: Computations for m of genus $0$}

In this section, we furnish computations that render %
\Cref{Theorem: asymptotic for twN(X) for m of genus $0$}, \Cref{Corollary: asymptotic for twN(X) for m of genus $0$}, \Cref{Theorem: asymptotic for NQ(X) for 5 < m <= 9}, and \Cref{Theorem: asymptotic for NQ(X) for m > 9}
completely explicit. These computations mirror those given in \cref{Section: Computations for m = 7} and \cref{Section: Computations for m = 10 and 25}.

\subsection*{Enumerating elliptic curves with a cyclic $m$-isogeny for \texorpdfstring{$m \in \set{6, 8, 9, 12, 16, 18}$}{m in 6, 8, 9, 12, 16, 18}}

We begin by outlining an algorithm for computing all elliptic curves equipped with (or admitting) a cyclic $m$-isogeny up to twist height $X$. Write $e_m$ for the maximum twist minimality defect of a pair $(\A m (a, b), \B m (a, b))$ with $\gcd(a, b) = 1$: these values are determined by Table \ref{table:Tm(2)} and Table \ref{table:Tm(3)}. For instance, $e_8 = 2^2$.

Recall from the proof of \Cref{Lemma: formula for R(X)} that 
\begin{equation}
\calR m (X) = \calR m (1) X^{1/(2 \degB m)};
\end{equation}
however, $\calR m (1)$ is a compact region, and thus is contained within a rectangular enveloping region $[a_{\textup{min}}, a_{\textup{max}}] \times [0, b_{\textup{max}}]$ (see Table \ref{table:Rm}). For all coprime pairs of integers $(a, b)$ within the enveloping region 
\begin{equation}
\parent{[a_{\textup{min}}, a_{\textup{max}}] \times [0, b_{\textup{max}}]} \cdot e_m^{3/\degB m} X^{1/(2 \degB m)},
\end{equation}
we compute $(\A m (a, b), \B m (a, b))$.
For each such pair, we compute the twist minimality defect and thence the twist height of $(\A m (a, b), \B m (a, b))$, and we report the result if this twist height is at most $X$.

This algorithm furnishes all elliptic curves equipped with a cyclic $m$-isogeny. To instead enumerate elliptic curves admitting a cyclic $m$-isogeny, we simply omit duplicate elliptic curves.

In Table \ref{table:enumerating}, we record a bound $X$ and the number of elliptic curves equipped with cyclic $m$-isogeny up to that bound $X$, as well as approximately how many CPU minutes it took us to compute that list of elliptic curves.

\jvtable{table:enumerating}{
\rowcolors{2}{white}{gray!10}
\begin{tabular}{c | c c c c }
$m$ & $X$ & $\twNeq m (X)$ & $\twNad m (X)$ & CPU minutes \\
\hline
\hline
$4$ & $10^{21}$ & $6\,299\,452$ & $3\,149\,720$ & $46$ \\ %
$6$ & $10^{42}$ & 7\,551\,963 & $7\,550\,700$ & $60$ \\ %
$8$ & $10^{42}$ & $5\,855\,992$ & $2\,927\,707$ & $27$ \\ %
$9$ & $10^{42}$ & $4\,671\,446$ & $4\,671\,446$ & $13$ \\ %
$12$ & $10^{84}$ & $10\,478\,972$ & $5\,239\,486$ & $74$ \\ %
$16$ & $10^{84}$ & $7\,836\,058$ & $3\,918\,029$ & $33$ \\%
$18$ & $10^{126}$ & $9\,730\,625$ & $9\,730\,625$ & $76$ %
\end{tabular}}{Enumerating elliptic curves up to quadratic twist for $m \in \set{4, 6, 8, 9, 12, 16, 18}$}

For $m \in \set{4, 6, 8, 9, 12, 16, 18}$, we list the first few twist minimal elliptic curves admitting a cyclic $m$-isogeny in Table \ref{table:firstfew4}, Table \ref{table:firstfew6}, Table \ref{table:firstfew8}, Table \ref{table:firstfew9}, Table \ref{table:firstfew12}, Table \ref{table:firstfew16}, and Table \ref{table:firstfew18} below.

\jvtable{table:firstfew4}{\small
\rowcolors{2}{white}{gray!10}
\begin{tabular}{c c c c}
$(A, B)$ & $(a, b)$ & $\twistheight(E)$ & $\twistdefect(E)$ \\
\hline\hline
$(1, 0)$ & $(0, 1)$ & $4$ & $3$ \\
$(-2, 1)$ & $(3, 1)$ or $(3, 5)$ & $32$ & $6$ \\
$(1, 2)$ & $(-3, 1)$ or $(3, 7)$ & $108$ & $6$ \\
$(6, 7)$ & $(-1, 1)$ or $(1, 3)$ & $1323$ & $2$ \\
$(-11, 14)$ & $(3, 2)$ or $(3, 4)$ & $5324$ & $3$ \\
$(-11, 6)$ & $(9, 1)$ or $(9, 17)$ & $5324$ & $6$ \\
$(13, 14)$ & $(-3, 5)$ or $(3, 11)$ & $8788$ & $6$ \\
$(-2, 21)$ & $(-9, 1)$ or $(9, 19)$ & $11907$ & $6$ \\
$(13, 34)$ & $(-3, 2)$ or $(3, 8)$ & $31212$ & $3$ \\
$(22, 23)$ & $(-3, 7)$ or $(3, 13)$ & $42592$ & $6$ \\
$(-23, 28)$ & $(6, 1)$ or $(6, 11)$ & $48668$ & $3$ \\
$(-23, 42)$ & $(9, 5)$ or $(9, 13)$ & $48668$ & $6$ \\
$(-26, 51)$ & $(9, 7)$ or $(9, 11)$ & $70304$ & $6$ \\
$(-26, 5)$ & $(15, 1)$ or $(15, 29)$ & $70304$ & $6$ \\
\end{tabular}
}{$E \in \twistE$ with a cyclic $4$-isogeny and $\twht E \leq 10^{5}$}

\jvtable{table:firstfew6}{\small
\rowcolors{2}{white}{gray!10}
\begin{tabular}{c c c c}
$(A, B)$ & $(a, b)$ & $\twistheight(E)$ & $\twistdefect(E)$ \\
\hline\hline
$(0, 1)$ & $(1, 1)$ & $27$ & $4$ \\
$(-12, 11)$ & $(9, 5)$ & $6912$ & $12$ \\
$(-15, 22)$ & $(2, 1)$ & $13500$ & $1$ \\
$(33, 74)$ & $(0, 1)$ & $147852$ & $3$ \\
$(60, 61)$ & $(-3, 1)$ & $864000$ & $12$ \\
$(93, 94)$ & $(-1, 1)$ & $3217428$ & $4$ \\
$(-120, 11)$ & $(21, 13)$ & $6912000$ & $12$ \\
$(-75, 506)$ & $(9, 7)$ & $6912972$ & $12$ \\
$(-132, 481)$ & $(11, 7)$ & $9199872$ & $4$ \\
$(-123, 598)$ & $(7, 5)$ & $9655308$ & $4$ \\
$(-255, 502)$ & $(12, 7)$ & $66325500$ & $3$ \\
$(-348, 2497)$ & $(5, 1)$ & $168576768$ & $4$ \\
$(-372, 2761)$ & $(7, 3)$ & $205915392$ & $4$ \\
$(-387, 766)$ & $(7, 4)$ & $231842412$ & $1$ \\
$(-408, 3107)$ & $(9, 1)$ & $271669248$ & $12$ \\
$(-423, 1342)$ & $(8, 5)$ & $302747868$ & $1$ \\
$(-327, 3454)$ & $(4, 3)$ & $322113132$ & $1$ \\
$(-435, 2162)$ & $(27, 17)$ & $329251500$ & $12$ \\
$(-372, 3611)$ & $(15, 11)$ & $352061667$ & $12$ \\
$(213, 3674)$ & $(3, 5)$ & $364453452$ & $12$ \\
$(-552, 4979)$ & $(15, 7)$ & $672786432$ & $12$ \\
\end{tabular}
}{$E \in \twistE$ with a cyclic $6$-isogeny and $\twht E \leq 10^{9}$}

\jvtable{table:firstfew8}{\small
\rowcolors{2}{white}{gray!10}
\begin{tabular}{c c c c}
$(A, B)$ & $(a, b)$ & $\twistheight(E)$ & $\twistdefect(E)$ \\
\hline\hline
$(6, 7)$ & $(0, 1)$ or $(4, 3)$ & $1323$ & $4$ \\
$(-3, 322)$ & $(2, 3)$ or $(6, 5)$ & $2799468$ & $4$ \\
$(-138, 623)$ & $(4, 1)$ or $(8, 5)$ & $10512288$ & $4$ \\
$(141, 142)$ & $(-2, 1)$ or $(10, 7)$ & $11212884$ & $4$ \\
$(-138, 2567)$ & $(4, 5)$ or $(8, 7)$ & $177916203$ & $4$ \\
$(141, 4718)$ & $(1, 2)$ or $(5, 4)$ & $601007148$ & $1$ \\
$(-579, 5362)$ & $(3, 1)$ or $(5, 3)$ & $776418156$ & $1$ \\
$(582, 4417)$ & $(-4, 1)$ or $(16, 11)$ & $788549472$ & $4$ \\
\end{tabular}
}{$E \in \twistE$ with a cyclic $8$-isogeny and $\twht E \leq 10^{9}$}

\jvtable{table:firstfew9}{\small
\rowcolors{2}{white}{gray!10}
\begin{tabular}{c c c c}
$(A, B)$ & $(a, b)$ & $\twistheight(E)$ & $\twistdefect(E)$ \\
\hline\hline
$(0, 2)$ & $(0, 1)$ & $108$ & $2$ \\
$(24, 2)$ & $(2, 1)$ & $55296$ & $6$ \\
$(-48, 142)$ & $(-2, 5)$ & $544428$ & $6$ \\
$(-51, 142)$ & $(-1, 2)$ & $544428$ & $1$ \\
$(69, 362)$ & $(1, 2)$ & $3538188$ & $3$ \\
$(-96, 362)$ & $(-2, 1)$ & $3538944$ & $2$ \\
$(-120, 502)$ & $(-4, 1)$ & $6912000$ & $6$ \\
$(-120, 506)$ & $(-2, 3)$ & $6912972$ & $2$ \\
$(-75, 506)$ & $(-1, 4)$ & $6912972$ & $3$ \\
$(165, 502)$ & $(1, 1)$ & $17968500$ & $1$ \\
$(-264, 1654)$ & $(-4, 7)$ & $73864332$ & $6$ \\
$(-219, 1654)$ & $(-1, 3)$ & $73864332$ & $1$ \\
\end{tabular}
}{$E \in \twistE$ with a cyclic $9$-isogeny and $\twht E \leq 10^{9}$}

\jvtable{table:firstfew12}{\small
\rowcolors{2}{white}{gray!10}
\begin{tabular}{c c c c}
$(A, B)$ & $(a, b)$ & $\twistheight(E)$ & $\twistdefect(E)$ \\
\hline\hline
$(213, 3674)$ & $(-9, 5)$ or $(3, 1)$ & $364453452$ & $108$ \\
$(-1947, 108214)$ & $(-5, 3)$ or $(1, 1)$ & $316177284492$ & $4$ \\
$(-5907, 61486)$ & $(-9, 7)$ or $(-3, 5)$ & $824443510572$ & $108$ \\
$(-9867, 324934)$ & $(-6, 5)$ or $(-3, 4)$ & $3842513269452$ & $27$ \\
$(-41547, 3259514)$ & $(-9, 4)$ or $(-6, 1)$ & $286865949497292$ & $27$ \\
$(-65307, 874294)$ & $(-4, 3)$ or $(-1, 2)$ & $1114138529957772$ & $1$ \\
$(-71643, 7378058)$ & $(-15, 7)$ or $(-9, 1)$ & $1470893677938828$ & $108$ \\
$(-11667, 11349074)$ & $(-12, 7)$ or $(3, 2)$ & $3477639977751852$ & $27$ \\
$(-168843, 12140858)$ & $(-15, 11)$ or $(-3, 7)$ & $19253477048692428$ & $108$ \\
$(-212043, 28562182)$ & $(-7, 5)$ or $(-1, 3)$ & $38135707808174028$ & $4$ \\
$(228813, 5274866)$ & $(-21, 11)$ or $(9, 1)$ & $47918374464655188$ & $108$ \\
$(276837, 35589962)$ & $(-15, 8)$ or $(6, 1)$ & $84865738374033012$ & $27$ \\
$(-386067, 92329774)$ & $(-7, 3)$ or $(-5, 1)$ & $230169637578251052$ & $4$ \\
$(-627483, 187952182)$ & $(-9, 8)$ or $(-6, 7)$ & $988247863401150348$ & $27$ \\
\end{tabular}
}{$E \in \twistE$ with a cyclic $12$-isogeny and $\twht E \leq 10^{18}$}

\jvtable{table:firstfew16}{\small
\rowcolors{2}{white}{gray!10}
\begin{tabular}{c c c c}
$(A, B)$ & $(a, b)$ & $\twistheight(E)$ & $\twistdefect(E)$ \\
\hline\hline
$(-3, 322)$ & $(0, 1)$ or $(4, 3)$ & $2799468$ & $4$ \\
$(-11523, 476098)$ & $(3, 1)$ or $(5, 3)$ & $6120074050668$ & $1$ \\
$(-11523, 584962)$ & $(2, 3)$ or $(6, 5)$ & $9238874618988$ & $4$ \\
$(-15843, 767522)$ & $(4, 1)$ or $(8, 5)$ & $15906413128428$ & $4$ \\
$(-15843, 1441118)$ & $(1, 2)$ or $(5, 4)$ & $56074169427948$ & $1$ \\
$(30237, 1524962)$ & $(-2, 1)$ or $(10, 7)$ & $110579874088212$ & $4$ \\
$(30237, 3904418)$ & $(-1, 1)$ or $(7, 5)$ & $411600957805548$ & $1$ \\
$(-311043, 66769598)$ & $(5, 2)$ or $(7, 4)$ & $120370838936786028$ & $1$ \\
$(-311043, 69595202)$ & $(4, 5)$ or $(8, 7)$ & $130774287818361708$ & $4$ \\
\end{tabular}
}{$E \in \twistE$ with a cyclic $16$-isogeny and $\twht E \leq 10^{18}$}

\jvtable{table:firstfew18}{\small
\rowcolors{2}{white}{gray!10}
\begin{tabular}{c c c c}
$(A, B)$ & $(a, b)$ & $\twistheight(E)$ & $\twistdefect(E)$ \\
\hline\hline
$(-75, 506)$ & $(-4, 1)$ & $6912972$ & $3888$ \\
$(-1515, 22682)$ & $(1, 2)$ & $13909063500$ & $243$ \\
$(-24555, 1485286)$ & $(-2, 1)$ & $59564011548492$ & $16$ \\
$(-172227, 27405506)$ & $(-2, 5)$ & $20434485218644332$ & $3888$ \\
$(-393195, 94898662)$ & $(1, 1)$ & $243155418015559500$ & $1$ \\
$(-1123275, 458221178)$ & $(8, 1)$ & $5669154212905687500$ & $3888$ \\
$(-1143435, 440919866)$ & $(-1, 4)$ & $5979907087519351500$ & $243$ \\
$(1324653, 1127890514)$ & $(-10, 1)$ & $34347699312421973292$ & $3888$ \\
$(-2065467, 1142549354)$ & $(4, 5)$ & $35246400621552810252$ & $3888$ \\
$(-2046747, 1164275786)$ & $(-5, 2)$ & $36599528858379780492$ & $243$ \\
$(-2752707, 1757875394)$ & $(5, 1)$ & $83433402148362948972$ & $243$ \\
$(-4906515, 4040728274)$ & $(2, 7)$ & $472475598685352563500$ & $3888$ \\
$(2097645, 14640824018)$ & $(-7, 1)$ & $5787550654003232936748$ & $243$ \\
$(-11577747, 12814884434)$ & $(1, 5)$ & $6207720523046065646892$ & $243$ \\
$(-24565035, 46858172762)$ & $(-4, 7)$ & $59294191693335605671500$ & $3888$ \\
$(-26453307, 52372560746)$ & $(-8, 5)$ & $74057898215523422065932$ & $3888$ \\
$(-76172547, 255885590014)$ & $(4, 1)$ & $1767890750702694206045292$ & $16$ \\
$(-76182627, 255814479646)$ & $(-1, 2)$ & $1768592684578826203703532$ & $1$ \\
$(-88080555, 318174471718)$ & $(-2, 3)$ & $2733380669628994787815500$ & $16$ \\
$(-122727387, 523299579766)$ & $(2, 3)$ & $7394085267179261308598412$ & $16$ \\
$(-117935067, 566044198774)$ & $(-3, 1)$ & $8650962944073889823783052$ & $1$ \\
\end{tabular}
}{$E \in \twistE$ with a cyclic $18$-isogeny and $\twht E \leq 10^{25}$}

\subsection*{Computing \texorpdfstring{$\twconsteq m$ and $\twconstad m$}{tctwm and ctwm}}

In this subsection, we estimate the constants appearing in \Cref{Theorem: asymptotic for twN(X) for m of genus $0$} and \Cref{Corollary: asymptotic for twN(X) for m of genus $0$}. The constants $\Q m$ were computed already in Table \ref{table:Qm}.

We can estimate the area $\R m$ of the region $\calR m$ by performing rejection sampling on an enveloping rectangle. This computation proceeds, \textit{mutatis mutandis}, in exactly the same way that the analogous computations for $m \in \set{7, 10, 25}$ proceeded. We record these computations in Table \ref{table:Rm} below, along with how many CPU days it took to complete them.

\jvtable{table:Rm}{
\rowcolors{2}{white}{gray!10}
\begin{tabular}{c | c c c c}
$m$ & Enveloping region & $\#$ of trials & $\#$ of successes & CPU days \\
\hline\hline
$4$ & $[-0.4583, 0.4583] \times [0, 0.9166]$ & $319\,525\,000\,000$ & $65\,681\,836\,724$ & $24$ \\ %
$6$ & $[-0.677, 1.7036] \times [0, 1.0338]$ & $331\,210\,000\,000$ &  $55\,726\,701\,265$ & $34$ \\ %
$8$ & $[-0.677, 2.0309] \times [0, 1.3539]$ & $326\,270\,000\,000$ &  $42\,861\,204\,516$ & $34$ \\ %
$9$ & $[-0.677, 0.677] \times [0, 0.6801]$ & $316\,167\,000\,000$ & $131\,922\,088\,793$ & $34$ \\ %
$12$ & $[-1.6456, 0.8228] \times [0, 0.8228]$ & $224\,310\,000\,000$ & $69\,681\,481\,937$ & $34$ \\ %
$16$ & $[-0.8228, 2.4684] \times [0, 1.6456]$ & $110\,452\,000\,000$ & $19\,715\,984\,750$ & $14$ \\ %
$18$ & $[-0.8781, 0.8781] \times [0, 0.5532]$ & $130\,060\,000\,000$ & $61\,569\,780\,450$ & $23$ %
\end{tabular}
}{Approximating $\R m$ for $m \in \set{4, 6, 8, 9, 12, 16, 18}$}

Recall \eqref{Equation: twconst for m of genus $0$} and \eqref{Equation: twconsteq for m of genus $0$}. We have assembled everything we need to compute $\twconsteq m$ and $\twconstad m$. We report each constituent factor of $\twconsteq m$ in Table \ref{table:ingredients}. The quantity $\R m$, along with its standard error from the rejection sampling recorded in Table \ref{table:Rm}, is estimated there as well.

In Table \ref{table:twconsteq}, we then record $\twconsteq m$, with its propagated standard error, as well as the ratio
\begin{equation}
\frac{ \twistNeq {m}(X)}{\twconsteq {m}X^{1/\degB m}}
\end{equation}
for $m$ and $X$ as in Table \ref{table:enumerating}. In Table \ref{table:twconstad}, we perform analogous computations for $\twconstad m$.

\jvtable{table:ingredients}{
\rowcolors{2}{white}{gray!10}
\begin{tabular}{c | c c c c}
$m$ & $\delta_m$ & $\Q m$ & $\R m$ & Standard error \\
\hline
\hline
$4$ & $2$ & $6$ & $0.172\,703\,107$ & $6.1 \cdot 10^{-7}$ \\
$6$ & $1$ & $3$ & $0.414\,078\,663$ & $1.6 \cdot 10^{-6}$ \\
$8$ & $2$ & $2$ & $0.481\,622\,136$ & $2.2 \cdot 10^{-6}$ \\
$9$ & $1$ & $2$ & $0.384\,231\,017$ & $8.1 \cdot 10^{-7}$ \\
$12$ & $2$ & $1 + \sqrt 3$ & $0.630\,926\,202$ & $2.0 \cdot 10^{-6}$ \\
$16$ & $2$ & $4/3$ & $0.966\,770\,617$ & $6.3 \cdot 10^{-6}$  \\
$18$ & $1$ & $(1 + 2^{1/3})(1 + 3^{2/3})/2$ & $0.459\,917\,568$ & $1.4 \cdot 10^{-6}$ 
\end{tabular}}{Ingredients to compute $\twconsteq m$ and $\twconstad m$ for $m \in \set{4, 6, 8, 9, 12, 16, 18}$}

\jvtable{table:twconsteq}{
\rowcolors{2}{white}{gray!10}
\begin{tabular}{c | c c c}
$m$ & $\twconsteq m$ & Standard error & $\twistNeq {m}(X)/\twconsteq {m}X^{1/\degB m}$ \\
\hline
\hline
$4$ & $0.629\,945\,396$ & $2.2 \cdot 10^{-6}$ & $0.999\,999\,688\ldots$ \\
$6$ & $0.755\,188\,924$ & $3.0 \cdot 10^{-6}$ & $1.000\,009\,767\ldots$  \\
$8$ & $0.585\,582\,298$ & $2.7 \cdot 10^{-6}$ & $1.000\,028\,863\ldots$ \\
$9$ & $0.467\,168\,897$ & $9.9 \cdot 10^{-7}$ & $0.999\,947\,990\ldots$ \\
$12$ & $1.047\,897\,587$ & $3.3 \cdot 10^{-6}$ & $0.999\,999\,630\ldots$ \\
$16$ & $0.783\,634\,746$ & $5.1 \cdot 10^{-6}$ & $0.999\,963\,062\ldots$ \\
$18$ & $0.973\,099\,640$ & $2.9 \cdot 10^{-6}$ & $0.999\,961\,832\ldots$
\end{tabular}}{The constant $\twconsteq m$, its standard error, and a related ratio}

\jvtable{table:twconstad}{
\rowcolors{2}{white}{gray!10}
\begin{tabular}{c | c c c}
$m$ & $\twconstad m$ & Standard error & $\twistNad {m}(X)/\twconstad {m}X^{1/\degB m}$ \\
\hline
\hline
$4$ & $0.314\,972\,698$ & $1.1 \cdot 10^{-6}$ & $0.999\,997\,783\ldots$ \\
$6$ & $0.755\,188\,924$ & $3.0 \cdot 10^{-6}$ & $0.999\,842\,524\ldots$  \\
$8$ & $0.292\,791\,149$ & $1.4 \cdot 10^{-6}$ & $0.999\,930\,158\ldots$ \\
$9$ & $0.467\,168\,897$ & $9.9 \cdot 10^{-7}$ & $0.999\,947\,990\ldots$ \\
$12$ & $0.523\,948\,794$ & $1.7 \cdot 10^{-6}$ & $0.999\,999\,630\ldots$ \\
$16$ & $0.391\,817\,373$ & $2.6 \cdot 10^{-6}$ & $0.999\,963\,062\ldots$ \\
$18$ & $0.973\,099\,640$ & $2.9 \cdot 10^{-6}$ & $0.999\,961\,832\ldots$
\end{tabular}}{The constant $\twconstad m$, its standard error, and a related ratio}

\subsection*{Computing \texorpdfstring{$\consteq {m}$, $\consteqprm m$, $\constad m$, $\constadprm m$ for $m \in \set{6, 8, 9}$}{tcm, tcm' and cm, cm' for m in 6, 8, 9}}

In this subsection, for $m \in \set{6, 8, 9}$, we estimate the constants $\consteq {m}$, $\constad {m}$, $\consteqprm m$ and $\constadprm m$, 
 which are defined in \eqref{Equation: constants for 5 < m <= 9} and used in \Cref{Theorem: asymptotic for NQ(X) for 5 < m <= 9}.  Of course, $\consteq {9} = \constad {9}$ and $\consteqprm 9 = \constadprm 9$ by \Cref{Lemma: Difference between twNeq and twNad}. 

We follow the same strategy as in \cref{Section: Computations for m = 7}. For $m \in \set{6, 8, 9}$, we have the identities $\consteq m = \twconsteq m / (3 \zeta(2))$ and $\constad m = \twconstad m / (3 \zeta(2))$, whence we obtain the following table.

\jvtable{table:constform<=9}{
\rowcolors{2}{white}{gray!10}
\begin{tabular}{c | c c c c c}
$m$ & $\consteq m$ & Standard error & $\constad m$ & Standard error \\
\hline
\hline
$6$ & $0.153\,033\,271$ & $6.1 \cdot 10^{-7}$ & $0.153\,033\,271$ & $6.1 \cdot 10^{-7}$ \\
$8$ & $0.118\,663\,783$ & $5.5 \cdot 10^{-7}$ & $0.059\,331\,892$ & $2.9 \cdot 10^{-7}$\\
$9$ & $0.094\,668\,211$ & $2.1 \cdot 10^{-7}$ & $0.094\,668\,211$ & $2.1 \cdot 10^{-7}$
\end{tabular}}{The constants $\consteq m$ and $\constad m$ and their standard errors for $m \in \set{6, 8, 9}$}

We can approximate $\widetilde{\ell}_{m, 0}$ and $\ell_{m, 0}$ by truncating the integrals \eqref{equation: Definition of elleq0 for m of genus $0$} and \eqref{equation: Definition of ell0 for m of genus $0$} and using our approximations for $\twconsteq m$ and $\twconstad m$. Truncating this integral at the $X$ recorded in Table \ref{table:enumerating} yields estimates for $\widetilde{\ell}_{m, 0}$ and $\ell_{m, 0}$ which we record in Table \ref{table:ell}. 

We now assess the error in these estimates. In \Cref{Theorem: asymptotic for twN(X) for m of genus $0$}, we have shown that for some $M > 0$ and for all $u > X$, we have
\begin{equation}
\abs{\twistNad m(u) - \twconst m \floor{u^{1/6}}} < M u^{1/12} \log u.
\end{equation}
Thus
\begin{equation}
\begin{aligned}
&\abs{\int_X^\infty \parent{\twistNad m(u) - \twconst m \floor{u^{1/6}}} u^{-7/6} \,\mathrm{d}u} \\
&\qquad\qquad < M \int_X^\infty u^{-13/12} \log u\,\mathrm{d}u \\
&\qquad\qquad = 12 M X^{-1/12} (\log X + 12);
\end{aligned}
\end{equation}
this gives us a bound on our truncation error. 
We do not know the exact value for $M$, but empirically, we find that throughout the range $10^{19} \leq u \leq 10^{42}$, we have
\begin{equation}
\begin{aligned}
    -1.75 \cdot 10^{-1} &\leq \frac{\twistNeq 6(u) - \twconsteq 6 \floor{u^{1/6}}}{u^{1/12} \log u} \leq 6.52 \cdot 10^{-3}, \\
    -1.75 \cdot 10^{-1} &\leq \frac{\twistNad 6(u) - \twconstad 6 \floor{u^{1/6}}}{u^{1/12} \log u} \leq 0, \\
    -1.35 \cdot 10^{-1} &\leq \frac{\twistNeq 8(u) - \twconsteq 8 \floor{u^{1/6}}}{u^{1/12} \log u} \leq 9.80 \cdot 10^{-3}, \\
    -6.72 \cdot 10^{-2} &\leq \frac{\twistNad 8(u) - \twconstad 8 \floor{u^{1/6}}}{u^{1/12} \log u} \leq 4.59 \cdot 10^{-3}, \\
    -1.36 \cdot 10^{-1} &\leq \frac{\twistNeq 9(u) - \twconsteq 9 \floor{u^{1/6}}}{u^{1/12} \log u} \leq 3.30 \cdot 10^{-2}, \\
    -1.36 \cdot 10^{-1} &\leq \frac{\twistNad 9(u) - \twconstad 9 \floor{u^{1/6}}}{u^{1/12} \log u} \leq 3.30 \cdot 10^{-2}, \\
\end{aligned}
\end{equation}
If we assume these bounds continue to hold for larger $u$, and combine the resulting truncation error with the standard errors in our approximations for $\twconsteq m$ and $\twconstad m$, we obtain the uncertainties in Table \ref{table:ell} and Table \ref{table:constm}.

\jvtable{table:ell}{
\rowcolors{2}{white}{gray!10}
\begin{tabular}{c | c c c c}
$m$ & $\widetilde{\ell}_{m, 0}$ & Uncertainty & $\ell_{m, 0}$ & Uncertainty \\
\hline
\hline
$6$ & $-0.572\,182\,786$ & $1.21 \cdot 10^{-2}$ & $-0.677\,298\,843$ & $1.21 \cdot 10^{-2}$ \\
$8$ & $-0.742\,599\,375$ & $9.33 \cdot 10^{-3}$ & $-0.387\,595\,153$ & $4.65 \cdot 10^{-3}$ \\
$9$ & $-0.231\,206\,039$ & $9.37 \cdot 10^{-3}$ & $-0.231\,206\,039$ & $9.37 \cdot 10^{-3}$ \\
\end{tabular}}{The constants $\widetilde{\ell}_{m, 0}$ and $\ell_{m, 0}$ and their numerical uncertainties for $m \in \set{6, 8, 9}$}

\jvtable{table:constm}{
\rowcolors{2}{white}{gray!10}
\begin{tabular}{c | c c c c}
$m$ & $\consteqprm m$ & Uncertainty & $\constadprm m$ & Uncertainty \\
\hline
\hline
$6$ & $-0.037\,215\,321$ & $1.47 \cdot 10^{-2}$ & $-0.165\,021\,121$ & $1.47 \cdot 10^{-2}$ \\
$8$ & $-0.392\,302\,971$ & $1.14 \cdot 10^{-2}$ & $-0.215\,964\,396$ & $5.65 \cdot 10^{-3}$ \\
$9$ & $0.126\,227\,997$ & $1.14 \cdot 10^{-2}$ & $0.126\,227\,997$ & $1.14 \cdot 10^{-2}$ \\
\end{tabular}}{The constants $\consteqprm m$ and $\constadprm m$ and their numerical uncertainties for $m \in \set{6, 8, 9}$}

As a sanity check, we now compute $\NQeq m (X)$ and $\NQad m (X)$ for $X = 10^{42}$, and verify 
\begin{equation}
\frac{\NQeq m(X)}{X^{1/6}} - \consteq m \log X \approx \consteqprm m,
\end{equation}
and
\begin{equation}
\frac{\NQad m(X)}{X^{1/6}} - \constad m \log X \approx \constadprm m,
\end{equation}
as shown in Table \ref{table:enumeratingoverQ} below.

\jvtable{table:enumeratingoverQ}{
\rowcolors{2}{white}{gray!10}
\begin{tabular}{c | c c c c c}
$m$ & $X$ & $\NQeq m (X)$ & $\frac{\NQeq m(X)}{X^{1/6}} - \consteq m \log X$ & $\NQad m (X)$ & $\frac{\NQad m(X)}{X^{1/6}} - \constad m \log X$ \\
\hline
\hline
$6$ & $10^{42}$ & $147\,624\,808$ & $-0.037\,148\,635$ & $146\,343\,642$ & $-0.165\,265\,235$ \\ %
$8$ & $10^{42}$ & $110\,837\,024$ & $-0.392\,102\,846$ & $55\,219\,676$ & $-0.215\,935\,023$ \\ %
$9$ & $10^{42}$ & $92\,813\,182$ & $0.126\,090\,504$ & $92\,813\,182$ & $0.126\,090\,504$ \\ %
\end{tabular}}{Enumerating elliptic curves up to $\bbQ$-isomorphism for $m \in \set{6, 8, 9}$}

We emphasize that for $m \in \set{6, 8, 9}$, the estimates in Table \ref{table:enumeratingoverQ} for $\consteqprm m$ and $\constadprm m$ depend on empirical rather than theoretical estimates for the implicit constant in the error term in the asymptotics of $\twistNeq m(X)$ and $\twistNad m (X)$.

\subsection*{Computing \texorpdfstring{$\consteq {m}$ and $\constad {m}$ for $m \in \set{12, 16, 18}$}{tcm and cm for m in 12, 16, 18}}

In this subsection, for $m \in \set{12, 16, 18}$, we estimate 
\begin{equation}
\consteq {m} = 2 \twistLeq {m}(1/6)/\zeta(2) \ \text{and} \ \constad {m} = 2 \twistLad {m}(1/6)/\zeta(2),
\end{equation}
the constants which appear in \Cref{Theorem: asymptotic for NQ(X) for m > 9}. Of course, $\consteq {18} = \constad {18}$ by \Cref{Lemma: Difference between twNeq and twNad}. 

We follow the same strategy as in \cref{Section: Computations for m = 10 and 25}. We first compute the partial sums
\begin{equation}
\begin{aligned}
	\sum_{n \leq 10^{84}} \frac{\twistheq {12}(n)}{n^{1/6}} &= 0.212\,842\,775\,719\,150\,50, \\
 \sum_{n \leq 10^{84}} \frac{\twisthad {12}(n)}{n^{1/6}} &= 0.106\,421\,387\,859\,575\,25, \\
 \sum_{n \leq 10^{84}} \frac{\twistheq {16}(n)}{n^{1/6}} &= 0.269\,169\,745\,679\,529\,87, \\
 \sum_{n \leq 10^{84}} \frac{\twisthad {16}(n)}{n^{1/6}} &= 0.134\,584\,872\,839\,764\,93, \ \text{and} \\
\sum_{n \leq 10^{126}} \frac{\twisthad {18}(n)}{n^{1/6}} &= 0.107\,025\,809\,032\,944\,72.
 \end{aligned}
\end{equation}
We empirically confirm that 
\begin{equation}
\begin{aligned}
\twistNeq {12} (X) &< 1.090\,472 X^{1/12} \ \text{for} \ X \leq 10^{84}, \\
\twistNad {12} (X) &< 0.545\,236 X^{1/12} \ \text{for} \ X \leq 10^{84}, \\
\twistNeq {16} (X) &< 0.847\,726 X^{1/12} \ \text{for} \ X \leq 10^{84}, \\
\twistNad {16} (X) &< 0.423\,863 X^{1/12} \ \text{for} \ X \leq 10^{84}, \ \text{and} \\
\twistNad {18} (X) &< 1.007\,095 X^{1/18} \ \text{for} \ X \leq 10^{126}.
 \end{aligned}
\end{equation}
If these bounds continue to hold for larger $X$, then partial summation, using the endpoint values from Table \ref{table:enumerating}, gives
\begin{equation}\label{Equation: tail of twL for m > 9}
\begin{aligned}
\sum_{n > 10^{84}} \frac{\twistheq {12}(n)}{n^{1/6}} &< 1.133\,046\,8 \cdot 10^{-7}, \\
\sum_{n > 10^{84}} \frac{\twisthad {12}(n)}{n^{1/6}} &< 5.665\,234 \cdot 10^{-8}, \\
\sum_{n > 10^{84}} \frac{\twistheq {16}(n)}{n^{1/6}} &< 9.118\,462 \cdot 10^{-8}, \\
\sum_{n > 10^{84}} \frac{\twisthad {16}(n)}{n^{1/6}} &< 4.559\,231 \cdot 10^{-8}, \ \text{and} \\
\sum_{n > 10^{126}} \frac{\twisthad {18}(n)}{n^{1/6}} &< 5.375\,8 \cdot 10^{-15}. \\
\end{aligned}
\end{equation}
Assuming \eqref{Equation: tail of twL for m > 9}, for $m \in \set{12, 16, 18}$, we therefore have the following estimates for $\consteq m$ and $\constad m$.

\jvtable{table:constform>9}{
\rowcolors{2}{white}{gray!10}
\begin{tabular}{c | c c c c c}
$m$ & $\consteq m$ & Error bound & $\constad m$ & Error bound \\
\hline
\hline
$12$ & $0.258\,785\,784$ & $1.4 \cdot 10^{-7}$ & $0.129\,392\,892$ & $6.9 \cdot 10^{-8}$ \\
$16$ & $0.327\,271\,167$ & $1.2 \cdot 10^{-7}$ & $0.163\,635\,583$ & $5.6 \cdot 10^{-8}$ \\
$18$ & $0.130\,127\,779\,817\,961\,2$ & $6.6 \cdot 10^{-15}$ & $0.130\,127\,779\,817\,961\,2$ & $6.6 \cdot 10^{-15}$
\end{tabular}}{The constants $\consteq m$ and $\constad m$ and their conditional error bounds for $m \in \set{12, 16, 18}$}

We emphasize that for $m \in \set{12, 16, 18}$, the estimates in Table \ref{table:constform>9} for $\consteq m$ and $\constad m$ depend on empirical rather than theoretical estimates for the implicit constant in the error term in the asymptotics of $\twistNeq m(X)$ and $\twistNad m (X)$.

%% file: Ch8_m11_14_15_17_19_21_27_37_43_67_163.tex
\chapter{Counting elliptic curves with a cyclic \texorpdfstring{$m$}{m}-isogeny when \texorpdfstring{$X_0(m)$}{X0(m)} has nonzero genus}\label{Chapter: m of nonzero genus}

For completeness, we prove \Cref{Intro Theorem: asymptotic for NQ(X) for m of nonzero genus} (\Cref{Theorem: asymptotic for NQ(X) for m of nonzero genus} and \Cref{Theorem: asymptotic for NQ(X) for m of nonzero genus given RH}) and \Cref{Intro Theorem: asymptotic for twN(X) for m of nonzero genus} (\Cref{Theorem: asymptotic for twN(X) for m of nonzero genus}), giving asymptotic counts for the number of elliptic curves over $\bbQ$ admitting (equivalently, equipped with) a cyclic $m$-isogeny when
\[
m \in \set{11, 14, 15, 17, 19, 21, 27, 37, 43, 67, 163}.
\]
In each of these cases, the compactified moduli space $X_0(m)$ has genus $g > 0$ and only finitely many rational points. These points were enumerated classically \cite{Mazur2}, and summing over their quadratic twists gives us our desired asymptotics. Although we are unaware of any reference for \Cref{Intro Theorem: asymptotic for NQ(X) for m of nonzero genus} in the literature, we expect that the claims of this chapter are familiar to experts in the theory of arithmetic statistics.

Throughout the remainder of this chapter, we assume 
\[
m \in \set{11, 14, 15, 17, 19, 21, 27, 37, 43, 67, 163}.
\]
By \Cref{Corollary: Elliptic curves for which twNeq(X) = twNad(X)}, we have
\begin{equation}
\twNeq m (X) = \twNad m (X) \ \text{and} \ \NQeq m (X) = \NQad m (X)
\end{equation}
for all $X > 0$, so we may use either notation interchangeably. We opt to work with $\twNad m (X)$ and related functions.

In \cref{Section: Counts for twist classes}, we record all elliptic curves admitting a cyclic $m$-isogeny up to twist equivalence, write $\twistNad m (X)$ in terms of the Heaviside step function, and point out its long-run asymptotics. In \cref{Section: Working over the rationals for m of nonzero genus}, we leverage Walfisz's and Liu's asymptotics for the count of squarefree integers to give asymptotics for the number of elliptic curves with a cyclic $m$-isogeny up to $\bbQ$-isomorphism.

\section{Counts for twist classes when \texorpdfstring{$X_0(m)$}{X0(m)} has nonzero genus}\label{Section: Counts for twist classes}

In this section, we write down all elliptic curves admitting a cyclic $m$-isogeny, up to quadratic twist, and use this explicit enumeration to describe $\twistNad m (X)$ for 
\[
m \in \set{11, 14, 15, 17, 19, 21, 27, 37, 43, 67, 163}.
\]

The twist classes and their $j$-invariants are classical; see \cite[page 1]{Mazur2} and the references there. We use twist-minimal short Weierstrass representatives with $B>0$.

\jvtable{table:nonzerogenus}{
\rowcolors{2}{white}{gray!10}
{\footnotesize%
\begin{tabular}{c c|c c}
$m$ & $(A, B)$ & $\twistheight(E)$ & $j(E)$ \\
\hline\hline
$11$ & $(-264, 1\,694)$ & $77\,480\,172$ & $-32768$ \\
$11$ & $(-363, 10\,406)$ & $2\,923\,690\,572$ & $-121$ \\
$11$ & $(-4\,323, 109\,406)$ & $323\,181\,166\,572$ & $-24729001$ \\
\hline
$14$ & $(-35, 98)$ & $259\,308$ & $-3375$ \\
$14$ & $(-595, 5\,586)$ & $842\,579\,500$ & $16581375$ \\
\hline
$15$ & $(-75, 2\,950)$ & $234\,967\,500$ & $-25/2$ \\
$15$ & $(-435, 4\,210)$ & $478\,550\,700$ & $-121945/32$ \\
$15$ & $(3\,165, 31\,070)$ & $126\,818\,068\,500$ & $46969655/32768$ \\
$15$ & $(-18\,075, 935\,350)$ & $23\,621\,749\,807\,500$ & $-349938025/8$ \\
\hline
$17$ & $(-95\,115, 12\,657\,350)$ & $4\,325\,629\,743\,607\,500$ & $-882216989/131072$ \\
$17$ & $(-437\,835, 111\,510\,650)$ & $335\,734\,876\,712\,407\,500$ & $-297756989/2$ \\
\hline
$19$ & $(-152, 722)$ & $14\,074\,668$ & $-884736$ \\
\hline
$21$ & $(45, 18)$ & $364\,500$ & $3375/2$ \\
$21$ & $(-75, 262)$ & $1\,853\,388$ & $-140625/8$ \\
$21$ & $(-1\,515, 46\,106)$ & $57\,395\,607\,372$ & $-1159088625/2097152$ \\
$21$ & $(-17\,235, 870\,894)$ & $20\,478\,321\,699\,372$ & $-189613868625/128$ \\
\hline
$27$ & $(-120, 506)$ & $6\,912\,972$ & $-12288000$ \\
\hline
$37$ & $(-1\,155, 16\,450)$ & $7\,306\,267\,500$ & $-9317$ \\
$37$ & $(-29\,963\,955, 63\,131\,603\,150)$ & $107\,611\,181\,539\,805\,427\,907\,500$ & $-162677523113838677$ \\
\hline
$43$ & $(-3\,440, 77\,658)$ & $162\,830\,654\,028$ & $-884736000$ \\
\hline
$67$ & $(-29\,480, 1\,948\,226)$ & $102\,480\,782\,771\,052$ & $-147197952000$ \\
\hline
$163$ & $(-8\,697\,680, 9\,873\,093\,538)$ & $2\,631\,905\,352\,272\,628\,650\,988$ & $-262537412640768000$
\end{tabular}
}}{$E \in \twistE$ with a cyclic $m$-isogeny when $X_0(m)$ has nonzero genus}

By Mazur's theorem on isogenies (\Cref{Theorem: Kenku-Mazer Theorem}), Table \ref{table:nonzerogenus} is an exhaustive list of elliptic curves admitting a cyclic $m$-isogeny over $\bbQ$, up to quadratic twist, for 
\begin{equation}
m \not\in \set{1, 2, 3, 4, 5, 6, 7, 8, 9, 10, 12, 13, 16, 18, 25}.
\end{equation}
We emphasize that for each $m$ in Table \ref{table:nonzerogenus}, each elliptic curve associated to this $m$ admits exactly one unsigned cyclic $m$-isogeny.

Recall that the \defi{Heaviside step function} $\heaviside : \bbR \to \bbR$ is given by
\begin{equation}\label{Equation: heaviside step function}
	\heaviside(X) \coloneqq \begin{cases} 1 & \text{if} \ X \geq 0, \\
	0 & \text{if} \ X < 0.
	\end{cases}
\end{equation}

The following lemma is immediate from Table \ref{table:nonzerogenus}.

\begin{lemma}\label{Lemma: twN(X) as a sum of heaviside functions for m of nonzero genus}
	Let $X > 0$ be arbitrary. We have the following identities:
		\begin{equation} 	\begin{aligned}
		\twistNad {11} (X) =& \heaviside(X - 77\,480\,172) + \heaviside(X - 2\,923\,690\,572) + \heaviside(X - 323\,181\,166\,572), \\
		\twistNad {14} (X) =& \heaviside(X - 259\,308) + \heaviside(X - 842\,579\,500), \\
		\twistNad {15}(X) =& \heaviside(X - 234\,967\,500) + \heaviside(X - 478\,550\,700) \\
		&+ \heaviside(X - 126\,818\,068\,500) + \heaviside(X - 23\,621\,749\,807\,500), \\ 
		\twistNad {17} (X) =& \heaviside(X - 4\,325\,629\,743\,607\,500) + \heaviside(X - 335\,734\,876\,712\,407\,500), \\
		\twistNad {19} (X) =& \heaviside(X - 14\,074\,668), \\
		\twistNad {21}(X) =& \heaviside(X - 364\,500) + \heaviside(X - 1\,853\,388)\\
		&+ \heaviside(X - 57\,395\,607\,372) + \heaviside(X - 20\,478\,321\,699\,372), \\ 
		\twistNad {27} (X) =& \heaviside(X - 6\,912\,972), \\		
		\twistNad {37} (X) =& \heaviside(X - 7\,306\,267\,500) + \heaviside(X - 107\,611\,181\,539\,805\,427\,907\,500), \\				
		\twistNad {43} (X) =& \heaviside(X - 162\,830\,654\,028), \\		
		\twistNad {67} (X) =& \heaviside(X - 102\,480\,782\,771\,052), \\		
		\twistNad {163} (X) =& \heaviside(X - 2\,631\,905\,352\,272\,628\,650\,988).
		\end{aligned}\end{equation} 	
\end{lemma}

\begin{proof}
    We inspect the twist heights in Table \ref{table:nonzerogenus}; as this table is exhaustive, the lemma follows.
\end{proof}

We recover \Cref{Intro Theorem: asymptotic for twN(X) for m of nonzero genus}, which reports the asymptotic behavior of $\twistNad m (X)$, from \Cref{Lemma: twN(X) as a sum of heaviside functions for m of nonzero genus}.

\begin{theorem}\label{Theorem: asymptotic for twN(X) for m of nonzero genus}
	For $X$ sufficiently large, we have the following identities:
		\begin{equation} 
  \begin{aligned}
		\twistNad {11} (X) =& 3, \ \twistNad {14} (X) = 2, \ \twistNad {15}(X) = 4, \ \twistNad {17} (X) = 2, \\
    \twistNad {19} (X) =& 1, \twistNad {21}(X) = 4, \twistNad {27} (X) = 1, \ \twistNad {37} (X) = 2, \\ 
    \twistNad {43} (X) =& 1, \ \twistNad {67} (X) = 1, \ \twistNad {163} (X) = 1.
		\end{aligned}
  \end{equation} 
\end{theorem}

\begin{proof}
    We take the limits of the identities listed in \Cref{Lemma: twN(X) as a sum of heaviside functions for m of nonzero genus}. Alternately, examine the first table in \cite{Mazur2}.
\end{proof}

\section{Estimates for rational isomorphism classes when \texorpdfstring{$X_0(m)$}{X0(m)} has nonzero genus}\label{Section: Working over the rationals for m of nonzero genus}

In this section, we use the values we read off of Table \ref{table:nonzerogenus}, together with Walfisz's and Liu's asymptotics (\Cref{Theorem: Estimating squarefree(X)} and \Cref{Theorem: Estimating squarefree(X) given RH}) to prove \Cref{Intro Theorem: asymptotic for NQ(X) for m of nonzero genus}.

\Cref{Lemma: twN(X) as a sum of heaviside functions for m of nonzero genus} implies that $\twistLad m (s)$, as a finite sum of terms of the form $n^{-s}$, is holomorphic on $\bbC$. For each elliptic curve $E$ occurring in Table \ref{table:nonzerogenus}, we have $j(E) \neq 0, 1728$, so no additional casework is necessary. It would therefore be straightforward to apply \Cref{Theorem: Landau's Tauberian theorem} to 
\begin{equation}
\LQad m (s) = \frac{2 \zeta(6s) \twistLad m (s)}{\zeta(12s)},
\end{equation}
and obtain results akin to \Cref{Theorem: asymptotic for NQ(X) for m = 7}. For 
\[
m \in \set{11, 14, 15, 17, 19, 21, 27, 37, 43, 67, 163},
\]
write
\begin{equation}\label{Equation: Defining const m for m of nonzero genus}
    \constad m \colonequals \frac{2 \twistLad m (1/6)}{\zeta(2)}.
\end{equation}
These constants are computed numerically in Table \ref{table:constantsfornonzerogenus} below.
For any $\epsilon > 0$, an argument along the lines we have sketched yields the asymptotic
\begin{equation}
\NQad m (X) = \constad m X^{1/6} + O(X^{1/12 + \epsilon})
\end{equation}
for $X \geq 1$.

However, we can do better than this.

Recall \eqref{Equation: Definition of squarefree count}, which defines
\begin{equation}
\kfree 2 (X) = \# \set{n \in \bbZ_{>0} : n \leq X, \ n \ \textup{squarefree}}
\end{equation}
to be the number of squarefree integers with size at most $X$.

\begin{lemma}\label{Lemma: NQ(X) as a sum of squarefree(X) for m of nonzero genus}
	Let $X > 0$ be arbitrary. We have the following identities:
		\begin{equation} 	\begin{aligned}
		\NQeq {11} (X) =& 2 \kfree 2\parent{\parent{X/77\,480\,172}^{1/6}} + 2 \kfree 2\parent{\parent{X/2\,923\,690\,572}^{1/6}} \\
        &+ 2 \kfree 2\parent{\parent{X/323\,181\,166\,572}^{1/6}}, \\
		\NQeq {14} (X) =& 2 \kfree 2\parent{\parent{X/259\,308}^{1/6}} + 2 \kfree 2\parent{\parent{X/842\,579\,500}^{1/6}}, \\
		\NQeq {15}(X) =& 2 \kfree 2\parent{\parent{X/234\,967\,500}^{1/6}} + 2 \kfree 2\parent{\parent{X/478\,550\,700}^{1/6}} \\
		&+ 2 \kfree 2\parent{\parent{X/126\,818\,068\,500}^{1/6}} + 2 \kfree 2\parent{\parent{X/23\,621\,749\,807\,500}^{1/6}}, \\ 
		\NQeq {17} (X) =& 2 \kfree 2\parent{\parent{X/4\,325\,629\,743\,607\,500}^{1/6}} + 2 \kfree 2\parent{\parent{X/335\,734\,876\,712\,407\,500}^{1/6}}, \\
		\NQeq {19} (X) =& 2 \kfree 2\parent{\parent{X/14\,074\,668}^{1/6}}, \\
		\NQeq {21}(X) =& 2 \kfree 2\parent{\parent{X/364\,500}^{1/6}} + 2 \kfree 2\parent{\parent{X/1\,853\,388}^{1/6}}\\
		&+ 2 \kfree 2\parent{\parent{X/57\,395\,607\,372}^{1/6}} + 2 \kfree 2\parent{\parent{X/20\,478\,321\,699\,372}^{1/6}}, \\ 
		\NQeq {27} (X) =& 2 \kfree 2\parent{\parent{X/6\,912\,972}^{1/6}}, \\		
		\NQeq {37} (X) =& 2 \kfree 2\parent{\parent{X/7\,306\,267\,500}^{1/6}} + 2 \kfree 2\parent{\parent{X/107\,611\,181\,539\,805\,427\,907\,500}^{1/6}}, \\				
		\NQeq {43} (X) =& 2 \kfree 2\parent{\parent{X/162\,830\,654\,028}^{1/6}}, \\		
		\NQeq {67} (X) =& 2 \kfree 2\parent{\parent{X/102\,480\,782\,771\,052}^{1/6}}, \\		
		\NQeq {163} (X) =& 2 \kfree 2\parent{\parent{X/2\,631\,905\,352\,272\,628\,650\,988}^{1/6}}.
		\end{aligned}\end{equation} 	
\end{lemma}

\begin{proof}
	Let $E$ be an elliptic curve from Table \ref{table:nonzerogenus}. Recalling \eqref{Equation: quadratic twists multiply height by c^6}, we have
    \begin{equation}
    \begin{aligned}
        \# \set{E^{(c)} : c \in \bbZ, \ c \ \text{squarefree}, \ \hht(E^{(c)}) \leq X} %
        =& 2 \kfree 2 ((X/\hht(E))^{1/6}).
    \end{aligned}\end{equation} 
    Summing over the elliptic curves $E$ associated to each $m$, our claim follows.
\end{proof}

We are ready to prove \Cref{Intro Theorem: asymptotic for NQ(X) for m of nonzero genus}, with a modestly improved error term.

\begin{theorem}\label{Theorem: asymptotic for NQ(X) for m of nonzero genus}
	Let $m \in \set{11, 14, 15, 17, 19, 21, 27, 37, 43, 67, 163}$, and let $\constad m$ be given by \eqref{Equation: Defining const m for m of nonzero genus}. Then for $\kappa$ sufficiently small, we have
	\begin{equation}
	\NQeq m (X) = \constad m X^{1/6} + O\parent{X^{1/12} e^{-\kappa \frac{\log^{3/5} X}{\log^{1/5} \log X}}}
	\end{equation}
	for $X \geq 2$. The implicit constant depends on $\kappa$ and $m$.
\end{theorem}

\begin{proof}
    We substitute the asymptotic for $\kfree 2(X)$ given by \Cref{Theorem: Estimating squarefree(X)} into the identities of \Cref{Lemma: NQ(X) as a sum of squarefree(X) for m of nonzero genus}.
\end{proof}

In the presence of the Riemann hypothesis, we can use \Cref{Theorem: Estimating squarefree(X) given RH} to say even more.

\begin{theorem}\label{Theorem: asymptotic for NQ(X) for m of nonzero genus given RH}
	Let $m \in \set{11, 14, 15, 17, 19, 21, 27, 37, 43, 67, 163}$, and let $\constad m$ be given by \eqref{Equation: Defining const m for m of nonzero genus}. Then for any $\epsilon > 0$, we have
	\begin{equation}
	\NQeq m (X) = \constad m X^{1/6} + O\parent{X^{11/210 + \epsilon}}
	\end{equation}
	for $X \geq 1$. The implicit constant depends on $\epsilon$ and $m$.
\end{theorem}

\begin{proof}
    We substitute the asymptotic for $\kfree 2(X)$ given by \Cref{Theorem: Estimating squarefree(X) given RH} into the identities of \Cref{Lemma: NQ(X) as a sum of squarefree(X) for m of nonzero genus}.
\end{proof}

\jvtable{table:constantsfornonzerogenus}{
\rowcolors{2}{white}{gray!10}
\begin{tabular}{c c}
$m$ & $\constad m$ \\
$11$ & $0.10571705075608457\ldots$ \\
$14$ & $0.19181968564707050\ldots$ \\
$15$ & $0.11675063269362475\ldots$ \\
$17$ & $0.0044705452369290256\ldots$ \\
$19$ & $0.07824834140601367\ldots$ \\
$21$ & $0.28048805576004344\ldots$ \\
$27$ & $0.08809218206211060\ldots$ \\
$37$ & $0.027777766140563242\ldots$ \\
$43$ & $0.01645351934345495\ldots$ \\
$67$ & $0.005620493220876168\ldots$ \\
$163$ & $0.00032721745018530919\ldots$ \\
\end{tabular}
}{$\constad m$ when $X_0(m)$ has nonzero genus}